\documentclass{acta_2024}
\usepackage{amssymb,upquote,enumitem}
\usepackage{amsmath}  

\def\ccite#1{(\cite{#1})}
\def\ODEs{ODE\kern .4pt s}
\def\PDEs{PDE\kern .4pt s}
\def\RBFs{RBF\kern .4pt s}
\def\BVPs{BVP\kern .4pt s}
\def\real{\mathbb{R}}
\def\complex{\mathbb{C}}
\def\ep{\varepsilon}
\def\pO{\partial\kern .7pt \Omega}
\def\half{{\textstyle{1\over 2}}}
\def\C{{\kern .7pt\cal C}}
\def\Cp{\complex_+}
\def\Cpp{\overline{\Cp}}
\def\Cm{\complex_-}
\def\Cmp{\overline{\Cm}}
\def\In{I_n^{}}
\def\Rn{{\cal R}_n}
\def\signEF{\hbox{sign}_{E/F}^{}}
\def\Rmn{{\cal R}_{mn}^{}}
\def\xz{\vec x_0^{}}
\def\zk{\zeta_{\kern .5pt k}}
\def\zn{\zeta_{\kern .5pt n}}
\def\znn{\zeta_{\kern .5pt 2n}}
\def\zz{\zeta_{\kern .5pt 0}}
\def\pk{\pi_k^{}}
\def\fj{f_j^{}}
\def\fk{f_k^{}}
\def\wk{w_k^{}}

\def\In{I_n^{}}
\def\zj{z_j^{}}
\def\tk{t_k^{}}
\def\tj{t_j^{}}
\def\ak{\alpha_k^{}}
\def\bk{\beta_k^{}}
\def\bj{\beta_j^{}}
\def\bb{\beta}
\def\rnh{\hat r_n^{}}
\def\rns{r_n^*}
\def\Re{\hbox{\rm Re\kern .5pt}}\def\Im{\hbox{\rm Im\kern .5pt}}
\def\pput(#1,#2)#3{\noindent\smash{\raise#2pt\hbox to 0pt
   {\kern #1pt #3\hss}}\ignorespaces}

   \usepackage[longnamesfirst]{natbib} 
   \setcitestyle{comma,aysep={}} 
   \shortcites{} 
   
\DeclareSymbolFont{CMlargesymbols}{OMX}{cmex}{m}{n} 
\DeclareMathDelimiter{(}{\mathopen} {operators}{"28}{CMlargesymbols}{"00}
\DeclareMathDelimiter{)}{\mathclose}{operators}{"29}{CMlargesymbols}{"01}
\DeclareMathAlphabet\mathcal{OMS}{cmsy}{m}{n} 
\SetMathAlphabet\mathcal{bold}{OMS}{cmsy}{b}{n} 
\pagerange{\pageref{firstpage}--??}
\doi{XXXXXXXX}  
\usepackage[twoside,
	paperheight=247mm,
	paperwidth=174mm,
	marginratio={3:5,2:3},
	left=18mm,
	top=20mm]{geometry}
\usepackage{amssymb,amsmath}
\numberwithin{figure}{section}
\numberwithin{table}{section}
\usepackage{graphics,graphicx}
\usepackage[hang]{subfigure}
\usepackage[UKenglish]{babel}
\newcommand{\ignore}[1]{}
\allowdisplaybreaks

\newtheorem{theorem}{Theorem}[section]

\title[Applications of AAA rational approximation]{Applications of AAA rational approximation}

\author[Nakatsukasa and Trefethen]{Yuji Nakatsukasa\\
Mathematical Institute, University of Oxford,\\
Woodstock Road, Oxford, OX2 6GG, UK\\
E-mail: {nakatsukasa@maths.ox.ac.uk}\and
Lloyd N. Trefethen\\
School of Engineering and Applied Sciences, Harvard University,\\
29 Oxford St., Cambridge, MA 02138, USA\\
E-mail: {trefethen@seas.harvard.edu}}

\begin{document}

\label{firstpage}
\maketitle

\pput(280,150){\small\sl Draft \today}

\begin{abstract}
The AAA algorithm for rational approximation is employed to illustrate
applications of rational functions all across numerical
analysis.
\end{abstract}

\tableofcontents 

\section{\label{sec:intro}Introduction}

Suppose we had an algorithm that could quickly compute a rational
function $r$ to approximate almost any function $f$ on almost any
real or complex domain $Z$ to machine precision.  What could we do
with such a tool?  That is the topic of this review.

Such an algorithm, the AAA algorithm, was introduced in 2018
\ccite{aaa}.  Section \ref{sec:AAA} describes how it works, with
extensions presented in sections \ref{sec:lawson}--\ref{sec:sets},
but our focus will be on applications, not the algorithm itself.

In most areas of numerical analysis, the starting point is
polynomials.  We see this for example with series expansions
(Gregory, Newton, Taylor), interpolation (Harriot, Gregory, Newton,
Waring), quadrature (Gregory, Newton, Cotes, Simpson, Gauss,
Jacobi), rootfinding (Newton, Raphson, Simpson, Good), splines
(Schoenberg, de Casteljau, Bezier, de Boor), matrix iterations
(Lanczos, Hestenes, Stiefel, Saad), solution of \ODEs\ (Adams,
Runge, Heun, Kutta) and solution of \PDEs\ by finite difference
(Richardson, Southwell, Courant, Friedrichs, Lewy), finite element
(Courant, Argyris, Clough) and spectral methods (Lanczos, Clenshaw,
Elliott, Orszag, Gottlieb).  The names listed are some early figures
in each area, not a systematic list.

Rational functions, which are quotients of the form $r(z)
= p(z)/q(z)$ where $p(z)$ and $q(z)$ are polynomials, are
generalisations of polynomials since polynomials correspond to
the case $q=1$.  The reason for the power of rational functions
is clearer, however, if we think of polynomials as the special
case of rational functions whose poles are constrained all to lie
at $z=\infty$.  For approximating analytic functions far from any
singularities, this constraint is mild, and polynomials are useful
for all kinds of problems, giving fast exponential convergence.
Bulletproof algorithms are available for working with them in many
settings, often based on Chebyshev polynomials, Chebyshev points and
Chebyshev interpolants \ccite{atap}.  Many of these reliable methods
are brought together in the Chebfun software system \ccite{chebfun}.
For robust polynomial computation on regions other than intervals
and disks, see \ccite{VA}.

In other circumstances, however, polynomials lose their power.  This
is especially true in approximation of functions with singularities
on or near the approximation domain, on unbounded domains, or in
applications requiring extrapolation off the domain of approximation.
Here rational functions, whose poles may go wherever they are needed,
have great advantages.  Nevertheless, they have not played a large
role in numerical analysis, and the biggest reason for this has
been the lack of algorithms for working with them easily.

Our aim in this paper is to show how rational functions can change
the game in 26 different application areas, the topics of sections
\ref{sec:approx}--\ref{sec:hilbert}.  We present AAA code snippets
in MATLAB with many figures to make it easy for readers to adapt
rational functions to their own problems.  We are equally interested
in outlining the state of theoretical understanding of these methods,
and in each section we strive to point to results that are known---or
not known---to delineate the possibilities.  (Twenty open problems
are listed in the Discussion at the end.)  We also cite a good
deal of related literature on numerical methods involving rational
functions, though it is impossible to be exhaustive.

Historically, the AAA algorithm represents a break with the past:
what used to be difficult is now usually easy.  Two essential
features of its success are that it is a greedy descent algorithm
rather than aiming to enforce optimality conditions (thereby making
it somewhat analogous to gradient descent methods in machine
learning) and that it employs a barycentric representation of a
rational function rather than the exponentially unstable quotient
representation $r(z) = p(z)/q(z)$.  However, we do not claim
that AAA is the only method that can be effective for numerical
rational approximation.  Recently, Salazar Celis and Driscoll have
shown that similar speed and accuracy can be achieved in many
cases by a greedy interpolation algorithm based on an entirely
different representation, namely Thiele continued fractions
\ccite{salazar,driscollzhou,driscolljuliacon}.  We do not know if
AAA will remain the best approach for these problems in the decades
to come, but it seems clear that it has ushered in a new era for
rational functions.

Throughout the paper we use this notation:
\smallskip

\noindent\kern 1.5cm  \hbox to 1.47in{$\{\pk\}$: poles,\hfill}
$\{\zk\}$: interpolation points,
\vskip 2pt

\noindent\kern 1.5cm \hbox to 1.47in{$\{z_j^{}\}$: sample points,\hfill}
$\{\tk\}$: barycentric support points,
\vskip 2pt

\noindent\kern 1.5cm \hbox to 1.47in{$\{f_j^{}\}$: sample values,\hfill}
$\{\bk\}$: barycentric weights.

\section{\label{sec:AAA}The AAA algorithm}

We start with a vector $Z\in \complex^m$ of $m$
distinct real or complex numbers and a vector
$F\in\complex^m$ of associated function values.  The AAA algorithm aims to find
a rational function $r$ such that the maximum error
$\|r(Z)- F\|$ is small: our standard convergence condition
is $\|r(Z)- F\|/\|F\| \le 10^{-13}$.
(Throughout this paper, $\|\cdot\|$
without a subscript is the $\infty$-norm for a vector or
function, depending on context.)
The approximation is represented in barycentric form,
\begin{equation}
r(z) = {N(z)\over D(z)} = \sum_{k=0}^n {\fk\kern .5pt \bk\over z-\tk} \left/
\sum_{k=0}^n {\bk\over z-\tk} \right.,
\label{aaa-baryrep}
\end{equation}
where $N$ and $D$ stand for numerator and denominator and normally $n\le (m-1)/2$.
The numbers $\{\kern .7pt \bk\}$ are {\em barycentric weights}, and
the numbers $\{\tk\}$, which are $n+1$ distinct entries of
$Z$, are {\em support points} or {\em nodes}.
If we assume that all the barycentric weights are nonzero, then
each $\tk$ is a pole
of the numerator in (\ref{aaa-baryrep}) (assuming $\fk\ne 0$)
and also a pole of the denominator.  The poles
cancel out in the quotient, giving the interpolatory property
\begin{equation}
\lim_{z\to\tk} r(z) = \fk.
\end{equation}
Thus $\tk$ is a removable singularity of (\ref{aaa-baryrep}),
and we define $r(\tk) = \fk$ accordingly.  As can be verified
by multiplying $N(z)$ and $D(z)$ by $\prod_{k=0}^n(z-\tk)$, $r$
belongs to $\Rn$, the set of rational functions of degree $n$, i.e.,
the set of quotients $p(z)/q(z)$ where $p$ and $q$ are polynomials
of degree $\le n$.

The use of a representation of rational functions different from the
obvious choice $p(z)/q(z)$ is crucial to the success of AAA.\ \ Good
rational approximations need to have poles and zeros clustering at
singularities, with the consequence that both $|\kern .7pt p(z)|$ and
$|q(z)|$ in such a representation will have to vary widely with $z$.
Such variation quickly leads to loss of accuracy in floating point
arithmetic if $p(z)$ and $q(z)$ are expressed in a global basis,
such as monomials or Chebyshev polynomials.  In other words, the
$p(z)/q(z)$ representation of rational functions is unstable.
This is the reason why Varga and his collaborators had to use
200-digit extended precision arithmetic in the 1990\kern .4pt s to
compute best approximations accurate to 16 digits \ccite{varga}.
By contrast, with the help of the barycentric representation,
the AAA algorithm reliably achieves close to machine accuracy in
standard 16-digit arithmetic.

There remains the problem of how to find good
values of $\{\tk\}$ and $\{\kern .7pt \bk\}$.  Before AAA, most methods for
rational approximation followed nonlinear optimisation strategies motivated by
an attempt to satisfy optimality criteria, often involving an iterative adjustment
of poles.
The principal methods we are aware of are presented in
\ccite{sk},
\ccite{osborne}, 
\ccite{ew},
\ccite{istace},
\ccite{vf} (``Vector fitting''),
\ccite{agh},
\ccite{irka} (``IRKA''),
and~\ccite{rkfit} (``RKFIT'').
Crucially, AAA is not based on optimality criteria, and it does
not work with poles.  It is a greedy
algorithm, which just aims locally for descent.  It needs no initial
guess to get started, and no specification of the degree of the
approximation (though a degree can be specified if desired).

Let the nonlinear approximation problem $F\approx N(Z)/D(Z)$ be linearised
to $F\cdot D(Z) \approx N(Z)$, that is,
\begin{equation}
f(z) \sum_{k=0}^n {\beta_k^{}\over z-\tk} 
\,\approx\, \sum_{k=0}^n {\fk\kern .7pt \beta_k^{}\over z-\tk} , \quad z\in Z
\label{aaa-AAAeq}
\end{equation}
together with the normalisation
\begin{equation}
\sum_{k=0}^n |\kern .5pt \beta_k|^2 = 1.
\label{aaa-normalize}
\end{equation}
The AAA algorithm solves a sequence of problems 
(\ref{aaa-AAAeq})--(\ref{aaa-normalize})
with $n = 0, 1, 2,\dots.$  The support points $t_0^{}, t_1^{},\dots$ are chosen along
the way, one by one, from $Z$.  This means that the least-squares problem at step $n$
actually involves sampling in the remaining $m-n-1$ points.  To be precise, at step
$n$, we compute weights $\{\kern .7pt \beta_k^{}\}$ satisfying (\ref{aaa-normalize})
such that
\begin{equation}
\sum_{j=1}^m{}' \left[ \fj \sum_{k=0}^n {\beta_k^{}\over \zj-\tk} 
- \sum_{k=0}^n {\fk\kern .7pt \beta_k^{}\over \zj-\tk} \right]^2 = \,
\hbox{minimum},
\label{aaa-AAA2}
\end{equation}
where the prime on the summation sign indicates that the summation extends over
the points $Z\backslash \{t_0^{},\dots,t_n^{}\}$.
This is a problem of numerical linear algebra involving
an $(m-n-1)\times (n+1)$ matrix $A$ with entries
\begin{equation}
a_{jk} = {\fj-\fk \over \zj-\tk} ,
\label{aaa-Aentries}
\end{equation}
known as a {\em Loewner matrix\/} \ccite{abg}.
If $\bb$ denotes the
column vector with entries $\beta_0,\dots,\beta_n$,
the problem (\ref{aaa-normalize})--(\ref{aaa-AAA2}) is to find a vector $\bb$ satisfying
\begin{equation}
\|A\bb \|_2^{} = \hbox{minimum}, \quad \|\bb\|_2^{}=1,
\label{aaa-AAAprob}
\end{equation}
where $\|\cdot\|_2^{}$ is the vector 2-norm.  It is routine
matter to solve (\ref{aaa-AAAprob}) by computing the singular
value decomposition (SVD) of $A$ and setting $\bb$ equal to the
right singular vector corresponding to the smallest singular value
$\sigma_{n+1}$.  (In section~\ref{sec:zolo} we will modify
(\ref{aaa-AAAprob}) in certain contexts and also mention
an adjustment introduced in 2025 to improve accuracy in
cases where $A$ is highly ill-conditioned.)

The other half of the algorithm is the greedy choice of
support points.  At step $n$,
the next support point $t_n$ is selected as a sample point $z_j$ where
the approximation error $|\kern .4pt \fj - r(\zj)|$ is maximal, and then the next SVD calculation
is carried out.  We can write the full algorithm like this:
\vskip 1pt
\begin{displaymath}
\parskip=2pt
\framebox[4.4in][c]{\parbox{4.2in}{\vspace{2pt}\em
{\sc AAA Algorithm\vrule width 0pt depth 8pt}\hfill\break
Given: $m$-vectors $Z$ of sample points
and\/ $F$ of \kern 1pt function values.\hfill\break
$r := \hbox{the zero function}$.\hfill\break
{\bf for} $n = 0,1,2,\dots$:\hfill\break
\indent ~~~~~$t_n^{} := \hbox{a sample point } \zj \hbox{ where }
|\kern .4pt \fj - r(\zj)| \hbox{ is maximal}$;\hfill\break
\indent ~~~~~Construct the $(m-n-1)\times (n+1)$ matrix $A$ from $(\ref{aaa-Aentries})$;\hfill\break
\indent ~~~~~Compute the vector $\bb$ of $(\ref{aaa-AAAprob})$ from the SVD of $A$;\hfill\break
\indent ~~~~~Define\/ $r$ by $(\ref{aaa-baryrep})$.
\vspace{2pt}}}
\end{displaymath}
\vskip 5pt

\noindent Depending on one's purposes, the iteration may be run for
a fixed number of steps or until a convergence criterion is met.
In the Chebfun code {\tt aaa.m}, the default condition is that
the maximal error should be less than $10^{-13}$ times the maximal
value $|\kern .4pt \fk|$.  There is no guarantee that this will be
achieved in the presence of rounding errors, but the reliability
in practice is striking.

The name AAA (``triple-A'') comes from ``adaptive
Antoulas-Anderson.''  Antoulas and Anderson \ccite{aa} proposed
the use of the barycentric formula for rational interpolation (see
section 11 of \ccite{aaa} for a discussion of the history), and the
AAA algorithm couples this with its adaptive selection of support
points and the linearised least-squares determination of barycentric
weights.  The original AAA paper includes a basic MATLAB code, which
evolved with various checks and extra features into the {\tt aaa.m}
command in Chebfun, which has been used for all the examples in this
paper except for the black curves in Figure~\ref{analcont-six}.
Different MATLAB implementations are available in the MATLAB RF,
Control System and System Identification toolboxes \ccite{rf},
and there are also open-source AAA implementations in Julia
\ccite{driscolljuliacon,macmillen} and Python/SciPy \ccite{virtanen}.

\label{symmetrypage}
Quite a few modifications to the basic AAA algorithm can be added to
introduce useful properties or enhance its reliability.  The most
important property not imposed by the algorithm as described above
is real symmetry.  If $Z$ is a complex set with $\bar Z = Z$, and if
$F$ takes conjugate values at conjugate points, then it is natural
to want $r$ to have the property $r(\bar z) = \overline{r(z)}$.
This is not achieved by the basic AAA algorithm, even in the absence
of rounding errors.  A number of authors, beginning with Hochman
even before the original AAA paper appeared in print \ccite{hochman}
and also including the MathWorks RF Toolbox \ccite{rf}, have fixed
that by introducing new support points as conjugate pairs when they
are complex, rather than one-by-one.  (This seems to bring a slight
increase in accuracy, too.)  Regarding reliability improvements,
as mentioned above, one that matters occasionally in practice
is to modify the calculation so that the barycentric weights are
constructed from several right singular vectors rather than just
the minimal one.  Another is to apply column scaling to the matrix
$A$ when it is highly ill-conditioned.  We shall discuss these
modifications in section \ref{sec:zolo}.

The computational cost of the AAA algorithm depends strongly on
the degree $n$.  The iteration involves SVD computations of sizes
$(m-1)\times 1, (m-2)\times 2, \dots, (m-n-1)\times (n+1)$, for
a total operation count of the order of $mn^3$, which amounts
to $n^4$ since $m$ will generally need to scale with $n$ for
effective resolution of a function.  Most applications we know
of involve reasonably small degrees and sample sets, say $n <
150$ and $m < 3000$, giving computation times usually under one
second on today's laptop and desktop machines.  In this regime
it is not clear that any good methods are available to speed up
the linear algebra significantly, though some proposals have been
put forward in \ccite{hochman,setAAA} (based on updating Cholesky
factorisations) and \ccite{park} (based on randomised sketching).
(Hochman's work reduces the complexity to $O(n^3)$ in theory, but
it is based on the Cholesky factorisation of the Gram matrix, which
has stability issues associated with squaring the condition number.)
The $O(n^4)$ operation count certainly stands out, however, and we
presume that as time passes and larger scale rational approximations
become more common, the challenge of speeding up the linear algebra
may grow in importance.

One may ask what theory there may be to guarantee the success of AAA.\ \ In
a certain superficial sense, the answer is easy, as pointed out
first in \ccite{hochman}:
\begin{theorem}
In exact arithmetic, AAA always converges to an approximation $r$
satisfying the prescribed accuracy criterion $\|r(Z)-F\|/\|F\| \le \hbox{tol}$.
\end{theorem}
The proof is merely the observation that in the worst case,
a rational function can simply interpolate all the data values.
This will happen at the latest at step $n = \lfloor m/2 \rfloor$ if
all the barycentric weights are nonzero at that step.  If there are
zero weights at this step and they keep reappearing at subsequent
steps, then one could continue all the way to $n=m-1$, at which
point all the sample points are support points and the weights can
be arbitrarily set to $1$ or to any other nonzero values.

This mathematically valid theorem, however, is a ``Yogiism''
\ccite{yogi} in that it misses the actual purpose of most
approximations.  Convergence on a discrete set $Z$ does not
imply convergence on a continuum that $E$ is intended to
approximate.\footnote{With polynomials, this is the familiar
lesson of the Runge phenomenon \ccite{atap}: interpolation of all
the data does not ensure that anything useful has been achieved.}
Specifically, the convergence theorem does not exclude ``bad poles''
in the approximant, lying in regions near $Z$ where one doesn't
want or expect them.  The most familiar instance is the situation
in which poles turn up between the sample points in approximation
of a real function on a discretised real interval---in odd-degree
approximations of $|x|$ on $[-1,1]$, for example, but also in
many other problems without such symmetries.  This is not just
a theoretical problem but a practical one.  AAA is at its least
reliable for real approximation on real intervals.  One way that
trouble may arise in AAA approximations is in the form of ``spurious
poles,'' typically with residues at the level of machine epsilon or
other noise in the data; one also speaks of ``Froissart doublets''
since these are poles that are paired with zeros at almost but
not exactly the same position, cancelling to invisibility further
away \ccite{stahl98}.  Spurious poles or Froissart doublets will be
mentioned at various points in this paper, notably near the ends of
sections \ref{sec:poles}, \ref{sec:equi} and \ref{sec:potential}.
The code {\tt aaa.m} has a ``cleanup'' feature that can sometimes
remove spurious poles, first mentioned in \cite{aaa}---with also
a \verb|'cleanup',2| variant due to Stefano Costa---but we will
not discuss cleanup options in this paper, as they are not fully
understood and the experiments here do not take advantage of them.

Another case where the choice of $Z$ is an issue is when $F$ has one
or more singularities on the continuum $E$ of ultimate interest, so
that poles and zeros need to be exponentially clustered near there.
In such a situation it is important to make sure that $Z$ likewise
contains points exponentially clustered near the singularities.
Usually one would also loosen the tolerance, perhaps to $10^{-6}$
or $10^{-8}$, since resolution of a singularity to 13-digit accuracy
will usually require rational degrees in the hundreds.  Good rational
approximations of functions with branch point singularities converge
root-exponentially, i.e., at a rate $\exp(-C\sqrt n \kern 1pt)$
for some $C>0$ as $n\to\infty$ \ccite{newman,gt}.  With polynomial
approximations, by contrast, since the convergence is only algebraic,
one would typically need degrees in the millions to get 6 or 8
digits of accuracy.

To circumvent problems of discretisation of a complex
domain by a finite set $Z$, a continuum version of the AAA
algorithm was introduced in \ccite{continuum}; for software
see \ccite{driscolljuliacon}.   This applies {\em a priori\/}
to intervals, disks, circles, lines and half-planes, or to
more general domains via boundary parameterisation.  The method
appears very promising but has not yet become a standard, no doubt
in part because working with a discrete $Z$ is so straightforward
and usually successful.  In this review all our illustrations are
based on the discrete case, except that continuum AAA is discussed
and illustrated in section~\ref{sec:continuum}.

Returning to the matter of theoretical justification of AAA, suppose
we could guarantee that under suitable assumptions AAA (or some
modification of AAA) always converged throughout $E$ when applied on
a discrete set $Z$ approximating a continuum $E$ sufficiently finely.
This would be helpful, yet such a theorem would still miss the point
that usually what is wanted is not just ``any old approximant'' but
one that compresses the data. This could be for the simple purpose
of reducing the number of data values needed in a representation,
but more often and more interestingly, the point of an approximation
$r\approx f$ is to estimate information about $f$ such as its values,
or the locations of its poles or other singularities, beyond the data
domain.  For this to happen, an approximation needs to be near-best,
that is, not too far from the optimal approximation of its degree.
The true theoretical challenge of AAA approximation is to establish
a theorem to ensure that, under appropriate assumptions, AAA finds
not just accurate but near-best approximations.

An interesting analogy can be drawn between the AAA algorithm
for rational approximation and the QR algorithm for eigenvalues
of matrices.  (The analogy is mentioned in a footnote in
\ccite{continuum}.)  Both algorithms consist of an alternation
between an elementary nonlinear step (choice of the next support
point for AAA, shifting by a multiple of the identity for QR)
and a routine matrix computation (tall-skinny SVD for AAA, QR
factorisation and recombination for QR).\ \ Both turn problems that
were previously troublesome into black-box solvers that almost always
get the desired answer fast.  In both cases the basic algorithm is
extraordinarily simple to describe and highly reliable in practice,
but not quite infallible.  In both cases fully reliable performance
requires a few extra tricks that users are likely to be unaware of
(for QR, ``balancing'' the matrix is in this category), yet rigorous
guarantees of success may still be elusive.

The rest of this paper is mostly about AAA applications, not about
the algorithm, though some algorithmic details will be discussed
in sections~\ref{sec:poles} (computation of poles and zeros),
\ref{sec:diff} (computation of derivatives),
and \ref{sec:zolo} (the \verb|'sign'| and
\verb|'damping'| features).
Generalisations of AAA to best approximation, approximation on a continuum,
periodic functions, type $(m,n)$ approximation, and vector- or matrix-valued functions are
discussed in sections \ref{sec:lawson}--\ref{sec:sets}.

\section{\label{sec:approx}Approximating functions}

Our first application of rational approximation is to the
approximation of functions. That may sound tautological, but there are
things to be said.  A key distinction is between what we might call
{\em off\kern .5pt line\/} or {\em optimised\/} approximation and
{\em online\/} or {\em on-the-fly\/} approximation.  AAA contributes
to both, but it is particularly striking for on-the-fly applications,
so we begin there.  As in every section of this paper, our aim
is to provide ideas and references and illustrations that may be
useful for the reader, not to attempt a comprehensive account of
all the possibilities.

The key features are speed, reliability and flexibility about
sample sets.  For example, suppose the values of a meromorphic
function $f$ are known at 50 scattered sample points in $\complex$.
(A meromorphic function is one that is
analytic apart from poles.)  For illustration we take $f(z) = \tan(z)/\tan(2)$
and a set of 50 complex normally distributed random points,
\begin{figure}
\vskip 1pt
\begin{center}
\includegraphics[trim=40 170 215 10,clip,scale=1.2]{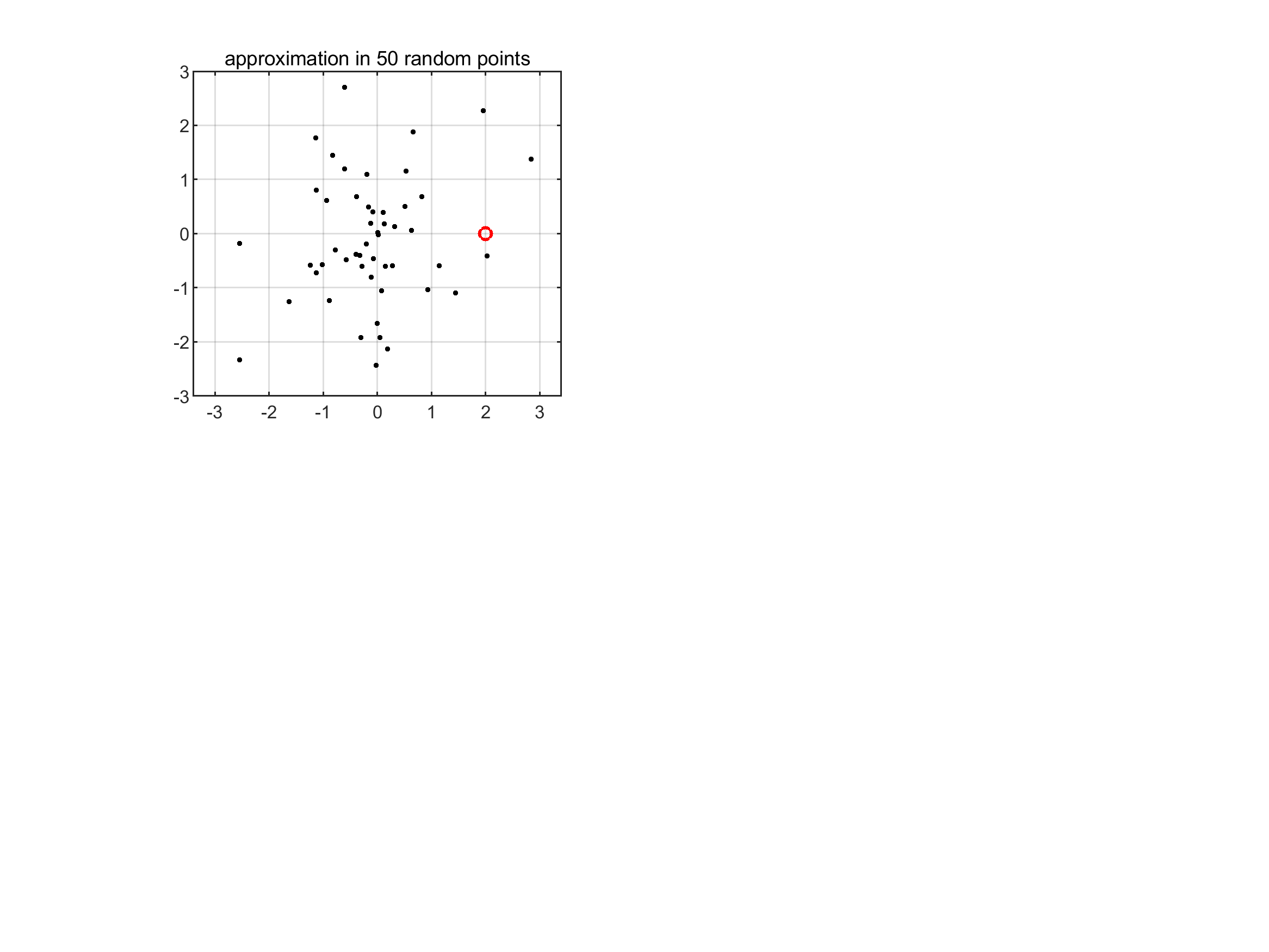}
\end{center}
\vskip -12pt
\caption{\label{fig3.1}Fifty scattered points in the complex plane
at which a meromorphic function $f(z)$ is sampled.  The aim is
to estimate the value at $z=2$ (the red circle).  
For $f(z) = \tan(z)/\tan(2)$, AAA approximation gets a result
accurate to 11 digits in a few milliseconds.}
\end{figure}

{\small
\begin{verbatim}
      rng(2)
      Z = randn(50,1) + 1i*randn(50,1);
      F = tan(Z)/tan(2);
\end{verbatim}
\par}

\noindent
as plotted in Figure~\ref{fig3.1}.
Suppose we want to evaluate $f$ at an arbitrary point---for
example, at $z=2$, where the exact value is 1.  The command

{\small
\begin{verbatim}
      r = aaa(F,Z)
\end{verbatim}
\par}

\noindent produces a degree $10$ rational approximation $r\approx f$ in a couple
of milliseconds on a laptop, and for this choice of $f$,
$r(2)$ matches the true value $f(2)=1$ to 11 digits:

{\small
\begin{verbatim}
      r(2)
      ans =
        1.000000000008511 - 0.000000000007166i
\end{verbatim}
\par}

\noindent
If the random number seed in this code is set to the values
$1,2,\dots,10$ to give an indication of how the results depend
on the set of sample points, the numbers of digits of accuracy
(rounded values of $-\log_{10}|r(2)-1|$) come out as 12, 11, 10,
13, 10, 13, 14, 10, 14, 11 (mean: $11.8$).  If the number of sample
points is increased to $100$, the numbers of digits go up to 14,
15, 13, 15, 15, 13, 14, 12, 15, 13 (mean:$13.9$).

Rational approximation is not the only way in which one might
approximate $f$, but its simplicity and accuracy are remarkable. The
reader might find it interesting to consider the alternatives.
What computations could be carried out to obtain an analogous
approximation by a polynomial rather than a rational function, say,
even if $f$ didn't happen to have any poles?  (A well-conditioned
basis for the polynomial can be constructed by Vandermonde with
Arnoldi orthogonalisation~\ccite{VA}.)  Or perhaps by some other
fitting method?

\begin{figure}
\vskip 7pt
\begin{center}
\includegraphics[trim=50 140 20 0, clip, scale=0.84]{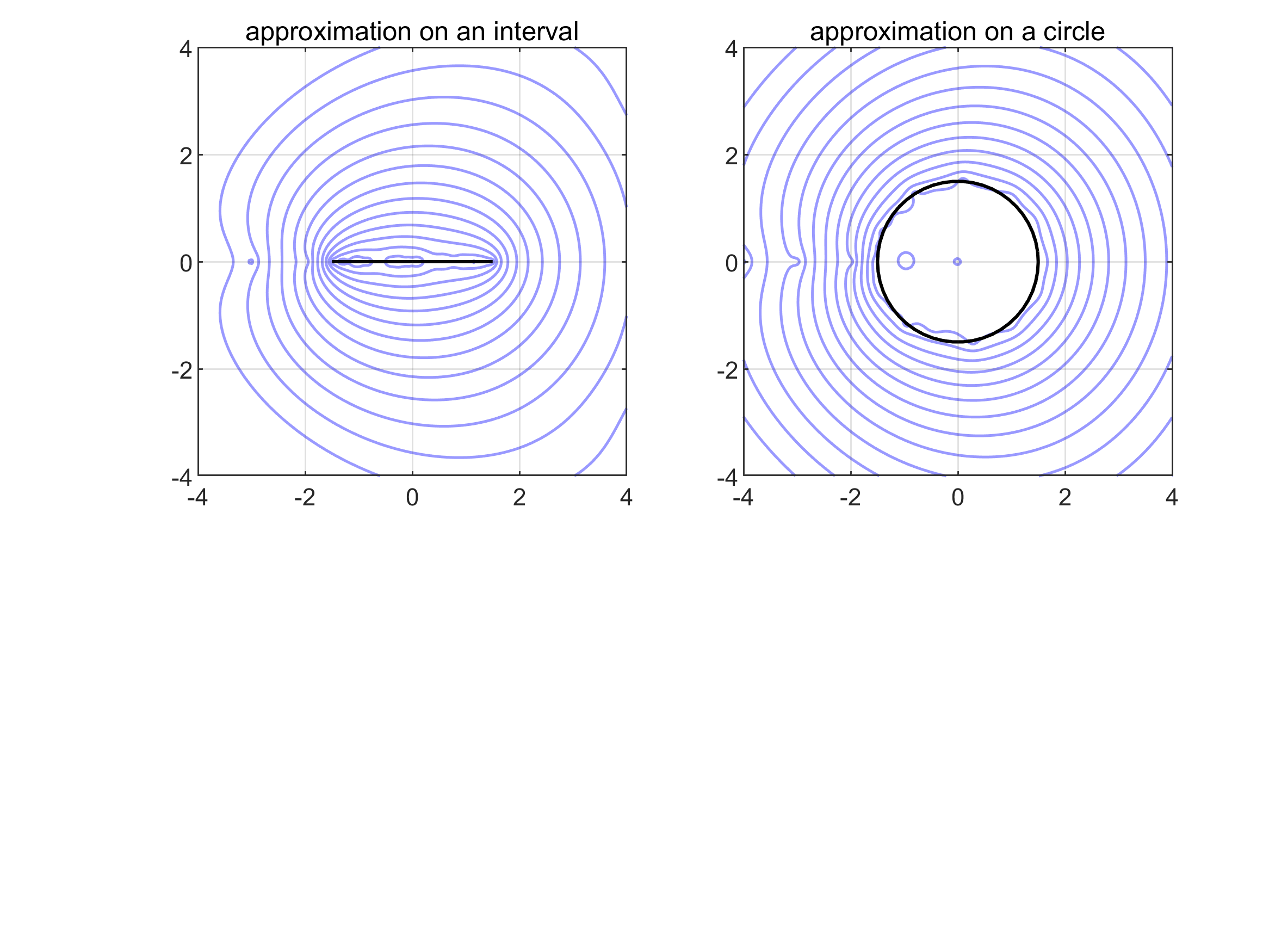}
\end{center}
\vskip -11pt
\caption{\label{fig3.2}AAA approximation of $\Gamma(z)$ based
on samples at $50$ Chebyshev points in $[-1.5,1.5]$ (left) and
$50$ equispaced points on $|z|=1.5$ (right).  From inside out,
the contours represent relative errors $|r(z)-\Gamma(z)|/|\Gamma(z)|
= 10^{-13}, 10^{-12}, \dots,10^{-1}$.}
\end{figure}

Other on-the-fly function approximation problems may start
from data on a continuum, such as an interval or a circle.
In Figure~\ref{fig3.2}, we suppose that a meromorphic function $f(z)$
is known on the interval $[-1.5,1.5]$ (left) or the circle $|z|=1.5$
and we wish to find a way to evaluate it at other points nearby.
For this illustration, we take $f$ to be the gamma function,
$f(z) = \Gamma(z)$.  The MATLAB {\tt gamma} command only evaluates
$\Gamma(z)$ for real arguments, so the left-hand part of the figure
represents a computation that a MATLAB user might find helpful:
a quick method for evaluating $\Gamma(z)$ for complex values of $z$
near the origin.  Again in a few milliseconds, a call to AAA computes
a rational approximation of degree $9$ (left) or $12$ (right) with a
relative error on the approximation domain of less than $10^{-13}$.
To evaluate $\Gamma(z)$ to full accuracy for complex $z$ for these
illustrations, we have used the code of \ccite{godfrey}.

Classically, the best known use of rational functions for
approximating functions would be in {\em off\kern .5pt line\/} mode,
notably in the design of a special functions library.  Here it may be
worthwhile to make greater efforts to obtain an approximation that
is truly optimal in some sense.  For example, the \hbox{FUNPACK}
\ccite{funpack} and \hbox{SPECFUN} \ccite{specfun} projects
of the last century used rational minimax (best $\infty$-norm)
approximations to evaluate Bessel functions, error functions, gamma
functions and exponential integrals.  A more recent discussion of
such methods can be found in \ccite{muller}.

For a continuous real function on a real interval, a minimax
approximation can be obtained with the Remez algorithm based on
equioscillation, whose most robust implementation that we are aware
of is that of the Chebfun {\tt minimax} command \ccite{minimax}.
Though much slower to compute than a AAA approximation, this
might typically have accuracy about one digit better, or it might
have comparable accuracy but with a rational degree about $10\%$
lower.  On domains other than an interval, the Remez algorithm is
not applicable.  In the last section we cited various published
algorithms for minimax rational approximation, but the best
available technique in practice may often be the extension of AAA
known as the AAA-Lawson algorithm \ccite{aaaL}, to be discussed in
section \ref{sec:lawson}.  Here, AAA approximation is run as usual
(perhaps with a specified degree rather than a specified convergence
tolerance) and then the approximation is improved by a subsequent
nonlinear iteratively reweighted least-squares iteration.

\begin{figure}
\vskip 7pt
\begin{center}
\includegraphics[trim=0 150 90 10, clip, scale=1]{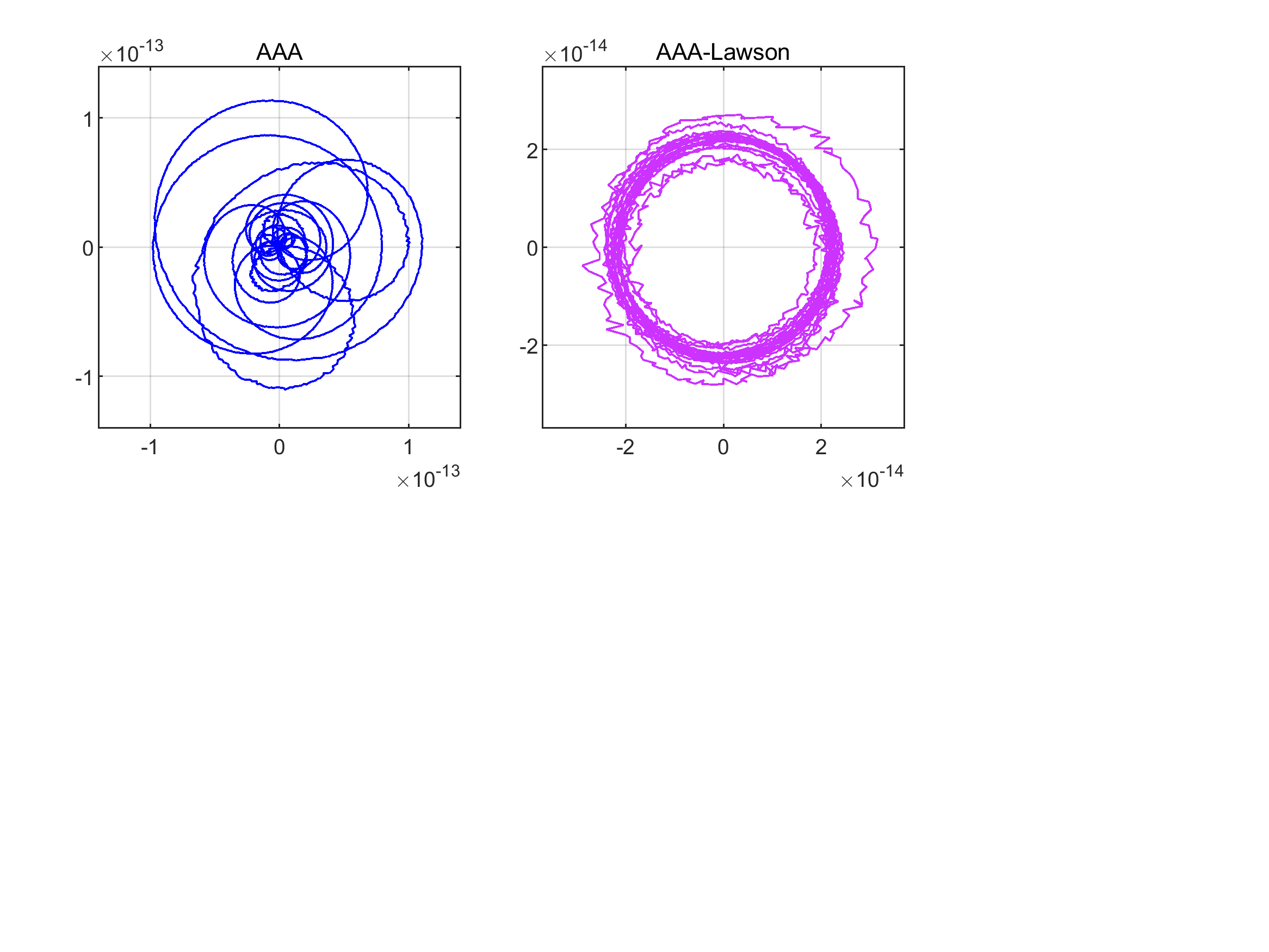}
\end{center}
\vskip -10pt
\caption{\label{fig3.3}Error curves for
AAA approximation of $\Gamma(z)$ on the circle $|z|=1.5$.
On the left, straight AAA approximation.  On the right,
AAA-Lawson approximation with 20 Lawson steps.  The error curve
has winding number $21$ but is ragged because the
errors are close to the level of machine precision (compare
Figure~\ref{gamma2fig}).
In this paper we use magenta to distinguish
error curves from AAA-Lawson best approximation.}
\end{figure}

Figure~\ref{fig3.3} illustrates the use of AAA-Lawson approximation
for the problem of approximating $\Gamma(z)$ on $|z|=1.5$ as
in Figure~\ref{fig3.2}.  Both images in Figure~\ref{fig3.3} show
error curves $e(z) = \Gamma(z) - r(z)$, where $z$ ranges over 5000
equispaced points on the circle of approximation.  On the left,
with standard AAA approximation, the error curve is irregular and the
absolute error is about $1.2\cdot 10^{-13}$, which corresponds to a
relative error less than $10^{-13}$ since $\Gamma(-1.5) \approx 2.4$.
This is not a best approximation of the given degree $12$, but a
near-best approximation.  On the right, 20 iterative Lawson steps
have been appended with a call to Chebfun {\tt aaa.m} of the form

{\small
\begin{verbatim}
      r = aaa(F,Z,'lawson',20)
\end{verbatim}
\par}

\noindent and the result is an error curve that is closer to
optimal, with the error reduced by a factor of about 4, to
$3\cdot 10^{-14}$.  This is a genuine improvement, at least in
the $\infty$-norm.  (The root-mean-square error has roughly
halved, from $4.9\cdot 10^{-14}$ to $2.2\cdot 10^{-14}$.)
We can see that the approximation is nearly optimal from the
fact that the error is approximately a circle of winding number
21 \ccite{nearcirc}.\label{footpage}\footnote{Near-optimality
follows by an argument related to Rouch\'e's theorem, since the
winding number $21$ is $\ge 2\kern .3pt n + 1 - 2\kern .3pt p$,
where $n=12$ is the degree of $r$ and $p=2$ is the number of poles
of $\Gamma$ within the circle of approximation.  To be precise in
this case, suppose $r=p/q$ with $p$ and $q$ of degree 12 and there
were another approximation $\tilde r = \tilde p/ \tilde q$ with
$\max_{|z|=1.5} |f(z)-\tilde r(z)|< \min_{|z|=1.5} |f(z)-r(z)|$.
Then $r-\tilde r = p/q - \tilde p /\tilde q$ would also have
winding number 21, so assuming $q$ and $\tilde q$ each have 2
zeros with $|z|<1.5$ corresponding to the poles of $\Gamma(z)$ at
$z=-1$ and $z=0$, $p\kern .3pt \tilde q - \tilde p\kern .3pt  q$
would have winding number $25$.  This means $p\kern .3pt \tilde q -
\tilde p\kern .3pt  q$ would have $25$ zeros with $|z|<1.5$, which
is not possible since it is a nonzero polynomial of degree $24$.}
The true best approximation would have an error curve circular
to plotting accuracy and a maximum error perhaps $10\%$ smaller,
but we cannot find it in this case because we are too close to
machine precision, as is revealed by the jaggedness of the plot.
Section \ref{sec:lawson} will show examples further from machine
precision, hence with cleaner behaviour.

The reason a discussion of AAA-Lawson approximation is deferred to
near the end of this article is that we do not recommend it as the
method to turn to initially for applications, because it is slower
and less reliable than standard AAA.\ \ It almost always runs into
difficulty at accuracies close to machine precision, and even at
accuracies well above that level, it does not achieve the ``99\%
reliability'' (not a precise figure) of AAA.\ \ Sometimes better
performance can be obtained by running the AAA-Lawson algorithm
with the ``sign'' and/or ``damping'' improvements described in
\ccite{zolo}; see sections~\ref{sec:zolo} and~\ref{sec:lawson}.

Here at the end of this section on approximating functions,
let us give a brief and non-rigorous summary of some theorems
about convergence of polynomial vs. rational best approximations
to an analytic function $f$ on a compact continuum $E\subset
\complex$, asymptotically as the degree $n$ approaches~$\infty$.
The convergence rates in these theorems depend on how $f$ can be
analytically continued to the rest of $\complex$, and they are
proved by arguments of potential theory to be outlined in section
\ref{sec:potential}.  If $f$ is an entire function---analytic
throughout $\complex$---then both polynomial and rational best
approximations converge superexponentially, i.e.\ with error
$O(C^{-n})$ for any $C>1$, though possibly at very different
rates despite that (Figure~\ref{potential-aaashowold}).  If $f$
is analytic apart from poles, whether finite or countably infinite
in number, then polynomials converge exponentially, i.e.\ at
the rate $O(C^{-n})$ for some $C>1$, whereas rational functions
still converge superexponentially (Figure~\ref{potential-pot2}).
The same applies if there are essential singularities as well as
poles---but in this case the superexponential convergence of rational
approximations is likely to look just exponential for the ranges
of $n$ accessible computationally.  If $f$ has one or more branch
points, then both polynomial and rational approximations converge
only exponentially---though again perhaps at very different rates
(Figure~\ref{potential-pot3}).  If $E$ is a nonconvex domain and
$f$ has branch points in typical locations \ccite{schwarzfun} in
reentrant ``inlets'' of $E$, then both polynomials and rational
functions converge exponentially, but the convergence constant
for the polynomials is exponentially worse as a function of
the length-to-width ratio of the inlet---i.e., $O(C^{-n})$ with
$C$ exponentially close to $1$ (Figure~\ref{potential-pot4}).
These results are mostly due to the Moscow school associated
with Gonchar, except the last from \ccite{inlets}; for surveys,
see (\cite{inlets,analcont}).  Finally, if $f$ is not analytic on
the boundary of $E$ but has one or more branch points there, then
polynomials converge just algebraically, e.g.\ at the rate $O(C/n)$,
whereas rational functions converge {\em root-exponentially,}
at the rate $O(C^{-\sqrt n})$ for some $C>1$ (\cite{newman,gt}).
This root-exponential convergence is made possible by exponential
clustering of poles and zeros near the singularities, a phenomenon
investigated in detail in \ccite{clustering} and visible in this
paper in Figures \ref{schwarz-squiggles3}, \ref{wienerhopffig},
\ref{bigaaals}, \ref{helmdaisy}, \ref{aaawedges} and \ref{hilbfig}.

\section{\label{sec:poles}Locating poles and zeros}

One of the long-established applications of rational functions is
the determination of poles and other singularities of analytic
functions in the complex plane.  (By ``analytic'' here we mean
analytic away from those singularities.)  The version of this
technique that is best known is based on Pad\'e approximation,
in which a function $f$ is approximated from knowledge of some of
its Taylor coefficients at a point \ccite{bgm,stahl}.  With AAA
approximation, we can pursue similar ideas when the data consists
of function values rather than Taylor coefficients.

Before considering such applications, we explain how poles and zeros are
computed by the AAA algorithm.  Poles and zeros play no
part in the construction of the AAA approximant, which adjusts
barycentric weights, not pole locations.  Only at the end, after
the approximant has been found, are they computed if desired.
This is done by solving a matrix generalised eigenvalue
problem involving the support points and barycentric coefficients.
We first note that the zeros of $r(z) = n(z)/d(z)$ are those of 
$$
n(z)=\sum_{k=0}^n  \frac{ f_k^{}\beta_k}{z- t_k^{}},
$$
and the poles are (generically) the zeros of
$$
d(z)=\sum_{k=0}^n  \frac{ \beta_k}{z- t_k^{}}.
$$
Both of these expressions are in partial fractions form,
and their zeros can be found
by solving the $(n+2)\times (n+2)$ generalised eigenvalue problem 
\begin{equation}\label{ch2-gep}
\def\crr{\cr\noalign{\vskip 2pt}}
\begin{pmatrix}
0& \beta_0^{} & \beta_1^{} & \cdots & \beta_n^{} \crr
1 & t_0^{} \crr 1 & & t_1^{} \crr \vdots & & & \ddots \crr
1 & & & & t_n^{}\end{pmatrix}x 
= \lambda 
\begin{pmatrix}
0\crr &1 \crr & & 1 \crr & & & \ddots \crr & & & & 1   
\end{pmatrix} x
\end{equation}
(\cite{corless},\cite[chap.~2]{lawrence}).
Indeed, by setting
\begin{equation}
x=\left(1,\frac{1}{\lambda-t_0},\frac{1}{\lambda-t_1},\ldots,\frac{1}{\lambda-t_n}\right)^T,
\end{equation}
one can verify that (\ref{ch2-gep}) holds for any zero
$\lambda$ of $d$. 
The generalised eigenvalue problem~(\ref{ch2-gep}) has at least
two eigenvalues at infinity, and the remaining $n$ are the zeros
of $d$.  (An alternative formulation without the eigenvalues at
infinity has been proposed in \cite[sec.\ 3]{deckers}.)  The zeros
of $n(z)$ are found analogously, with $\beta_k^{}$ replaced by
$f_k^{}\beta_k^{}$. Computing these eigenvalues using the standard
QZ algorithm~\ccite{molerstewart} requires $O(n^3)$ operations,
and is usually a negligible part of the AAA computation.

Concerning the accuracy of these computations of poles and zeros, one
can say that they are backward stable in the usual sense of numerical
linear algebra, but we do not know what theoretical guarantees there
may be concerning their accuracy as poles and zeros of the rational
approximation $r(z)$.\footnote{Nian Shao (personal communication)
has pointed out that the backward stability depends on the poles
being not too close to the sample points.}  In no applications that
we are aware of has accuracy of computation of AAA poles or zeros
appeared to be a concern.

The locations of the AAA poles are interesting in almost every
application, as this paper will illustrate, and the ability to
calculate them easily means that with any approximation, one can
perform an {\em a posteriori\/} check to make sure there are no
``bad poles'' in unwanted regions of the complex plane---e.g.,
no poles in $[a,b\kern .5pt]$ if the aim is approximation on that
interval.  If there are bad poles, the standard fix is to replace
the AAA approximation~$r$ in its barycentric representation by a new
rational function $r$ in partial fractions form, whose poles are the
acceptable ones from AAA and whose coefficients are determined by
solving a matrix least-squares problem.  This may worsen the error,
but often it actually improves it, since this least-squares problem
involves the errors of the approximation one ultimately cares
about rather than the errors of the linearised barycentric form
(\ref{aaa-AAA2}).  See \ccite{aaals} and \ccite{continuum}.

In ordinary AAA operations, the computation of poles is carried
out in this {\em a posteriori\/} fashion, at the end of the
approximation process.  However, in the continuum AAA algorithm
mentioned above \ccite{continuum}, poles are calculated at every
step of the iteration.

Now we look at applications.  To begin, let us return to the
example of Figure~\ref{fig3.1}, in which the function $f(z) =
\tan(z)/\tan(2)$ was approximated from its values at 50 scattered
points in the $z$-plane.  Instead of {\tt r = aaa(F,Z)}, we can
execute

\begin{figure}
\vskip 10pt
\begin{center}
\includegraphics[trim=25 190 25 35, clip,scale=.9]{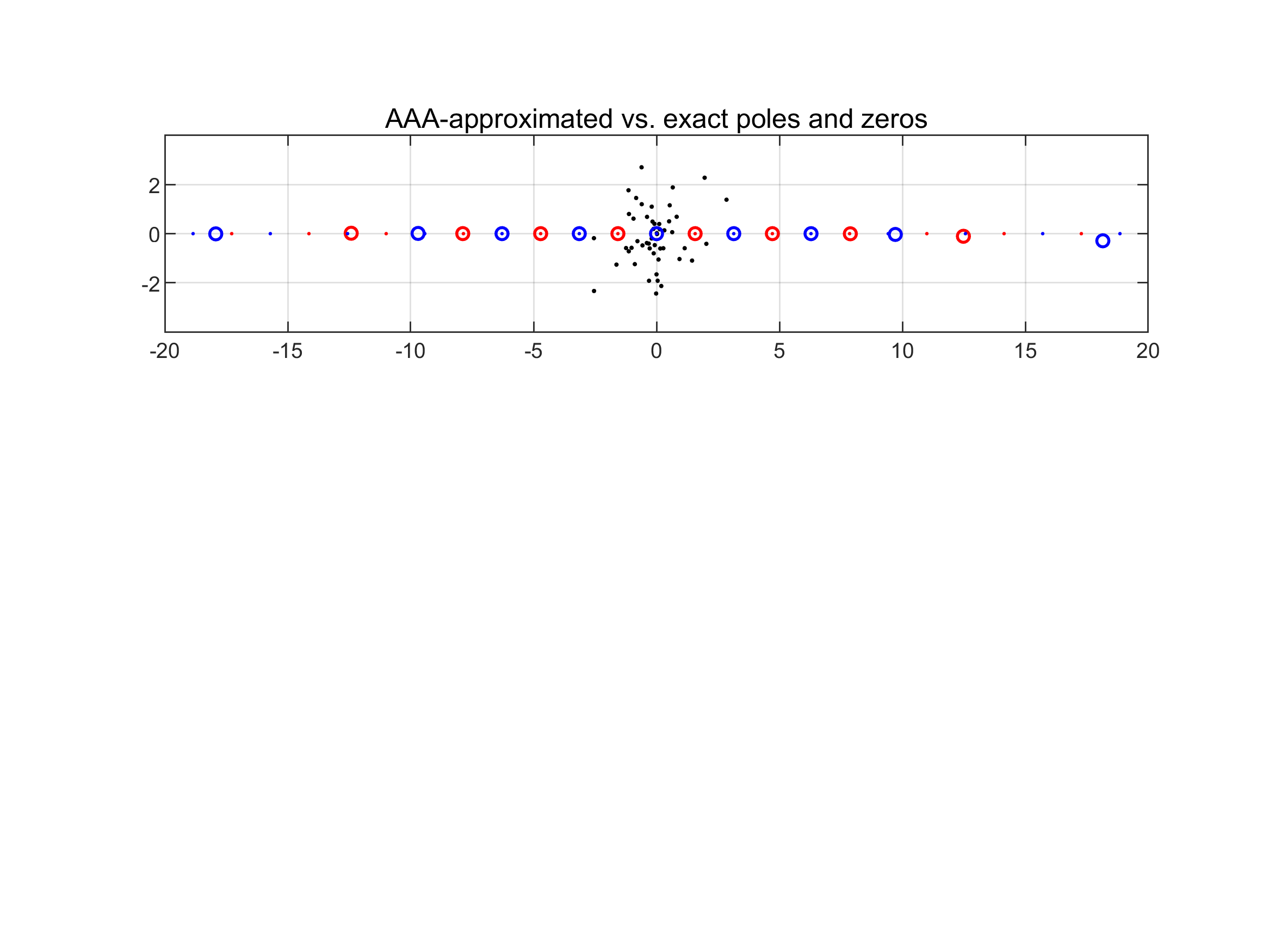}
\vskip -15pt
\end{center}
\caption{\label{fig4.1}Repetition of the experiment of
Figure~\ref{fig3.1}, now with the poles and zeros of $r(z)$ shown
as red and blue circles, respectively.  The inner poles and zeros
closely match those of the target function $f(z) = \tan(z)/\tan(2)$
(faint red and blue dots).}
\end{figure}

{\small
\begin{verbatim}
      [r,pol,res,zer] = aaa(F,Z)
\end{verbatim}
\par}

\noindent to obtain additionally the poles, residues and zeros of
the approximation $r(z)\approx f(z)$.  Again the computation takes
a few milliseconds on a laptop.  Figure~\ref{fig4.1} plots the
sample points again, as in Figure~\ref{fig3.1}, together with the
poles (red circles) and zeros (blue circles) of $r$ and the poles
and zeros of $f$ (faint red and blue dots).  We see immediately that
the poles and zeros of $r$ are approximations to those of $\tan(z)$,
which are the odd and even multiples of $\pi/2$, respectively.
This rational function $r$ is of degree $10$, and in the figure
there are 8 poles and 9 zeros.  One more pair of poles lies at
approximately $z= \pm 35$, and as for the tenth zero, it is at
$\approx -800 - 1000\kern .5pt i$.  We can think of this
as an approximation to $\infty$, so that 
$r$ is in some sense approximately a rational function of type
$(9,10)$ with a numerator of degree $9$ and a denominator of degree
$10$---a consequence of $f$ being an odd function.%

To see how closely the poles and zeros of $r$ match those of $f$, it
is convenient to divide these numbers by $\pi/2$ so that the targets
for comparison are integers.  Here is what we find for the poles,
divided by $\pi/2$.  (All these numbers would be purely real if
computed with a real-symmetric set of sample points by a version of
AAA maintaining real symmetry as discussed in section~\ref{sec:AAA}.)

{\footnotesize
\begin{verbatim}
                1.00000000000033 + 0.00000000000003i
               -1.00000000000002 - 0.00000000000027i
                3.00000332651624 - 0.00000243518063i
               -3.00000327305622 + 0.00000129315010i
                5.01823394381568 - 0.00346572096760i
               -5.01704794958420 + 0.00126227829889i
                7.94615053339014 - 0.06783633480729i
               -7.90540076322207 + 0.01056856361015i
               22.84117829937205 - 0.68627425575586i
              -22.13232732826279 - 0.22693756275246i
\end{verbatim}
\par}

\noindent Here are the zeros, divided by $\pi/2$:

{\footnotesize
\begin{verbatim}
              -0.00000000000000  +   0.00000000000000i
               2.00000000074146  -   0.00000000265986i
              -2.00000000080136  +   0.00000000129131i
               4.00058578212117  -   0.00020264321204i
              -4.00054629016288  +   0.00009068387390i
               6.18331670334207  -   0.02050187149938i
              -6.17371491450402  +   0.00563688314931i
              11.55660055153871  -   0.18489790406735i
             -11.41334018853083  -   0.00529302069495i
            -515.84554454068018  - 630.11666599322939i
\end{verbatim}
\par}

\noindent Evidently the poles and zeros within the cloud of
data are correct to about 13 digits, which is the tolerance of
the AAA approximation.  Moving outside the data set, the poles of
$f(z)/(\pi/2)$ at $\pm 3$, $\pm 5$ and $\pm 7$ are captured to about
6, 2 and 1 digits, respectively, and the zeros of $f(z)/(\pi/2)$
at $\pm 2$, $\pm 4$ and $\pm 6$ are captured to about 9, 4 and 2
digits.  This is typical of rational approximations: highly accurate
poles and zeros at distances surprisingly far from the data set.
One may ask, why do rational approximations locate poles and
zeros so well?  A rigorous answer is not fully established, though
most experts would connect the matter with Hermite integrals and
potential theory, to be discussed in section~\ref{sec:potential}.
What's clear is that poles and zeros are merely by-products of
the approximation process---surprisingly accurate by-products.
Similar effects are the basis of numerical methods for computing
matrix eigenvalues.  The Arnoldi and Lanczos iterations, for example,
minimise a certain vector $2$-norm at each iterative step, and when
the roots of the polynomials involved approximate eigenvalues of
a matrix, that is because placing the roots like that is a good
strategy for the minimisation problem~\ccite{greenbaum,kuijlaars}.

When poles are captured accurately, the same will normally be true of
residues.\footnote{The method for computing residues in Chebfun {\tt
aaa.m} changed in 2025.  Previously, residues were computed from
(\ref{aaa-baryrep}) by an algebraic formula.  In the new method,
they are obtained by solving a linear least-squares problem to fit
the sample values as well as possible by a linear combination of
the poles.  For applications like those of section~\ref{sec:quad},
this sometimes contributes several digits of improved accuracy.}
Here for example are the residues of the ten poles of $r(z)$ listed
above, all multiplied by $\tan(2)$ so that we recognise the target
value of $-1$ corresponding to the residues of $\tan(z)$:

{\footnotesize
\begin{verbatim}
                -1.00000000000249 + 0.00000000000049i
                -0.99999999999953 - 0.00000000000229i
                -1.00002115229468 + 0.00001322035816i
                -1.00002043469684 + 0.00000683956587i
                -1.05083731528699 + 0.00770105605601i
                -1.04773789789134 + 0.00251606979538i
                -2.35467344729027 + 0.07037511036527i
                -2.30241052490136 + 0.00113182260999i
               -21.40951807722749 - 0.67668458618553i
               -22.91889560491100 + 1.25209694878431i
\end{verbatim}
\par}

\noindent
These match the target value to about 12, 5 and 2 digits for the
first, second and third pairs out from the origin.

Along with the approximation of $\tan(z)/\tan(2)$ at scattered sample
points, the last section considered the approximation of $\Gamma(z)$
at 50 points of the interval $[-1.5,1.5\kern .5pt ]$ or the circle
$|z|=1.5$ (Figure \ref{fig3.2}).  The gamma function has poles at
$0,-1,-2,-3,-4,\dots,$ and the poles of the rational approximations
computed there match these to about $15, 15, 7, 3, 1,\dots$ digits
(for approximation on the interval) and $15,14,11,5,3,3,\dots$
digits (for approximation on the circle).

For another illustration of locating poles and zeros by rational
approximation, we turn to the Riemann zeta function $\zeta(z)$, an
example presented with a colourful figure in \ccite{jjiam}.  The function is defined
by the Dirichlet series
\begin{equation}
\zeta(z) = \sum_{k=1}^\infty k^{-z},
\end{equation}
which converges just for $\Re z > 1$ even though $\zeta(z)$ is analytic
throughout the complex plane apart from a simple pole at $z=1$.
With $\Re z = 4$, we can evaluate $\zeta(z)$ to 13-digit accuracy
quickly by summing $20{,}000$ terms of the series.  Suppose we evaluate
$\zeta(z)$ in this way at 100 points along the
line segment from $4-40\kern .7pt i$ to $4+40\kern .7pt i$
and then calculate a AAA approximant $r(z)$ to these data values:
\begin{figure}
\vskip 10pt
\begin{center}
\includegraphics[trim=60 15 63 8, clip,scale=.55]{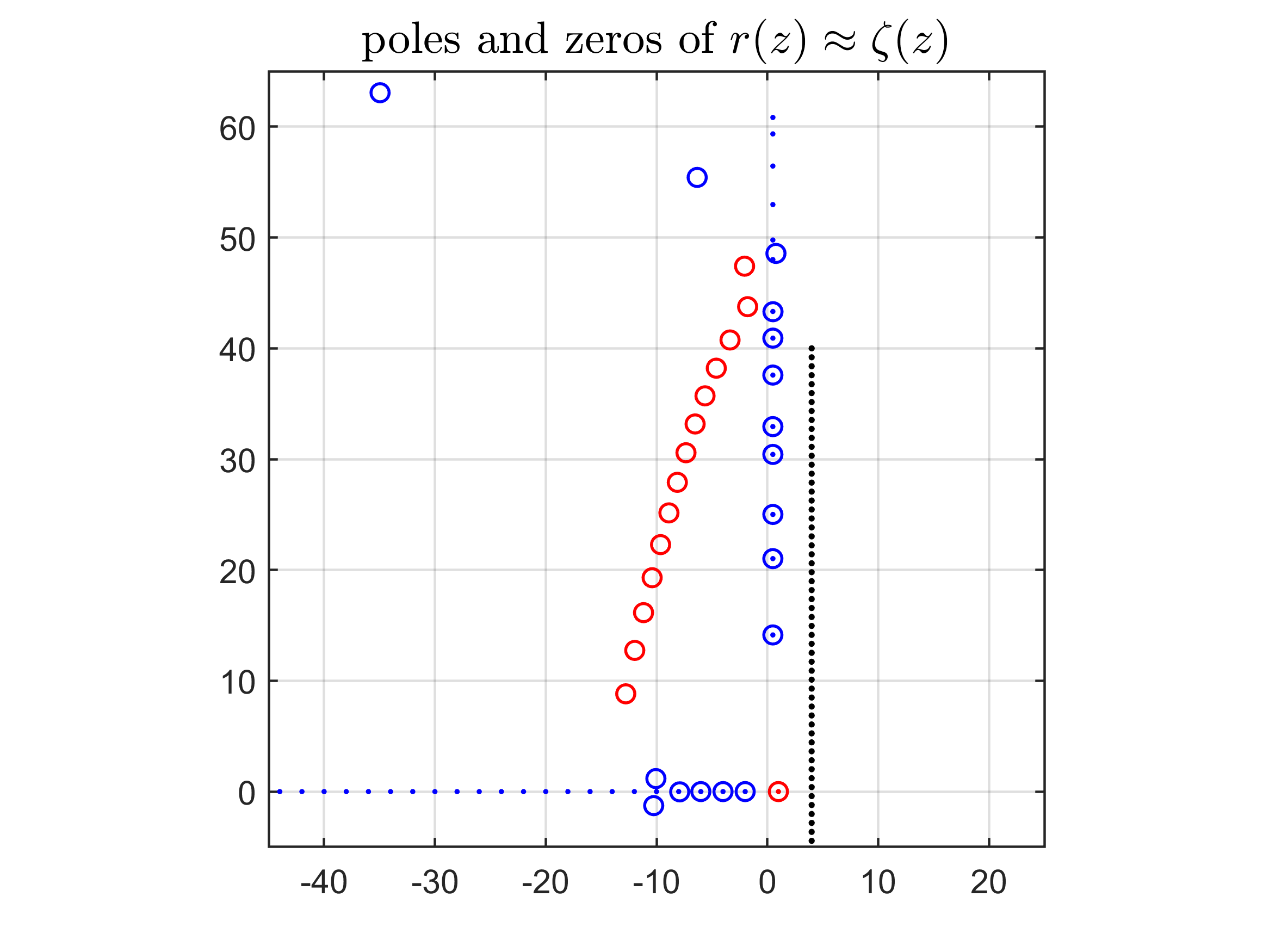}~~~~
\includegraphics[trim=63 15 60 8, clip,scale=.55]{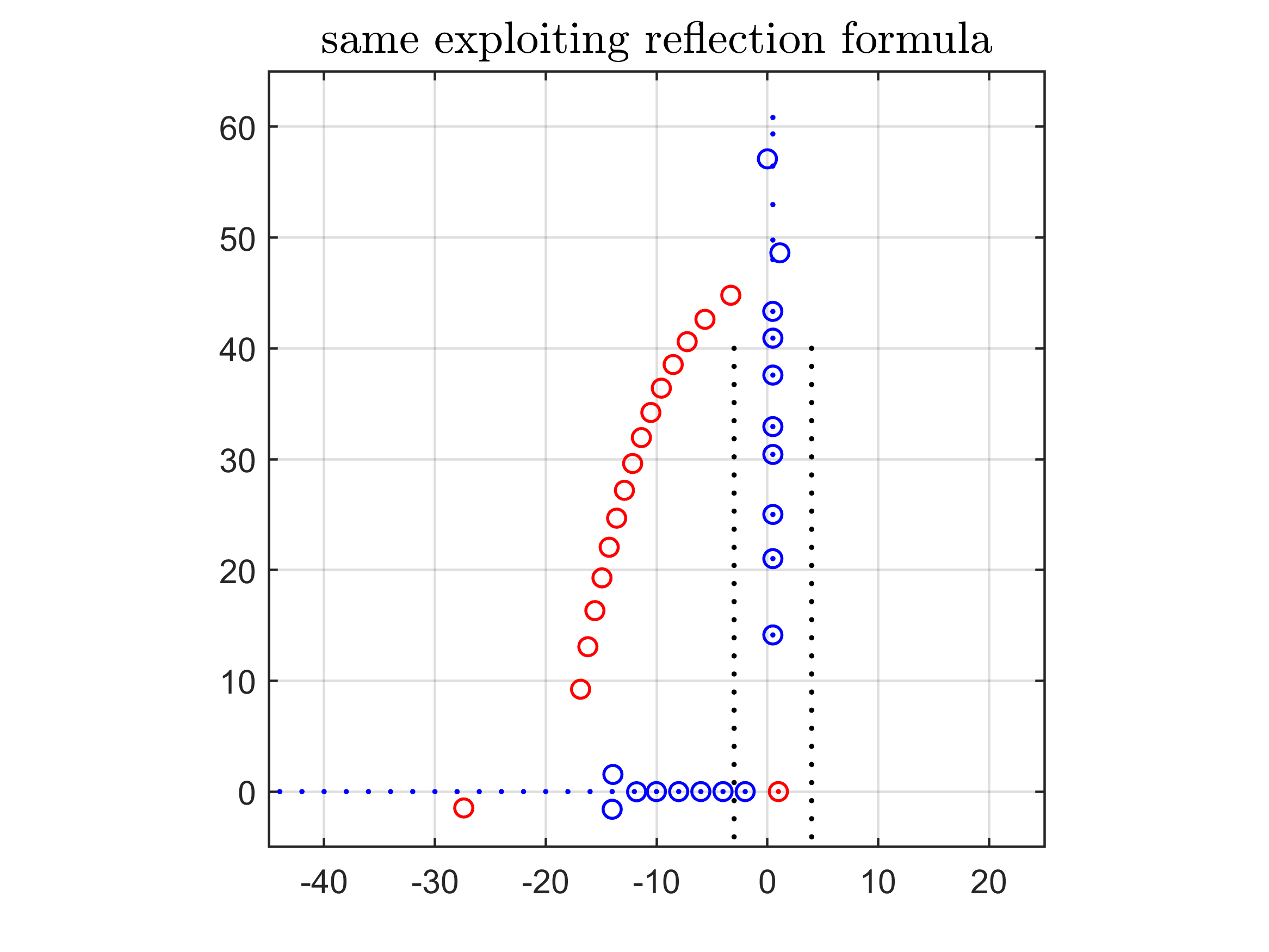}
\end{center}
\caption{\label{fig4.2}On the left, poles (red circles) and
zeros (blue circles) of the degree $29$ AAA
approximation of the Riemann zeta function
$\zeta(z)$ based on sample values at 100 points with 
$\Re z = 4$ (black dots).  The pole and the zeros of
$\zeta(z)$ are marked by red and blue dots, respectively.
The pole of $\zeta(z)$ is captured to 11 digits
of accuracy, whereas the other poles of $r(z)$ are artifacts of the approximation
(red circles).  The first few zeros of $\zeta(z)$ on the critical line
$\Re z = 1/2$ are also highly accurate.  This figure shows the upper
half of a symmetric configuration in the complex plane.
On the right, the same except with $50$ sample points in the right
half-plane and also their reflections by (\ref{poles-reflec})
in the left half-plane.
The rational approximation is now of degree $32$; the accuracy is similar.}
\end{figure}

{\small
\begin{verbatim}
      Z = linspace(4-40i,4+40i).';
      zeta = @(z) sum((2e4:-1:1).^(-z),2);
      [r,pol,res,zer] = aaa(zeta(Z),Z);
\end{verbatim}
\par}

\noindent The left side of Figure~\ref{fig4.2} shows the resulting poles
(red circles) and zeros (blue circles) of $r(z)$.  For comparison,
the pole and the zeros of $\zeta(z)$ are shown by red
and blue dots.  The agreements are very good.  On the right side
of the figure, the experiment is repeated but now with 50 points with
$\Re z = 4$ and $50$ more points at their reflections in the left
half-plane with $\Re z = -3$, where
$\zeta(z)$ is evaluated by the reflection formula
\begin{equation}
\zeta(1-z) = 2 (2\pi)^{-z} \cos(\pi z/2) \,\Gamma(z) \,\zeta(z).
\label{poles-reflec}
\end{equation}
In both experiments the pole matches the correct value $z=1$ to $11$
digits, and as for the zeros, here are a pair of columns showing
the absolute errors in the zeros in the upper half-plane.
The reason for the blank in the final row of the first column of
data is that this approximation is of degree 1 less than the other.

{\footnotesize
\begin{verbatim}
                               error                     error
    imag part of zero  100 samples on one side   50 samples on each side
         14.135            0.000000000080            0.000000000011
         21.022            0.000000000316            0.000000000027
         25.011            0.00000000082             0.000000000088
         30.425            0.0000000074              0.00000000069
         32.935            0.000000034               0.0000000029
         37.586            0.0000030                 0.00000040
         40.919            0.00089                   0.00020
         43.327            0.026                     0.012
         48.005            0.62                      0.88
         49.774            8.9                       7.3
         52.970           37.                       28.
         56.446                                    164.
\end{verbatim}
\par}

The examples of this section so far have dealt with meromorphic
functions, which have poles as well as zeros (actually, the gamma
function has no zeros).  Since rational functions can approximate
poles, this is an obvious type of problem to attempt. But the AAA
algorithm can be a powerful and exceptionally simple rootfinding
tool even for functions that are analytic, when there are no poles.
Indeed, AAA may often outperform more complicated bespoke rootfinding
methods.  To illustrate how simple its use can be, suppose we want
to find the roots of the function $f(x) = \sin(5\pi x)$ in $[-1,1]$.
The commands

{\small
\begin{verbatim}
      X = chebpts(40);
      F = sin(5*pi*X);
      [~,~,~,zer] = aaa(F,X);
      zeros = sort(zer(imag(zer)==0 & abs(zer)<=1.1))
\end{verbatim}
\par}

\noindent based on just 40 function values produce
this output in a few milliseconds, with all the zeros captured
to accuracy better than $10^{-13}$:

{\footnotesize
\begin{verbatim}
        zeros =
          -1.000000000000001
          -0.800000000000001
          -0.599999999999981
          -0.399999999999969
          -0.200000000000000
           0.000000000000039
           0.200000000000001
           0.399999999999939
           0.599999999999971
           0.800000000000001
           1.000000000000002
\end{verbatim}
\par}

\noindent A general term for this kind of algorithm, introduced by
Boyd \ccite{boydbook}, is {\em proxy rootfinding\/}.  The idea
is to approximate a general function $f$ by a function with
simple enough algebraic structure that its roots can be determined
easily---typically via a matrix eigenvalue problem.  The best-known
version of this idea involves polynomial approximation, as
suggested first by Good in his paper on ``colleague matrices''
\ccite{good}, then developed further by Boyd and used as the
main rootfinding engine for polynomial approximations in Chebfun
\ccite{battles,chebfun}.  Examples like this one illustrate that
rational functions can do the job too, and with the AAA algorithm,
they can be applied to virtually arbitrary rootfinding problems
on arbitrary domains.  This method is advocated by Costa in one
of the online Chebfun examples \ccite{costaexample} and developed
more fully in \ccite{bowhay}.  For rational rootfinding before AAA,
see \cite[section 6]{akt}.

In this section we have considered how AAA may be used to find
both poles and zeros.  It is natural to ask, are these two
problems symmetric, or might AAA perhaps be more accurate, say,
as a polefinder than as a zerofinder?  If this were true, then
there might be an advantage sometimes in replacing data values $F$
by their reciprocals $1/F$.

So far as we are aware, the answer to this question is perhaps
surprising, and surprisingly elementary.  It seems that polefinding
and zerofinding are at bottom essentially equivalent, but that
polefinding is sometimes slightly more accurate for the accidental
reason that a function tends to be bigger near its poles than near
its zeros!  Since AAA is designed to approximately minimise the
$\infty$-norm of the absolute error, this means that in a relative
sense, it sometimes gives more weight to approximation near a pole
than near a zero.  Perhaps the balance would be more even if AAA
measured the $\infty$-norm of errors via a metric on the Riemann
sphere rather than the complex plane, but so far as we know, no
AAA codes have tried this.

AAA approximation of poles is the foundation of several of the
topics discussed later in this article, notably in sections
\ref{sec:solution}, \ref{sec:res}--\ref{sec:nonlin} and \ref{sec:quad}.
See the
final paragraph of section~\ref{sec:potential} for an informal
classification of computed poles into {\em approximation poles\/},
{\em pole poles\/}, {\em branch cut poles\/} and {\em spurious
poles\/}.

\section{\label{sec:diff}Computing derivatives and integrals}

Differentiation and integration are two of the basic problems of
numerical analysis, and rational approximations sometimes lead to
good methods for them.

Let us begin with the problem of calculating the derivative of
\begin{equation}
f(x) = \tanh(8\kern .5pt x)
\label{diff-tanh}
\end{equation}
for $x\in [-1,1]$.  One approach would be by polynomial
approximation.  We can compute a polynomial approximant to $f$
accurate to 10 digits and differentiate it with the Chebfun commands

{\small
\begin{verbatim}
      p = chebfun('tanh(8*x)','eps',1e-10);
      dp = diff(p);
\end{verbatim}
\par}

\begin{figure}
\begin{center}
\includegraphics[trim=50 20 50 0, clip,scale=1]{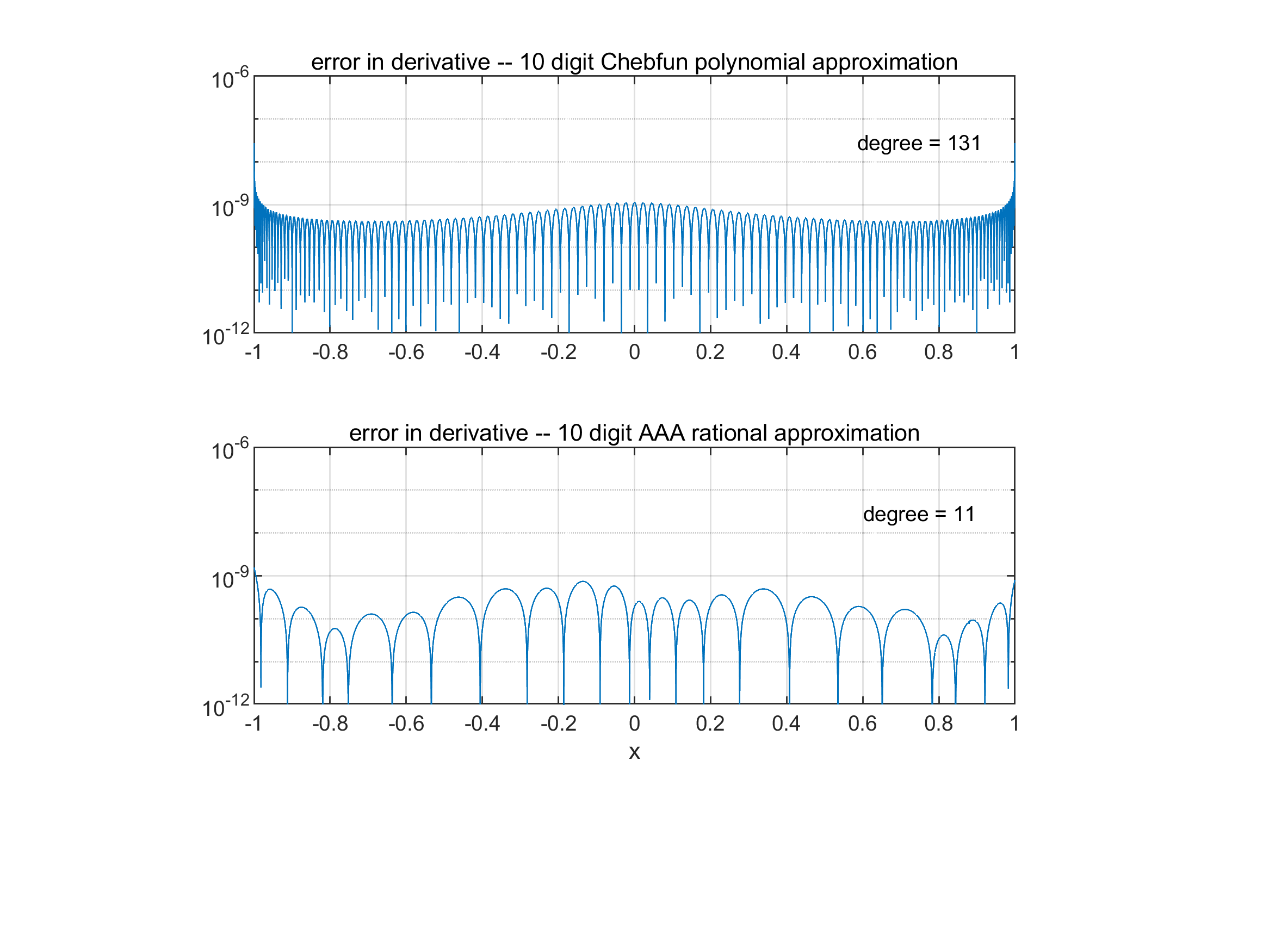}
\end{center}
\vspace{-44pt}
\caption{\label{diff-figdtanh}Numerical computation of the derivative
of $f(x) = \tanh(8\kern .4pt x)$ for $x\in [-1,1]$ based on 
polynomial and rational approximations of $f$ accurate to $10$ digits.
For $f(x) = \tanh(80\kern .4pt x)$, the degrees increase to
1201 and 22.}
\end{figure}

\noindent
Executing {\tt length(p)} reveals that the approximation has
degree 131.\ \  (This is a polynomial interpolant in an adaptively
determined number of Chebyshev points.)  As shown in the upper
half of Figure \ref{diff-figdtanh}, the accuracy of the derivative
is about $10^{-9}$ (hence a relative error of order $10^{-10}$).
The plot shows the absolute error $|f(x) - p(x)|$ on a log scale,
evaluated on a 10{,}000-point equispaced grid in $[-1,1]$.

Another approach to differentiating (\ref{diff-tanh}) would be by
rational approximation.  We can compute a AAA approximation to 10
digit accuracy and differentiate it with the commands

{\small
\begin{verbatim}
      X = linspace(-1,1,200)';
      F = tanh(8*X);
      [rr,pol] = aaa(F,X,'tol',1e-10,'deriv_deg',1);
      r = rr{1}; dr = rr{2};
\end{verbatim}
\par}

\noindent
Executing {\tt length(pol)} reveals that the rational approximation
has degree~11.  The lower half of Figure \ref{diff-figdtanh}
shows that the accuracy is about the same as for the polynomial
approximation, though the degree is ten times lower.  For other
functions, the ratio can be more extreme.

Computing derivatives is a step toward computing minima and maxima.
For example, where does the function $f$ of (\ref{diff-tanh})
attain its maximum slope, and what is that value?  The computation

{\small
\begin{verbatim}
      [rr,pol] = aaa(F,X,'tol',1e-10,'deriv_deg',2);
      ddr = rr{3};
      [r,pol,res,zer] = aaa(ddr(X),X);
      pos = zer(imag(zer)==0 & abs(zer)<1)
      val = dr(pos)

      pos = -9.87e-12
      val = 7.99999999981
\end{verbatim}
\par}

\noindent gives answers accurate to 10 or more digits.

These examples illustrate results of computing a derivative by
rational approximation, but not the mechanics of how it is done.  
Formulas for this
were introduced in \ccite{schneider}, where it was shown that if\/ $r(z)$ is
represented in the barycentric representation (\ref{aaa-baryrep}), 
\begin{equation}
r(z) = \sum_{k=0}^n {\fk\kern .5pt \bk\over z-\tk} \left/
\sum_{k=0}^n {\bk\over z-\tk} \right.,
\label{diff-baryrep}
\end{equation}
then the derivative $r'(z)$ can be written in the same way, except with
the data $\fk$ replaced by the divided differences
$(r(z)-\fk)/(z-\tk)$:
\begin{equation}
r'(z) = \sum_{k=0}^n {(r(z) - \fk)\kern 1pt \bk\over (z-\tk)^2}
\kern -2pt \left/
\sum_{k=0}^n {\bk\over z-\tk} \right..
\label{diff-baryrep1}
\end{equation}
As usual, these formulas apply when $z$ is not one of the
support points; of course $z$ should also not be one of the
poles of $r$.  If $z=\tj$ for some $j$, we have the special formulas
\begin{equation}
r(\tj) = \fj,\quad
r'(\tj) = - {1\over \bj}\sum_{k\ne j} {\bk(\kern .7pt\fj-\fk)\over \tj-\tk}.
\label{diff-baryrepsupp}
\end{equation}
Equations (\ref{diff-baryrep1})--(\ref{diff-baryrepsupp})
cover the computation of the first derivative, and for higher derivatives, 
Schneider and Werner present a generalised formula
involving higher-order divided differences and show how it
can be evaluated by a simple recursion.

Barycentric evaluations of polynomials and rational functions
are famously stable even at points $z$ close to the support
points $\{\tk\}$, where one comes close to ``dividing infinity by
infinity'', and theorems to this effect for the polynomial case were
presented by \ccite{higham}.  Barycentric evaluations of first and
higher-order derivatives, however, are another matter.  The formula
(\ref{diff-baryrep1}) fails numerically as $z$ approaches a support
point, and for accurate results, other methods are needed such as
a second AAA fit to data at support and/or other points defining
the derivative rather than the function value itself.  At present,
such refinements are not included in Chebfun {\tt aaa.m}, nor so
far as we know in other AAA software.

Rational approximations can also be used for numerical integration,
also known as quadrature: if $f(x) \approx r(x)$, then $\int \kern
-1pt f(x) \kern .7pt dx \approx \int \kern -1pt r(x) \kern .7pt dx$.
We are unaware of formulas that evaluate $\int \kern -1pt r(x)
\kern .7pt dx$ directly in the barycentric representation, so as an
alternative, one can convert a rational approximation to partial
fractions and then integrate these term by term.  A variation on
this theme has been proposed in \ccite{brunosant} to cope with
problems in which one or more singularities lie very close to the
integration contour.  After finding these singularities and their
residues with AAA, one can subtract them off and integrate the
smoother integrand that remains by AAA or other methods.

An interesting special case of AAA-based integration
has been suggested
to us by Dave Darrow, in which one wishes to
calculate an integral over the whole real line,
\begin{equation}
I = \int_{-\infty}^\infty \kern -1pt f(x) \kern 1pt dx.
\label{diff-real}
\end{equation}
Now there is the prospect of evaluating $I$ entirely by residue calculus
applied to an approximation $r\approx f$.
Specifically, if\/ $r$ satisfies
$r(z) = O(|z|^{-2})$ as $|z|\to\infty$ in the upper half-plane, 
then the contribution to a contour integral over a semicircle in
the upper half-plane with
$|z| = R$ will vanish as $R\to\infty$, and by the residue theorem we are left with 
\begin{equation}
I_r^{} = \int_{-\infty}^\infty r(x) \kern 1pt dx
= 2\pi\kern .7pt i \sum_{k=1}^n a_k^{},
\label{diff-rescalc}
\end{equation}
where $\{ a_k^{} \}$ are the residues of\/ $r$ in the upper
half-plane.
(In practice, instead of 
$r(z) = O(|z|^{-2})$ we are likely to have 
$|r(z)|\le \varepsilon/|z|$ as $|z|\to\infty$ with $\varepsilon\ll 1$, 
perturbing the result by at most $2\pi \varepsilon$.)
As an example, suppose we want to evaluate numerically
$$
{1\over \sqrt \pi\kern 1pt} \int_{-\infty}^\infty e^{-(x-1)^2} dx,
$$
whose exact value is $1$.
We use the tangent function to discretise
$\real$ by a grid with points extending to large values and approximate
the integrand there by a rational function $r$, which has
degree $30$.

{\small
\begin{verbatim}
      X = tan(.99*linspace(-pi/2,pi/2,200))'; 
      F = exp(-(X-1).^2)/sqrt(pi);
      [r,pol,res] = aaa(F,X);
\end{verbatim}
\par}

\noindent
Approximating the integral of\/ $r$
by residue calculus gives 12 digits of accuracy:

{\small
\begin{verbatim}
      ii = find(imag(pol)>0);
      I = -2*pi*imag(sum(res(ii)))
        
      I = 1.000000000003281
\end{verbatim}
\par}

\noindent
Applying the same method to
$$
\int_{-\infty}^\infty e^{-(x-1)^2} \sqrt{0.001 + x^2}\kern 1pt  dx
$$
gives a rational approximant of degree $52$ and an estimated integral
matching the exact value $1.863662433495\dots$ to 11 digits.

Another application that can be approached by AAA
approximation is the computation of complex Cauchy integrals, and a successful
example of this kind can be found in \ccite{kisil3}.  On the other hand,
sometimes there may be better approaches.  Cauchy integrals typically arise in
applications related to residue calculus, and it may be easier
and more accurate
to proceed directly from a rational approximation rather than by means
of an integral.  For example, suppose we want to compute the zeros and
poles of the function
\begin{equation}
f(z) = {\sin(z)\over z^2-{1\over 4}}
\end{equation}
in the unit disk.  There are methods for such computations based on Cauchy
integrals \ccite{akt}, but it is simpler to follow the methods of the
last section and execute

{\small
\begin{verbatim}
      Z = exp(2i*pi*(1:50)/50);
      F = sin(Z)./(Z.^2-1/4);
      [r,pol,res,zer] = aaa(F,Z);
      pol = pol(abs(pol)<1)
      zer = zer(abs(zer)<1)
        
      pol =
        -0.499999999999999 - 0.000000000000000i
         0.500000000000000 - 0.000000000000000i
      zer =
        -3.83e-16 + 1.69e-16i
\end{verbatim}
\par}

The relationship between Cauchy integrals and rational approximations
will reappear in sections \ref{sec:res}--\ref{sec:nonlin} and
\ref{sec:quad}--\ref{sec:cauchy}, and section \ref{sec:quad}
will present some surprisingly deep connections between rational
approximations (of a Cauchy transform) and quadrature formulas
(over an associated contour).

\section{\label{sec:branch}Analysis of branch points}
It is perhaps not surprising that rational approximations are good
at finding poles and zeros, as we saw in section \ref{sec:poles}.
A deeper challenge arises if a function has branch points.
Rational functions can be powerful here too, and their role for
such computations became well known in condensed matter physics
in the use of Pad\'e approximations for the numerical analysis of
critical points related to phase transitions \ccite{baker,guttmann}.

For an example of AAA rather than Pad\'e type, suppose $f$ is the function
\begin{equation}
f(z) = (2-z)^{-0.4}.
\label{branch-branchfun}
\end{equation}
This function has a branch point at $z_0=2$, which means that not
only is it not analytic at $z_0$, but it is not analytic in any
punctured neighborhood of $z_0$.  Now suppose we know some function
values for $f$ well away from the singularity, say, at 50 points
in $[-1,1]$:

{\small
\begin{verbatim}
      f = @(z) (2-z).^(-0.4);
      Z = linspace(-1,1,50)';
      F = f(Z);
\end{verbatim}
\par}

\noindent  We can use these values to approximate $f(z)$ by AAA,
and the approximation will be excellent for many values of $z$.  For
example, at $z=-4$, which is well beyond the radius of convergence
of $2$ of the Taylor series, we get six digits of accuracy:

{\small
\begin{verbatim}
      [r,pol,res,zer] = aaa(F,Z);
      f(-4), r(-4)

      ans = 0.488359342
      ans = 0.488359411
\end{verbatim}
\par}

Suppose, on the other hand, our goal is not just to evaluate $f$ but
to determine the location and the nature of its singularity.  We find
that there are six zeros and poles of $r$, lining up in interlacing
order on the real axis to the right of $z_0$, clustering at that
point.  The configuration is plotted in Figure~\ref{sec:branch}.1,
and here are the numbers.  In this example the zeros and poles
happen to be real, but in general there may also be complex conjugate
pairs even in approximation of a real function.

{\footnotesize
\begin{verbatim}
                     zeros         poles        residues
                   2.0929234     2.0308549    -0.1350795
                   2.4436865     2.2630564    -0.1938906
                   3.2970933     2.8646639    -0.2887175
                   5.5180478     4.3453159    -0.5017359
                  13.5608306     8.8550598    -1.1596502
                 100.5398233    34.7222948    -5.3256391
\end{verbatim}
\par}

It is clear from the figure and the data that there is some kind
of singularity near $z=2$.  It must be a branch point, since it
is not approximated simply by a pole.  But how could we determine
its location to more than, say, a couple of digits of accuracy?
And how could we determine the nature of the singularity?

\begin{figure}
\vskip 5pt
\begin{center}
\includegraphics[trim=30 235 0 0, clip,scale=1]{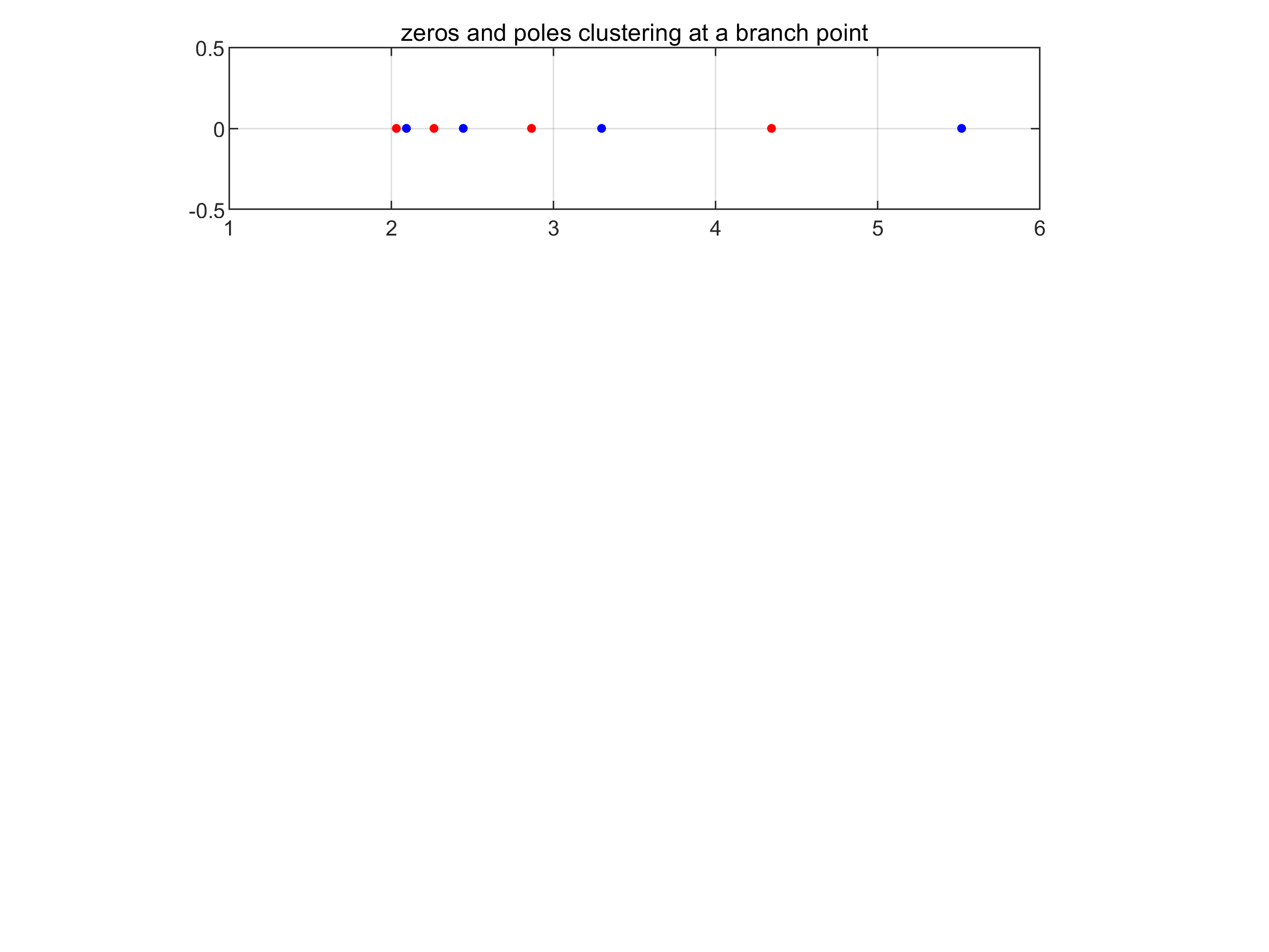}~~~~
\vskip -10pt
\end{center}
\caption{\label{fig5.1}Zeros (blue) and poles (red) of a AAA
approximation to the function (\ref{branch-branchfun}) with a branch point at
$z_0=2$.  The zeros and poles cluster near the singularity.  By rational
approximation of the logarithmic derivative (\ref{branch-logdiff}), one
can determine the singularity location and exponent to $10$ digits
of accuracy.}
\end{figure}

A powerful method for some such problems is to work with the
{\em logarithmic derivative\/} of $f$,
\begin{equation}
{f'(z)\over f(z)} = (\log f)'.
\label{branch-logdiff}
\end{equation}
Suppose $f$ is a function of the form
\begin{equation}
f(z) = g(z) (z-z_0)^\alpha,
\label{branch-alphaform}
\end{equation}
where $z_0$ is an unknown complex number, $g$ is an unknown
analytic function, and $\alpha$ is an unknown real number that is not
an integer.  Then we have
\begin{equation}
f'(z) = g'(z) (z-z_0)^\alpha + \alpha g(z) (z-z_0)^{\alpha -1},
\end{equation}
which implies
\begin{equation}
{f'(z)\over f(z)} = {g'(z)\over g(z)} + \alpha (z-z_0)^{-1}.
\end{equation}
The logarithmic derivative has converted the branch point to a simple pole.
If we now approximate $f'/f$ by a AAA rational function $r$, we can expect to
find a pole of $r$ very close to $z_0$ with a residue that is very close
to $\alpha$.

To apply this method we need the derivative $f'$.  Working from
function values but not derivatives, we can estimate
$f'$ from a preliminary rational approximation as in the last section.
Thus the natural AAA method for this
problem consists of three steps, applied at a suitable
set of sample points $Z$.  So far as we know, this method
has not been published previously.
\medskip

~~1.  Approximate $f(z)$ for $z\in Z$ by a rational function $r_1^{}(z)$;
\smallskip

~~2.  Approximate $r_1'(z)/r_1^{}(z)$ on $z\in Z$ by another rational \rlap{function $r(z)$;}
\smallskip

~~3.  Estimate $z_0$ and $\alpha$ by a pole and a residue of $r(z)$.
\medskip

\noindent 
For step (2), we need to calculate $r_1'(z)$, where $r_1^{}$ is a
AAA rational function in the barycentric representation.  In this
computation, since one rational approximation is being built from
another, 13-digit precision is not achievable.  We proceed with
the tolerance of the second approximation weakened to 11 digits:

{\small
\begin{verbatim}
      rr = aaa(F,Z,'deriv_deg',1);
      f2 = @(z) rr{2}(z)./rr{1}(z);
      [r,pol,res] = aaa(f2(Z),Z,'tol',1e-11);
\end{verbatim}
\par}

\noindent The computed pole location and exponent come out
with about ten correct digits:

{\small
\begin{verbatim}
      pol =  2.000000000024
      res = -0.400000000012
\end{verbatim}
\par}

As pointed out in \ccite{bowhay}, this method for analysis of
branch points also works for the analysis of poles (or zeros) of
order higher than~$1$: where $\alpha$ in (\ref{branch-alphaform})
is an integer $\le -2$ (or $\ge 2$).  For example, if AAA is
applied to the function $f(z) = e^z/(z-z_0^{})^3$ with $z_0^{} =
1+0.5\kern .5pt i$, the triple pole is approximated by three poles
in a near-circle about $z_0^{}$  of radius $0.000134$, with residues
of magnitude about $5.03\cdot 10^7$.  By contrast, a repetition of
the logarithmic derivative calculation above finds the triple pole
location with accuracy $5\cdot 10^{-13}$ and its multiplicity (which
of course we know to be an integer) with accuracy $4\cdot 10^{-12}$.

The method we have discussed also works if $f$ has several branch points
related multiplicatively.  For example, suppose we have the function
\begin{equation}
f(z) = g(z) (z-z_1^{})^\alpha (z-z_2^{})^{\kern 1pt \beta}
\end{equation}
and we want to find $z_1^{}, z_2^{}, \alpha,$ and $\beta$.
The logarithmic derivative is
\begin{equation}
{f'(z)\over f(z)}  = {g'(z)\over g(z)} + {\alpha\over z-z_1^{}}
+ {\beta\over z-z_2^{}},
\end{equation}
indicating that a rational approximation can find both branch points and exponents.
To illustrate, consider the function
\begin{displaymath}
f(z) = {e^z (2\kern .6pt i-z)^{1/2}\over (1.5-z)^{3/4}}.
\end{displaymath}
If we sample in 50 roots of unity and follow the procedure above,
we find these poles and residues also correct to 10 places:

{\small
\begin{verbatim}
      pol =
         0.000000000105 + 2.000000000013i
         1.499999999999 + 0.000000000008i
        
      res =
         0.500000000017 - 0.000000000085i
        -0.749999999996 - 0.000000000024i
\end{verbatim}
\par}

A greater challenge arises if we have a function with
several branch points related additively.  For example,
the logarithmic derivative of
\begin{equation}
f(z) = (z-z_1^{})^\alpha + (z-z_2^{})^{\kern 1pt \beta}
\end{equation}
is
\begin{equation}
{f'(z)\over f(z)} = {\alpha(z-z_1^{})^{\alpha-1} + \beta (z-z_2^{})^{\kern 1pt \beta-1}\over
(z-z_1^{})^\alpha + (z-z_2^{})^{\kern 1pt \beta}}.
\end{equation}
Like $f$ itself, this function has branch points at 
$z_1^{}$ and $z_2^{}$, so looking at the logarithmic
derivative will not obviously be helpful.
Indeed, the same challenge applies for an additive combination involving
just one branch point, such as this seemingly minor modification of (\ref{branch-branchfun}),
\begin{equation}
f(z) = 1 + (2-z)^{-0.4},
\label{branch-branchfun2}
\end{equation}
with logarithmic derivative
\begin{equation}
{f'(z)\over f(z)} = {0.4(2-z)^{-1.4} \over
1 + (2-z)^{-0.4}}.
\end{equation}

\section{\label{sec:inv}Inverse functions}
If $f: z\mapsto w$ is a one-to-one function of a real or complex
variable $z$, then the inverse function $f^{-1}$ is defined by
the condition $f^{-1}(f(z)) = z$.  (As often in this paper there
is a matrix analogue, but our focus is on the scalar case.)
Some inverse functions are known analytically, but others must
be found numerically, and this can be challenging since $f^{-1}$
may have singularities even though $f$ does not.  In this section
we show how rational approximations can be used to compute inverse
functions numerically.

To start with an easy example, suppose we didn't know about $\log(t)$
and wanted to find an inverse of $f(x) = e^x$ for $x\in [-1,1]$,
which it maps into the interval $[e^{-1},e]$.  The code

{\small
\begin{verbatim}
      X = linspace(-1,1)';
      f = @(x) exp(x);
      [finv,pol] = aaa(X,f(X));
\end{verbatim}
\par}

\noindent returns a function handle {\tt finv} in about 5
ms that matches $\log(t)$ on $[e^{-1},e]$ with maximum error
$1.6\cdot 10^{-14}$ (on all of the interval, not just at the sample
points).  The rational approximation is of degree 8, with poles in
$(-\infty,0\kern .5pt ]$ clustering at $0$, where the true inverse
function has a logarithmic branch point.

As a more challenging example, again with an exact solution available
for comparison, suppose we want the inverse of $\sin(x)$ for $x\in
[-\pi/2,\pi/2]$.  Since $f(x)$ has derivative zero at $x=\pm \pi/2$,
$f^{-1}(t) = \sin^{-1}(t)$ has square root singularities at $t =
\pm 1$.  Because of ill-conditioning, we cannot expect a global
approximation to $f^{-1}(t)$ on $[-1,1]$ to have accuracy better than
about the square root of machine precision.  With this in mind we
compute a AAA approximant with tolerance $10^{-7}$, using a sample
grid of 400 points exponentially clustered near $\pi/2$ and $-\pi/2$:

{\small
\begin{verbatim}
      X = (pi/2)*tanh(linspace(-14,14,400)');
      f = @(x) sin(x);
      [finv,pol] = aaa(X,f(X),'tol',1e-7);
\end{verbatim}
\par}

\noindent This time the computation takes about 20 ms, giving a
rational function of degree 40, with 20 poles in $(-\infty,-1]$
and 20 in $[1,\infty)$.  The maximum error is $1.4\cdot 10^{-7}$.

Of course, inverses can also be computed of functions tht one does not know 
analytically.  For example, the ODE initial-value problem
\begin{equation}
x' = \exp(\cos(t^2 x)), \quad u(0) = 0, \quad t\in [\kern .5pt 0,2\kern .4pt]
\label{invode}
\end{equation}
can be solved numerically by the Chebfun commands
\begin{figure}[t]
\vskip 5pt
\begin{center}
\includegraphics[trim=0 150 0 0, clip,scale=.82]{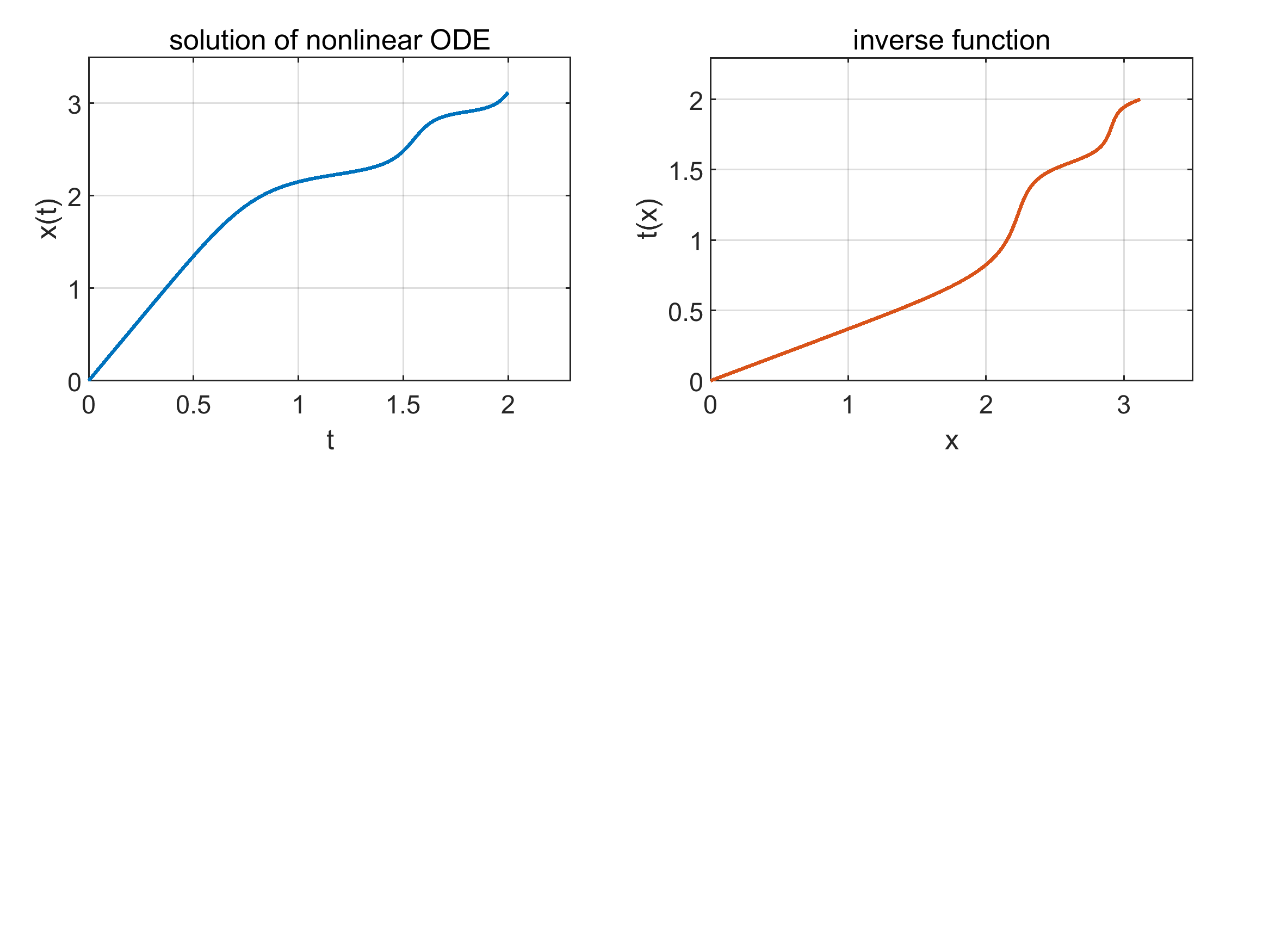}~~~~
\vspace{-20pt}
\end{center}
\caption{\label{odefig}Solution of the ODE (\ref{invode}) and its
inverse function computed with AAA.}
\end{figure}

{\small
\begin{verbatim}
      N = chebop([0 2]);
      N.op = @(t,x) diff(x)-exp(cos(t^2*x));
      N.lbc = 0;
      x = N\0;
\end{verbatim}
\par}

\noindent
Figure \ref{odefig} shows the result together with the inverse function
(approximated by a rational function of degree 34) computed by the sequence

{\small
\begin{verbatim}
      T = linspace(0,2,500)';
      [t,pol] = aaa(T,x(T),'tol',1e-10);
\end{verbatim}
\par}

Our next problem comes from a Chebfun example posted by Kuan Xu in
2012 \ccite{kuanxu}.  The planar curve shown in Figure~\ref{arcfig}
is parametrised by $x\in [\kern .5pt 0,1]$ and defined by these
Chebfun commands, using complex arithmetic for convenience:

{\small
\begin{verbatim}
      x = chebfun('x',[0 1]); 
      f = exp(1i*2*pi*x).*(0.5*sin(8*pi*x).^2+0.5);
\end{verbatim}
\par}

\noindent  For any open or closed curve like this, one can
define a function $s(x)$ representing the arc length
over the parameter interval $[\kern .5pt 0,x]$:
\begin{equation}
s(x) = \int_0^x | f'(t) | \kern .7pt dt.
\label{arcformula}
\end{equation}
Chebfun realises this definition in one line:

{\small
\begin{verbatim}
      s = cumsum(abs(diff(f)));
\end{verbatim}
\begin{figure}[t]
\vskip 5pt
\begin{center}
\includegraphics[trim=0 45 0 0, clip,scale=.7]{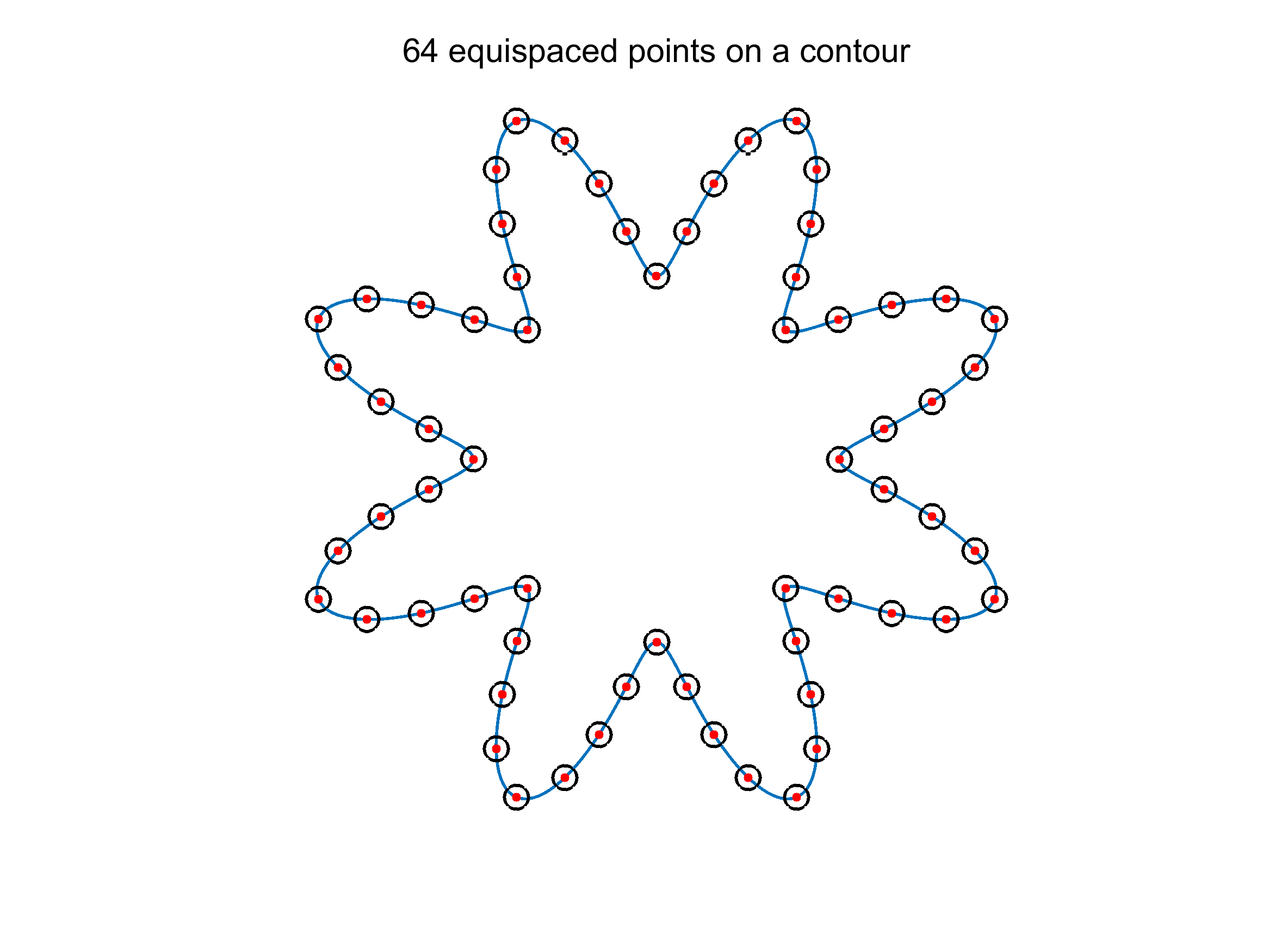}~~~~
\vspace{-10pt}
\end{center}
\caption{\label{arcfig}64 equispaced points on a curve determined by
a numerical computation of the inverse of an arc length function.
Circles come from a polynomial approximation in Chebfun, dots from a rational
approximation with AAA.}
\end{figure}

\noindent Now we ask, how can points be placed equidistantly along
the curve?  This is an inverse problem: we want values of $f(x)$
corresponding to equispaced values of $s$.  In \ccite{kuanxu}, the
result is computed by means of the {\tt inv} command in Chebfun,
generating the black circles plotted in the figure:

{\small
\begin{verbatim}
      chebfuninv = inv(s);
      distances = s(1)*(1:64)/64;
      chebfunpts = f(chebfuninv(distances));
\end{verbatim}
\par}
\par}

\noindent With AAA, we can obtain the same result by rational
approximation.  The following computation gives the red dots in
the figure, whose maximum difference from the Chebfun circles is
$1.5\cdot 10^{-5}$:

{\small
\begin{verbatim}
      T = linspace(0,1,800)';
      AAAinv = aaa(T,s(T),'tol',1e-5,'mmax',200);
      AAApts = f(AAAinv(distances));
\end{verbatim}
\par}

\noindent The timings for this experiment on our laptop are about 1.9
s for Chebfun and 1.1 s for AAA.\ \ This is not much of a speedup,
because a function like this, which is smooth and has many points of
complexity, does not play to the strengths of rational approximation.
It is interesting though how simple the inverse function calculation
is conceptually and computationally: just apply AAA in the reverse
of the usual direction by means of a command {\tt aaa(T,s(T))}
rather than the more familiar {\tt aaa(s(T),T)}.

For another example with a similar flavor but a clearer benefit
from rational approximation,
consider this $5\times 5$ matrix $A$ explored in \ccite{maxwell}:
\vspace{2pt}

\scalebox{.71}{\tt
~0.2560+0.0573i~~~0.0568+0.0800i~~~0.1597+0.2204i~~-0.1649+0.1315i~~-0.3639+0.0091i}

\vspace{-4pt}
\scalebox{.71}{\tt
~0.4733+0.2805i~~-0.3192+0.1267i~~~0.0810+0.0687i~~~0.5213+0.1574i~~-0.0596+0.2879i}

\vspace{-4pt}
\scalebox{.71}{\tt
~0.1447+0.3037i~~~0.2942+0.1844i~~-0.2918+0.0364i~~-0.2714+0.0265i~~-0.0849+0.2264i}

\vspace{-4pt}
\scalebox{.71}{\tt
-0.0650+0.1360i~~~0.0952+0.0813i~~-0.0503+0.0920i~~-0.1500+0.0814i~~~0.4742+0.1514i}

\vspace{-4pt}
\scalebox{.71}{\tt
~0.1938+0.0344i~~~0.0419+0.1868i~~-0.0453+0.0988i~~-0.2207+0.2483i~~-0.0772+0.1793i}
\vspace{4pt}

\begin{figure}
\vskip 5pt
\begin{center}
\includegraphics[trim=0 0 0 0, clip,scale=.66]{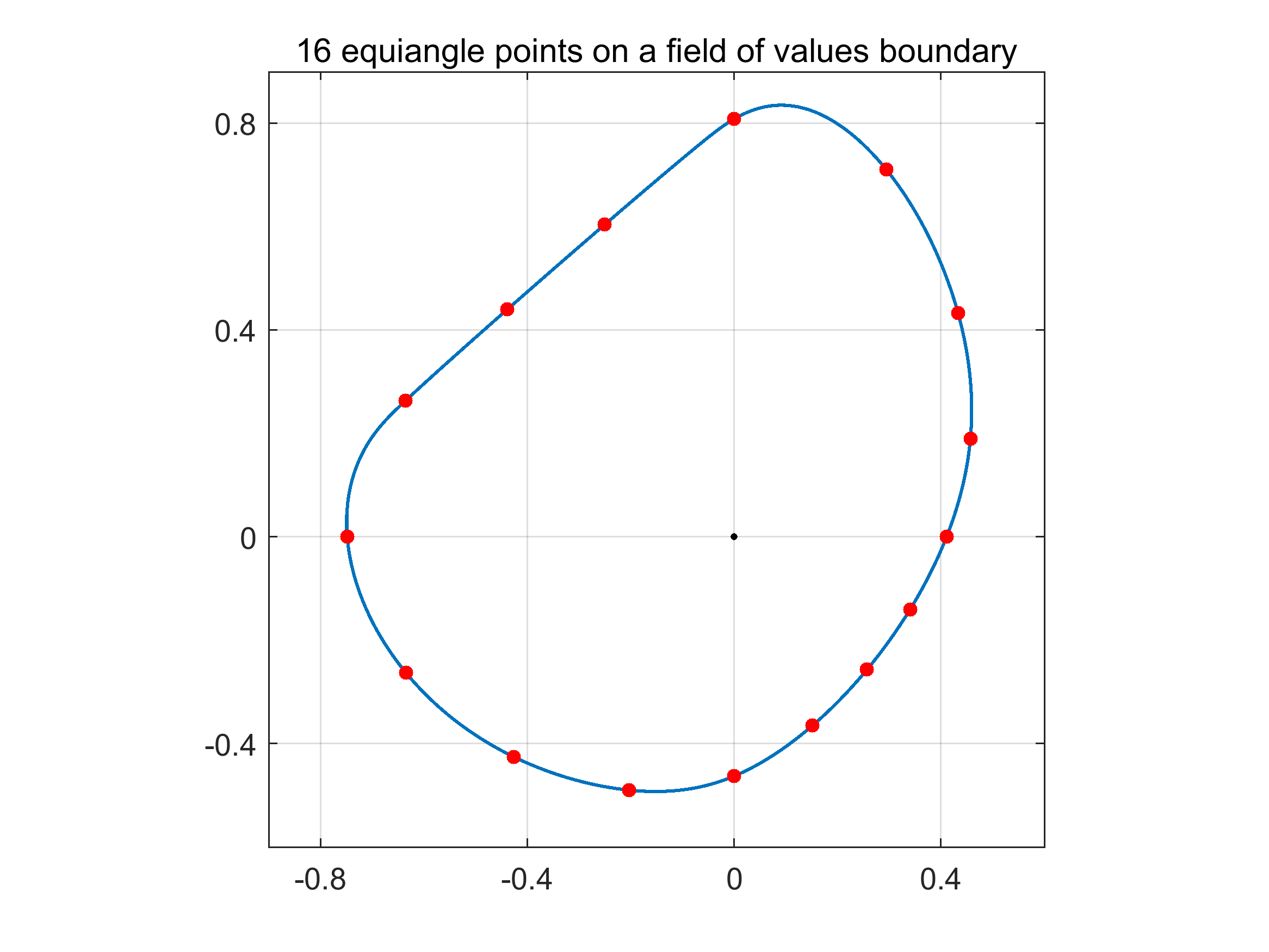}~~~~
\vspace{-17pt}
\end{center}
\caption{\label{equianglefig}Sixteen points on the boundary of
the field of values of a $5\times 5$ matrix spaced at
equal angles with respect to the origin.  AAA approximation of an inverse
function simplifies the representation of this curve from a degree 5703 polynomial
to a degree 68 rational function.}
\end{figure}

\noindent 
The boundary of the field of values (= numerical range) of $A$
can be computed by the Chebfun command

{\small
\begin{verbatim}
      C = fov(A)
\end{verbatim}
\par}

\noindent and is plotted in Figure~\ref{equianglefig}.  This is a
smooth curve, but as discussed in \ccite{fovexample}, that doesn't
mean it is obvious how to find a representation of it with a small
number of parameters.  In fact, the Chebfun representation is by a
polynomial of degree $5703$.  Large polynomial degrees are typical
in Chebfun representations of fields of values, and in this case
the behaviour is particularly extreme because a matrix has been
chosen whose field of values boundary includes a nearly straight
segment \ccite{maxwell}.

With a AAA rational approximation of degree just $68$, we can convert the curve
to a representation as a function of angle with respect to the origin:

{\small
\begin{verbatim}
      a = -unwrap(angle(C));
      tt = chebpts(2000,[0 2*pi]);
      [r,pol] = aaa(tt,a(tt),'tol',1e-12);
\end{verbatim}
\par}

\noindent  The commands

{\small
\begin{verbatim}
      angles = (1:16)*2*pi/16;
      angles2 = mod(angles-a(0),2*pi) + a(0) - 2*pi;
      plot(C(r(angles2)),'.')
\end{verbatim}
\par}

\noindent  give the red dots shown in the figure, 16 points equally
spaced with respect to angle.

The application that first made us aware of the power of rational
approximation for computing inverse functions was numerical conformal
mapping, to be discussed in section \ref{sec:conformal}.

\section{\label{sec:equi}Equispaced interpolation}
\begin{figure}
\vskip 7pt
\begin{center}
\includegraphics[trim = 0 70 0 0, clip, scale=.77]{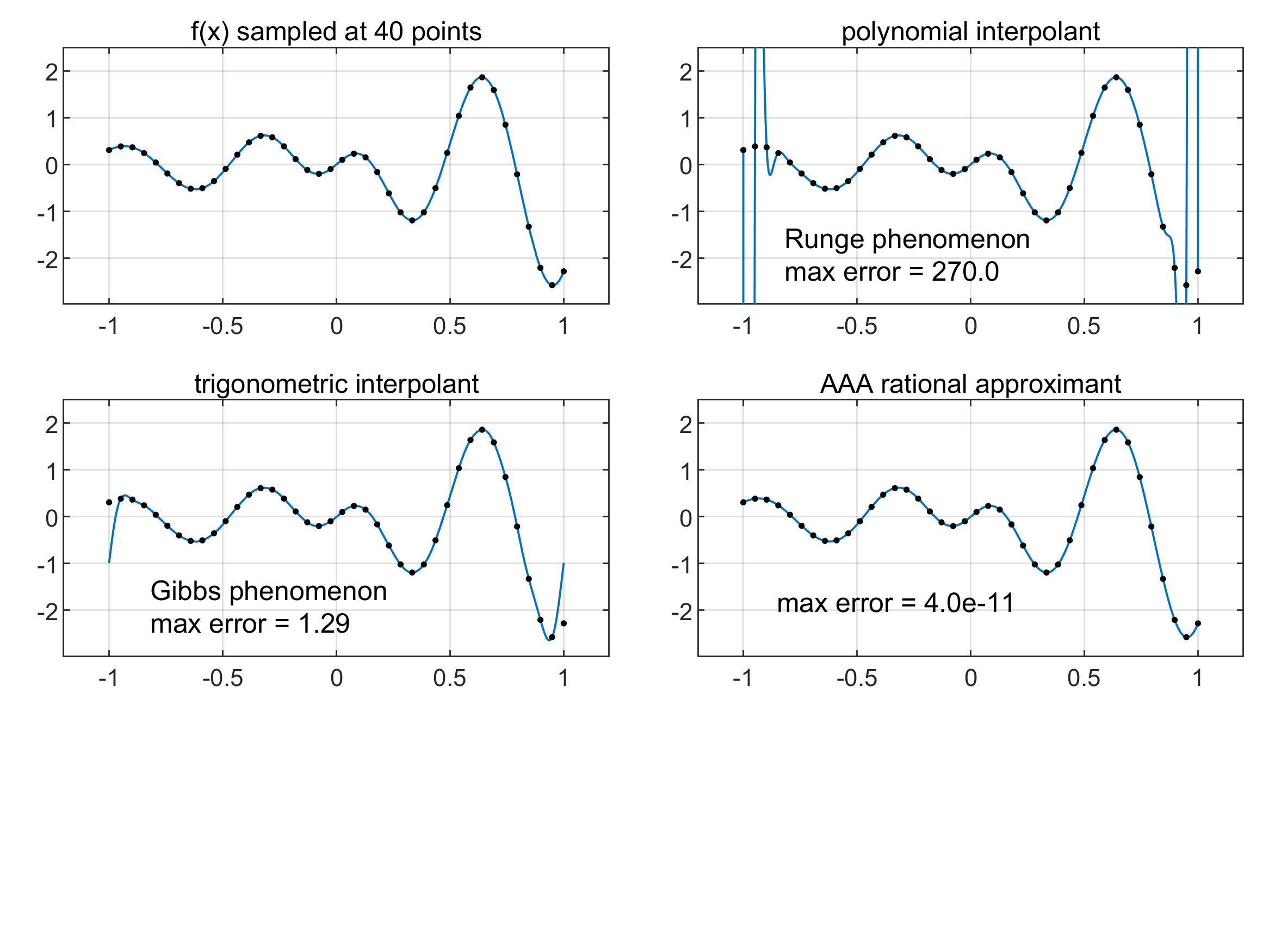}
\end{center}
\vskip -16pt
\caption{\label{equi-equifig} For interpolation of nonperiodic equispaced data,
polynomials suffer from the Runge phenomenon and trigonometric polynomials
from the Gibbs phenomenon.
The best general-purpose method appears to be AAA approximation.}
\end{figure}

One of the oldest computational challenges is the interpolation of
data in equally spaced points.  To illustrate, Figure~\ref{equi-equifig}
considers samples at
40 equispaced points from $-1$ to $1$ of an analytic function,
\begin{equation}
f(x) = e^x \kern -.5pt \cos(10\kern .5pt x)\kern .5pt  \tanh(4\kern .4pt x).
\end{equation}
From these data, how might we recover $f$ to good accuracy between
the samples?

One idea would be polynomial interpolation.  This method
was discredited by Runge, who showed that degree $n-1$
polynomial interpolants of data at $n$ equispaced points don't
generally converge as $n\to\infty$, even when the underlying
function is analytic \ccite{runge}.  The second panel of
Figure~\ref{equi-equifig} shows that for this problem, all
accuracy is lost for $|x|>0.8$, and in fact, the error near the
endpoints for this computation with $n=40$ is as large as $270$.
(With 60 equispaced interpolation points, it would be $6450$.)
This effect is called the {\em Runge phenomenon\/} \ccite{atap}.

Another idea would be Fourier or trigonometric interpolation, i.e.,
interpolation by a minimal-degree trigonometric polynomial composed
of sines and cosines scaled to the given interval.  This would be
excellent if $f$ were periodic, but in the non-periodic case the
idea fails for a reason also going back to around 1900, the {\em
Gibbs phenomenon\/} \ccite{atap}.  For the variant of trigonometric
interpolation shown in the figure, the value interpolated at $x=\pm
1$ is the mean of $f(-1)$ and $f(1)$, leading to an error at that
point of $|f(1)-f(-1)|/2\approx 1.29$.  There are large errors also
nearby, though this is not so visible in the plot, and even in the
middle of the interval, the accuracy is no better than $0.003$.

In the face of these long-recognised failures of the two most
obvious methods of interpolation, a large literature has grown
up of papers with titles like ``Defeating the Runge phenomenon''
and ``On the Gibbs phenomenon and its resolution.''  Many methods
have been proposed, some of them rather complicated, including
polynomial least-squares fitting, Fourier extensions, Fourier
series with boundary corrections, multi-domain approximations,
splines, conformal maps, Gegenbauer reconstruction, regularisation,
and Floater-Hormann interpolation of adaptively determined degree.
This list comes from \ccite{equi}, where references can be found;
see also the earlier review in \ccite{impossibility}.  In the latter
paper the ``impossibility theorem'' is proved, which asserts that no
method for approximating analytic functions from equispaced samples
can be exponentially convergent unless it is exponentially unstable.
This establishes in a rigorous sense that a comprehensive defeat
of the Runge and Gibbs effects is not possible.

Examples can be found in which almost any of these methods
outperforms the others.  The best general purpose method for
equispaced interpolation, however, appears to be simply to
construct a AAA fit to the data.  (Strictly speaking the result
is then an approximation, not an interpolant.)  The final part of
Figure~\ref{equi-equifig} shows that for this example, the accuracy
is ten digits as delivered by the code segment

{\small
\begin{verbatim} 
      f = @(x) exp(x).*cos(10*x).*tanh(4*x);
      X = -1+2*(0:39)'/39;
      r = aaa(f(X),X);
\end{verbatim}  
\par}

\begin{figure}
\vskip 7pt
\begin{center}
\includegraphics[trim=0 0 130 201, clip, scale=1.14]{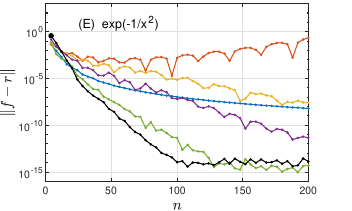}
\includegraphics[trim=190 0 0 201, clip, scale=1.14]{htfig-eps-converted-to.pdf}
\end{center}
\vskip -8pt
\caption{\label{equi-htfig} Illustration from \ccite{equi}
of six different methods of equispaced interpolation
applied to the function $f(x) = \exp(-1/x^2)$ for $x\in[-1,1]$.
AAA approximation shows steady convergence down to accuracy $10^{-13}$ as
the number of sample points $n$ increases.
The other methods converge less quickly.  Examples
vary, but this behaviour is, broadly speaking, typical.}
\end{figure}

\noindent The degree of the AAA rational approximant is
$18$.  For an extensive discussion with
many further examples and comparisons with other methods, see \ccite{equi}.
To illustrate one such comparison, Figure~\ref{equi-htfig} reproduces
an image from \ccite{equi} showing the behaviour as a function of the number of
sample points in approximation of the $C^\infty$ function $f(x) = \exp(-1/x^2)$
on $[-1,1]$.

Perhaps this is a good place to emphasise a consideration that
applies throughout this article and to AAA approximation generally.
AAA can be extremely accurate for smooth data, but in the presence
of noise, it typically fails.  If data are known to be corrupted
by noise at a level $\varepsilon$, one can often do quite well by
running AAA with a tolerance set above that level.  This is why
the default AAA tolerance is as high as $10^{-13}$ even though IEEE
machine epsilon is about $10^{-16}$.  Attempting to fit closer to
the noise level is likely to lead to loss of accuracy and spurious
poles, and for noisy problems in general, and noisy equispaced
interpolation in particular, one may be well advised to rely on
methods that are more robust---more linear---even if they will
usually be less accurate when the problem is smooth.

\section{\label{sec:imputing}Imputing missing data}
A variation on the last section's theme of interpolation between
equispaced data points is the problem of filling in larger stretches
of a smooth function that may be absent from a dataset.  An example
of this kind (for trigonometric as opposed to ordinary rational
approximation) appeared in Heather Wilber's 2021 PhD
thesis \ccite{wilber} and was reproduced in \ccite{wdt}.
\begin{figure}
\vskip 0pt
\begin{center}
\includegraphics[trim=30 100 20 0, clip, scale=.82]{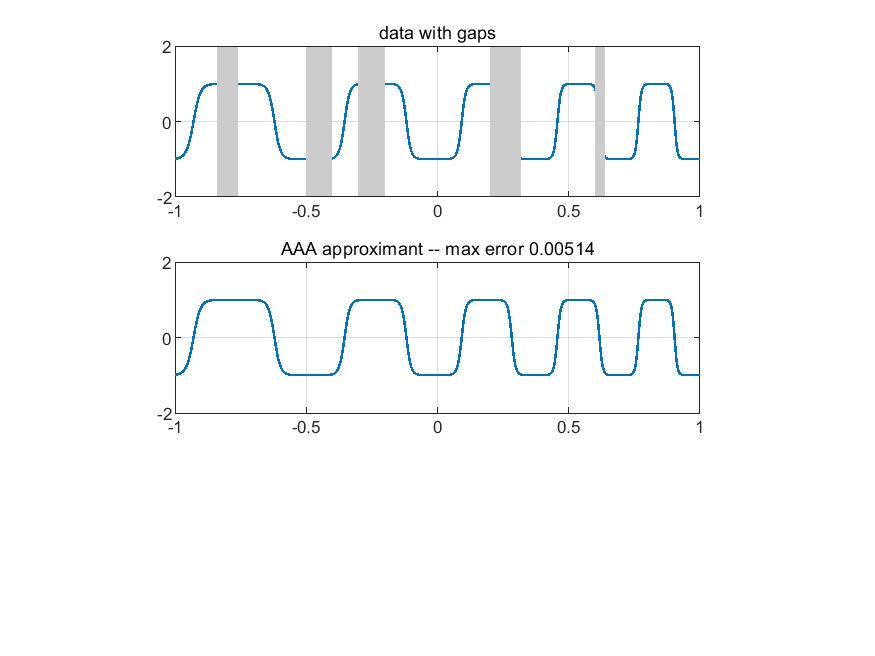}
\end{center}
\vskip -11pt
\caption{\label{imputing-gap1} Use of AAA approximation to 
fill in gaps in an analytic function.}
\end{figure}

\begin{figure}[t]
\vskip 0pt
\begin{center}
\includegraphics[trim=30 100 20 0, clip, scale=0.82]{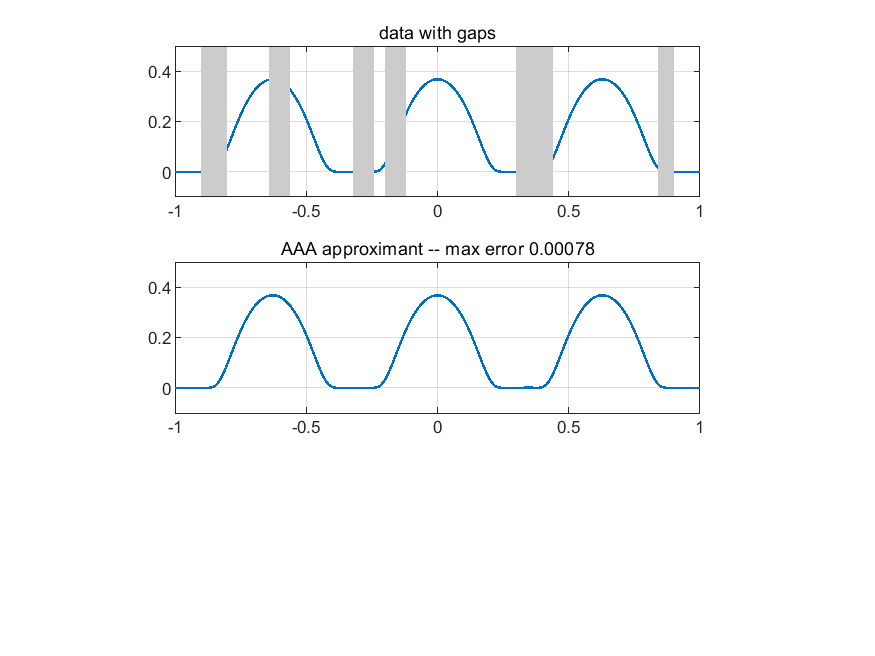}
\end{center}
\vskip -11pt
\caption{\label{imputing-gap2} A similar experiment to that
of Figure \ref{imputing-gap1} but now for a $C^\infty$ function.}
\end{figure}

Figure~\ref{imputing-gap1} illustrates the process.
The function 
$f(x) = \tanh(4\sin(30\kern .5pt e^{x/2}))$
has been sampled in 1000 equally spaced points in
$[-1,1]$, and then 275 of the samples, marked by
the grey stripes, have been removed.  When AAA approximation
is applied to the remaining data, $f(x)$ is
recovered in about $0.1$ s on a laptop
by a rational function of degree $90$ that matches it with
a maximum error of $0.8\cdot 10^{-5}$ on $[-1,1]$.
For this experiment the standard tolerance of $10^{-13}$ leads
to spurious poles with very small residues on $[-1,1]$, so 
it was loosened to $10^{-10}$.
Figure~\ref{imputing-gap2} shows a similar experiment but
with a function that is $C^\infty$ rather than analytic, 
$f(x) = \exp(-1/\cos(5\kern .5pt x)^2).$
This time the rational approximation has degree $72$.

\begin{figure}
\vskip 10pt
\begin{center}
\includegraphics[trim=50 175 20 0, clip, scale=0.89]{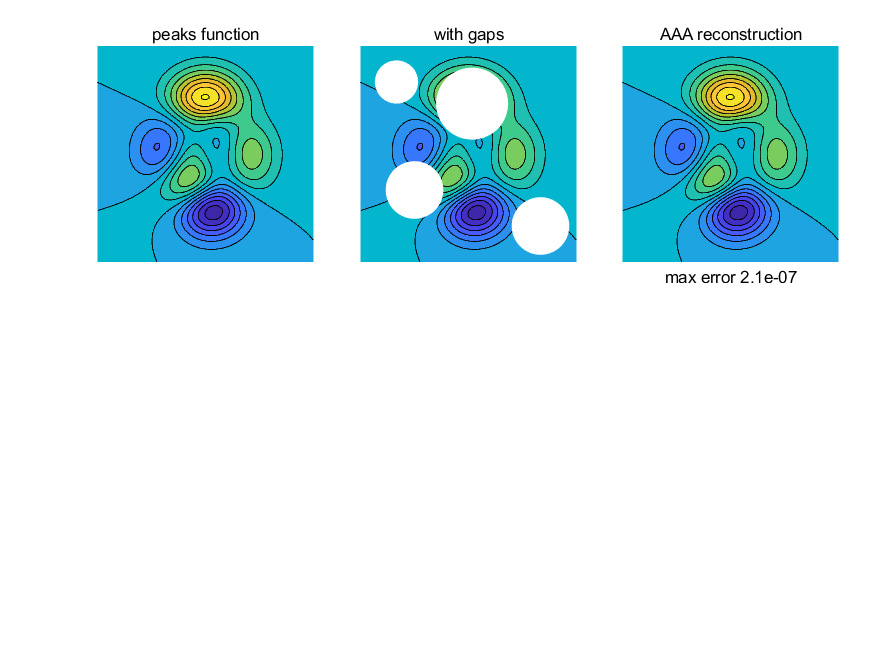}
\end{center}
\vskip -6pt
\caption{\label{imputing-gap3} Filling in missing data in the bivariate
MATLAB ``peaks'' function (\ref{imputing-peaks}).
Here 100 univariate AAA approximations are
computed along horizontal slices of the image.  The total computing time is
still less than a second.}
\end{figure}

A generalisation of this idea can be applied in multiple dimensions, as
illustrated in Figure~\ref{imputing-gap3}.
Here we start from the MATLAB ``peaks'' function, which is the
bivariate real analytic function
\begin{equation}
\begin{split}
f(x&,y) = 3(1-x)^2e^{-x^2 - (y+1)^2} \\
   &- 10(x/5 - x^3 - y^5)\kern .7pt e^{-x^2-y^2}
   - e^{-(x+1)^2 - y^2}\kern -2pt /3.
\end{split}
\label{imputing-peaks}
\end{equation}
Ideally, we would like to find a near-best approximation of $f$
by a bivariate rational function $r(x,y)$.  However, no method
for computing multivariate rational approximations is known
that has anything like the power of AAA in the univariate case
(see section \ref{sec:bivariate}).  Instead, to fill in the gaps,
we can use a sequence of univariate approximations.  (This idea
originates in \ccite{jjiam}.)  In the example, $f$ is sampled on
a $100\times 100$ regular grid, and then on each horizontal line,
a AAA approximation is computed from the data after eliminating
the values in the gaps.   In other words, we compute 100 1D
approximations like those illustrated in Figures~\ref{imputing-gap1}
and~\ref{imputing-gap2}.  The whole computation takes $0.7$ s on
a laptop and obtains excellent accuracy.

In image analysis, the process of filling in missing data is known
as {\em inpainting\/}, and this is an area of engineering on which
many methods have been brought to bear, including machine learning.
The strength of rational approximation is its ability to deal with
smooth functions to high accuracy, and we would not claim that this
regime is typical of inpainting problems.

\section{\label{sec:analcont}Analytic continuation}
In one way or another, half the applications of rational
approximation and half the sections of this paper might be said to
have something to do with analytic continuation.  This is the process
whereby information about an analytic function in some domain of
the complex plane is used to infer information about it elsewhere.
The idea is fundamental in theory, but problematic in practice.
For an extended discussion, see \ccite{jjiam}.

In theory, analytic continuation can be carried out by the
``chain-of-disks'' method of Weierstrass.  If $f$ is analytic at a
point $z_1^{}$, then its Taylor series converges in a disk about
$z_1^{}$.  At any other point $z_2^{}$ in that disk, $f$ is also
analytic, so one can construct a new Taylor series about $z_2^{}$.
By moving along in this way one can in principle extend an analytic
function along an arbitrary path in the complex plane so long as no
singularities are encountered.  Poles and essential singularities
are isolated, and one can analytically continue around them.
Branch points are not isolated, and lead to analytic continuations
having different values depending on which way you pass around
them---``multivalued analytic functions'' with different branches,
but still analytic.  These are the singularities occurring most
commonly in practice, but another possibility is that $f$ might have
a {\em natural boundary\/} beyond which no analytic continuation
at all is possible.  For example, a Taylor series with random
coefficients from the standard normal distribution will have the unit
circle as a natural boundary with probability~$1$ \ccite{kahane}.

The trouble is that these operations assume one has exact information
about a function.  If $f$ is only known to accuracy $\ep$, all
is lost.  Here is how it is put in Theorem 5.1 of \ccite{analcont}.
An earlier paper that expresses the same result is \ccite{miller}.

\begin{theorem}
\label{analcont-thm1}
Let $\Omega$ be a connected open region of $\complex$ and let
$E$ be a bounded nonempty continuum in $\overline{\Omega}$
whose closure $\overline{E}$ does not enclose any points
of $\Omega\backslash \overline{E}$.  Let $f$ be an analytic function in
$\Omega\cup E$ satisfying $\|f\|_E^{} \le\ep$ for some $\ep>0$.
Then no bounds whatsoever can be inferred on the value of
$f$ at any point $z\in\Omega\backslash\overline{E}$.
\end{theorem}

This theorem would appear to suggest that numerical analytic
continuation is impossible, but in fact, much can often be achieved
for reasons whose theoretical basis is not always clear.  One case
that can be analysed mathematically, going back to the Hadamard
three-circles theorem of complex analysis, is the situation where
it is known that $f$ is bounded in $\Omega$.  Here is how it is
put in Theorem 5.2 of \ccite{analcont}.

\begin{theorem}
Let $\Omega$, $E$, $\ep$, and $f$ be as in Theorem \ref{analcont-thm1},
but now with\/ $f$ additionally satisfying
$\|f\|_\Omega^{} \le 1$, and let $z$ be a point in
$\Omega\backslash \overline{E}$.  Assume that the boundary of
$E$ is piecewise smooth (a finite union of smooth Jordan arcs).  Then
there is a number $\alpha\in(0,1)$, independent of $f$ though not of
$z$, such that for all\/ $\ep>0$,
\begin{equation}
|f(z)| \le \ep^{\alpha(z)}.
\label{analcont-alpha}
\end{equation}
\end{theorem}

The fractional power $\alpha(z)$ in this theorem has a simple
interpretation.  {\em If we know $f$ to $d$ digits of accuracy on
$E$, then it is determined to $\alpha(z)\kern .7pt d$ digits of
accuracy at $z\in\Omega\backslash \overline{E}$.}  For a bounded
analytic function $f$ in the annulus $r_1<|z|<r_2$, for example, the case of
the three-circles theorem, if $d$ digits are given on the inner
circle, then the number of digits falls off smoothly to~$0$ on
the outer circle according to the formula
\begin{equation}
\alpha(z) = 1 - {\log |z/r_1|\over \log (r_2/r_1)}.
\label{analcont-falloff}
\end{equation}
Results for other domains can be obtained by conformal
transplantation, and are generally more adverse.  In an infinite
half-strip of half-width 1, for example, if one knows $f$ to $d$
digits on the end segment, the number of digits falls off by a
factor of $10$ with each translation along the strip by a distance
$(2/\pi)\log 10 \approx 1.47$.  At a distance of 6 into the strip,
therefore, which is just three times the width, the number of digits
has been divided by more than $10{,}000$.  Similarly if one tracks
a function like $\sqrt z$ known to be bounded and analytic for $0
< |z| < 2$ from $z=1$ around the origin and back to $z=1$ again,
the number of accurate digits may be divided by $(\pi/4)\exp(\pi^2)
\approx 15{,}000$.  You'd need to start out with $15{,}000$ digits
of accuracy to end up with $1$---or $225{,}000{,}000$ digits for
two trips around the origin.  All these results are discussed in
\ccite{analcont}.

It is interesting to note that (\ref{analcont-alpha}) indicates
that analytic continuation, under a boundedness assumption, is an
example of a problem that is well-posed (continuous dependence on
data) but with infinite condition number (not Lipschitz continuous).

In the face of such pessimistic bounds, it is often surprising how successfully
analytic continuation can be carried out in practice by various means.
A chapter of the textbook \ccite{fpiret} describes nine methods that are
often useful (several of them analytic or semi-analytic
as opposed to strictly numerical):
(1) Circle-chain method, (2) Schwarz reflection principle,
(3) Use of a functional equation,
(4) Partitioning of an integration interval,
(5) Replace Taylor coefficients by integrals or sums,
(6) Subtraction of a similar series or integral,
(7) Borel summation, (8) Ramanujan’s
formula, and (9) Pad\'e approximations.

Here, we will focus on a tenth approach (or one could regard it
as an extension of the ninth), namely rational approximation.
So far as we are aware, this is the best method available for
general-purpose strictly numerical analytic continuation.

\begin{figure}
\begin{center}
\includegraphics[trim=0 158 0 0, clip, scale=0.89]{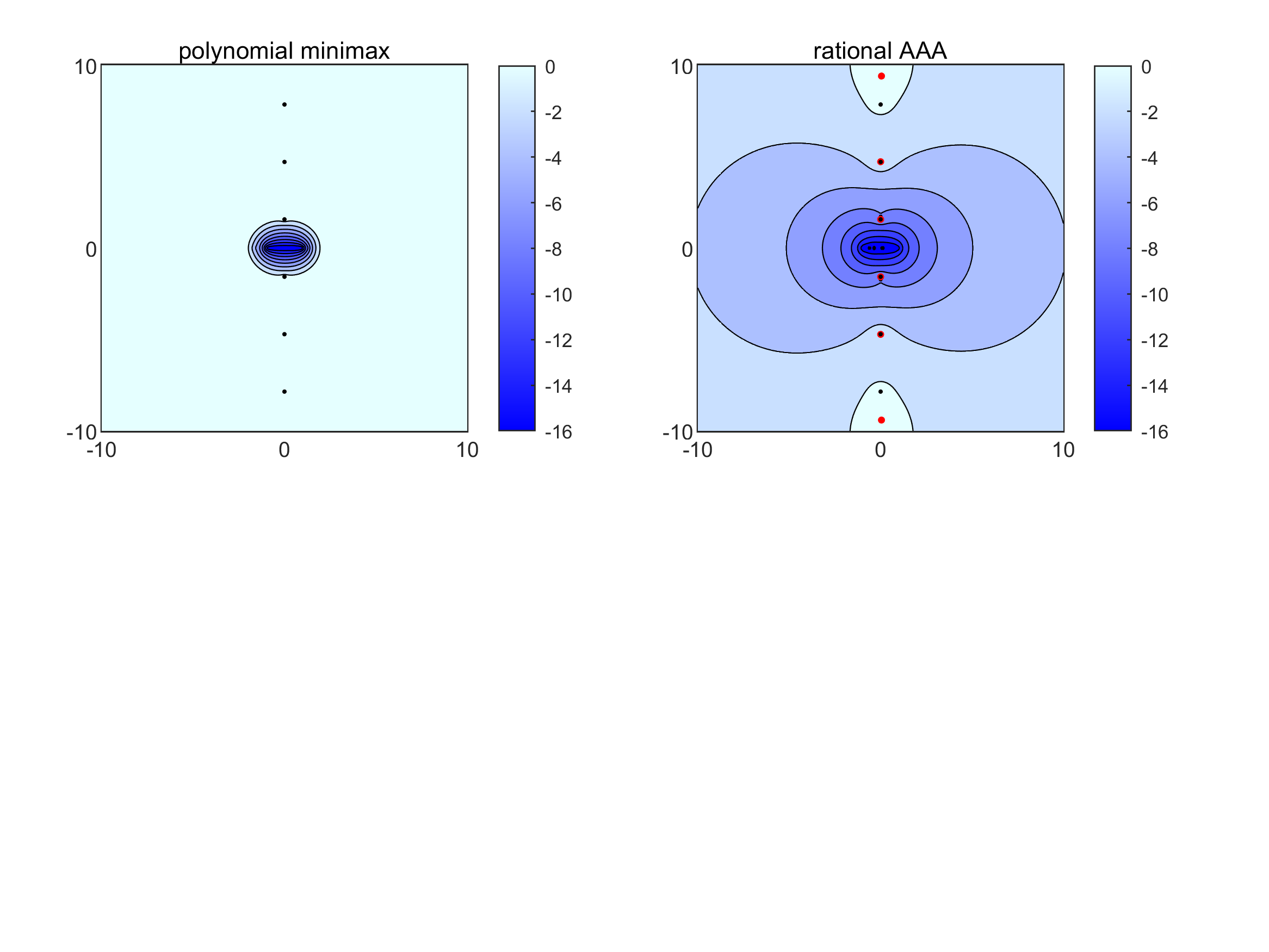}~~
\end{center}
\vskip -6pt
\caption{\label{analcont-contours} Numerical analytic continuation
of $\tanh(z)$ from data on $[-1,1]$ by degree $30$ polynomial
minimax fitting (left) and degree $7$ AAA rational fitting (right)
(error contours of $\log_{10}^{} |f(z)-r(z)|$).
Black dots mark poles of $f$, and red dots mark poles of the
rational approximation $r$.  As is well known, polynomials have
no accuracy beyond the first singularity---with ``beyond'' and ``first'' for this
case of approximation on the interval defined by Bernstein ellipses
surrounding $[-1,1]$.  Figure from \ccite{jjiam}.}
\end{figure}

Any experiment in which a AAA approximant is evaluated outside
the approximation domain can be taken as an illustration of
numerical analytic continuation by rational approximation.
Figure~\ref{analcont-contours}, taken from \ccite{jjiam}, shows
one such.  Note that polynomial approximation, shown on the left,
cannot get past the closest singularity (``closest'' as defined in
a fashion related to potential theory), whereas rational functions
can do much better.  We have already seen these effects earlier in
the paper.  For example, the contour lines of Figure \ref{fig3.2}
show considerable accuracy in analytic continuation of $\Gamma(z)$
from the unit interval or disk out to values of $|z|$ as great as
$3$ or $4$, despite the poles at $x = 0, -1, -2, \dots.$ Similarly,
the zeta function example of Figure \ref{fig4.2}, though the images
focus on zeros and poles, reflect accurate analytic continuation
from the line $\Re(z) = 4$ well across the imaginary axis and into
the left half-plane.

Though such experiments are encouraging, it is hard to know what
can be said precisely.  As discussed at length in \ccite{jjiam},
the main analysis available in this area is derived from
potential theory following ideas of Walsh, Gonchar, Stahl and
others which we will present here in section~\ref{sec:potential}.
Approximately speaking, this theory suggests that numerical analytic
continuation by near-best rational approximations will typically
achieve accuracies related to a potential function defined by $n$
poles $\{\pk\}$ of $r$ and $2\kern .3pt n+1$ interpolation points
$\{\zk\}$ where $r(\zk) = f(\zk)$, and that the positions of these
poles and interpolation points will be determined by singularities
of~$f$.  However, we are not aware if arguments of this kind have
been carried very far to attempt to quantify convergence rates.
There is also the problem that convergence at particular points
$z\in \complex$ can never be guaranteed because of difficulties
related to ``Froissart doublets''---poles with near-zero residues
at seemingly arbitrary locations in $\complex$ that appear spurious
from an approximation point of view and yet may be unavoidable
even in exact arithmetic (and are almost invariably worse with
rounding errors).  Such effects have been studied extensively in
connection with Pad\'e approximation \ccite{bgm,stahl98}.

\begin{figure}
\vskip 12pt
\begin{center}
\includegraphics[trim=0 0 0 0, clip, scale=1]{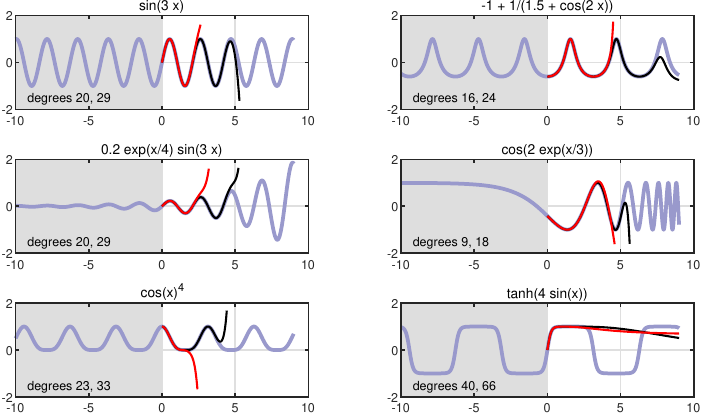}~~
\end{center}
\vskip -3pt
\caption{\label{analcont-six}
Illustration of the ``one-wavelength
principle'' by six examples in which an analytic function
$f(x)$ is analytically continued by AAA approximation from data for
$x\in [-10,0\kern .5pt]$
to $x>0$.  The thick blue-grey curves show $f(x)$,
the red curves show approximations
$r(x)$ with $13$-digit accuracy, and 
the black curves show approximations
$r(x)$ with $29$-digit accuracy.
Approximately one and two wavelengths of successful continuation,
respectively, appear to be achieved at least in the first four cases.}
\end{figure}

Figure \ref{analcont-six} illustrates an idea introduced in
\ccite{jjiam}, the {\em one-wavelength principle\/}.  This is a
rule of thumb, not a theorem.  The observation is that in many
cases in which AAA is applied to a function which has some kind of
oscillatory behaviour, about one wavelength of successful analytic
continuation beyond the data domain is observed, after which accuracy
is lost.  One wavelength is not an absolute figure, but is tied
to the accuracy of the rational approximation, which by default
in the Chebfun AAA code is 13 digits.  Figure \ref{analcont-six}
shows that if the AAA accuracy is increased to 29 digits (thanks to
Driscoll for 128-bit computations in Julia \ccite{driscolljuliacon}),
the number of wavelengths approximately doubles.  Continuing such
experiments to higher precision, Figure~7 of \ccite{jjiam}, based
on data from Daan Huybrechs, shows a roughly linear relationship
of the number of wavelengths with the number of digits of accuracy
of the rational approximation.  Thus about two wavelengths are
possible with 26-digit approximations, three wavelengths with
39-digit approximations, and so on.

It is not clear how general the one wavelength principle is, or how
it might be formulated precisely.  Before \ccite{jjiam}, we have
found just one publication that seems to speak to the topic: an
article by Henry Landau in 1986 written in the context of sampling
theory.  See Theorem 2 of \ccite{landau}.  Unfortunately, the route
this suggests to try to make the principle precise, by restricting
attention to band-limited functions of a certain class, has been
falsified by a theorem proved by Alex Cohen of MIT \ccite{cohen}.
So the problem remains open of finding a theorem to explain the
one-wavelength effect.  We shall discuss the related literature of
``function extension'' in section~\ref{sec:bivariate}.

\section{\label{sec:schwarz}Computing the Schwarz function} The last
section mentioned a list of nine methods of analytic continuation
discussed in the textbook \ccite{fpiret}, of which the second was the
Schwarz reflection principle.  This is one of the fundamental results
of complex analysis, going back to Schwarz long ago \ccite{schwarz}.
The present section on the generalisation known as the Schwarz
function is both interesting in its own right and important for
the theoretical understanding of the PDE solution techniques to be
discussed in sections \ref{sec:aaals}--\ref{sec:hilbert}.

The simplest version of the reflection principle runs as follows.
For a proof and further details, see \ccite{schwarzfun}.

\begin{figure}
\begin{center}
\vspace{10pt}
\includegraphics[clip, scale=.80]{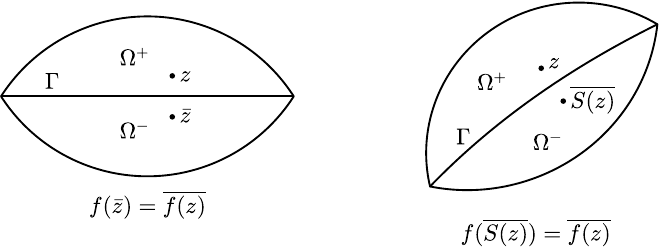}
\vspace{-6pt}
\end{center}
\caption{\label{schwarz-schematic}On the left, the Schwarz reflection principle
(Theorem \ref{schwarz-thm1}) extends an analytic function that takes real values
on a real interval across the interval from one half-plane to
the other.  On the right, the generalisation by the Schwarz function
$S(z)$ for reflection of an analytic function
across an analytic arc $\Gamma$, again assuming it takes real
values on $\Gamma$.  Figure from \ccite{schwarzfun}.}
\end{figure}

As sketched on the left in Figure \ref{schwarz-schematic}, let
$\Omega\subseteq \complex$ be a connected open set that is symmetric
about $\real$, $\overline{\Omega} = \Omega$, and whose intersection
with $\real$ is a finite or infinite open interval $\Gamma$.
(The symbol $\overline{\Omega}$ denotes the complex conjugate of
$\Omega$, in contrast to the last section, where it represented
the closure of $\Omega$.)  Define $\Omega^+ = \{z\in \Omega: \Im
(z) > 0\}$ and $\Omega^- = \{z\in \Omega: \Im (z) < 0\}$.

\begin{theorem}
\label{schwarz-thm1}
{\bf Schwarz reflection principle.} Let\/ $\Omega$, $\Gamma$,
$\Omega^+$ and\/ $\Omega^-$ be as defined above, and let $f$
be analytic in the interior of\/ $\Omega^+$ and continuous
on $\Omega^+\cup \Gamma$, taking real values on $\Gamma$.
Then the formula \begin{equation} f(\bar z ) = \overline{f(z)}
\label{analcont} \end{equation} defines an analytic continuation of\/
$f$ to all of\/ $\Omega$.
\end{theorem}

The Schwarz function generalises this result to cases where
$\Gamma$ is not a real segment but an analytic curve,
as sketched on the right in Figure~\ref{schwarz-schematic}.
Let $\Gamma$ be an analytic Jordan arc or Jordan curve,
that is,
the image in $\complex$ of the unit circle or an arc of the unit circle
under a one-to-one analytic function with nonvanishing derivative.
Consider the function
\begin{equation}
S(z) = \bar z,  \quad z\in \Gamma.
\label{Sdef}
\end{equation}
Like any other analytic function defined on a curve in $\complex$,
$S$ can be analytically continued to a complex neighborhood of
$\Gamma$, and this is the Schwarz function.  (It may be startling to
read that $z\mapsto \bar z$ is an analytic function, and this would
not be true on a domain with interior, but here we are dealing with just
a curve.)  The point of $S$ is that
$z\mapsto \overline{S(z)}$ is the generalisation to a general arc
$\Gamma$ of the reflection $z\mapsto \bar z$ when $\Gamma$ is real.
As the simplest example, the Schwarz function of the unit circle---or
of any arc of the unit circle---is $S(z) = z^{-1}$.

The analytic continuation works as follows.  Since $S$ is analytic,
it must extend to an analytic function on some neighborhood
of $\Gamma$.  Since $|d\kern 1pt \overline{S(z)}/dz| = 1$ on $\Gamma$, it follows
that in a sufficiently small neighborhood of $\Gamma$, $\overline{S(z)}$ maps a point
$z$ close to $\Gamma$ on one side to a point close to $\Gamma$ on
the other side.  This function
$\overline{S(z)}$ is co-analytic, not analytic, but if we compose it with
itself we get an analytic function, the identity:
\begin{equation}
\overline{S( \overline{S(z)} )} = z.
\label{schwarz-ident}
\end{equation}
With this property in hand, as sketched on the right side
of Figure \ref{schwarz-schematic}, we define
a {\em reflection domain} for $\Gamma$
to be an open set $\Omega\subseteq\complex$ containing $\Gamma$
that is divided into two disjoint open sets $\Omega^-$
and $\Omega^+$ by $\Gamma$ and in which a single-valued
analytic Schwarz function $S$ for $\Gamma$ can
be defined that satisfies
\begin{equation}
\overline{S(\Omega^+)} = \Omega^-, \quad
\overline{S(\Omega^-)} = \Omega^+,
\label{defn1}
\end{equation}
and
\begin{equation}
\overline{S( \overline{S(z)} )} = z, \quad z\in \Omega.
\label{defn2}
\end{equation}
With this definition, the generalisation of the reflection principle
is as follows.  See \ccite{schwarzfun} for a proof.

\begin{theorem}
{\bf Reflection across an analytic curve.}
Let $\Omega$ be a reflection domain for an analytic Jordan curve
or Jordan arc $\Gamma$ and let 
$f$ be analytic in the interior
of\/ $\Omega^+$ and continuous on
$\Omega^+\cup \Gamma$, taking real values on\/ $\Gamma$.
Then the formula
\begin{equation}
f(\overline{S(z)}) = \overline{f(z)}
\label{Scont}
\end{equation}
defines an analytic continuation of\/ $f$ to all of\/ $\Omega$.
\end{theorem}

Although some of these ideas go back as far as Schwarz
\ccite{schwarz}, the Schwarz function was not defined and named until
almost a century later \ccite{dp}.  Two books have been published
about it \ccite{davis,shapiro}, but they are nonnumerical.  So far
as we are aware, no numerical method was proposed for computing
the Schwarz function before the arrival of AAA, when the paper
\ccite{schwarzfun} has appeared in elaboration of earlier experiments
presented in \ccite{jjiam,inlets}.

\begin{figure}[t]
\begin{center}
\vspace{10pt}
\includegraphics[trim=5 134 0 4, clip, scale=.85]{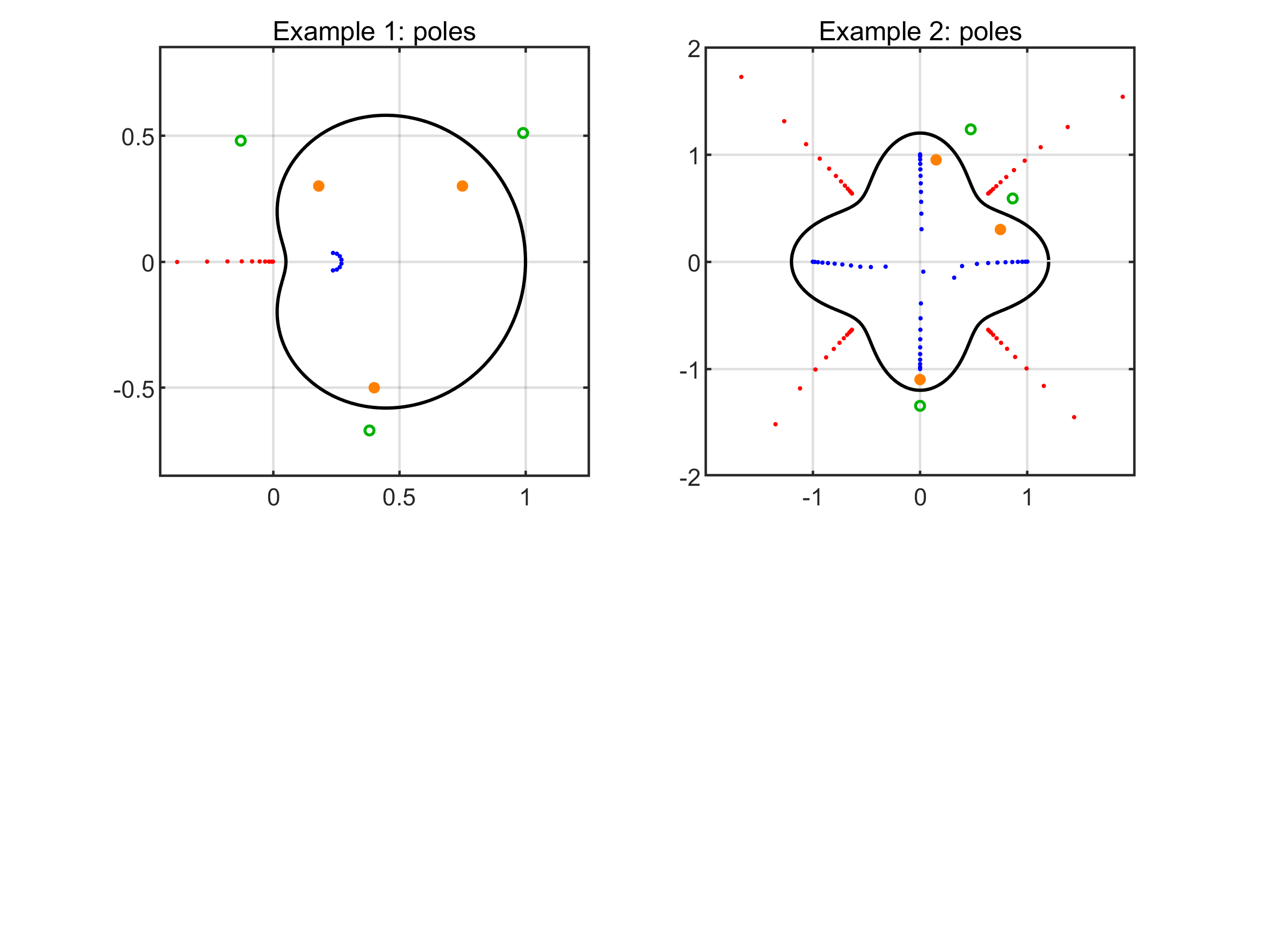}~~~~
\end{center}
\begin{center}
\includegraphics[trim=8 144 0 5, clip, scale=.85]{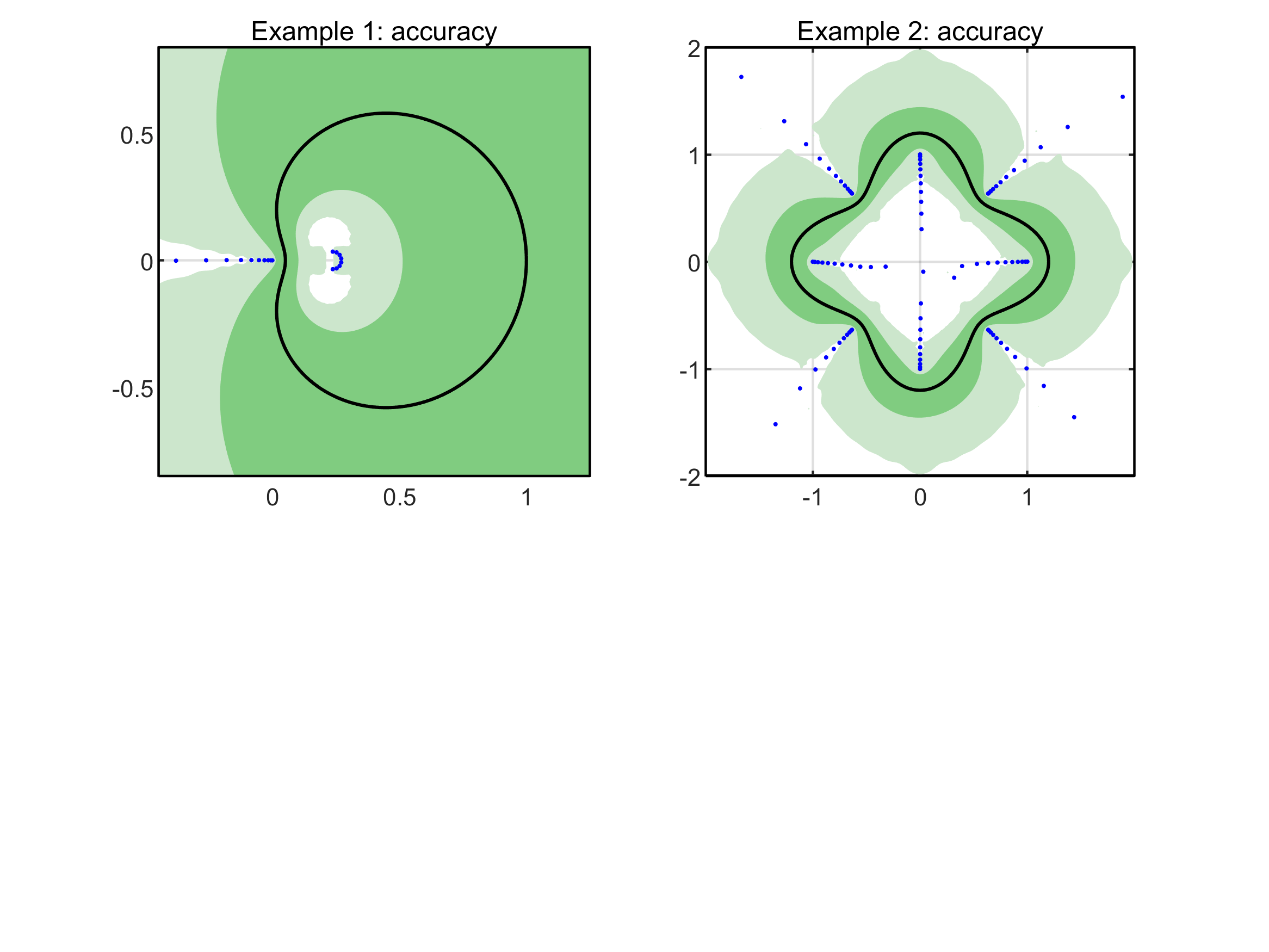}~~
\end{center}
\vskip -6pt
\caption{\label{schwarz-squiggles}Numerically computed Schwarz functions for
closed curves $\Gamma$ obtained as images of the unit circle under
maps $f$.  On the left, $f(z) = \exp(1.5(z-1))$, and on the right,
$f(z) = z(1+0.2\cos(4\kern .9pt \hbox{arg}(z)))$.  In the first row of plots,
the dots show poles of the AAA approximants, blue for poles inside $\Gamma$
and red for poles outside.  The
orange/green dots illustrate three arbitrary pairs of
points mapped to each other by 
$\overline{S(z)}$.  In the second row, as a test of (\ref{schwarz-ident}),
dark and light green shading
indicates regions where $|\overline{S(\overline{S(z)})} - z|$ is
less than $10^{-8}$ and $10^{-1}$, respectively.}
\end{figure}

The AAA method for computing the Schwarz function could
not be simpler: just calculate a rational approximation of
$\kern 1pt \overline z$ on $\Gamma$.  To illustrate, Figure
\ref{schwarz-squiggles} shows the poles of AAA approximations to
the Schwarz functions of two curves defined by:

{\small
\begin{verbatim}
      T = 2*pi*(1:500)'/500;     % angles
      C = exp(1i*T);             % unit circle
      Z = exp(1.5*(C-1));        % curve 1
      Z = C.*(1+.2*cos(4*T));    % curve 2
\end{verbatim}
\par}

\noindent
The AAA approximations and their poles are computed by the
command

{\small
\begin{verbatim}
      [r,poles] = aaa(conj(Z),Z);
\end{verbatim}
\par}

\noindent with poles plotted inside $\Gamma$ as blue dots and outside
as red dots.  To the experienced eye, these poles show a great deal.
Evidently the Schwarz function $S$ is well-behaved close to $\Gamma$,
but has branch points inside concave regions, whether interior or exterior to
the curve.  One cannot see the whole Riemann surface of $S(z)$ from such a plot,
but one can guess where a reflection domain for each curve appears to lie.
Such guesses are supported by the second row of the figure, in which the
identity (\ref{schwarz-ident}) for the AAA approximation $r$ to $S$
is verified by plotting level curves of the discrepancy
\begin{equation}
|\;\overline{r( \overline{r(z)} )} - z\,|.
\label{ident}
\end{equation}
(We thank Keaton Burns for suggesting this type of plot.)
The dark-green region is bounded by the curve where this measure
is $10^{-8}$, i.e., 8 digits of accurate approximation to the
Schwarz function in a reflection domain, and the light green region
corresponds to $10^{-1}$, one digit of accuracy.

\begin{figure}
\begin{center}
\vspace{10pt}
\includegraphics[trim=8 143 0 5, clip, scale=.85]{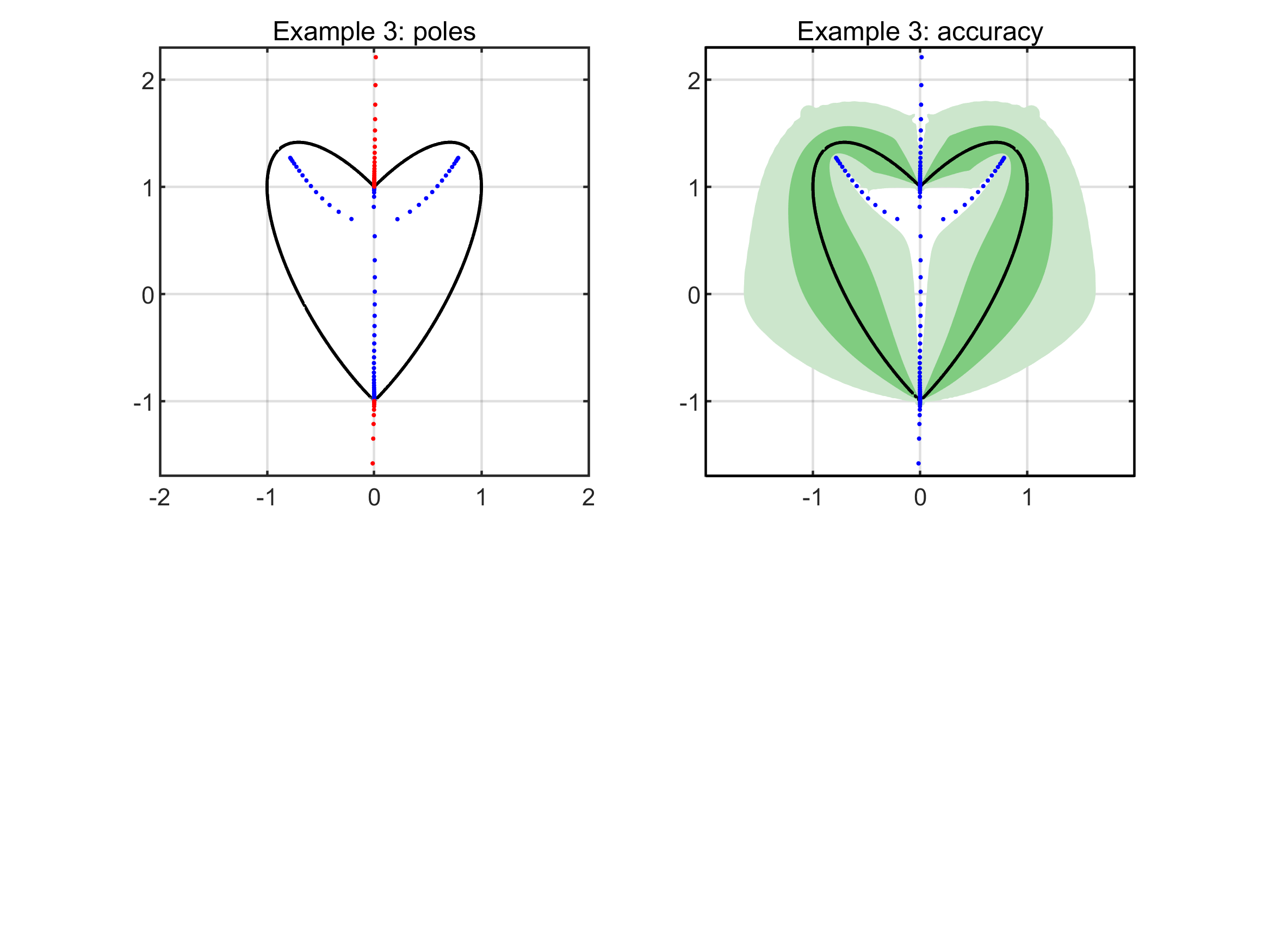}~~
\end{center}
\vskip -6pt
\caption{\label{schwarz-squiggles3}Similar plots for a curve $\Gamma$ consisting
of two pieces joined at branch points at $z= \pm\kern .4pt  i$.  Mathematically, we now have
a Schwarz function with unrelated branches associated with the left and right halves.}
\end{figure}

The Schwarz function is defined for open arcs $\Gamma$ as well as
closed ones.  In fact, if $\kern 1pt\Gamma$ is a closed analytic arc,
then in theory any open segment of it defines the full structure
of $S$ by analytic continuation, though this will be far from true
in numerical practice.  See \ccite{schwarzfun} for illustrations
of Schwarz functions of open arcs.

It is also interesting to consider what happens with non-smooth
contours, as in Figure \ref{schwarz-squiggles3}.  Here $\Gamma$
is a heart-shaped region consisting of $z+i\kern .7pt \Re(z)$ in
the right half-plane, where $z$ ranges over the right half of the
unit circle, together with its reflection in the left half-plane.
Mathematically, $\Gamma$ is here not a single analytic curve but two
unrelated components connected at a pair of branch points at $z =
\pm \kern .4pt i$.  Each of these pieces would have its own Schwarz
function, with poles presumably configured approximately like those
in the figure that are not along the imaginary axis.  The rest of
the poles in the figure, along the axis, are from this point of view
artifacts of constructing a global rational approximation to data
that consists of two unrelated branches (see sections \ref{sec:zolo}
and \ref{sec:quad}).  However, this strict interpretation of the
situation misses the deeper reality that the picture would look
much the same if, instead of breaks of analyticity at $z=\pm \kern
.4pt i$, $\Gamma$ merely had high curvature there.  For this reason
among others, computing a single Schwarz function for a contour
consisting of distinct branches often makes good sense in practice.
Examples will appear in sections \ref{sec:aaals}--\ref{sec:hilbert},
which concern applications of rational approximation to the solution
of \PDEs\ in the plane.  Here branch points (corners) and nonconvex
points (points of high curvature) arise frequently and are precisely
the situations where rational approximation shows its greatest power,
since poles can be placed in inlets nearby.  The Schwarz function
is a fundamental tool for explaining the success of these methods,
as discussed in \ccite{barnettb} and \ccite{inlets}.

\section{\label{sec:potential}Illustrating potential theory}

The theory of rational approximation has developed greatly since
Chebyshev first worked on such problems in the 19th century.
The central tool is the Hermite integral formula, whose roots go
back to Cauchy in the 1820\kern .4pt s (at \'Ecole Polytechnique)
and Hermite in 1878 (\'Ecole Polytechnique, University of Paris).
M\'eray (Dijon) and Runge (Hannover) showed that Hermite's integral
could be applied to bound the error of polynomial interpolants,
and beginning in the 1920\kern .4pt s, Walsh (Harvard) investigated
analogous results for rational interpolants.  This was the subject
of many papers by Walsh in the ensuing forty years and the central
topic of his magnum opus {\em Interpolation and Approximation by
Rational Functions in the Complex Domain\/}, which appeared in
five editions from 1935 to 1969 \ccite{walsh}.  Beginning in the
1960\kern .4pt s, Gonchar and Rakhmanov (Moscow) and Stahl (Berlin)
and others then laid the foundation for analysis of asymptotic
rates of best degree $n$ rational approximation as $n\to\infty$.
They showed that in this limit, the Hermite integrals lead to
equilibrium PDE problems of potential theory.

In its most familiar form as presented in Theorem 8.2 of \ccite{walsh},
the Hermite integral for a rational approximation $r\approx f$ is
\begin{equation}
f(z) - r(z) = {1\over 2\pi i}
\int_\Gamma {\phi(z)\over \phi(t)}\kern 1pt {f(t) \over t-z} \kern
1pt dt,
\label{potential-hermite}
\end{equation}
where $\phi(z)$ is a complex potential function, which
measures the product of the distances of a point $z$ from interpolation points
where $r=f$ relative to the product of its distances from the poles of $r$. 
The standard definition of $\phi$ is
\begin{equation}
\phi(z)
= \prod_{k=0}^n (z-\zk) \! \left/ \prod_{k=1}^n (z-\pk)\right.,
\label{potential-potential}
\end{equation}
where $\pi_1^{}, \dots, \pi_n^{}$ are the poles of $r$, assumed
at the outset to be simple and finite, and $\zz,\dots,\zn$ are
any $n+1$ interpolation points satisfying $r(\zk) = f(\zk)$.
The contour $\Gamma$ must lie in a domain of analyticity of $f$
and enclose the interpolation points and also the point $z$.
Normally we would choose $\Gamma$ to exclude the poles $\{\pk\}$
to get a good estimate from (\ref{potential-hermite}), but this
is not necessary for the validity of the equation.  Equations
(\ref{potential-hermite})--(\ref{potential-potential}) also apply
if some of the points are confluent (i.e., higher-order poles
or interpolation points), and if there are poles at $\infty$,
these can simply be dropped from the denominator product in
(\ref{potential-potential}).

What makes rational approximations so powerful, and
(\ref{potential-hermite}) so effective at revealing their power,
is that in many situations, $\phi$ is much smaller at points $z$
on a set $Z$ where $f$ is to be approximated than at points $t$ on
the contour $\Gamma$.  If $f$ is analytic in a neighborhood of $Z$,
so that $Z$ and $\Gamma$ can be taken to be disjoint, then suitable
choices of $\{\zk\}$ and $\{\pk\}$ will make $\phi(z)/\phi(t)$
exponentially small as $n\to \infty$, leading to exponential
convergence---and AAA reliably finds such choices.  If $f$ is entire
or meromorphic, or even if it has essential singularities, one is
guaranteed superexponential convergence, as mentioned at
the end of section~\ref{sec:approx}.
These are among the fundamental results that Gonchar, Stahl and
others established in the past 50 years.  By contrast, polynomials
can never converge superexponentially unless $f$ is entire.

AAA offers a beautiful method of illustrating the mathematics
of (\ref{potential-hermite})--(\ref{potential-potential}).
At step $n$, the algorithm has determined $n+1$ support points
(explicitly) and $n$ poles (implicitly).  These can be fed into
(\ref{potential-hermite}) and (\ref{potential-potential}) to give
plots of a kind first published in \ccite{notices}, which the figures
of this section are modeled on.  In all the experiments here, the
relative error tolerance has been set to $10^{-10}$ rather than
the usual $10^{-13}$, because we need a few extra digits to play
with for a reason to be explained on p.~\pageref{computepage}.

\begin{figure}
\begin{center}
\vspace*{8pt}
\includegraphics [trim=10 98 27 19, clip, width=4.1in]{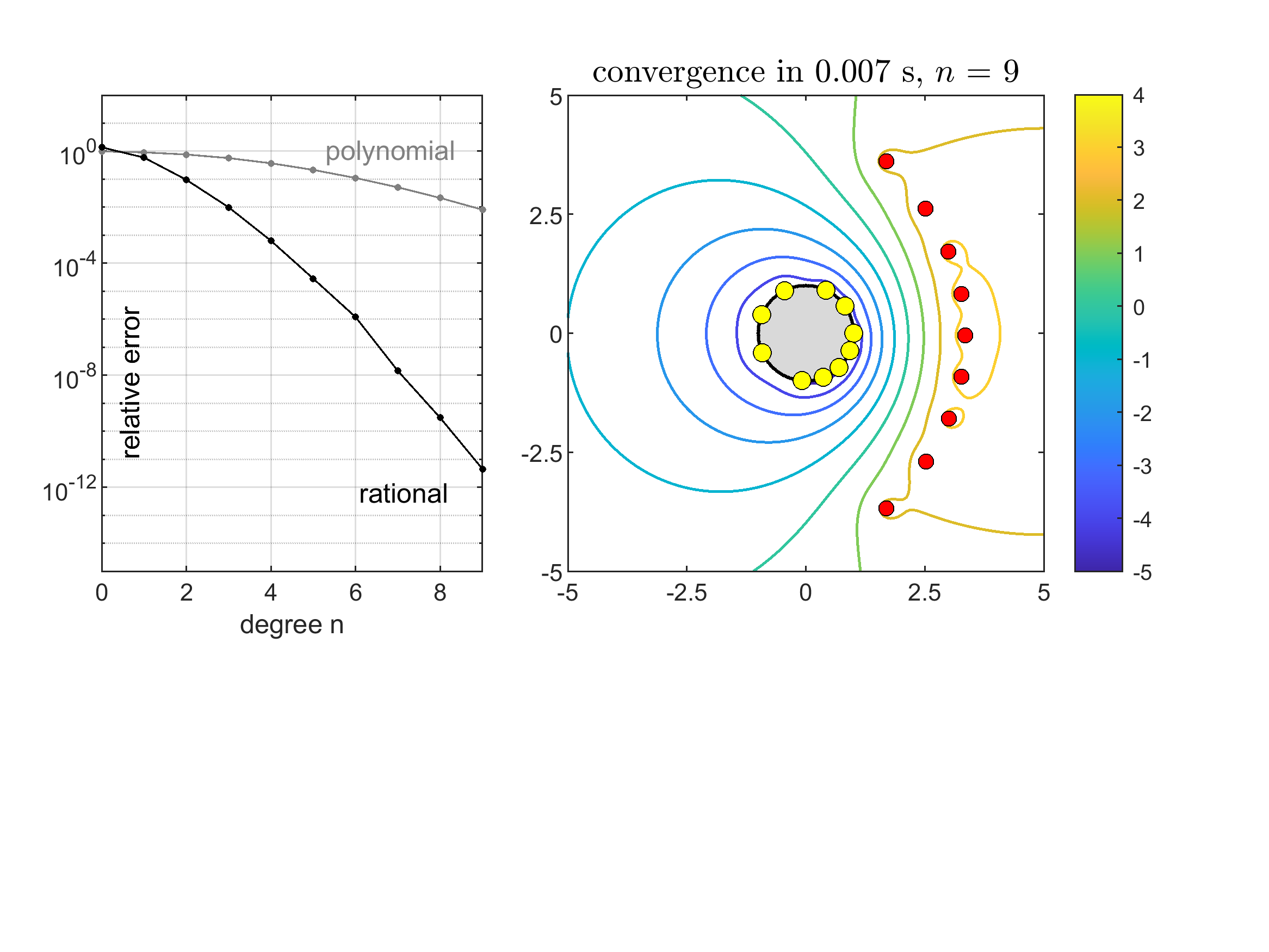}
\vspace*{-10pt}
\end{center}
\caption{\label{potential-aaashowold}Approximation of $f(z) = e^{4z}$
on the unit circle.  The plot on the left shows superexponential
convergence of AAA approximations 
as a function of degree $n$ down to tolerance $10^{-10}$ with $n=9$, with
slower superexponential convergence for polynomials.
On the right, level curves of the potential
function $(\ref{potential-potential})$ for $n=9$ based on the
$n+1$ AAA support points (yellow dots) and $n$ poles (red dots).  With
the Hermite integral $(\ref{potential-hermite})$, such plots explain at
least in part the high accuracy of rational approximants.}
\end{figure}

Consider Figure \ref{potential-aaashowold}, corresponding to
approximation of $e^{4z}$ in 100 roots of unity on the unit circle.
Accuracy $10^{-10}$ is achieved at degree $n=9$, as shown in
the convergence curve on the left.  The curve bends downward,
indicating superexponential convergence. The similar curve for
polynomial least-squares fitting on the same plot confirms that
for this entire function, polynomial approximants also converge
superexponentially, though at a lower rate.  In the plot on the
right, one sees the poles of the degree 9 rational approximant
as red dots in the complex plane, lining up along a curve whose
asymptotic shape and position as $n\to\infty$ were analysed by Saff
and Varga \ccite{saffvarga} (for Pad\'e approximation at the origin
rather than best approximation on the disk, but the asymptotics will
presumably be the same).  The 10 AAA support points are plotted as
yellow dots on the unit circle, somewhat denser on the right side
than the left.  Between the poles and the interpolation points,
level curves of $|\phi(z)|$ are plotted, with the numbers indicated
on the colourbar corresponding to base 10 logarithms of $|\phi(z)|$.
One sees that $|\phi(z)|$ decreases exponentially between the
poles and the unit disk, implying by (\ref{potential-hermite})
that $|f(z)-r(z)|$ will be extremely small for any $z$ there.

Images like Figure \ref{potential-aaashowold}, which are easily
produced after any AAA computation, explain how a rational
approximation achieves its accuracy.  It is remarkable how AAA
has detected all this structure in the complex plane, though it
is just a greedy descent algorithm working with function values
at the sample points.  Anyone familiar with potential theory will
recognise the flavor of these plots.  If we think of the yellow dots
as positive point charges generating potentials $\log|z-\zk|$ and the
red dots as negative point charges with potentials $-\log|z-\pk|$,
then the level curves are associated equipotential lines, and
it is evident that AAA has distributed the points so that the
interpolation points lie approximately on a low equipotential
curve and the poles on a high one.  In other words, the charges
lie in an approximately minimal-energy configuration.  This is the
effect exploited by Gonchar, Stahl and others in their analysis of
convergence of approximations as $n\to\infty$.  They showed that
the solution of a PDE problem involving continuous distributions
of signed charges provides rigorous upper bounds on best rational
approximation convergence rates in this limit.

All this may seem like a complete story, but it isn't.  If we look
at the numbers, we see that Figure \ref{potential-aaashowold} does
not give a quantitative explanation of the ten-digit accuracy of
this approximation.  The colourbar spans just about eight orders of
magnitude, not ten, and moreover, digits will be lost in an estimate
of $(\ref{potential-hermite})$ since this function $f(z) = e^{4z}$
has magnitude as large as $10^4$ near the poles.  For these reasons,
applying (\ref{potential-hermite})--(\ref{potential-potential})
based on what we see in Figure~\ref{potential-aaashowold} only
explains about five of the ten digits of accuracy.  From this point
of view one would have to explain the other five as a result of
cancellation in the integrand of (\ref{potential-hermite}) along a
contour $\Gamma$.  But why should there be five orders of magnitude
of cancellation?  It seems strangely fortuitous.

\begin{figure}
\begin{center}
\vspace*{8pt}

\includegraphics [trim=170 98 27 19, clip, width=2.35in]{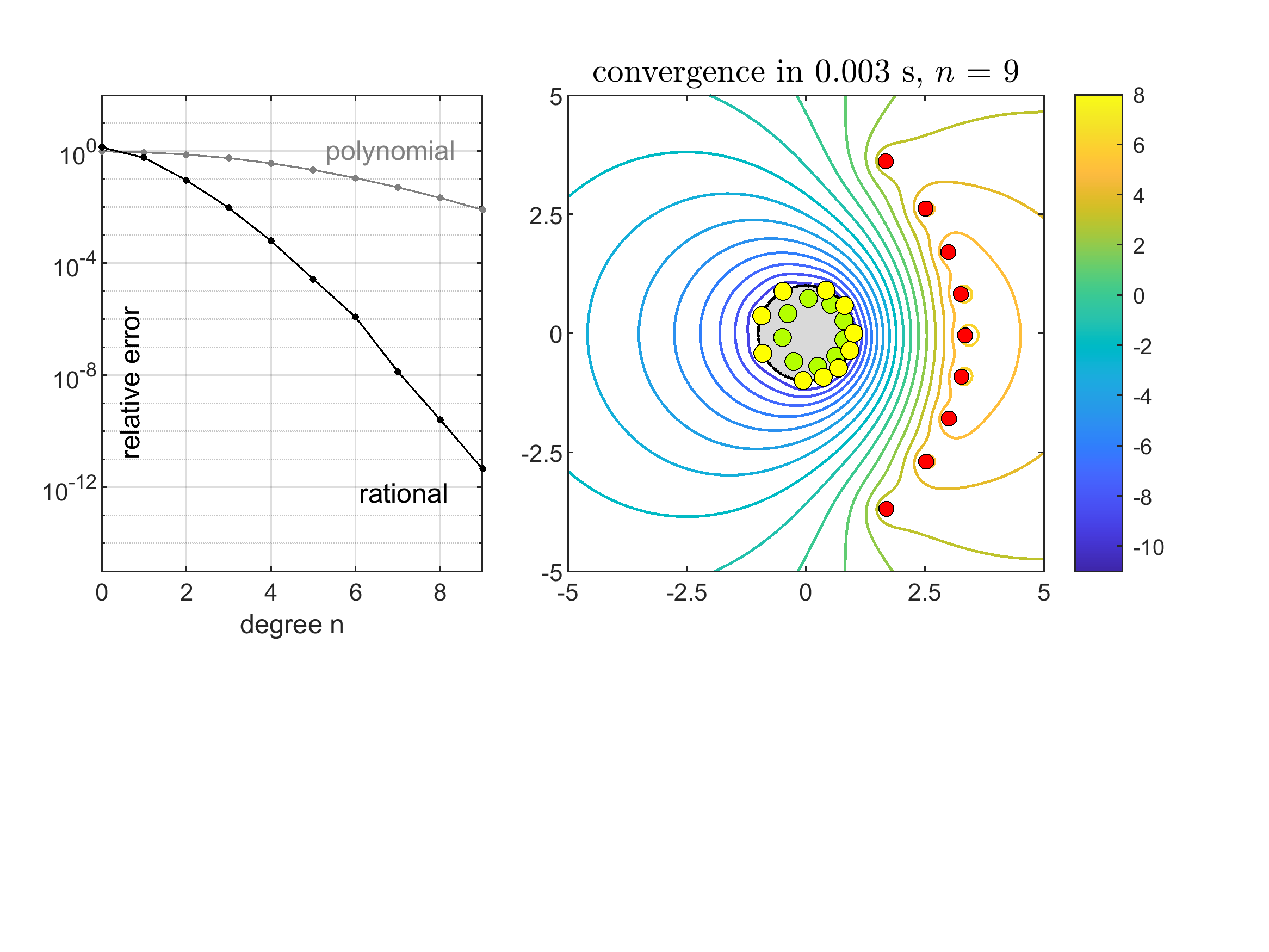}
\vspace*{-15pt}

\end{center}
\caption{\label{potential-aaashownew}Approximation of $f(z) = e^{4z}$ on
the unit circle again.  Now $n$ more interpolation points have been added to the image,
shown as green dots, and the potential function has been approximately squared
from $(\ref{potential-potential})$ to
$(\ref{potential-potential2})$.  This roughly doubles the span of amplitudes shown
in the colourbar, enabling $(\ref{potential-hermite})$ to fully explain the $10$-digit
accuracy of this approximation.}
\end{figure}

Happily, there is a modification of
(\ref{potential-hermite})--(\ref{potential-potential}) that changes
the figure in such a way as to yield
a wider range of contour levels and thereby
a quantitative as well as qualitative explanation.
We mentioned that the integral (\ref{potential-hermite})
is valid for any choice of $n+1$ interpolation points $\{\zk\}$,
but in many
best and near-best approximations, the total number of such points
in the approximation region of interest
will actually be $2\kern .3pt n+1$.  (One familiar case of this is in real
approximation on an interval, with equioscillation generically between
$2\kern .3pt n+2$ extrema.  Another is in degree $n$ Pad\'e approximation, where
the approximant interpolates the function generically to order $2\kern .3pt n+1$.)
In this case, (\ref{potential-hermite}) remains valid with a new
and stronger definition of $\phi$ based on all the interpolation
points $\zz,\dots,\znn$,
\begin{equation}
\phi(z)
= \prod_{k=0}^{2n} (z-\zk) \! \left/ \prod_{k=1}^n (z-\pk)^2\right..
\label{potential-potential2}
\end{equation}
Roughly speaking, both the numerator and the denominator
of (\ref{potential-potential}) have been squared!
Figure \ref{potential-aaashowold} improves to 
Figure \ref{potential-aaashownew}, with green dots now added to
show the $n=9$ additional points of interpolation of this approximation
in the unit disk.  Notice how much denser the contour lines
are than before.  Now the contours span a range of 14 orders
of magnitude, enough to explain
how the approximation achieved 10 digits of accuracy, not just 5.

We learned of the improvement
(\ref{potential-potential})$\rightarrow$(\ref{potential-potential2})
from Maxim Yattselev, whose assistance we gratefully acknowledge
(private communication, May 2025).  We state the result as a theorem,
a special case of \cite[Lemma 2]{stahl89}.  This result has roots
in earlier work by Gonchar for approximation of Markov functions
\ccite{gonchar75,goncharlopez}.

\begin{theorem}
{\bf Hermite integral formula with \boldmath$2\kern .3pt n+1$ interpolation points.}
Let $f$ be analytic in the closure of a
Jordan region $\Omega$ bounded by a Jordan contour $\Gamma$, and let 
the degree $n$ rational function $r(z)$ interpolate $f$ in $2\kern .3pt n+1$ points
$\zz,\dots,\znn\in \Omega$.   Let $r$ have $n$ finite poles
$\pi_1^{},\dots,\pi_n^{}$, and define
\begin{equation}
\phi(z)
= \prod_{k=0}^{2n} (z-\zk) \! \left/ \prod_{k=1}^n (z-\pk)^2\right..
\label{potential-potential3}
\end{equation}
Then for any $z\in\Omega$, 
\begin{equation}
f(z) - r(z) = {1\over 2\pi i}
\int_\Gamma {\phi(z)\over \phi(t)}\kern 1pt {f(t) \over t-z} \kern
1pt dt.
\label{potential-hermite2}
\end{equation}
\end{theorem}

The poles and interpolation points need not be simple.
For simplicity, we have limited the theorem statement to the generic case of
greatest interest, but it can be generalised to poles
at $\infty$ and numbers of interpolation points greater than or less than
$2\kern .3pt n+1$.
In fact, the denominator of
(\ref{potential-potential3}) does not have to be the square of that
of (\ref{potential-potential});
it can be the denominator of
(\ref{potential-potential}) multiplied by
any polynomial in $z$ of degree $n$.

This ``factor of 2'' improvement from (\ref{potential-potential})
to (\ref{potential-potential2}) was investigated by Stahl in the
paper just cited and others and is the subject of what Rakhmanov
calls the {\em Gonchar-Stahl $\rho^2$-theorem\/}; see \ccite{rakh}
and \cite[sec.\ 7.5]{jjiam}.  Not all best or near-best rational
approximations have the necessary doubling of the number of
interpolation points for this to apply, but in practice, many do.

\begin{figure}[p]
\begin{center}
\includegraphics [trim=10 98 27 19, clip, width=3.6in]{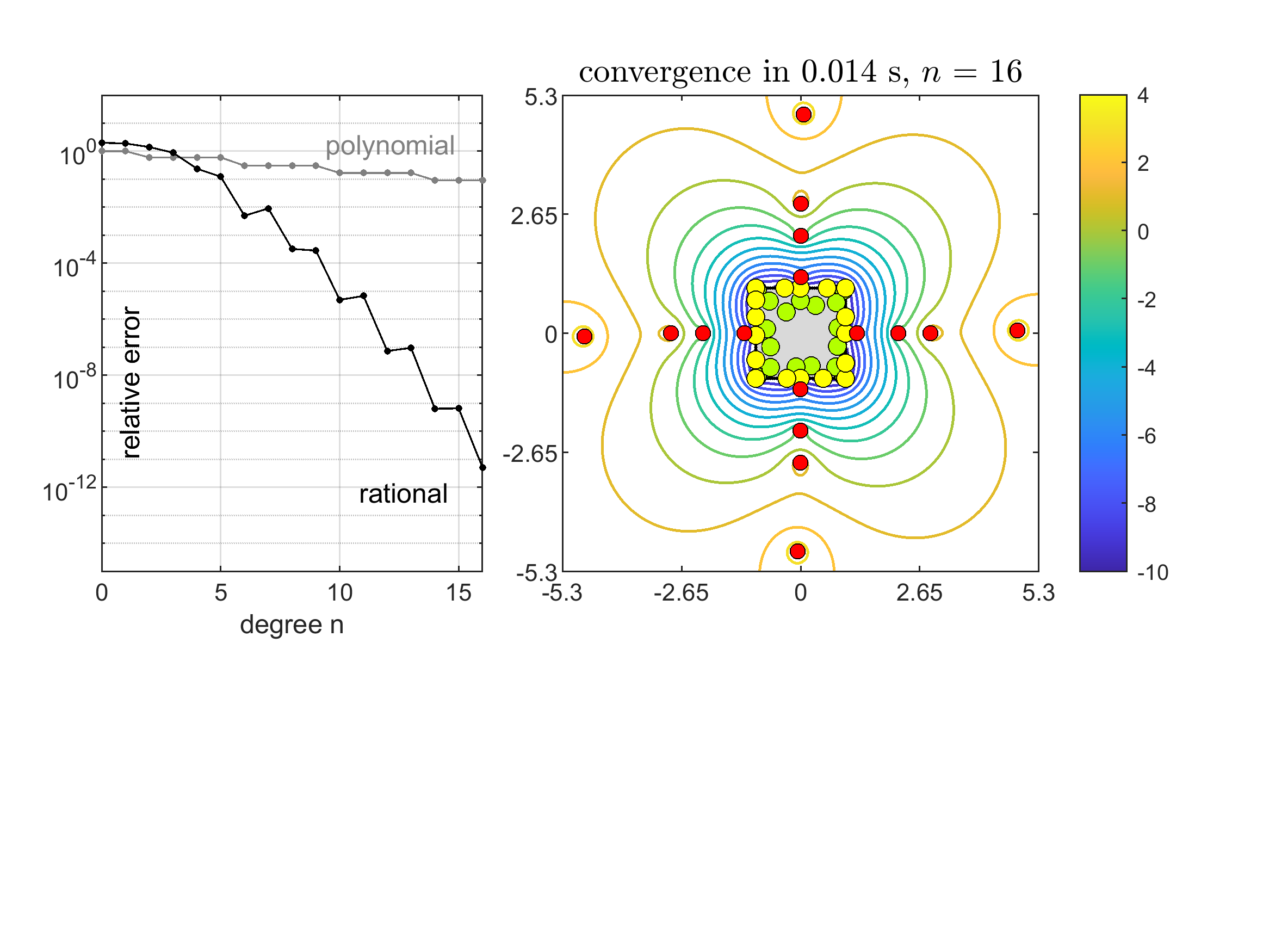}
\vspace{-13pt}
\end{center}
\caption{\label{potential-pot2}Approximation of $f(z) = \tan(z^2)$ on 
$[-1,1]^2$.  For a meromorphic function like this, polynomial approximations
converge exponentially whereas rational approximations converge superexponentially.}
\end{figure}

\begin{figure}[p]
\begin{center}
\includegraphics [trim=10 98 27 19, clip, width=3.6in]{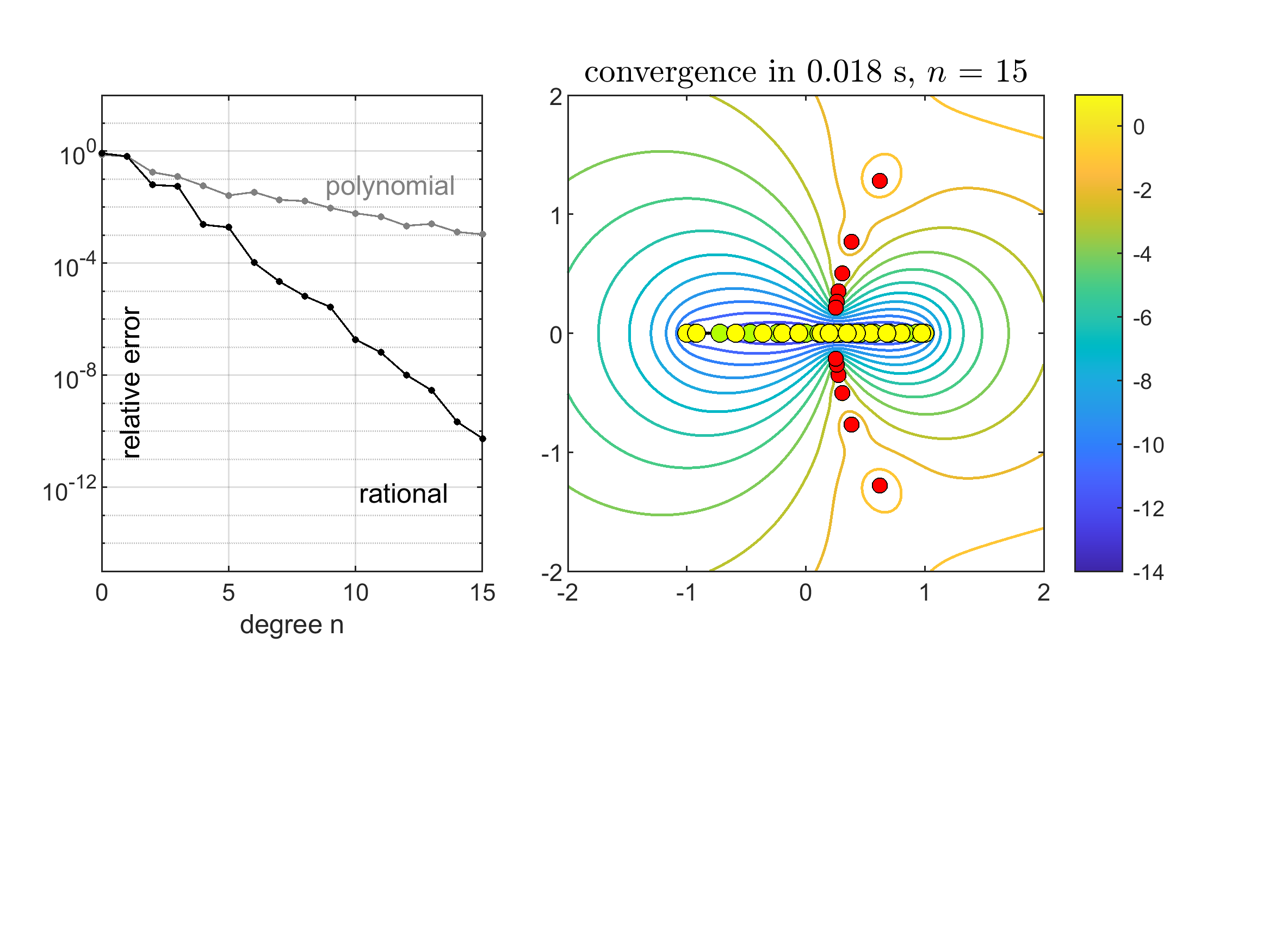}
\vspace{-13pt}
\end{center}
\caption{\label{potential-pot3}Approximation of $f(z) = 
\sqrt{1+25(z-0.25)^2}$ on $[-1,1]$.  For any function like this with branch points,
both polynomial and rational approximations converge exponentially.}
\end{figure}

\begin{figure}[p]
\begin{center}
\includegraphics [trim=10 98 27 19, clip, width=3.6in]{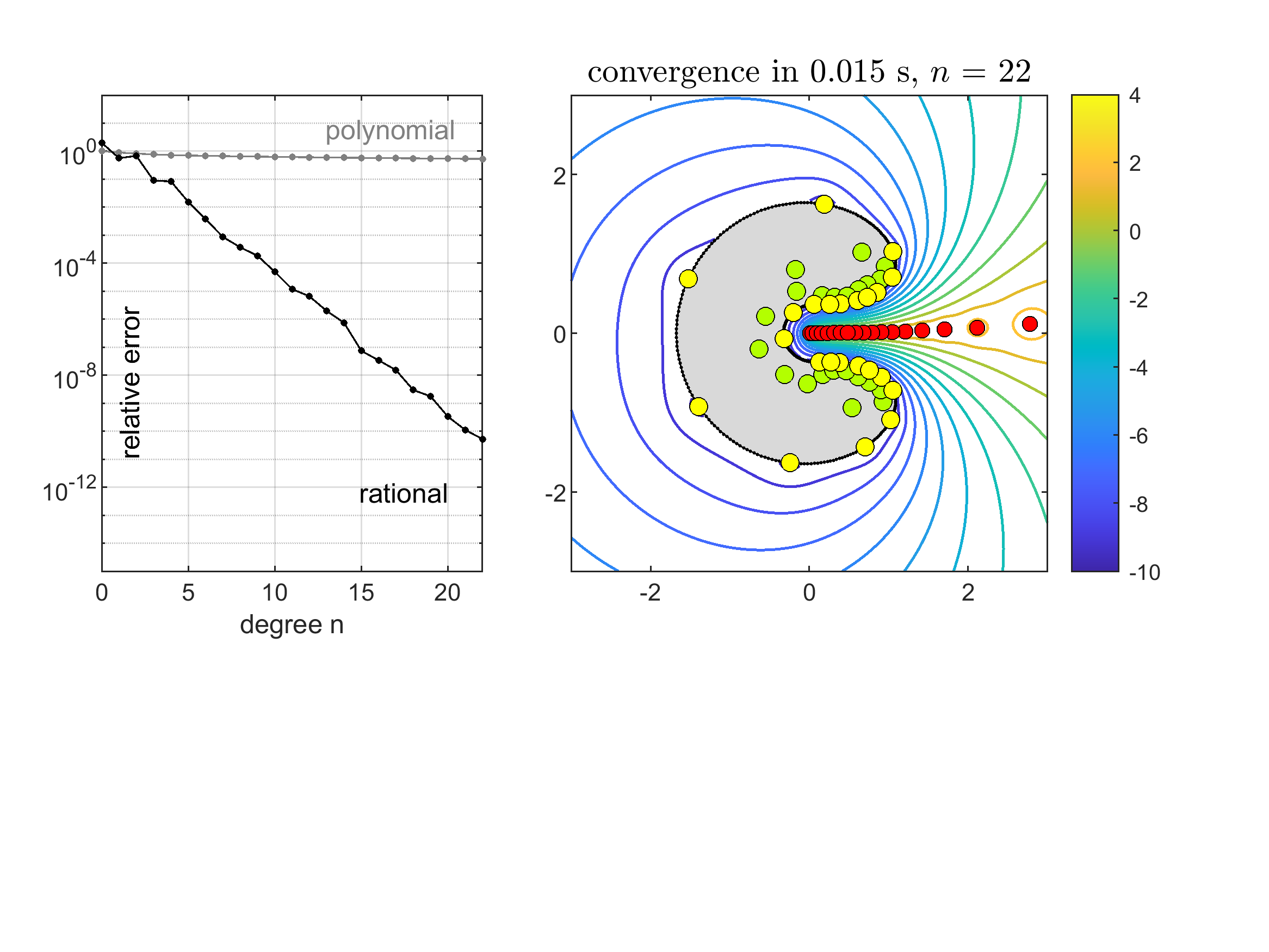}
\vspace{-13pt}
\end{center}
\caption{\label{potential-pot4}Approximation of
$f(z) = \log(-z)$ on a nonconvex region.  Both polynomial and rational
functions converge exponentially, but at vastly different rates
{\rm \ccite{inlets}}.}
\end{figure}

\label{computepage}
Producing a plot like Figure~\ref{potential-aaashowold} is
elementary once one has calculated the support points and poles
of a AAA approximation, but how did we find the additional
interpolation points---the green dots---needed to produce
Figure~\ref{potential-aaashownew}?  Following the zerofinding ideas
of section~\ref{sec:poles}, the method we have used is to
invoke AAA a second time to approximate the function $f(z) - r(z)$
by a new rational function $q(z)$, and take the zeros of $q$ as
numerically computed interpolation points of $f$ and $r$.  This is
where the tolerance $10^{-10}$ mentioned above comes into play.
We compute $r$ with that tolerance, and then $q$ with a tolerance
of just $10^{-3}$ on the theory that $f-r$ will have low relative
accuracy since it is smaller by a factor of $10^{10}$ than $f$
and $r$ individually.  This combination seems to reliably give us
the desired plots, although we have not investigated the tolerances
carefully.

Figures \ref{potential-pot2}--\ref{potential-pot4} show three
more examples, with some associated convergence rate results for
various classes of approximations stated in the captions; recall
the last paragraph of section \ref{sec:approx}.  An indication
of the implications of such results for rational vs.\ polynomial
approximation degrees can be seen in Table~\ref{speedtable}.

We have looked at hundreds of plots of rational approximations
in the past few years, perhaps thousands, and we note that it
often seems possible to classify the poles loosely into four
categories.  Sometimes, as in Figure~\ref{potential-aaashownew},
poles appear ``just'' for reasons of approximation---we
might call these {\em approximation poles\/}.  Sometimes,
as in Figure~\ref{potential-pot2}, they match poles of the
function being approximated---{\em pole poles\/}.  Sometimes,
as in Figures~\ref{potential-pot3} and~\ref{potential-pot4}, they
delineate branch cuts---{\em branch cut poles\/}.  Finally there
are {\em spurious poles\/} ({\em Froissart poles?\kern 1pt})
with negligible residues, often a result of rounding errors or
other noise \ccite{stahl98}.  We have not attempted to make such
distinctions precise, and of course, the boundaries between the
cases may be fuzzy.

\begin{table}
\label{speedtable}
\caption{Degrees of the rational approximations in
the examples of Figures 
\ref{potential-aaashownew}--\ref{potential-pot4}, and for
comparison, rough estimates of the degrees that would be needed for
polynomial approximations of the same accuracy.}
\begin{center}
\begin{tabular}{ccc}
Figure & Rational degree & Polynomial degree \\
\hline \\ [-8pt]
\ref{potential-aaashownew} & $7$ & $80$ \\
\ref{potential-pot2} & $16$ & $200$ \\
\ref{potential-pot3} & $15$ & $120$ \\
\ref{potential-pot4} & $22$ & $10{,}000$ \\
\end{tabular}
\end{center}
\end{table}

\section{\label{sec:solution}Solution of \ODEs\ and \PDEs}

We mentioned at the beginning that polynomials are the starting
point for algorithms in most areas of numerical analysis, but that
rational functions may bring advantages, and the previous sections
have shown examples in many areas.  Now we turn to the solution of
\ODEs\ and \PDEs.  Our view is that here, too, rational functions
should lead to new methods that are more powerful for certain
problems, notably those with singularities or near-singularities.
For example, Chebfun solves smooth ODE boundary-value problems by
means of polynomial and piecewise polynomial spectral approximations.
We are often asked, when is ``Ratfun'' coming along to
solve more difficult problems?

So far, there is no Ratfun.  Chebfun-style automated rational
approximation methods for \ODEs\ and \PDEs\ have not been developed,
but we hope that they will be in the next few years.  In this
section we outline what we think may be possible.

\begin{figure}
\vskip 7pt
\begin{center}
\includegraphics[trim=24 155 24 0, clip, scale=.46]{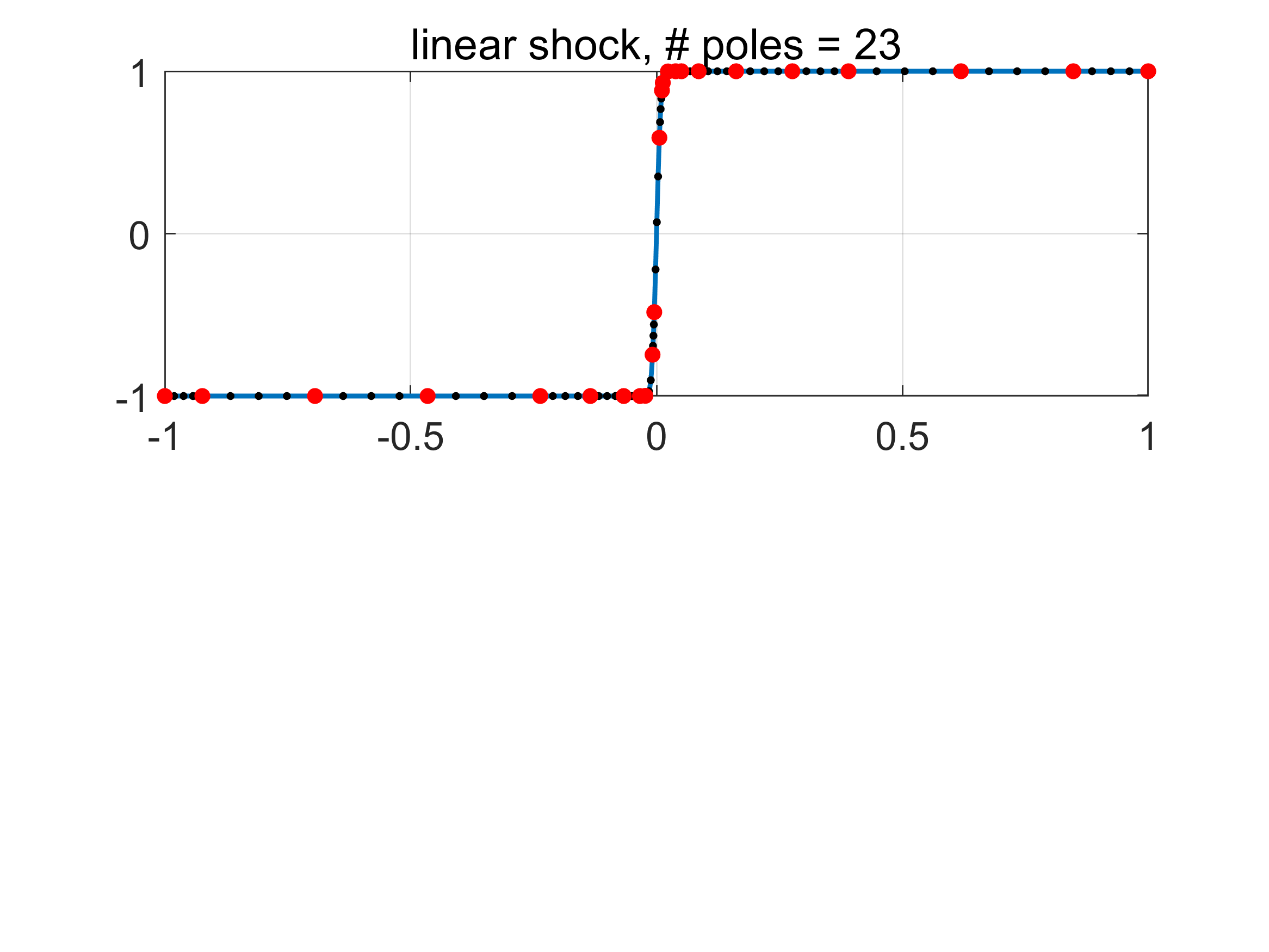}
\includegraphics[trim=24 155 24 0, clip, scale=.46]{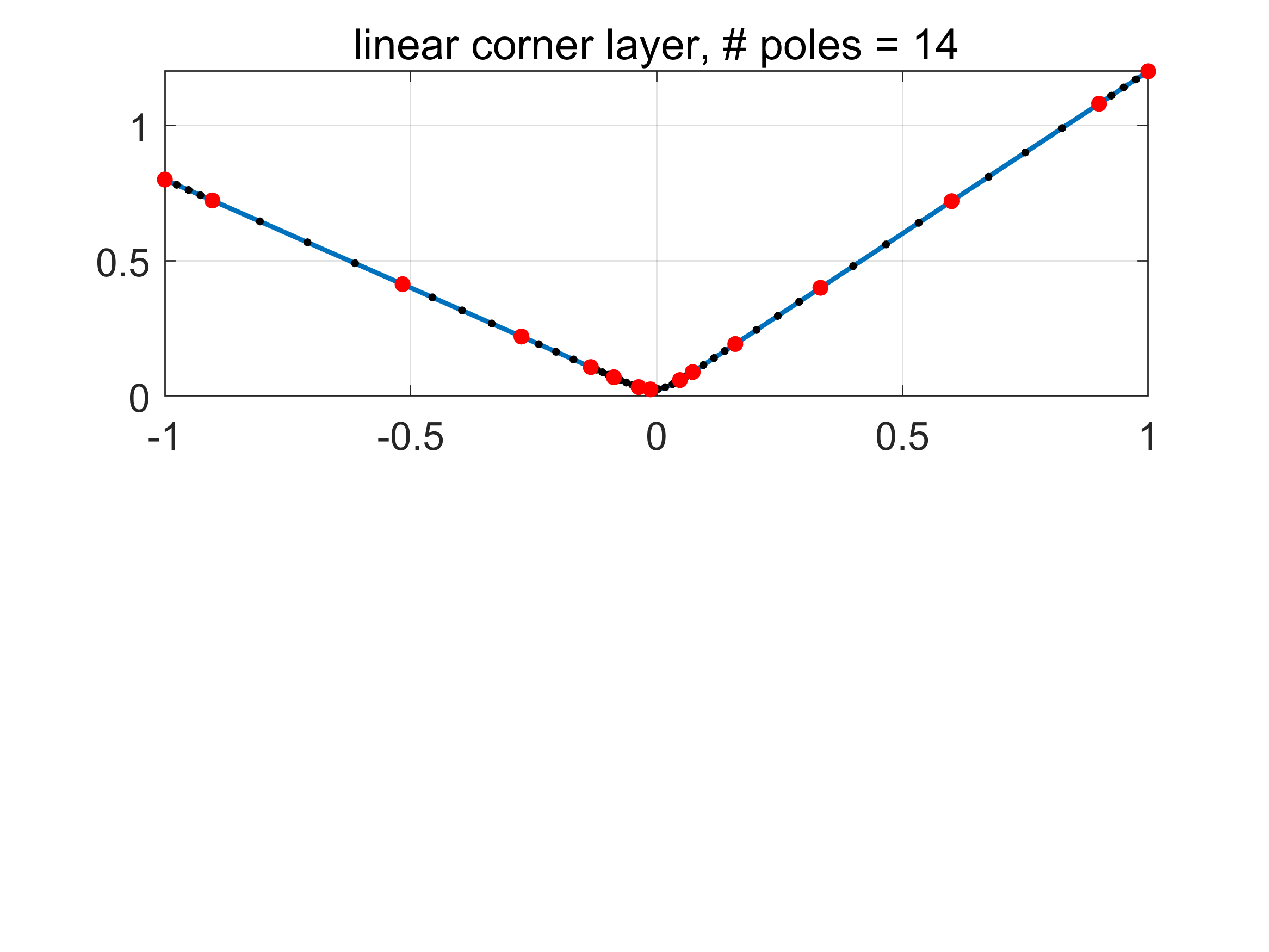}
\includegraphics[trim=24 155 24 0, clip, scale=.46]{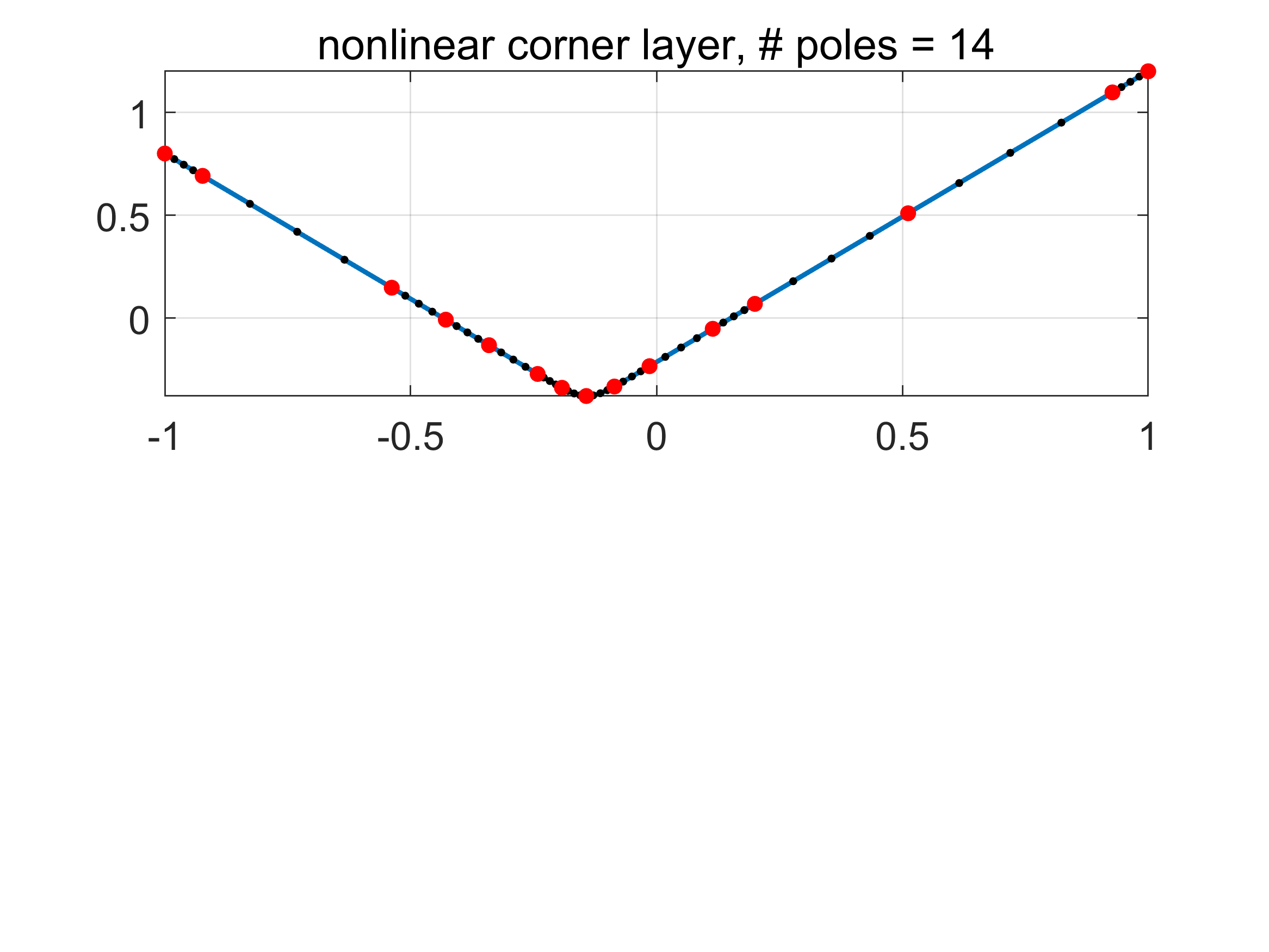}
\includegraphics[trim=24 155 24 0, clip, scale=.46]{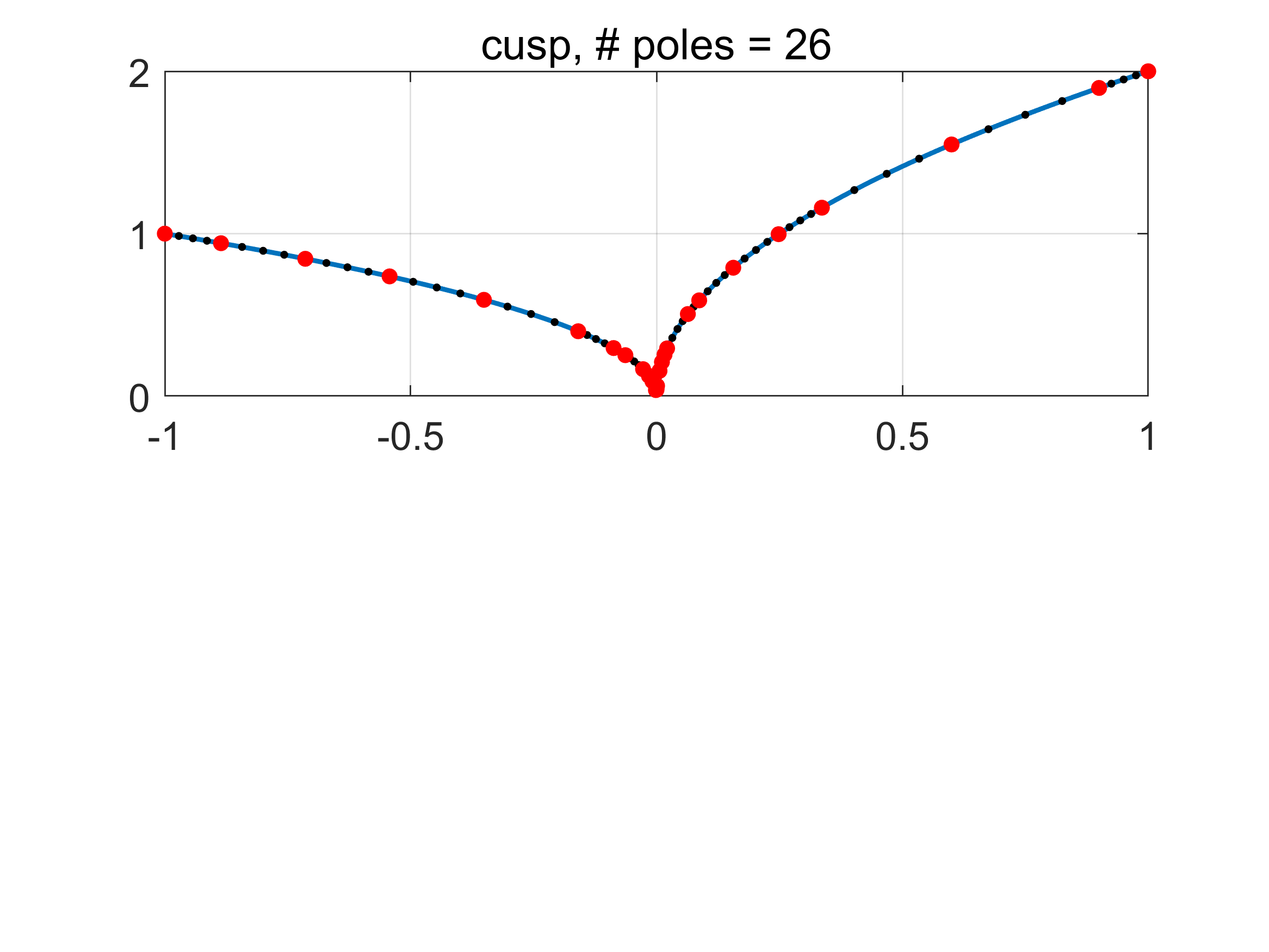}
\includegraphics[trim=24 155 24 0, clip, scale=.46]{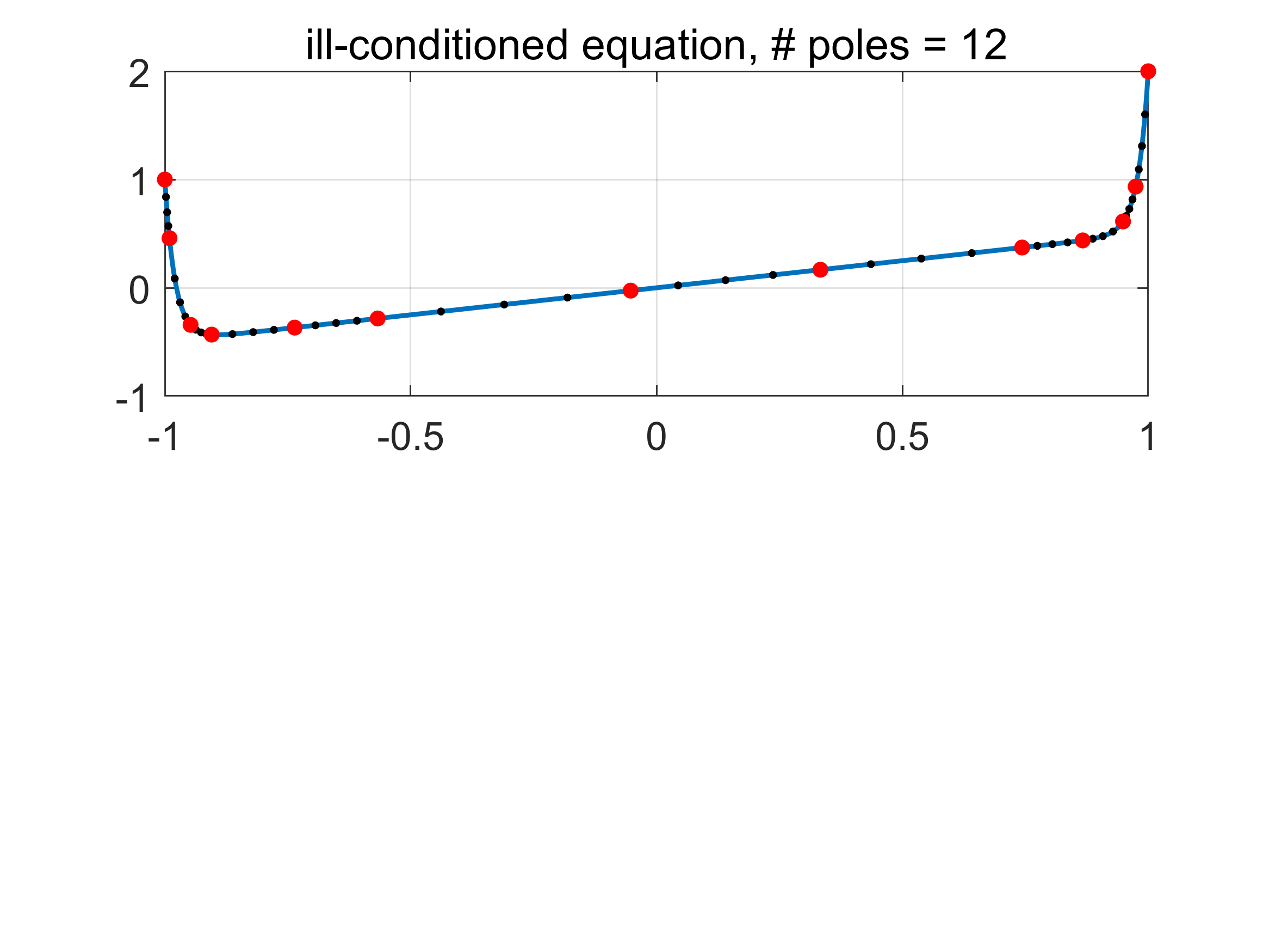}
\includegraphics[trim=24 155 24 0, clip, scale=.46]{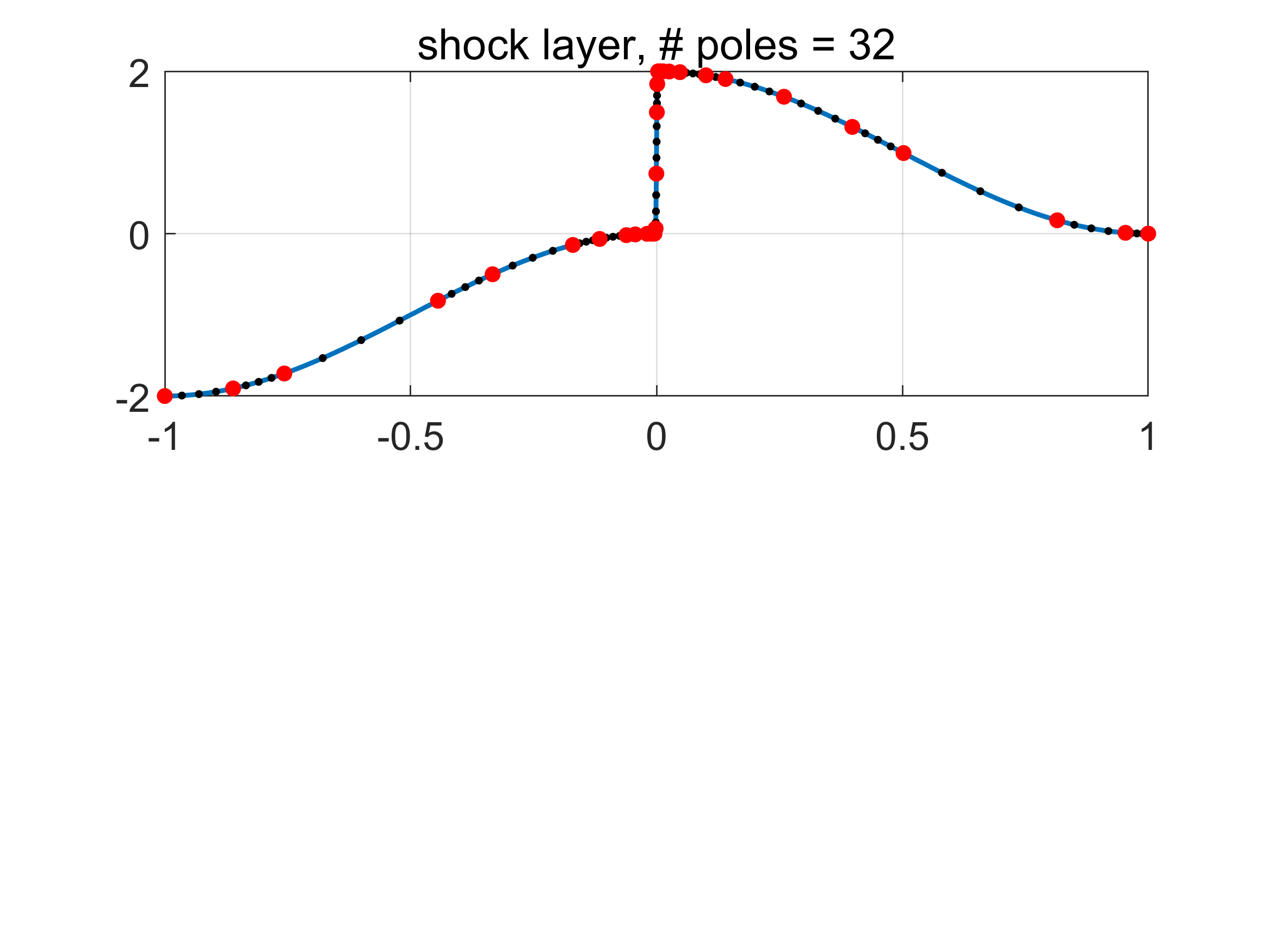}
\end{center}
\vskip -7pt
\caption{\label{sixfigs} Solutions to the six singularly
perturbed ODE boundary-value problems
from Examples 91--96 of Appendix~B of \ccite{ExplODE} approximated
to 6-digit accuracy by AAA rational approximations computed
by the continuum AAA algorithm of section~\ref{sec:continuum}.  Red dots mark
adaptively determined support points, and black dots mark
additional sample points.  Despite the narrow layers,
the fits are of low degrees, suggesting a strong potential for efficient rational function
solution of difficult \ODEs\ if a suitable algorithm can be developed.}
\end{figure}

To start with \ODEs, the natural terrain of rational approximation
methods will be problems with singularities or near-singularities,
such as the stiff boundary-value problems treated in \ccite{hemker}
and \ccite{leegreengard}.  Three examples from each of these sources
can be found as Examples 91--96 of Appendix~B of \ccite{ExplODE}.
The solutions to these six \BVPs, all involving boundary or interior
layers that can be made more challenging by decreasing further a
small parameter, are plotted in Figure~\ref{sixfigs}.

The computations for Figure \ref{sixfigs} are as follows.  First the
\BVPs\ are solved in Chebfun with these commands, which the reader
can decode to see the precise specifications of the \BVPs:

{\small
\begin{verbatim}
      % Example 91: Linear shock
      L = chebop([-1 0 1]);
      L.op = @(x,y) 1e-4*diff(y,2)+2*x*diff(y);
      L.lbc =-1; L.rbc = 1;
      y = L\0;

      % Example 92: Linear corner layer
      L = chebop(@(x,y) 1e-3*diff(y,2)+x*diff(y)-y);
      L.lbc = .8; L.rbc = 1.2;
      y = L\0;

      % Example 93: Nonlinear corner layer
      N = chebop(@(y) 0.05*diff(y,2)+diff(y)^2-1);
      N.lbc = .8; N.rbc = 1.2;
      y = N\1;
      
      % Example 94: Cusp
      L = chebop([-1 0 1]);
      L.op = @(x,y) 1e-6*diff(y,2)+x*diff(y)-.5*y;
      L.lbc = 1; L.rbc = 2;
      y = L\0;

      % Example 95: Ill-conditioned equation
      L = chebop(-1,1);
      L.op = @(x,y) 0.02*diff(y,2)-x*diff(y)+y;
      L.lbc = 1; L.rbc = 2;
      y = L\0;

      % Example 96: Shock layer
      L = chebop([-1 0 1]);
      L.op = @(x,y) 1e-6*diff(y,2)+x*diff(y);
      L.lbc =-2; L.rbc = 0;
      f = chebfun(@(x)-1e-6*pi^2*cos(pi*x)-pi*x*sin(pi*x));
      y = L\f; 
\end{verbatim}
\par}

\noindent Then each solution is approximated by a rational
function using the continuum AAA code {\tt aaax} to be
described in section \ref{sec:continuum} by executing

{\small
\begin{verbatim}
      deg = 150; nlawson = 0; tol = 1e-6; plt = 0;
      r = aaax(y,deg,nlawson,tol,plt);
\end{verbatim}
\par}

\noindent
Finally, the approximation is plotted with

{\small
\begin{verbatim}
      xx = linspace(-1,1,2000)';
      plot(xx,r(xx))
\end{verbatim}
\par}

\noindent The support points of the rational approximation are marked
by red dots, and the additional sample points, three between each
pair of support points, are marked by small black dots.

These computations show that rational functions computed by AAA can
represent solutions to certain difficult ODE \BVPs\ very efficiently.
The question is, how might AAA be used to actually compute such
solutions, rather than just represent them after they have been
computed by another tool?  Our vision for an algorithm of this kind
is inspired by an analogy with Chebfun.  In its early years, Chebfun
just approximated given functions by polynomials using the following
procedure: interpolate the function on a 17-point Chebyshev grid,
then a 33-point grid, then a 65-point grid, and so on until the
Chebyshev series of the interpolant goes down to machine precision.
Then in \ccite{bornemann} a method was introduced for generalising
this to solve ODE \BVPs.  Proceed just as before, but instead of {\em
sampling\/} a fixed function on each grid, {\em compute\/} a function
on each grid by a Chebyshev spectral method.  As before, stop when
the grid is fine enough that the solution is well resolved and one is
happy with the resulting Chebyshev series.  This is the basis of all
Chebfun solution of ODE \BVPs.  See Appendix A of \ccite{ExplODE}.

We suggest that a similar approach might work for solving ODE \BVPs\
such as those of Figure~\ref{sixfigs} by rational approximation.
Instead of Chebyshev grids of growing sizes, one might make use
of adaptively growing grids of support points and sample points
following the pattern of the continuum AAA algorithm to be described
in section \ref{sec:continuum}.  Perhaps such a method could be
based on finding a way to linearise the barycentric representation
at each step so that an appropriate residual could be minimised
in a linear least-squares problem.  As happened with Chebfun,
such an algorithm might be developed first for linear problems
and then extended to nonlinear ones via iteration.  At the moment,
all this is entirely speculative.

Rational functions have been advocated in the past for solving \ODEs,
for example in \ccite{mittelmann,mittelmann2}.  What is potentially
new now is the availability of the AAA algorithm to make the rational
approximations much more flexible and closer to optimal.

Moving from \ODEs\ to \PDEs, a setting in which univariate rational
approximations might be effective is in problems in one space
and one time variable with shocks, blowups, or other singular or
nearly-singular behaviour.  Previously, rational function methods
for such problems have been developed based on various strategies
(\cite{bberrut,berrutb}).  In \ccite{tee} and \ccite{haletee},
conformal maps were used to push singularities in the complex
plane further from the computational axis, yielding solutions to
difficult problems on grids with just a few points.  In this method,
positions of singularities are estimated numerically by Pad\'e or
Chebyshev-Pad\'e approximation, and then a highly nonuniform grid is
established by means of an explicit conformal map, leading to very
accurate spectral discretisation despite the singularities.  The new
vision is to find a way to achieve all this simply by rational
approximation, without the need for conformal maps.  In terms of
Ratfun, we are looking for an analogue of Chebfun {\tt pde15s}
as opposed to Chebfun backslash.  Again our view is that something
should be possible, but algorithms have not yet been developed.

\section{\label{sec:extending}Extending \ODEs\ and \PDEs\ into
\boldmath$\complex$} Often an ODE or PDE or related problem
is posed in real variables, yet for one reason or another,
there is interest in determining how the solution behaves for
complex arguments.  In particular, many authors have investigated
complex singularities of problems which may be singularity-free
for real arguments.  This kind of analysis has been pursued
for the Lorenz equations \ccite{taborweiss,vis}, the Burgers
equation \ccite{bessis,ssf,vanden,weideman03,weideman22}, the
nonlinear Schr\"odinger equation \ccite{ssf}, the Korteweg-de
Vries equation (\cite{kruskal,weideman03,weideman22,mccue}),
water waves \ccite{bakerxie,crew,lush,lustri}, vortex sheets and
layers \ccite{cgss,cafl}, Rayleigh-Taylor flow \ccite{tanveerRT},
Hele-Shaw flow \ccite{tanveerHS} and the Euler and Navier-Stokes
equations \ccite{siegel,cgss15}.

The direct approach to such analysis is simply to extend whatever
analytic or numerical method is being used for real variables to
complex variables.  In this fashion one can explore solutions of
problems in complex space and/or time with, typically, comparable
accuracy to what one has in the real case.  For obtaining particular
values at particular points, this can be an excellent strategy.
The drawback is that for a more global picture, for example if one
simply wants to know where singularities may be located, the work
involved may be considerable.

This leads to the alternative idea: use numerical analytic
continuation to move from real to complex variables.  With AAA
approximation this can be remarkably speedy and informative,
even if it cannot match the precise information of a true complex
calculation.  Numerical analytic continuation has also been used in
this manner for computational as opposed to analytical purposes,
to find nearby singularities of a numerically computed function
in order to construct a function-adapted computational grid
\ccite{tee,weideman03}.

There are just a handful of publications so far that have used
AAA for extending differential equations into the complex plane:
see \ccite{deng,jjiam,vanden,weideman22}.  More will presumably
appear in the years to come, as some of the older literature
relies on less effective methods such as Fourier series.  More
recently, some authors have used rational functions in the form of
Fourier-Pad\'e approximations \ccite{bgm} (sometimes called just
Pad\'e approximations---see section \ref{sec:trig}), and this can
be quite effective \ccite{cgss15,weideman03,weideman22}.

We will show two examples of AAA computations, both adapted from
\ccite{jjiam}, which differ mathematically in an interesting way.

\begin{figure}
\begin{center}
\includegraphics [trim=35 95 35 0, clip, width=4.2in]{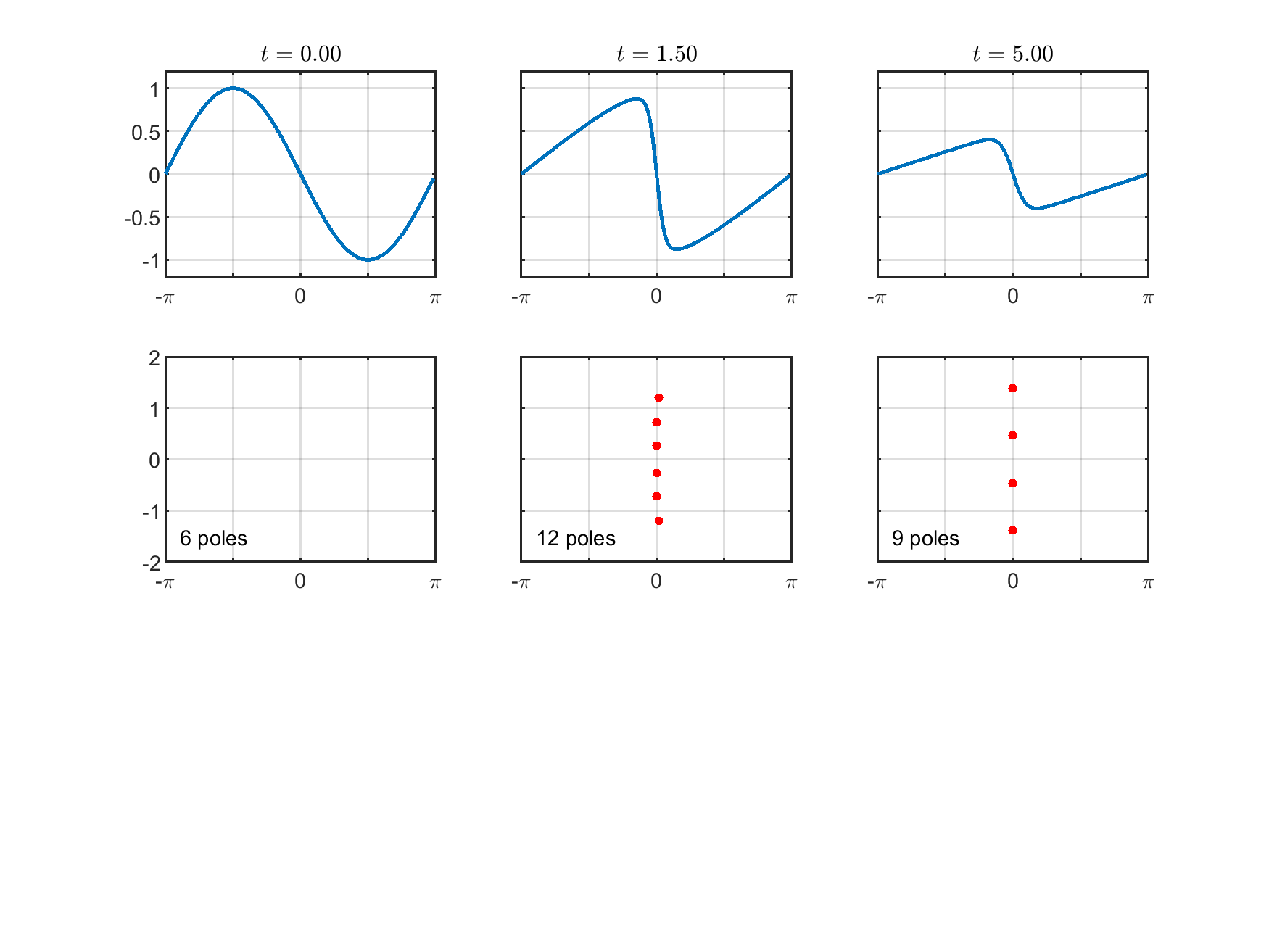}
\vspace{-25pt}
\end{center}
\caption{\label{DEs-burgersfig}Solutions to the periodic
Burgers equation (\ref{DEs-burgers}) at three times $t$
together with poles of associated AAA approximations in the
complex $x$-plane.  AAA is fast enough that it is feasible to compute
poles this way at every time step of a simulation.}
\end{figure}

The first example, illustrated in Figure \ref{DEs-burgersfig},
concerns the viscous Burgers equation on a periodic domain:
\begin{equation}
u_t^{} + u\kern .2pt u_x^{} = \nu\kern .5pt u_{xx}^{}, 
\quad
u(-\pi) = u(\pi), ~~ u_x^{}(-\pi) = u_x^{}(\pi).
\label{DEs-burgers}
\end{equation}
With $\nu = 0.075$ and initial condition $u(x,0) = -\sin(x)$, the
upper row of the figure shows the numerically computed solution at
times $t = 0,$ $1.5,$ and $5$, revealing an approximate shock that
steepens up at first and then loses amplitude because of dissipation.
Following \ccite{weideman22}, the lower row shows poles of a AAA
approximation to each solution in a nearby region of the complex
plane.  (Ordinary AAA approximation was used, though periodic AAA
would also be applicable here; see section \ref{sec:trig}.)  It is
known that for this problem the exact solution has an infinite train
of poles for each $t>0$ that move toward the real axis at first,
and then move out again as the amplitude dissipates.  Indeed, this
problem has the Painlev\'e property: the solutions are meromorphic
throughout the complex $x$-plane, so the only singularities are
poles \ccite{vanden}.

The AAA computation of poles in Figure \ref{DEs-burgersfig}
is quite accurate for the poles near the real axis and extremely
fast, so fast that it is easy to compute poles at every time step
in a numerical simulation.  The same would be true for many other
problems, indicating the exceptional convenience of estimating
complex behaviour by this method.

\begin{figure}
\begin{center}
\includegraphics [clip, width=4.2in]{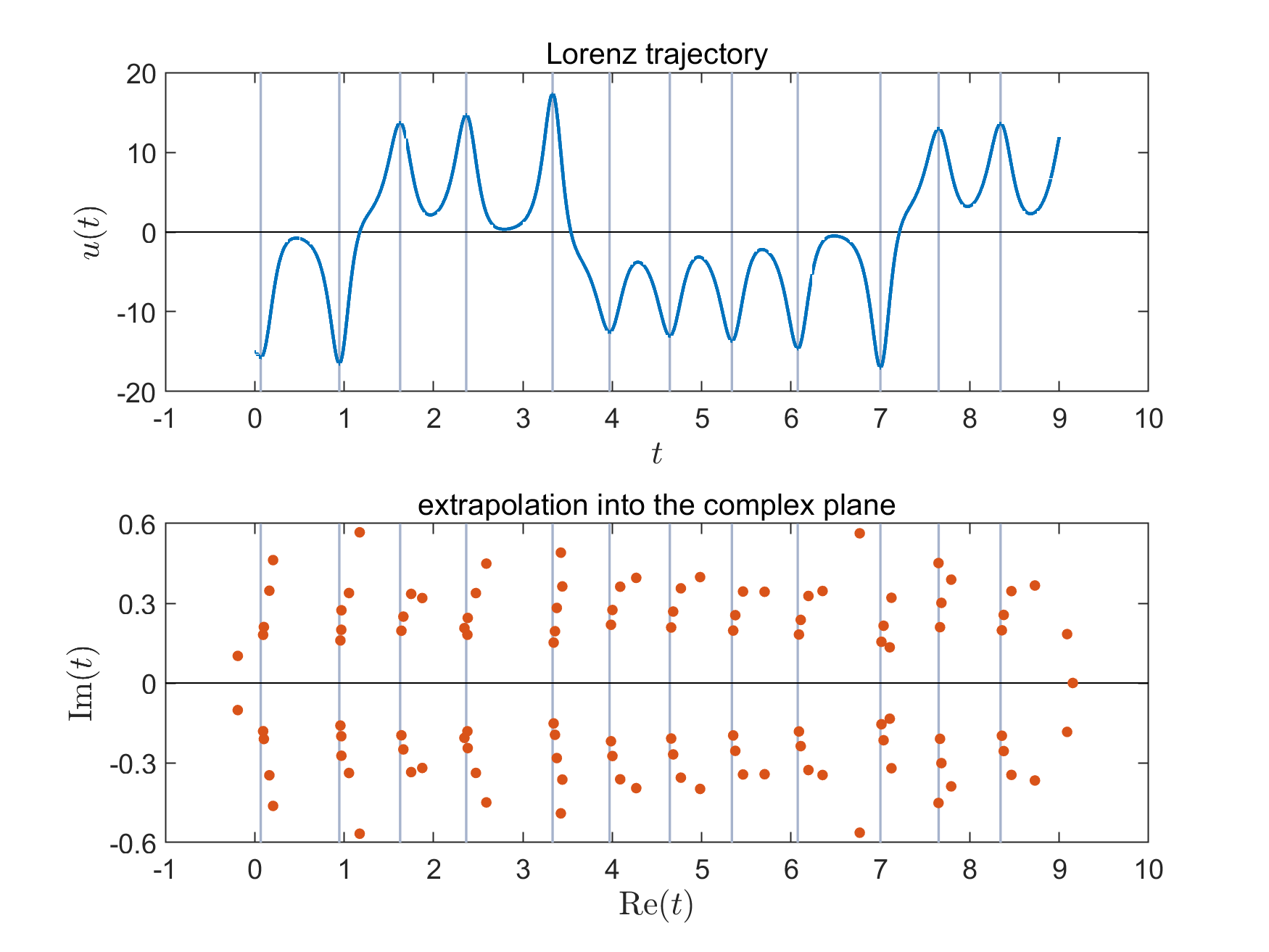}
\vspace{-15pt}
\end{center}
\caption{\label{DEs-lorenzfig}A trajectory $u(t)$
of the Lorenz equations
together with AAA poles in the complex $t$-plane.  The poles
line up along approximate branch cuts.}
\end{figure}

Our second example, illustrated in Figure \ref{DEs-lorenzfig},
shows a short stretch of a trajectory $u(t)$ of the Lorenz equations
(one of the three components $u,v,w$ in a numerical simulation)
together with, again, numerically computed poles in the complex
plane (now the complex $t$-plane rather than $x$-plane).  Once more
AAA makes it easy to calculate such estimates, and the vertical
gray lines show how the estimated poles match the local extrema
of the trajectory.  This problem has fundamentally different
behaviour from the Burgers example, however.  Here, it is known
that the analytic singularities are not poles but branch points
\ccite{taborweiss,vis}, and what Figure \ref{DEs-lorenzfig} is
displaying are trains of poles that approximate branch cuts near
those points, in the fashion we have seen in the figures of sections
\ref{sec:schwarz} and \ref{sec:potential}.

\begin{figure}[t]
\begin{center}
\vskip 0pt
\includegraphics [clip, width=4.6in]{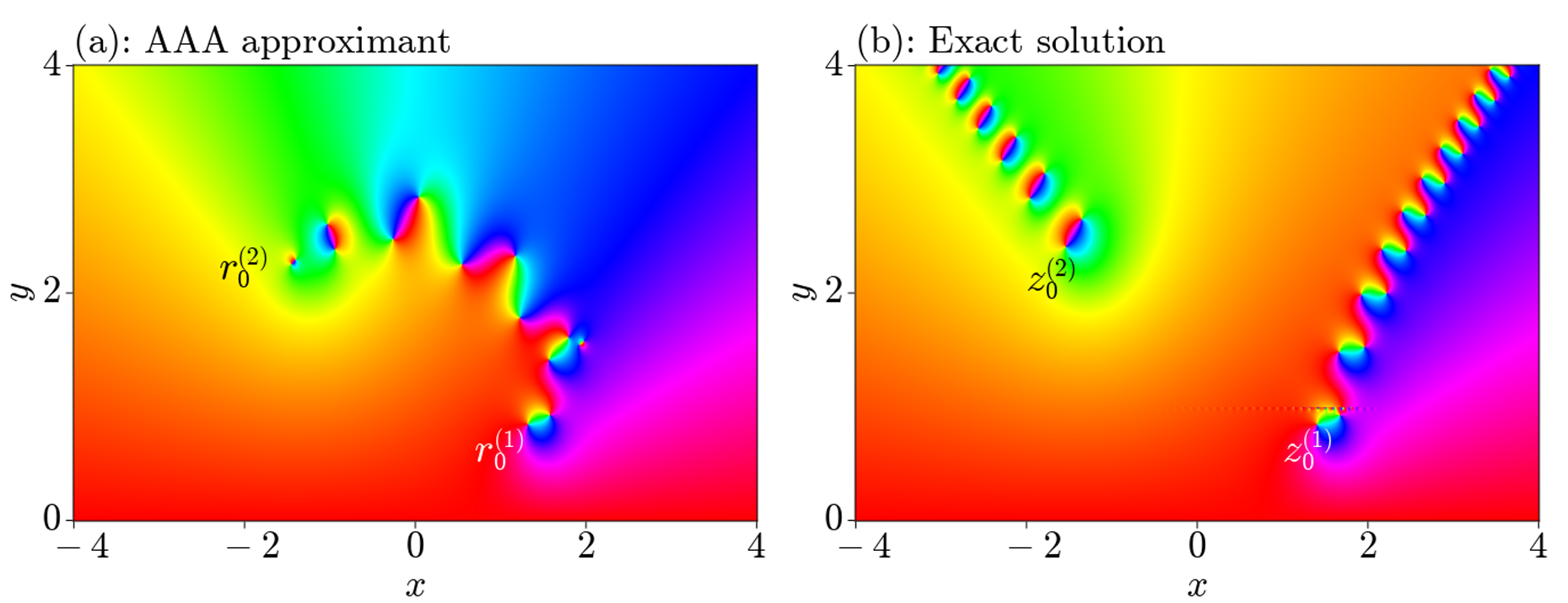}
\vspace{-8pt}
\end{center}
\caption{\label{DEs-vanden}A figure from \ccite{vanden} comparing
the exact pole structure of a Burgers solution with its AAA approximation from data
on the real axis.  The poles closest to the axis are
accurately captured, but not the poles further out.  (Note that the
real and imaginary axes in these images are not equally scaled.)}
\end{figure}

The distinction between poles and branch points, and between
resolution of exact singularities and mere function approximation by
convenient poles, seems unambiguous in principle but may be difficult
in practice.  The singularities of a complex function near the
$x$-axis will usually have a strong effect on its behaviour there,
but singularities further out may have very little influence.
Although mathematically well-defined in the analysis of an
exactly specified problem, such distant singularities may change
utterly when the problem is slightly modified.  Their locations
may be exponentially sensitive to perturbations, rendering their
mathematical significance questionable.  We close this section with a
beautiful figure from \ccite{vanden} showing a pronounced discrepancy
between the exactly known poles of a certain Burgers problem and
those computed by a AAA approximation.  The two sets of poles are
very different beyond the pair closest to the real axis, yet they
correspond to functions that agree to $10^{-13}$ on the real axis.

\section{\label{sec:bivariate}Multivariate analytic continuation}
Many applications involve functions of more than one variable that
one would like to extend beyond the domain where they are known a
priori.  In a typical case one might have a real analytic bivariate
or trivariate function $f(x,y)$ or $f(x,y,z)$ given in a domain in
$\real^2$ or $\real^3$, and one might wish to extend it to real
arguments outside this domain.  ``Real analytic'' is essentially
the same as analytic, being defined by the local convergence of
Taylor series (here these would be bivariate or trivariate series,
which would necessarily also converge in regions of $\complex^2$
or $\complex^3$).  What is distinctive about a function being real
as opposed to complex analytic is only that the focus is on real
values of the function taken on a real domain.

\begin{figure}
\begin{center}
\includegraphics[trim=10 106 0 90, clip, scale=.85]{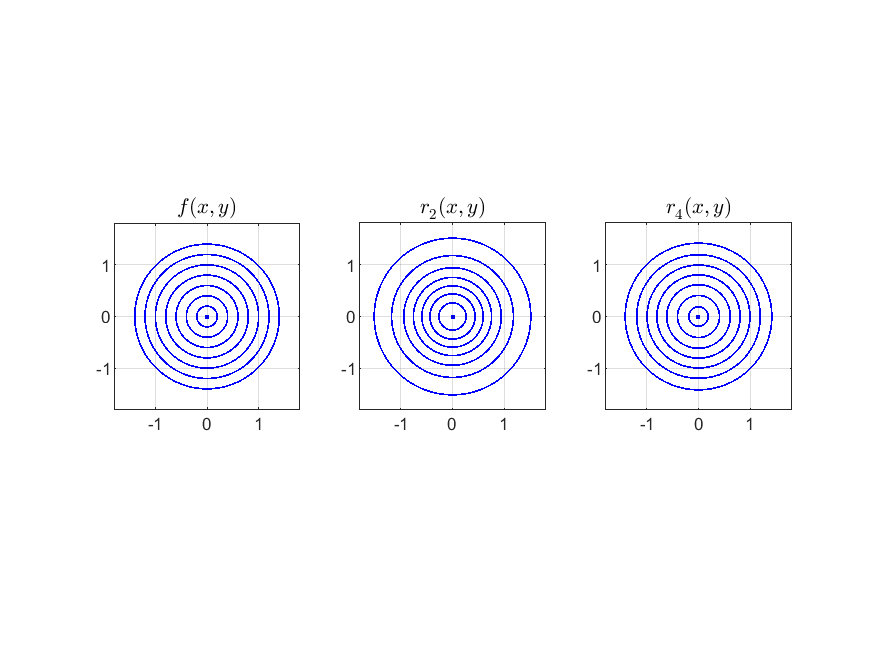}
\end{center}
\vskip -5pt
\caption{\label{bivariate-absvalfig}At the left, contours
$|f(x,y)| = 0,0.2,0.4,\dots, 1.4$
of the bivariate function (\ref{bivariate-absval}).  The other
two images show the same contours for the degree $2$ and $4$
rational approximations 
(\ref{bivariate-1}) and (\ref{bivariate-2}).  No numerical 
approximation algorithm is known for finding near-best
multivariate rational approximations like these.}
\end{figure}

For multivariate problems, no algorithm is known with the ``magic''
near-best approximation behaviour of AAA in the univariate case.  
At least that is how it appears if we consider the example
\begin{equation}
f(x,y) = \sqrt{x^2+y^2},  \quad -1\le x,y \le 1.
\label{bivariate-absval}
\end{equation}
This function is real analytic for all $(x,y) \in \real^2$ except
for the singularity at $(x,y) = (0,0)$.
There exist good low-degree bivariate rational
approximations of it, like 
\begin{equation}
r_0^{}(x,y) = 0.707, \quad 
r_2^{}(x,y) = {0.0618 + 2.257(x^2+y^2)\over
1 + 1.191(x^2+y^2)},
\label{bivariate-1}
\end{equation}
with maximum errors $0.707$ and $0.062$ on the given domain, respectively,
and
\begin{equation}
r_4^{}(x,y) = {0.0120 + 11.848(x^2+y^2) + 49.450(x^2+y^2)^2\over
1 + 40.910(x^2+y^2) + 18.788(x^2+y^2)^2},
\label{bivariate-2}
\end{equation}
with maximum error $0.012$.  The first image of
Figure \ref{bivariate-absvalfig} shows level curves of
(\ref{bivariate-absval}), which are equally spaced concentric circles
since $f$ is the $2$-norm of the vector $(x,y)^T$.  The next two
images show the corresponding contours for $r_2^{}$ and $r_4^{}$,
the latter almost indistinguishable from the contours for $f$.

The trouble is, no general algorithm is known for finding near-best
approximations like these.  In particular, no true multivariate
AAA algorithm has yet been discovered---seemingly because there
is no multivariate barycentric formula applicable with arbitrary
support points.  (The paper \ccite{pAAA} introduces a bivariate
approximation algorithm called p-AAA with support points restricted
to a tensor product grid.  For numerical multivariate approximation
by other methods, see \ccite{austin}.)  To find the approximations
(\ref{bivariate-1}) and (\ref{bivariate-2}), we did not use bivariate
numerical approximation but exploited the fact that this example
can be reduced to a univariate problem of approximation of $\sqrt s$
with $ s = x^2+y^2$ for nonnegative values of $s$.  (The reduction is
not exact, since the domain $-1\le x,y\le 1$ is not axisymmetric.)
This is Newman's problem \ccite{newman}, for which root-exponential
convergence as a function of the degree has been known since 1964.
To find the coefficients of $r_4^{}$ in (\ref{bivariate-2}), for
example, we computed minimax approximations of $\sqrt s$ for $ s\in
[\kern .5pt 0,2\kern .4pt ]$ with the Chebfun commands

{\small
\begin{verbatim}
      s = chebfun('s',[0 2]);
      [p,q,r,err] = minimax(sqrt(s),2,2);
      pcoeffs = poly(p)/q(0), qcoeffs = poly(q)/q(0)
\end{verbatim}
\par}

Is it reasonable to hope for a multivariate algorithm that could find
approximations like (\ref{bivariate-1}) and (\ref{bivariate-2})
as easily as AAA finds approximations in a single variable?
We do not know, and we regard this as an important open problem.
A pessimist (or an expert in several complex variables) may note
that the singularity of $f$ that is just the point $(0,0)$ for
$(x,y)\in \real^2$ becomes the pair of planes $y = \pm ix$ for $(x,y)
\in\complex^2$ (compare \ccite{bhh}).  An optimist may counter that
if good approximations exist, there ought to be ways to find them.

\begin{figure}
\begin{center}
\includegraphics[trim=10 106 107 90, clip, scale=.85]{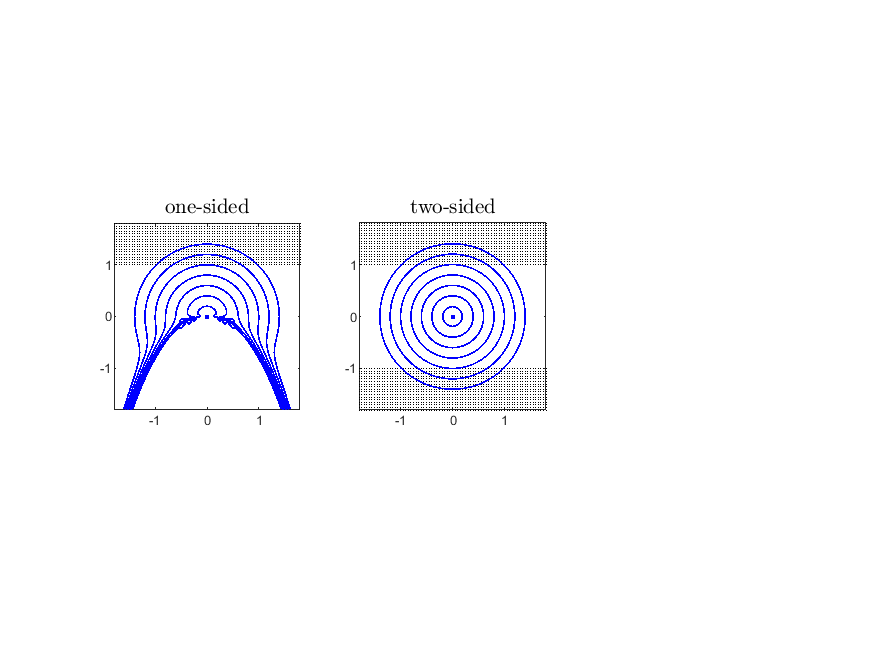}
\end{center}
\vskip -5pt
\caption{\label{bivariate-absval2fig}Bivariate approximation
of (\ref{bivariate-absval}) by
AAA approximation along vertical lines.  On the left, approximation along
vertical lines with data samples just above the singularity, marked
by black dots.
On the right, approximation along vertical
lines with samples on both sides of the singularity.}
\end{figure}

In the absence of a true multivariate AAA algorithm, an easy
alternative appears to be the method introduced in \ccite{jjiam}:
use a collection of univariate AAA approximations, each applied on
a different line or circle or other analytic arc.  For example,
Figure \ref{bivariate-absval2fig} shows the same problem as in
Figure \ref{bivariate-absvalfig}, but now with $f(x,y)$ approximated
numerically along vertical lines.  There are 73 vertical lines,
each with 73 sample points.  In the left image, the top 17 sample
points of each vertical line (i.e., those with $y\ge 1$) are used
to construct a AAA approximant, and the global result is plotted.
The computation time for all 73 approximations is small, about 1/3
seconds on our laptop.  In the right image a much better result is
obtained from AAA sample data both above and below the singularity,
the 34 samples with $|y|\ge 1$.

Since the function (\ref{bivariate-absval}) is elementary,
one can understand these approximations in detail.  Along
the vertical line in Figure \ref{bivariate-absval2fig} associated with
a particular value $x=x_0^{}$, the function being approximated is
\begin{displaymath}
f_{x_0^{}}^{}(y) = \sqrt{x_0^2 + y^2}.
\end{displaymath}
For $x_0^{} = 0$ this becomes simply $f_0^{}(y)=|y|$, whose rational
approximants are easy to understand---it will be $r(y) = y$ exactly
in the one-sided case and a good approximation $r(y) \approx |y|$
in the two-sided case (an approximation of $\hbox{sign}(x)$ flavor
with a branch cut, to be discussed in section \ref{sec:zolo}).

\begin{figure}
\begin{center}
\includegraphics[trim=150 165 15 0, clip, scale=1]{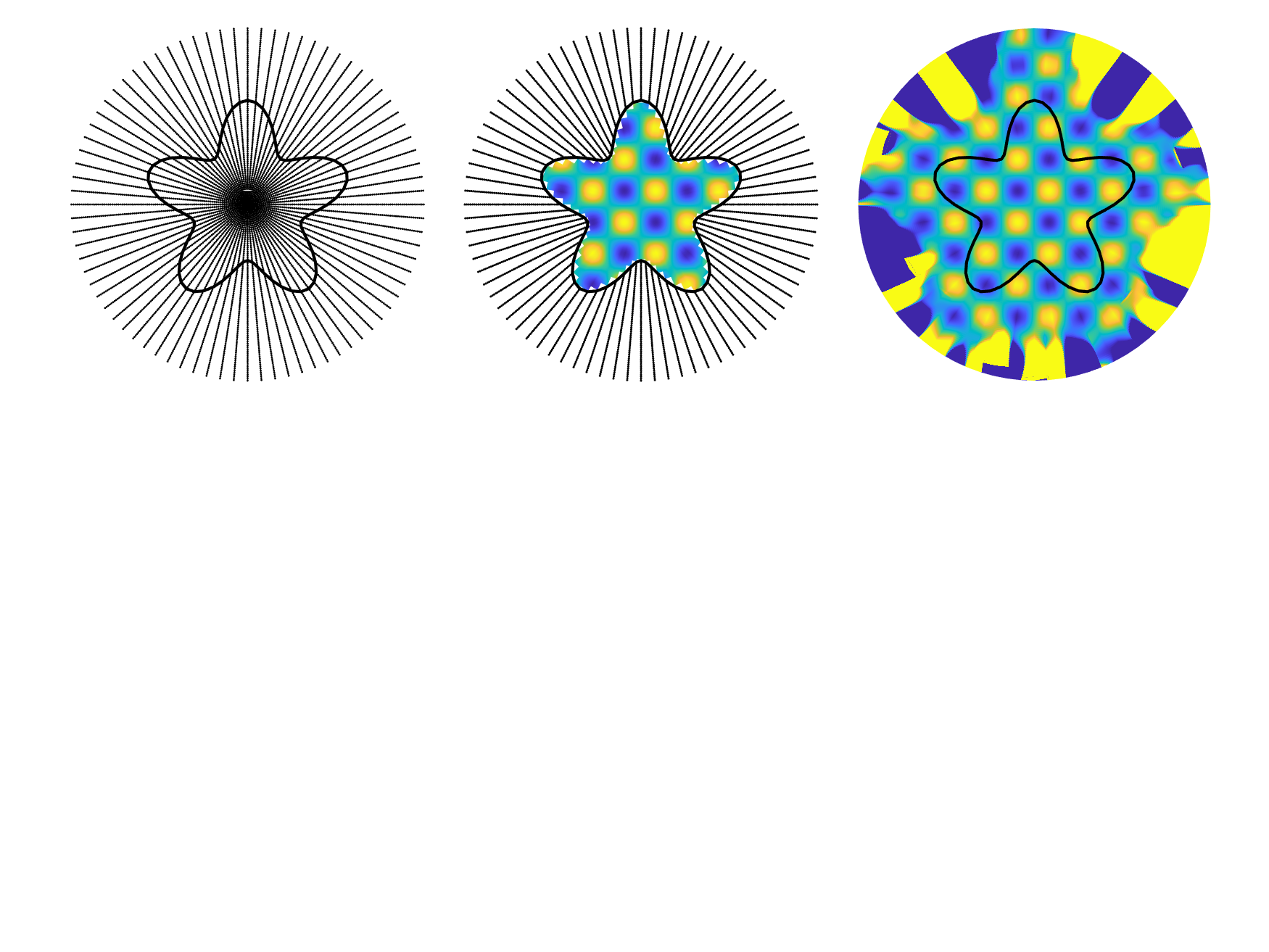}
\end{center}
\vskip -15pt
\caption{\label{bivariate-fryk}Reproduction of a figure from
\ccite{jjiam} showing continuation of a smooth function across
a boundary by univariate AAA approximation on 60 radial lines. As in
section \ref{sec:analcont}, about one wavelength of accurate extension is achieved.}
\end{figure}

In the literature of computational methods for \PDEs, multivariate
analytic continuation is considered in the context of what is called
{\em function extension\/}.  Here, for various algorithmic reasons,
it is desired to numerically continue a multivariate function
smoothly across a boundary.  Accuracy of analytic continuation may
not be the goal so much as smoothness.  Some contributions in this
area include \ccite{adcock,boyd,bruno,efj,fryklund,larsson}.

Figure \ref{bivariate-fryk}, which reproduces part of a figure from
\ccite{jjiam}, illustrates a computation of this kind.  Here AAA
approximations have been computed along each of 80 radial lines, and
the resulting data on a polar grid has then been fed to the MATLAB
contour plotter.  As discussed in section~\ref{sec:analcont},
about one wavelength of successful extension is
evidently achieved, with the number of accurate digits diminishing
from 13 on the inner boundary to 0 at the edge.

So far as we are aware, the other methods for function extension
proposed in the literature are not as accurate as the AAA approach.
For example, figures in \ccite{fryklund} and \ccite{efj} suggest
about half a wavelength of continuation rather than one wavelength.
On the other hand, these methods have advantages of speed,
simplicity and robustness related to their lack of dependence on
a process as nonlinear as rational approximation.  For example,
the method of \ccite{efj} involves an explicit formula that is the
same ``in all directions'' and therefore produces approximations
that are continuous, indeed analytic under appropriate hypotheses,
with respect to the direction tangential to the boundary.  The AAA
approach, by contrast, computes a different approximation along
each line and does not guarantee even tangential continuity, let
alone smoothness, though it often gives numerical continuity to
many digits of accuracy.

\section{\label{sec:res}Computing complex resonances} This section
presents an application of AAA approximation that was introduced in
\ccite{resonance}; a similar method was proposed in \ccite{betz}
and \ccite{bink}.  One can regard this as a special case of the
linear and nonlinear eigenvalue problems to be discussed in the
next two sections.

The physical problem is oscillation in a 2D or 3D cavity; in the 2D case
one often speaks of a ``drum.''  In the
standard and most familiar setting, the cavity is a closed planar domain bounded
by a curve $\Gamma$ or a closed spatial domain bounded by a surface
$\Gamma$, and we seek eigenvalues
and eigenfunctions for Laplace oscillations in this domain:
\begin{equation}
\Delta u +\kappa^2 u = 0 ,  \quad u = 0 \hbox{ on } \Gamma.
\label{res-eigprob}
\end{equation}
This is a self-adjoint problem, with negative real eigenvalues $-\kappa^2$ corresponding
to real frequencies $\kappa$.
\begin{figure}
\begin{center}
\vskip 10pt
\includegraphics[trim=0 0 0 0, clip, scale=.35]{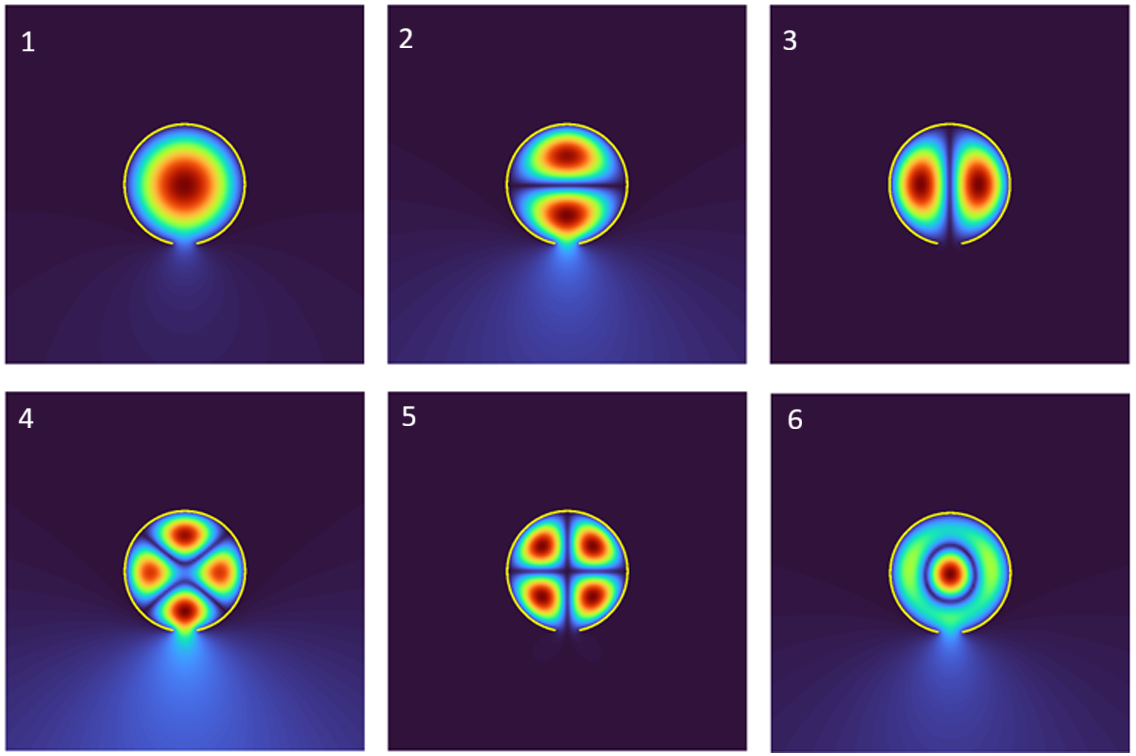}
\end{center}
\vskip -5pt
\caption{\label{res-disk} Images from \ccite{resonance} of the first
six complex 
resonances in the unit disk with a gap in the boundary of angle $\pi/8$.
Images 2--3 and 4--5 are near-degenerate pairs
that become exactly degenerate in the limit of zero gap size.}
\end{figure}

AAA can be applied to this well-known problem, as is shown in
section~\ref{sec:helmholtz} and
in Figures 2 and 6 of \ccite{resonance}, but this section
is devoted to a variation where AAA comes especially into its own:
when the boundary~$\Gamma$ is not fully closed and the eigenvalues
become slightly complex.  Various terms are in use for such
problems, including {\em complex resonances\/}, {\em scattering
resonances\/} and {\em scattering poles\/}; see \ccite{dyatlov}
for the mathematical foundations.  Figure \ref{res-disk} shows the
simplest possible example, a cavity bounded by a curve $\Gamma$
consisting of the unit circle with a gap in it, here a gap of angle $\pi/8$.
As suggested in the figure, we can imagine that this resonator
will have approximate eigenfunctions that oscillate like those of
a closed circular cavity but lose amplitude at a slow exponential
rate as energy is radiated out through the gap and away to infinity.

Equation (\ref{res-eigprob}) is still an eigenvalue problem, but
no longer self-adjoint, since energy is absorbed at $\infty$.
One way to interpret it as a scattering
problem might be as follows.  For each $\kappa$, let $G(\kappa)$ be the linear operator
that maps an incident wave at frequency $\kappa$ that hits $\Gamma$
to the corresponding wave at frequency $\kappa$ scattered by $\Gamma$.
If $\kappa$ is near a resonant frequency, $G(\kappa)$ will be large in norm, mapping
incident waves at this frequency to large scattered responses.
There will be no real values of $\kappa$ at which the size of $G(\kappa)$
is infinite, for that would correspond to nonzero radiation with
zero incident energy.  However, there can be complex eigenvalues.
A complex value of $\kappa$ with a small imaginary part will correspond
to an oscillation at the frequency $\Re (\kappa)$ with exponential decay
rate $\Im (\kappa)$.  If $\Im (\kappa)$ matches the decay rate associated
with radiation of energy, then $G(\kappa)$ will be infinite, mapping an
incident wave of amplitude zero to a finite response, oscillating
and decaying with time.

AAA works with scalar functions, not operators.  To turn the
concept just described into a scalar problem, we could make the input
a scalar amplitude by fixing the shape of the incident wave in
some arbitrary fashion: say, a plane wave of frequency $\kappa$ hitting
$\Gamma$ from a prescribed angle.  Likewise, we could make the output
a scalar by sampling the scattered wave at an arbitrary point $\xz =
(x_0^{},y_0^{})$.  The complex eigenvalue problem becomes: find a
complex frequency $\kappa$ at which the scattered amplitude at $\xz$
is infinitely larger than the incident amplitude.  We might be
unlucky, if $\xz$ is situated at a point where the scattered wave has
amplitude zero, but if $\xz$ is chosen at random, the probability
of this will be zero, and under reasonable assumptions
the probability of losing $d$ digits of
accuracy to a near-zero sampling point can be expected to be
of order $10^{-d}$.

This is the idea of the method presented in \ccite{resonance},
though the details are different.
The scattering problem is formulated in terms of an integral operator 
$F(k)$ that maps a distribution on $\Gamma$
of singular Hankel functions for frequency $\kappa$ 
to the outgoing wave field on $\Gamma$ generated
by these Hankel functions.  If the wave
field can be zero on $\Gamma$ with a nonzero distribution of Hankel functions, then
$\kappa$ is a complex
eigenvalue.  So we are looking for values of $\kappa$ for which $F(\kappa)$ is singular,
or equivalently, for which $(F(\kappa))^{-1}$ is unbounded.  We scalarise
the problem for AAA by considering
\begin{equation}
f(\kappa) = u^* (F(\kappa))^{-1} v,
\label{res-scalarized}
\end{equation}
where $u$ and $v$ are a pair of nonzero functions fixed at random (independently
of $\kappa$).
We speak of (\ref{res-scalarized}) as a {\em randomly scalarised
integral operator\/}.

What's needed is to find poles of $f(\kappa)$, and AAA does this just
as in section~\ref{sec:poles}, by computing approximations $r(\kappa)\approx f(\kappa)$.
The most straightforward way to proceed is to sample $f(\kappa)$ at
a number of real frequencies~$\kappa$.  Each sample evaluation will
require the numerical solution of an integral equation.  Then we construct
the approximation $r(\kappa)$ and consider poles of $r(\kappa)$ near the
sample set, which can be expected to be good approximations to
poles of $f(\kappa)$.  Details of the integral equation are given
in \ccite{resonance}, and its numerical discretisation comes
from \ccite{brunolint}.  The method is fast and accurate, giving
eigenvalues reported in \ccite{resonance} accurate to 13 digits.
Refinements of the method are discussed there involving adaptive
sampling of values of $\kappa$ to zero in on approximate poles as they
are identified.

Note that although Figure \ref{res-disk} shows a 2D example,
nothing in this method is restricted to 2D.\ \ For a 3D cavity
the discretisation of the integral equation is more complicated,
but the AAA approximation is the same.

To eight digits, the complex frequencies associated with Figure
\ref{res-disk} are as follows.  The third column shows the
corresponding roots of Bessel functions that are the eigenvalues
for the closed unit circle.

{\small
\begin{verbatim}

        Mode       Circle with gap pi/8      Closed circle
        ----       --------------------      -------------
\end{verbatim}
\begin{verbatim}
          1       2.3918509 - 0.0008668i       2.4048256
\end{verbatim}
\vspace{-7pt}
\begin{verbatim}
          2       3.7858514 - 0.0075513i       3.8317060
          3       3.8315198 - 0.0000008i          same    
\end{verbatim}
\vspace{-7pt}
\begin{verbatim}
          4       5.0664101 - 0.0227539i       5.1356223
          5       5.1345996 - 0.0000118i          same  
\end{verbatim}
\vspace{-7pt}
\begin{verbatim}
          6       5.4867988 - 0.0108397i       5.5200781

\end{verbatim}
\par}

There is some fascinating physics in these numbers, as
discussed in \ccite{resonance}.  Note that they are grouped as a
singleton, two pairs and another singleton.  Figure \ref{res-disk}
shows that these are associated with a lowest mode that would be
radially symmetric if there were no gap, a pair of modes that would
be degenerate if there were no gap, a second such near-degenerate
pair and then another singleton.  The gap in the boundary
has broken the degeneracies, and in each pair, we see that the
lower-frequency eigenfunction is nonzero near the gap and radiates
significantly outside (modes 2 and 4), whereas the higher-frequency
eigenfunction is near-zero near the gap and radiates so little
that the background appears a uniform deep blue (modes 3 and 5).
This is reflected in the eigenvalue imaginary parts being much
smaller and the real parts being much closer to the Bessel roots.
It is shown numerically in \ccite{resonance} based on further AAA computations
that if $\theta$ is the gap size, the real and imaginary parts of
the eigenvalues converge as $\theta\to 0$ at the rates $O(\theta^2)$
and $O(\theta^4)$ for the modes that are nonzero near the gaps,
and at the rates $O(\theta^4)$ and $O(\theta^8)$ for the modes that
are near zero there.  Thus if the gap size in Figure \ref{res-disk}
were reduced from $\pi/8$ to $\pi/16$, the decay rates for modes
2 and 3 would reduce from $8\cdot 10^{-3}$ to $5\cdot 10^{-4}$
and from $8\cdot 10^{-7}$ to $3\cdot 10^{-9}$, respectively.
Jon Chapman of the University of Oxford (private communication,
May 2025) has derived asymptotic formulas for modes 2 and 3 in the table
above that match the data beautifully.  Chapman finds
\begin{equation}
{\kappa\over \kappa_0^{}}  \sim 1 - {\theta^2\over 16}
- {\theta^4 \kern .5pt C \over 128\pi}\kern 1pt i
\quad \rlap{\hbox{(mode 2)}}
\label{chapman1}
\end{equation}
and 
\begin{equation}
{\kappa\over \kappa_0^{}}  \sim 1 - {\theta^4\over 512}
- {\theta^8 \kern .5pt C \over 2^{17}\pi}\kern 1pt i
\quad \rlap{\hbox{(mode 3),}}
\label{chapman2}
\end{equation}
where $\kappa_0^{} \approx 3.8317060$ is the first zero of $J_1$ and
\begin{equation}
C = \sum_{n=0}^\infty {n^2\over |H_n^{(1)}(\kappa_0^{})|^2},
\label{chapmanC}
\end{equation}
with $J_1$ and $H_n^{(1)}$ denoting the usual Bessel and Hankel functions.
All these observations and details have been stimulated by the
accurate calculations made possible by AAA.

See \ccite{betz,bink,resonance} for further examples of the use of
AAA to calculate complex resonances.

\section{\label{sec:eig}Matrix eigenvalue problems}
Let $A$ be an $N\times N$ matrix and let $I$ be the $N\times N$
identity.  Then
\begin{equation}
R(z) = (zI - A)^{-1},
\label{eig-res}
\end{equation}
known as the {\em resolvent\/} of $A$ (or more properly
the negative of the resolvent), is a matrix function of
$z$ defined for all $z\in\complex$ except the eigenvalues of $A$,
where it is, informally speaking, infinite.

One can consider $R(z)$ as a matrix rational function, but instead let
us scalarise it as in (\ref{res-scalarized})
by multiplying on the left and right by nonzero row and column vectors
$u^*$ and $v$:
\begin{equation}
R(z) = u^*(zI - A)^{-1} v.
\label{eig-res1}
\end{equation}
Here $u^*$ and $v$ are fixed, independent of $z$.  This scalarised
resolvent is a rational function of degree $N$.  This does
not mean that its degree is exactly $N$; it will be lower if $A$
has any eigenvalues with geometric multiplicity greater than $1$,
or if $u^*$ or $v$ happen to be orthogonal to one or more 
left or right eigenvectors of $A$, respectively.
For simplicity, let us suppose that $A$ has
distinct eigenvalues and that $u^*$ and $v$ are random vectors with
independent entries chosen from a continuous distribution such as the
normal distribution.  Then with probability $1$, $R$ is a rational
function of degree exactly $N$ having $N$ simple poles which are
the eigenvalues of $A$.  
This leads to the possibility that we can find the
eigenvalues by calculating a AAA approximation $r\approx R$ based
on samples of $R(z)$ at some points $z$ and then calculating the
poles of $r$, as described in section \ref{sec:poles}.

For example, suppose $A$ is the $20\times 20$ diagonal matrix
with diagonal entries $0,1,\dots,19$, which are accordingly its
eigenvalues.  (The fact that $A$ is diagonal will have little
effect on the results, since $u$ and $v$ are random.)  Here is
a computation of the kind just described based on samples of $R$
at $100$ equispaced points on the circle $|z|=1.05$ in the complex
plane.  Alternatively, for a real symmetric problem like this,
one could use real samples \ccite{filters}.

{\small
\begin{verbatim}
      Z = 1.05*exp(2i*pi*(1:100)/100);
      N = 20; A = diag(0:N-1); I = eye(N);
      rng(1), u = randn(N,1)+1i*randn(N,1);
              v = randn(N,1)+1i*randn(N,1);
      R = @(z) u'*((z*I-A)\v);
      for i = 1:length(Z), F(i) = R(Z(i)); end
      [~,pol] = aaa(F,Z); pol = sort(pol)
        
      pol =
        0.000000000000025 - 0.000000000000008i
        0.999999999999999 + 0.000000000000000i
        1.999999995063250 - 0.000000443219559i
        3.000344539241449 + 0.000522614895484i
        4.037538715031964 + 0.004841722473988i
        5.079074587752028 + 1.453783302951845i
        8.343483853643621 + 0.565153761906265i
       15.975301039683405 + 0.020589938026920i
\end{verbatim}
\par}

\noindent
The two eigenvalues $0$ and $1$ enclosed by the sample points have
been found to 13 digits or more, and the eigenvalues $2$, $3$,
and $4$ have been found to about 9, 4 and 2 digits, respectively.
After this there are three further poles with less connection
to individual eigenvalues.  (The nonzero imaginary parts in these
numbers result from the fact that {\tt aaa.m} does not maintain real
symmetry, as discussed on p.~\pageref{symmetrypage}.)
The results are similar, just slightly less accurate, if $A$
is made non-normal by the addition of random samples from $N(0,1)$ to
its strictly upper-triangular part.
If the diagonal matrix experiment is repeated but with dimension
$100$ or $1000$, the results are again not much different, with the number
of computed poles increasing in both cases from 8 to 9.

The fact that not all the eigenvalues of $A$ are computed is
the very reason why this method may be useful: it can treat
matrices of dimensions too large for standard methods such
as the QR algorithm.  Often $A$ will be a discretisation of a
linear operator of infinite dimension.  In such cases, a crucial
feature of eigenvalue calculation via rational approximation
is that attention is concentrated on a particular region of the
real line or complex plane.  An operator might be discretised by
a matrix of dimension $10^6$ with eigenvalues widely distributed
in $[\kern .5pt 0,\infty)$, for example.  If only the eigenvalues
in a subinterval $[a,b\kern .5pt]$ are of interest, then it may be
possible to determine these at the cost of a few dozen matrix solves,
one for each sample point over a set surrounding or approximating
$[a,b\kern .5pt]$.  These solves are trivially parallelisable.

\begin{figure}
\vspace*{-5pt}
\begin{center}
\includegraphics[trim=35 105 5 0, clip, scale=1.1]{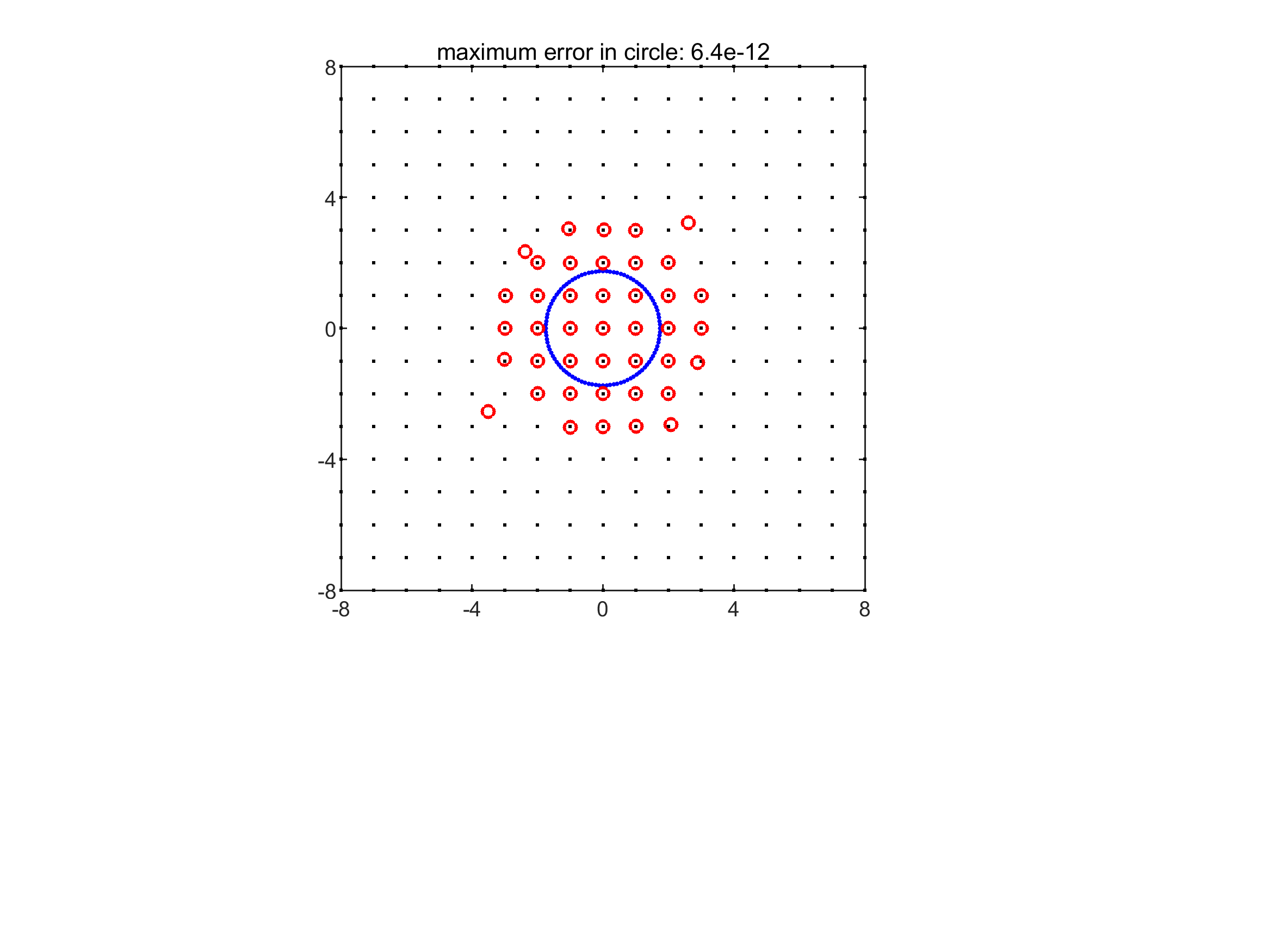}
\end{center}
\vspace*{-9pt}
\caption{\label{eig-gridfig} AAA calculation of the eigenvalues inside
the circle $|z| = 1.75$
of a $40401\times 40401$ diagonal matrix $A$ with diagonal entries
on the integer grid $j + i\kern .3pt k$, $-100\le j,k \le 100$.  Black dots mark
the eigenvalues of $A$ in the plotting region
and red circles show the poles of the approximation $r$ to
the scalarised resolvent.  The 100 AAA sample points are marked by blue dots.}
\end{figure}

Figure \ref{eig-gridfig} shows an example like the one
just displayed, except that $A$ is now the $40401\times 40401$ diagonal
matrix whose diagonal entries are the complex integers
$j + i\kern .3pt k$, $-100\le j,k \le 100$.
The scalarised resolvent
is evaluated at 100 equally spaced points on the circle $|z|=1.75$ (easy because of
the diagonal structure),
and AAA is applied as usual. The rational approximation $r(z)$ has
degree $42$, and one sees that its poles approximate a number of eigenvalues
of $A$ within and near the sample circle, with the approximations inside
the circle being accurate to 11 digits or more.

The idea of computing matrix eigenvalues via poles of rational
approximations seems to have been first proposed in \ccite{akt},
before the development of AAA.\ \ However, there were important contributions
in the preceding decade on the related method of computing
eigenvalues via {\em discretised contour integrals\/}.  This was
mainly the work of two schools: Sakurai and his collaborators
in Japan \ccite{ss03} and Polizzi and his collaborators in
Massachusetts, who call their method FEAST \ccite{feast,tang}.
An earlier
paper by Goedecker preceded both of these contributions \ccite{goedecker}.

To compute eigenvalues by discretised contour integrals, one
again starts from the idea of looking for poles of the scalarised
resolvent.  Mathematically, it turns out that these can be
represented in terms of Cauchy integrals of $R(z)$ and related
functions over a curve enclosing the desired poles; for details
see \ccite{akt} and \ccite{guetteltiss}.  When the integrals are
discretised, for example by the equispaced trapezoidal rule in
the case of a circular contour \ccite{trap}, one obtains numerical
approximations to the poles.  The contributions of Sakurai et al.\ and
Polizzi et al.\ are both of this kind, and have been shown to be
effective for large-scale problems.
These methods have been used a good deal in computational chemistry
and physics.

Although discretised contour integrals and rational approximations
may sound quite different, they are related.  As we will
discuss in section \ref{sec:quad},
every discrete quadrature formula is equivalent
to a rational approximation, with the nodes and weights of the
quadrature formula corresponding to the poles and residues of the
rational function \ccite{horning}.  So discretised Cauchy integrals give rational
approximations too, the difference being that they are constructed
by an {\em a priori\/} formula rather than by AAA.\ \ AAA approximation
is more flexible than contour integration, since the sample
points can be distributed arbitrarily rather than having to fit a
prescribed configuration, and it will often be more accurate since
the approximations are near-best.  Discretised Cauchy integrals
have the advantage of a stronger theoretical foundation, since
they are built on long-established theorems of complex analysis,
but difficulties may arise from the discretisation, especially if
there are eigenvalues near the contour.  For an analysis of the speed
advantage of approximation over discretised Cauchy integrals in a
polynomial setting, see Theorem 2.1 and Figure 3 of \ccite{akt}.
For a method of eliminating rounding error difficulties in
evaluations near poles, see \ccite{horningnak}.

All these methods tend to lose accuracy when there are multiple
eigenvalues in the region of interest.  For example, here is what
happens when five additional diagonal entries $0.1,0.2,\dots,0.5$
are inserted in the matrix $A$ from the beginning of this section,
increasing the dimension from 20 to 25.\ \ The accuracy of the
inner eigenvalues falls to just 8 or 9 digits.

{\small
\begin{verbatim}
      pol =
        0.000000000161100 - 0.000000000074178i
        0.099999997785597 + 0.000000008636569i
        0.199999996639445 - 0.000000001308943i
        0.300000007567588 + 0.000000004465621i
        0.399999999873358 - 0.000000000643139i
        0.500000000017551 + 0.000000000023886i
        1.000000000000000 - 0.000000000000000i
        1.999999999367556 + 0.000000000257540i
        2.999974537765880 - 0.000060424614904i
        3.998803671344900 - 0.008085233366773i
        4.995045433823007 + 0.044799185742410i
        6.059119288839810 - 0.010238777914851i
       12.218833692968895 - 0.344378150290857i
       18.120846774789243 + 0.801506642603402i
\end{verbatim}
\par}

\noindent Similarly, if the 100 sample points of Figure
\ref{eig-gridfig} with $|z|=1.75$ are replaced by 200 sample points
with $|z|=2.75$, so that there are 21 rather than 9 eigenvalues
within the circle, the maximum error deteriorates to $4.4\cdot
10^{-6}$.

To improve the accuracy in cases like this, three approaches have been considered:
\medskip

1. Zoom in near individual eigenvalues,

2. Sketch the resolvent with $k>1$ row and column vectors,

3. Iterate an approximation or integral as a ``filter''.
\medskip

\noindent 
The first method, zooming closer to individual eigenvalues, can
be done in many ways and was employed in \ccite{resonance} for the
resonance calculations described in the last section.  The second
method amounts to compressing the resolvent (\ref{eig-res})
to a $k\times k$ matrix function instead of a scalar function.
This function can then be approximated by the set-valued AAA
algorithm to be discussed in section \ref{sec:sets}, with considerable improvement in
accuracy though at the price of working with matrices with $k^2$
times as many rows (our experiments, not yet published).  The third
idea of filtering springs original from the work of Polizzi and
his collaborators \ccite{filters,brennan,guettelzolo,tang,vbk}.  An
investigation of the
properties and relative merits of these three possibilities in the
context of algorithms based on rational approximation has
yet to be carried out.

The discussion continues in the next section.

\section{\label{sec:nonlin}Nonlinear eigenvalue problems}
In a nonlinear (matrix) eigenvalue problem the resolvent
$(zI-A)^{-1}$ of (\ref{eig-res}) is generalised to $A(z)^{-1}$, where
$A(z)$ is an arbitrary (usually meromorphic) $N\times N$ function of
$z$.  If $A(\lambda)\kern .5pt v = 0$ for some $\lambda\in\complex$
and nonzero $v\in \complex^N$, then $\lambda$ is an eigenvalue of
the function $A(z)$ and $v$ is a corresponding right eigenvector.
Reviews of such problems can be found in
\ccite{mvoss,guetteltiss,brennan}, and there is an established test collection of
nonlinear eigenvalue problems known as the NLEVP collection \ccite{nlevp}.

\begin{figure}
\begin{center}
\includegraphics[trim=15 80 10 0, clip, scale=.95]{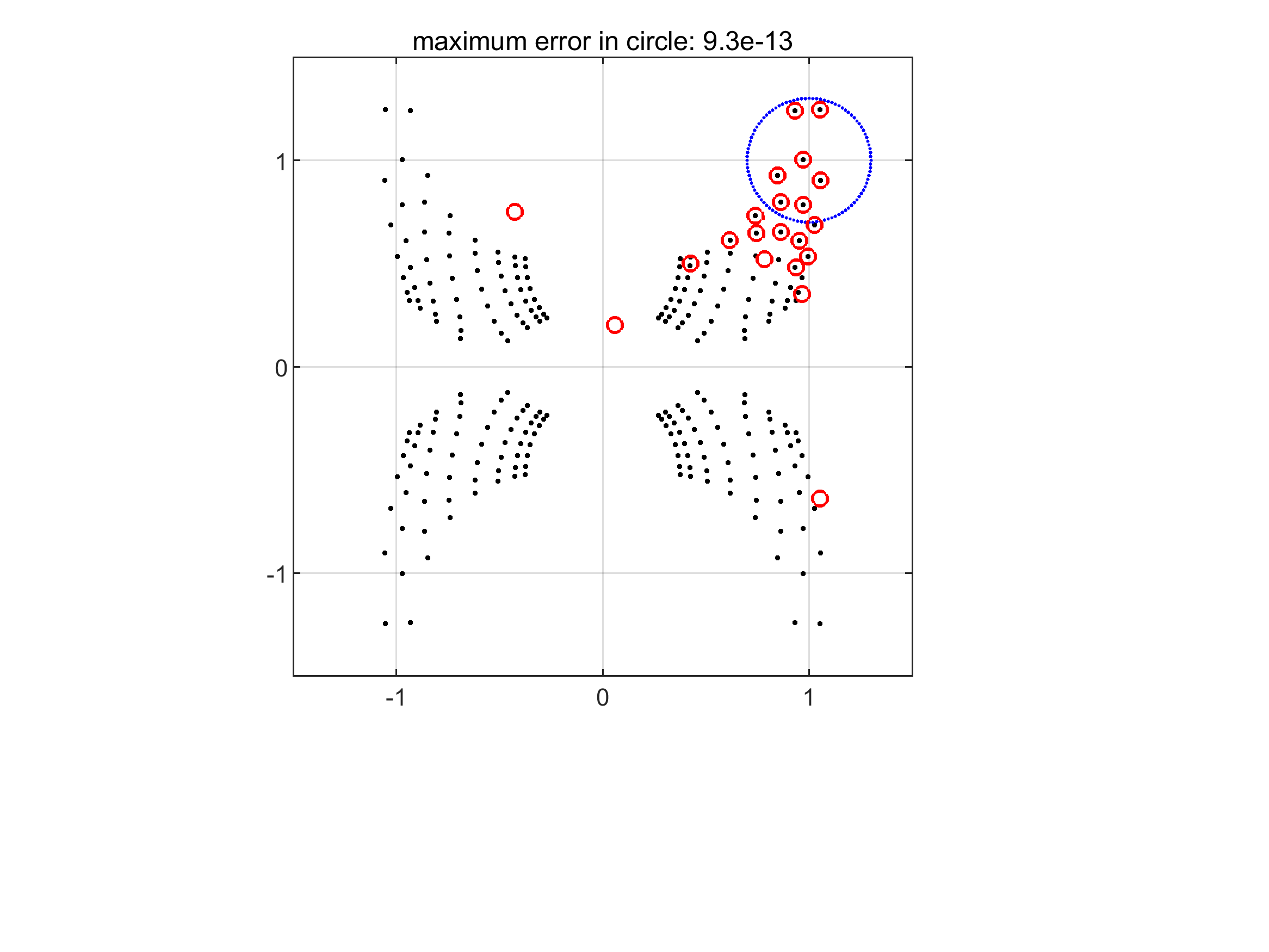}
\end{center}
\vspace*{-4pt}
\caption{\label{nep-butterfig} Calculation of 7 nonlinear eigenvalues of
the ``butterfly'' problem from the NLEVP collection
\ccite{nlevp}.  AAA approximation is based on 100 samples on the circle
$|z-(1+i)| = 0.3$, shown in blue.  Black dots mark all the eigenvalues for this problem,
and red circles show the AAA poles.}
\end{figure}

Historically, linear and nonlinear eigenvalue problems were
discussed separately, with linear cases founded in classical
tools of numerical linear algebra such as the QR algorithm.
Before the 21st century, there was not very much literature on
nonlinear generalisations.  When more literature began to appear,
the emphasis at first was on tools to approximate nonlinear problems
by linear ones.  For example, the approach known as {\em linearisation\/}
can reduce a nonlinear eigenvalue problem defined by a polynomial
$A(z) = A_0^{} + z A_1^{} + \cdots + z^k A_k^{}$ to an $N(k+1)\times
N(k+1)$ linear eigenvalue problem with block companion matrix
structure \ccite{mmmm}.  The special case of quadratic polynomials,
$k=2$, finds many applications \ccite{tissmeer}.
More generally, one approach to solving nonlinear eigenvalue
problems is to approximate the entries of $A(z)$ by rational
functions with shared poles, and then solve a rational
eigenvalue problem by a linearisation~\cite{setAAA,robust}.

In the era of contour integrals and rational approximations for
large-scale computations, however, the distinction between linear and
nonlinear has become less important.  Early nonlinear contributions based
on contour integrals were \ccite{beyn} and \ccite{asakura}.  Our view
is that rational approximation methods for nonlinear eigenvalue
problems are essentially the same as for linear ones, since both
involve a meromorphic function $A(z)^{-1}$ whose poles one would like
to locate.  For example, the resonance calculation of section \ref{sec:res}
involved a matrix $A(z)$ derived by numerical discretisation of
an integral equation \ccite{resonance}.  Note that $A(z)$ is a
linear operator for each $z$;
the nonlinearity lies in the dependence of $A(z)$ on $z$.

As in the last two sections, a powerful approach to finding the
poles of $A(z)^{-1}$ is to work with a
scalarised or ``sketched'' resolvent
\begin{equation}
R(z)=u^*\kern -1pt A(z)^{-1}v,
\label{nonlinsketched}
\end{equation}
where $u,v\in\complex^N$ are (typically) random vectors; compare
(\ref{res-scalarized}) and (\ref{eig-res1}).
The scalar function $R(z)$ has the same poles as $A(z)^{-1}$,
under the genericity condition that $u,v$ are not orthogonal to the 
left and right eigenvectors of $A(z)^{-1}$.
One can then find the poles of $R$ using AAA.\ \ As in the last
two sections, it is again
straightforward to focus on the eigenvalues in a given domain
$\Omega\subseteq\complex$ by taking sample points on the boundary of $\Omega$.

As a small-scale example of nonlinear eigenvalue calculation by AAA
approximation, Figure \ref{nep-butterfig} shows the ``butterfly''
example from the NLEVP collection, which serves as the logo of
that project \ccite{nlevp}, generated by this code calling the
NLEVP software:

{\small
\begin{verbatim}
      [c,~,A] = nlevp('butterfly');
      e = polyeig(c{:});
      plot(e,'.k'), axis equal, hold on
      Z = 1+1i + 0.3*exp(2i*pi*(1:300)/300); plot(Z,'.b')
      u = randn(64,1); v = randn(64,1);
      f = @(z) u'*(A(z)\v);
      for k = 1:length(Z), F(k) = f(Z(k)); end
      [r,pol] = aaa(F,Z); plot(pol,'or')
\end{verbatim}
\par}

\begin{figure}
\vspace*{6pt}
\begin{center}
\includegraphics[clip, scale=0.48]{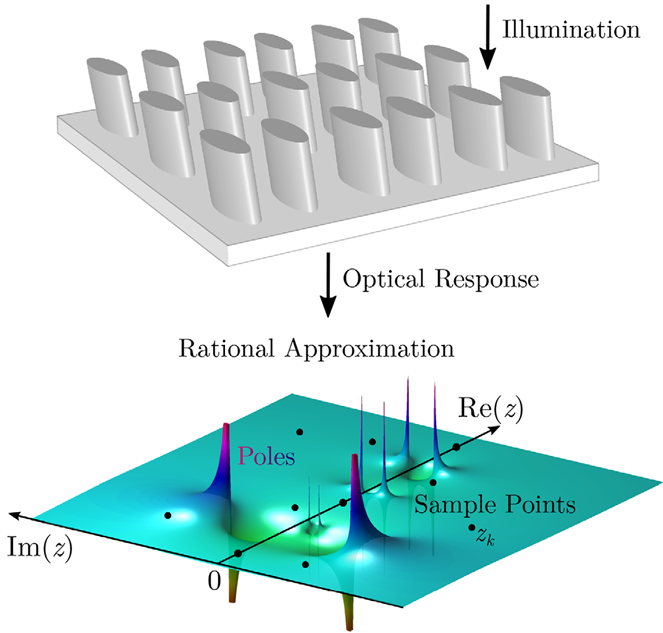}
\end{center}
\vspace*{-8pt}
\caption{\label{nep-betzfig} Image reproduced from \ccite{betz} suggesting
the use of AAA approximation to calculate nonlinear eigenvalues in
photonics applications.}
\end{figure}

At the time of this writing, there is not much experience with
using AAA to solve nonlinear eigenvalue problems via the poles
of the resolvent: the papers
available are \ccite{resonance}, which was highlighted in section
\ref{sec:res}, and \ccite{betz} and \ccite{bink}, which concern applications in photonics.
Figure \ref{nep-betzfig} reproduces a figure from \ccite{betz}
suggesting how rational approximation can be used in an application
of this kind.  In this image, in contrast to the application of
section \ref{sec:res}, one sees that the eigenvalues and poles lie in complex
conjugate pairs.  For a recent review of nonhermitian photonics
and wave physics, see \ccite{kim}.

\section{\label{sec:mor}Model order reduction}

One of the fundamental problems of {\em model order reduction (MOR)}
is to approximate a transfer function of the form
\begin{equation}
\label{eq:transfer_func}
H(s) = C(sE-A)^{-1}B+D,
\end{equation}
usually for $s$ on the imaginary axis.
Here $A$ and $E$ are large-scale $N\times N$ real matrices,
and $C\in\real^{p\times N}$ and $B\in\real^{N\times m}$ with typically $p,m\ll N$. 
This results from a linear dynamical system 
\begin{align*}
  E\dot x(t) &= Ax(t) + Bu(t),\\
 y(t) &= Cx(t) + Du(t) 
\end{align*}
and taking the Laplace transform~\cite[chap.~2]{abg}. 

Here we mainly focus on the Single Input Single Output (SISO) case $p=m=1$. 
The goal of MOR is to approximate $H$ by a simpler function,
typically one involving smaller matrices, i.e., 
\begin{equation}
  \label{eq:Happrox}
H(s)\approx \hat H(s) = \hat C(s\hat E-\hat A)^{-1}\hat B+\hat D,   
\end{equation}
where $\hat A, \hat E\in\real^{n\times n}$,
$\hat C\in\real^{p\times n}$ and $\hat B\in\real^{n\times m}$, with $n\ll N$.
One is usually interested in approximating
$H(s)\approx \hat H(s)$ on the imaginary axis $s\in {\rm i}\real$. 
Since the matrices are real, there is the real symmetry
$H(\bar s) = \overline{H(s)}$.  As discussed in section \ref{sec:AAA}, this
symmetry can be imposed exactly in a AAA code, but Chebfun {\tt aaa.m} does not
do this, so we work on a grid of imaginary points symmetric with respect
to the real axis, and the resulting approximation is nearly real-symmetric
but not exactly.

We have spoken as if the matrices $A, B, C, D, E$ are all known, and this is
the classical setting of MOR, but often this is not the case.  Often one
has an unknown system whose outputs
$H(s)$ can be sampled for various values of $s$, typically on the imaginary
axis, and the aim is to extract information about the system from
this limited information.  In this
context the alternative expressions {\em data-driven modelling\/} and
{\em Reduced Order Modelling (ROM)} are often
used.  For our purposes in this section, the distinctions between MOR and ROM will
not be important.

The transfer function $H$ in~\eqref{eq:transfer_func} has
$p\kern .6pt m$ entries, which are rational functions of degree $N$
sharing the same poles, equal to the eigenvalues of the
generalised eigenvalue problem $Ax=\lambda Ex$. 
The approximation $\hat H$ has entries that are rational functions of degree $n$.
Thus approximating $H$ by $\hat H$ is
an instance of a rational approximation problem,
where the original function is also rational (and matrix-valued), but of higher degree. 

When $p=m=1$, the problem reduces to a scalar rational approximation problem,
which AAA is well able to handle. 
One can sample $H$ at points on the imaginary axis and find an approximation 
\begin{equation}
H(s) \approx \hat H(s)=\sum_{k=0}^n 
\frac{ f_k^{}\kern .5pt \bk}{s- t_k^{}}\left/ \sum_{k=0}^n
\frac{ \bk}{s- t_k^{}}\right.
\end{equation}
as in (\ref{aaa-baryrep}).
One could take this to be the reduced transfer function.
Alternatively, one can use the poles and residues to
convert $\hat H(s)$ to transfer function form.   The matrices involved 
will be diagonal in the most obvious realisation, but an interesting
alternative is

\begin{equation}
\hat H(s) = [f_0\beta_0,\ f_1\beta_1,\ \ldots,\  f_n\beta_n]
\kern 2ptA^{-1}
\begin{bmatrix}
  1 \\ 0 \\0\\ \vdots \\0
\end{bmatrix}
\kern -1pt 
\end{equation}
with
\begin{equation}
A = \begin{bmatrix}
\beta_0&\beta_1&\ldots & \beta_{n-1} & \beta_n\\
s-t_0& t_1-s & \\
& s-t_1& t_2-s & \\
& & \ddots & \ddots  & \\
&&& s-t_{n-1}& t_{n}-s & 
\end{bmatrix},
\end{equation}
as proposed in \cite[Prop.~2]{setAAA}.

\begin{figure}
\begin{center}
\vspace*{8pt}
\includegraphics[trim=20 65 38 4, clip, scale=.48]{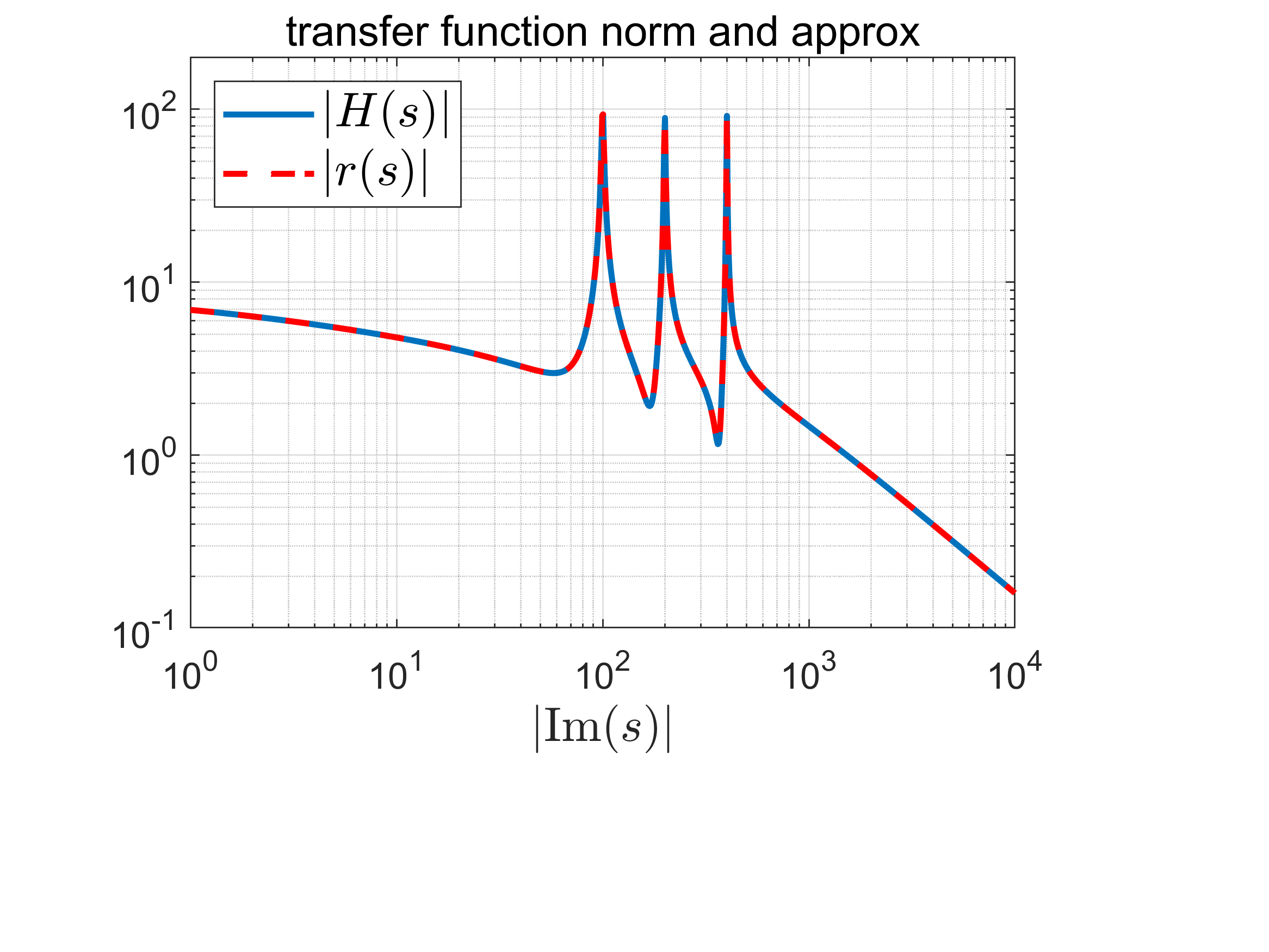}\includegraphics[trim=10 65 38 4, clip, scale=.48]{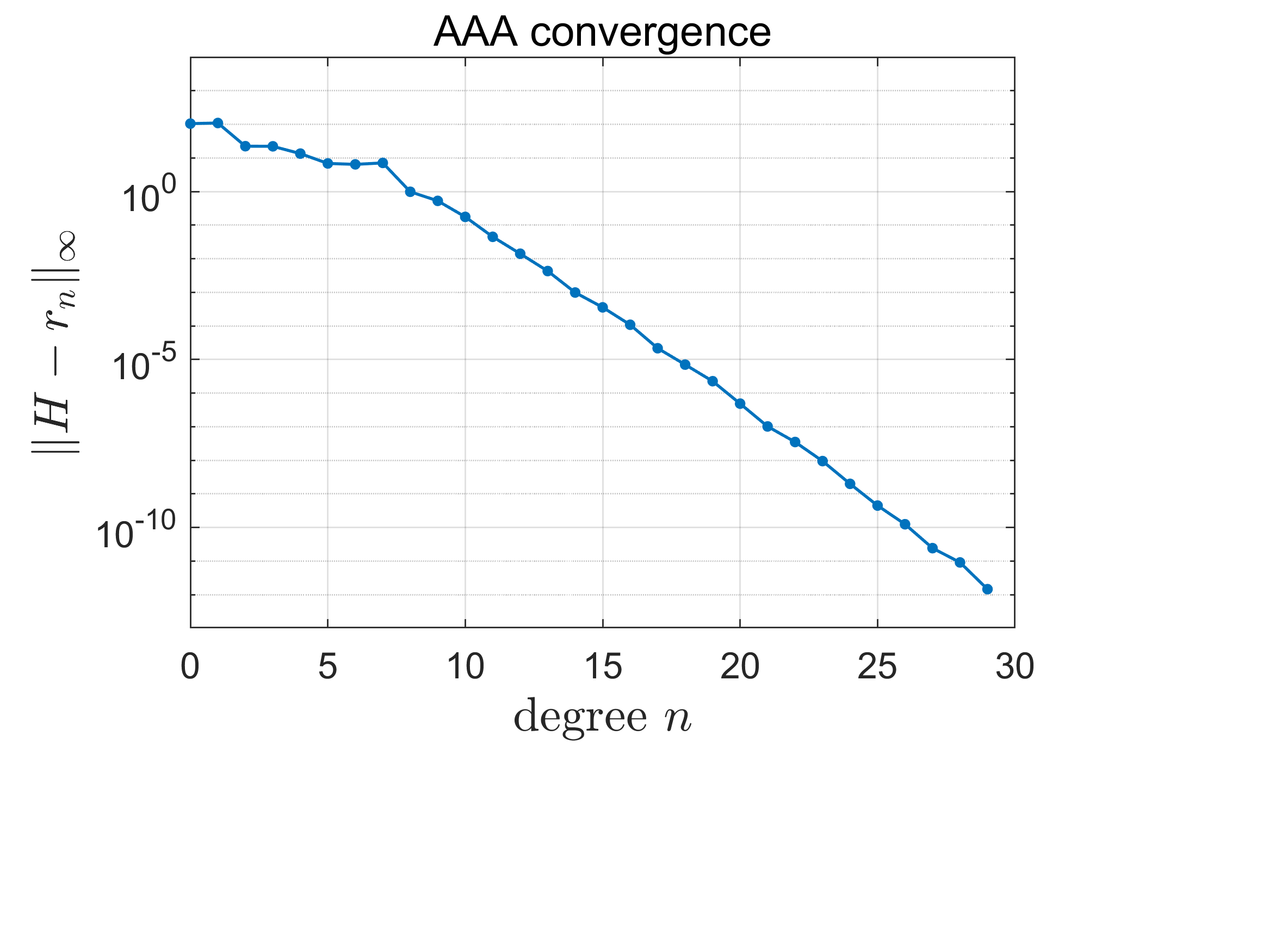}
\caption{\label{fig:mor}On the left, absolute values of the degree $1006$
transfer
function $H(s)$ and its degree 29 AAA approximant $r(s)$ for the FOM example
from SLICOT \ccite{slicot}.
On the right, AAA convergence for this example.  Rapid convergence begins at degree $n=6$, 
after the six poles at $-1\pm 100\kern .4pt i$, $-1\pm 200\kern .4pt i$,
and $-1\pm 400\kern .4pt i$ are captured.}
\vspace*{-6pt}
\end{center}
\end{figure}

To illustrate, we take Penzl's FOM model from the SLICOT
Benchmark examples collection~\ccite{slicot},
where $A\in\real^{1006\times 1006}$ and $p=m=1$,
so $H$ is a rational function of degree 1006.
We sample $H(s)$ at 100 logarithmically spaced points
in the interval $[\kern .3pt i,10^4\kern .3pt i]$ and their
complex conjugates and compute its
AAA approximation with the usual 13-digit relative accuracy,
which has degree 29.
Figure~\ref{fig:mor} shows $|H(s)|$ and $|r(s)|$ and the convergence
of $r$ to $H$ as the AAA iterations proceed. 
The two functions match to 10 digits of absolute accuracy on the whole imaginary
axis.

In the Multiple Input Multiple Output (MIMO) case,
with $p,m>1$, one can use the blockAAA or setAAA algorithms to
find approximations to all $p\kern .6pt m$ entries
simultaneously~\ccite{goseaalgs,setAAA}; see section \ref{sec:sets}. 

Classically, vector fitting~\ccite{vf}
was for many years among the most popular methods for
finding a rational approximation to a transfer function $H$.
Unlike AAA, in vector fitting the degree of $\hat H$ needs to be
specified in advance.  The convergence of vector fitting is not fully understood. 

So far as we are aware, compared with other MOR/ROM approaches,
AAA is often as efficient as any method for scalar problems with $p=m=1$,
and perhaps similarly for setAAA when $p,m$ are small. 
AAA has been applied in a variety of contexts in the scalar case,
including~\ccite{goseadata}
(higher-degree problems) and \ccite{pAAA} (parameterised problems). 
When $p$ and/or $m$ are large, however,
the situation can be different, and 
AAA-based methods in MOR/ROM are less mature. Studies in this direction
include~\ccite{aumann,goseaalgs,jonas}.

An important class of methods in MOR/ROM is the so-called {\em Loewner framework}.
In this method (as in AAA),
one only requires evaluations of the transfer function $H(s)$,
sometimes only {\em tangentially},
e.g.\ evaluations of $H(s_i)b_i$ and
$c_j^TH(s_j)$ for vectors $b_i$ and $c_j$.
Then one can use the sampled data to form matrices related
to the Loewner matrix (\ref{aaa-Aentries}) to
obtain a reduced model $\hat H$ that satisfies tangential interpolation,
i.e., $H(s_i)b_i=\hat H(s_i)b_i$, $c_j^TH(s_j)=c_j^T\hat H(s_j)$. 

Again, for scalar problems $p=m=1$,
AAA tends to be more efficient than Loewner,
and it benefits from the robustness brought by the least-squares part of 
AAA.\ \ By contrast, Loewner can be better for $p,m\gg 1$.
For example, with $p=m=n$ and $D=0$, the Loewner framework
is able to recover the system to obtain $\hat H = H$ by 
computing $H(s)$ at just two points,
$H(s_1)$ and $H(s_2)$\footnote{This may seem surprising at first,
but one can note that $H=C(sE-A)^{-1}B$
essentially contains $2\kern .3pt n^2$ degrees of freedom,
which are provided by $H(s_1)$ and $H(s_2)$.}.
AAA or setAAA would not be able to do this. 

When the matrices $A,B,C,D,E$ are known explicitly, 
a common approach in MOR for finding $\hat H$ from $H$
is {\em projection}: one finds orthonormal (tall-skinny)
matrices $V,W$ and sets $\hat E = W^T\!EV$, $\hat A = W^T\!EV$,
$\hat B = W^T\!B$, $\hat C = CV$ and $\hat D = D$. 
Remarkably, tangential interpolation ensues~\cite[chap.~3]{abg}.
For example,
for any fixed $b\in\real^m$,
$H(s_i)b=\hat H(s_i)b$ holds if $\mbox{Range}(V)$ includes
the vector $(s_iE-A)^{-1}Bb$. 

AAA (and AAA-Lawson) attempt to minimise the
$L_\infty$-norm of the error $H-\hat H$.
In the scalar case,
this is equivalent to the so-called $\mathcal{H}_\infty$ norm,
defined by $\|H\|_{\mathcal{H}_\infty}=\sup_{s\in \real}\|H({\rm i}s)\|_2$.
It is unclear how to do the same when $p,m>1$. 
In MOR, attention has also been paid to minimising the 
$\mathcal{H}_2$ norm of the error.
Extensive theory and optimality conditions have been developed,
and in particular IRKA (the Iterative Rational Krylov 
Algorithm)~\ccite{irka}
iteratively updates the tangential interpolation points and 
directions to find the best approximation $\hat H$ to $H$ in $\mathcal{H}_2$.
IRKA is applicable to MIMO systems with $p,m>1$
and converges rapidly for most problems,
although a full convergence proof is not available. 
For ROM as opposed to MOR, there is TF-IRKA (``transfer function IRKA'')
\ccite{TFIRKA}.

Approximating $A(z)$ can be more efficient than working
with the resolvent for problems where solving linear systems with respect to $A(z_i)$
is expensive.  On the other hand,
working with the resolvent has the advantage of simplicity
and having no distinction between linear and nonlinear problems.

For more details on MOR/ROM,
we refer to~\ccite{abg,benner} and the references therein. 
For practical information about AAA applied to ROM problems, see
the code {\tt rational} in the MathWorks RF Toolbox \ccite{rf},
which can handle multiple outputs as well as scalars.
It should be noted that applications such as those
targeted by the Toolbox often have significant amounts of noise, an issue
discussed here at the end of section~\ref{sec:equi}.  In such cases an
engineer may be pleased with two digits of accuracy, far less than the
ten digits of the SLICOT example we showed.
More recently, MathWorks has also introduced AAA fitting in its
Control System and System Identification toolboxes (version R2025a).
The applications here are closer to MOR as opposed to ROM.

\section{\label{sec:zolo}The Zolotarev sign problem}

Yegor Zolotarev was a student of Chebyshev who visited Berlin
in 1872 and learned about elliptic functions from Weierstrass.
He brought these ideas back to St.~Petersburg, where he published a
paper in 1877 on the theme of polynomial and rational approximations
on disjoint intervals \ccite{zolotarev}.  His work still resonates
today, and there is no knowing what more Zolotarev might have done
had he not been killed in a train accident in St.\ Petersburg in 1878,
at age 31.

Let $E$ and $F$ be disjoint closed sets in the complex plane,
and let $n\ge 0$ be an integer.  The {\em Zolotarev sign problem\/}
(also known as {\em Zolotarev's 4th problem\/}) is to find a rational
function $\rnh$ of degree $n$ that approximates $-1$ on $E$ and $+1$
on $F$ as closely as possible.
We denote this target function by
$\signEF(z)$ (not defined for values of $z$ outside $E$ and $F\kern
1pt$), so the problem can be stated as finding an optimal approximation
$\rnh(z)\approx\signEF(z)$ on $E\cup F$.
The minimum approximation error is denoted by $\tau_n^{}$:
\begin{equation}
\tau_n^{} = \|\rnh - \signEF\|_{E\cup F}^{} = \min_{r\in\Rn} \|r - \signEF\|_{E\cup F}^{},
\label{zolo1-tau}
\end{equation}
where $\|\cdot\|_{E\cup F}^{}$ denotes the $\infty$-norm on $E\cup F$.
A function $\rnh$ achieving this minimum
is an example of a rational function whose poles and zeros will delineate an
approximate branch cut, because no single analytic function
can take values $-1$ on $E$ and $+1$ on $F$.  (We assume that $E$
and $F$ are continua with no isolated points, though they need not
be connected.)

\begin{figure}
\vspace*{2pt}
\begin{center}
\includegraphics[trim=30 75 30 57, clip, scale=0.85]{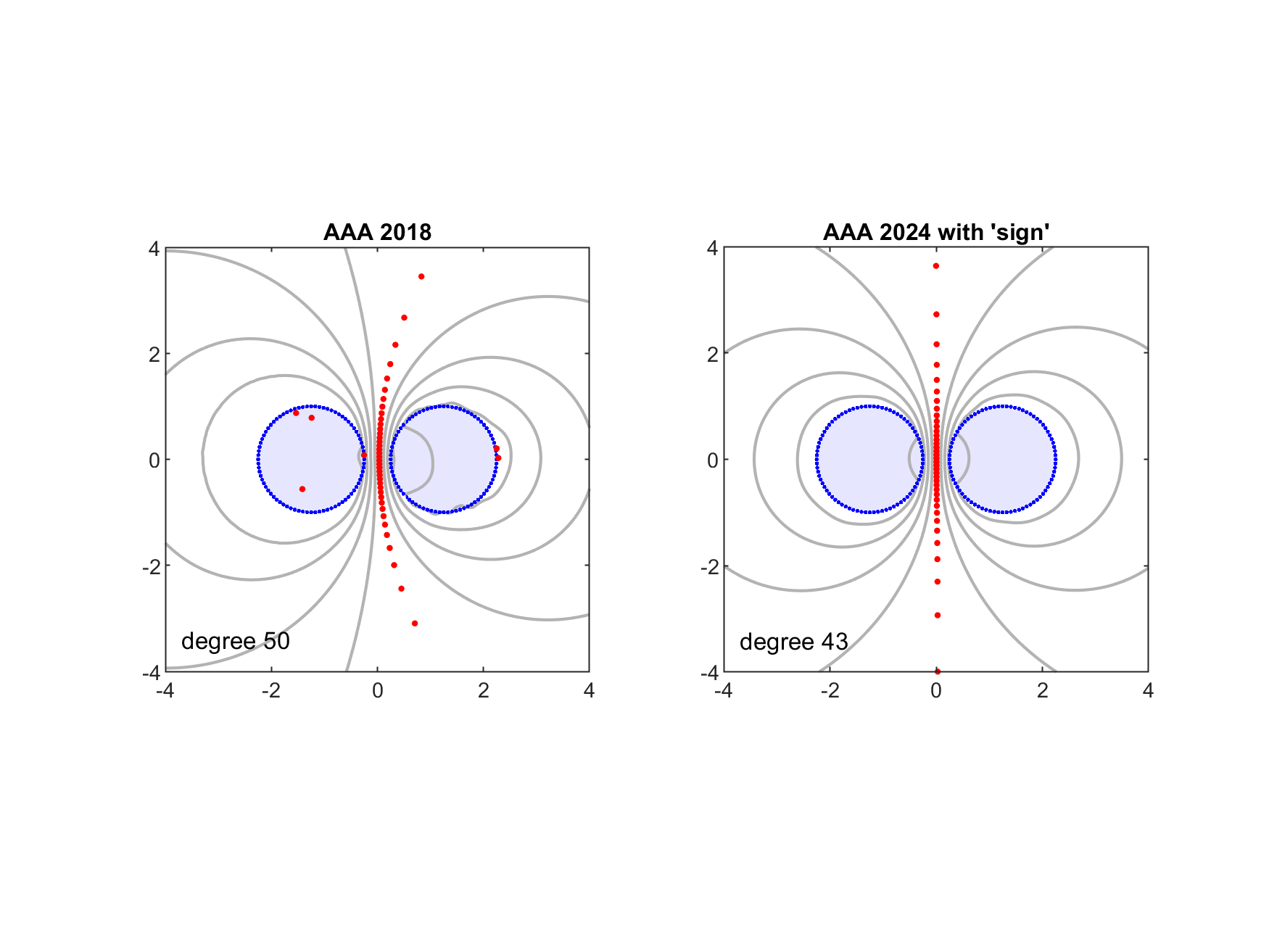}
\end{center}
\vspace*{-9pt}
\caption{\label{zolo-pairfig} Approximation of \smash{$\signEF(z)$} for
a pair of disks discretised by $80$ equispaced points on
each boundary.  Standard AAA does poorly on such problems (left),
whereas AAA with the \texttt{\textquotesingle sign\textquotesingle} modification introduced in
2024 and described in
\ccite{zolo} does much better (right).   The grey contours mark
deviations
$10^{-12}$, $10^{-9}$, $10^{-6}$ and $10^{-3}$ 
of $r(z)$ at each point $z$ from
$+1$ or $-1$, whichever is closer.}
\end{figure}

From the beginning of our AAA-related work, we tried experiments
with approximating sign functions, and the results appeared adequate
but not excellent.  The left image of Figure \ref{zolo-pairfig}
gives an illustration for a case where $E$ and $F$ are the circles
of radius $1$ about $z=-1.25$ and $z=+1.25$.  (This problem
has an analytic solution given below in 
(\ref{zolo2-exactrns})--(\ref{zolo2-exactrnh}).)  Two features of this
plot fall short of ideal.  One is the six spurious poles near
the boundary circles, but these can be argued away as artifacts
(Froissart doublets) since their residues are negligible, and they
can be removed numerically by a least-squares fit of $\signEF(z)$
in the space spanned by the remaining poles.  The other is the fact
that the remaining poles lie on a curve that bends away from the
imaginary axis, where by symmetry, the optimal poles should lie.
We interpreted this effect in the spirit of backward error analysis.
Many different pole configurations will give almost equally good
approximations on $E$ and $F$, we reasoned, so the pole locations
are exponentially ill-conditioned and can be expected to show some
oddities in a numerical computation.

Then, in the summer of 2024, together with Heather Wilber,
we investigated the situation more closely in order to pursue
applications to be described in the next section
\ccite{zolo}.  We found that the difficulties of our ``adequate''
sign function calculations were more fundamental than we had
realised.  In fact, approximation of sign functions using the AAA
algorithm as described in section~\ref{sec:AAA} almost invariably runs into
trouble.  The initial iterations tend to show no progress at all,
completely failing to find good approximations even though these
must exist even for low degrees like $1$ or $3$.

The details are not fully understood, but it became clear that
a source of the problem is that in approximating $\signEF(z)$,
AAA works with a Loewner matrix $A$ half of whose entries
(\ref{aaa-Aentries}) are zero.  This happens because $\fj$ and
$\fk$ are equal when $j$ and $k$ correspond both to $E$ or both to
$F$, and it makes $A$ a rectangular matrix with block $2\times 2$
diagonal structure.\footnote{We are grateful to Sam Hocking for
some analysis of this effect.} When AAA then computes a minimal
singular vector of $A$ to obtain barycentric weights, as described in
section~\ref{sec:AAA}, the blocks decouple and one obtains a weight
vector that more or less randomly attends to approximation on one
of $E$ and $F$ while ignoring the other.  Convergence is disrupted.

A modification to AAA was proposed in \ccite{zolo} that goes a
good way toward fixing this problem.  For calculating a new vector
{\tt wj} of barycentric weights at each step, we changed the MATLAB
lines in {\tt aaa.m}\medskip

{\small
\begin{verbatim}
      [~,S,V] = svd(A(J,:),0); 
      wj = V(:,end);
\end{verbatim}
\par}

\noindent to

{\small
\begin{verbatim}
      [~,S,V] = svd(A(J,:),0); 
      s = diag(S); wj = V*(1./s.^2); wj = wj/norm(wj);
\end{verbatim}
\par}
\medskip

\noindent
This ensures that {\tt wj} is a blend of all the singular vectors
rather than just one of them.  For most AAA computations, the
change has little effect, but for approximating sign functions,
the difference is considerable, as one sees in the right image of
Figure \ref{zolo-pairfig}.  As of this writing, the Chebfun code {\tt
aaa.m} offers \verb|'sign'| as an option (not turned on by default).

Another similar modification of AAA has been proposed in \ccite{mitchell}.

\begin{figure}
\vspace*{2pt}
\begin{center}
\includegraphics[trim=30 75 30 57, clip, scale=0.85]{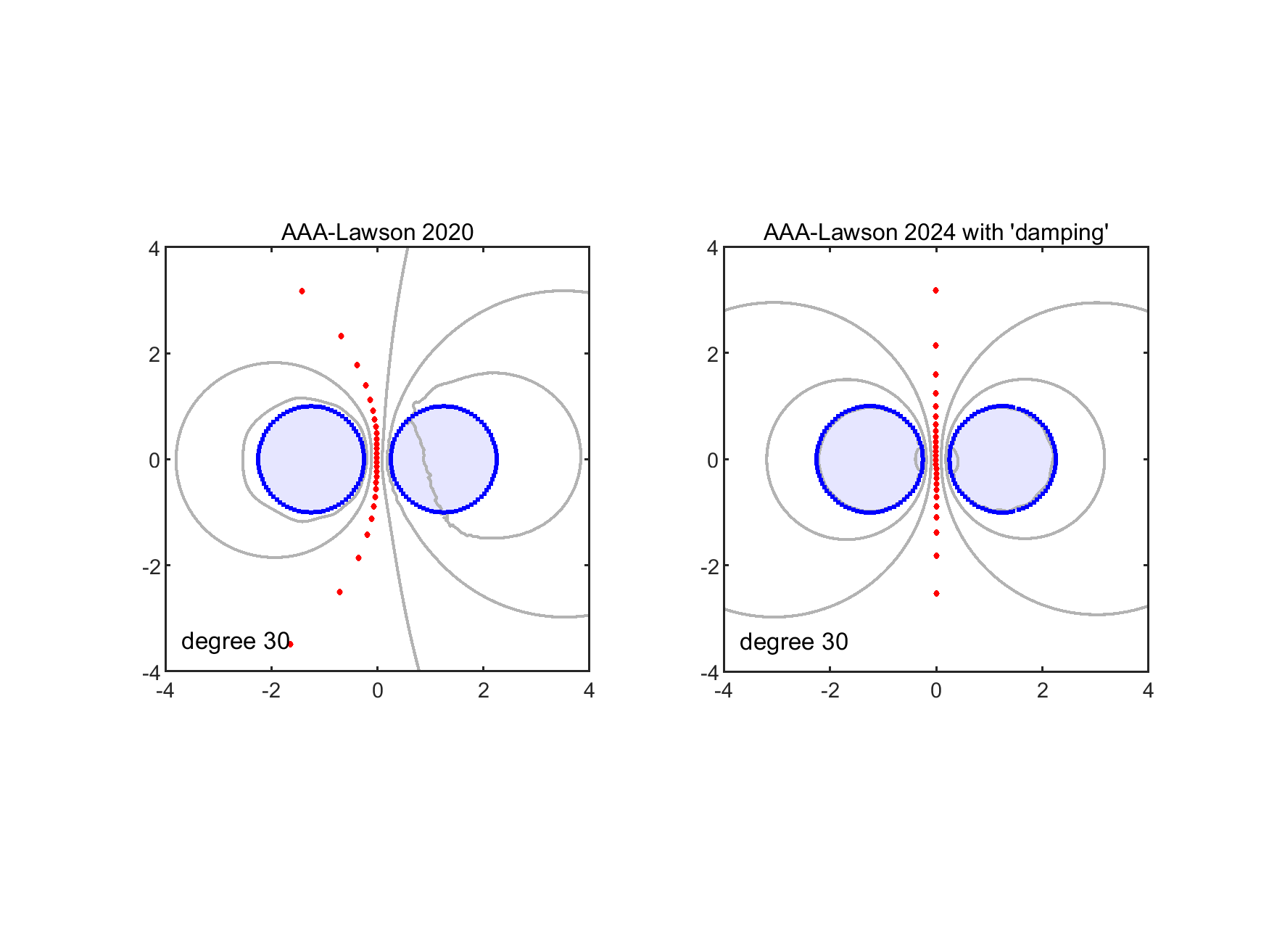}
\end{center}
\vspace*{-9pt}
\caption{\label{zolo-pairfiglawson} Repetition of
Figure \ref{zolo-pairfig} but now for degree $30$ approximations
with Lawson iteration turned on to attempt an improvement
from near-best to best approximation.
The contours mark accuracies $10^{-9}$, $10^{-6}$ and $10^{-3}$.
The plot on the right
shows the improvement both the {\tt \textquotesingle sign\textquotesingle } and 
the {\tt \textquotesingle damping\textquotesingle }
modifications introduced in 2024 \ccite{zolo}.}
\end{figure}

Fixing one problem in AAA sign function approximations unexpectedly
revealed another, illustrated in the left image of
Figure \ref{zolo-pairfiglawson}.  For the applications of the
next section, it was important for us
to be able to improve our Zolotarev approximations from near-best to 
best.  The standard algorithm for this is the AAA-Lawson method \ccite{aaaL}, 
which is a nonlinear
iteratively reweighted least-squares iteration in which at each
step, the weights of a least-squares problem are updated
by the formula
\begin{equation}
w_j^{\rm (new)} = w_j |e_j|,
\label{zolo-lawson}
\end{equation}
where $w_j$ denotes the current least-squares weight at sample
point $z_j$ and $e_j$ is the current rational approximation error
$e_j = r(z_j)-f_j$.  When the iteration is working well, the weights
converge steadily (linearly) to a configuration corresponding to the
best $L^\infty$ approximant, but we had long known that it sometimes
fails to converge, usually then exhibiting period-2 oscillations.
To our surprise, it emerged that $\signEF$ problems are a case
where these failures are exceptionally likely to occur.

And so we introduced a modification in the AAA-Lawson algorithm, which
usually (not always) fixes this second problem.  When the \verb|'damping'|
flag is specified, (\ref{zolo-lawson}) is changed to 
\begin{equation}
w_j^{\rm (new)} = \left((1-\delta) + {\delta \kern 1pt |e_j|\over \max_j |e_j|}\right) w_j,
\label{zolo-lawson2}
\end{equation}
where $\delta\in (0,1]$ is a user-supplied damping factor.
With $\delta =1$ we have standard AAA-Lawson, whereas for smaller
values it is a more robust iteration that converges more reliably,
if more slowly.  For details, see \ccite{zolo}.

As a typical example, the approximation plotted in the right image
of Figure \ref{zolo-pairfiglawson} can be generated by the code

{\small
\begin{verbatim}
     Z = exp(2i*pi*(1:80)'/80); Z = [1.25+Z; -1.25+Z];
     F = sign(real(Z));
     [r,pol] = aaa(F,Z,'degree',30,'sign',1,'damping',.95);
\end{verbatim}
\par}

Since the emphasis of this paper is on applications more than
algorithms, we will not say more about the \verb|'sign'| and
\verb|'damping'| modifications.  In both cases, there is a clear
improvement, but the modifications are not bulletproof,
and theoretical justifications are entirely lacking.  We
hope that these aspects of AAA will be made more solid by further
research.

Instead, in the rest of this section we will show four more examples
and then comment on applications.

\begin{figure}
\vspace*{-5pt}

\begin{center}
\includegraphics[trim=30 68 30 57, clip, scale=0.85]{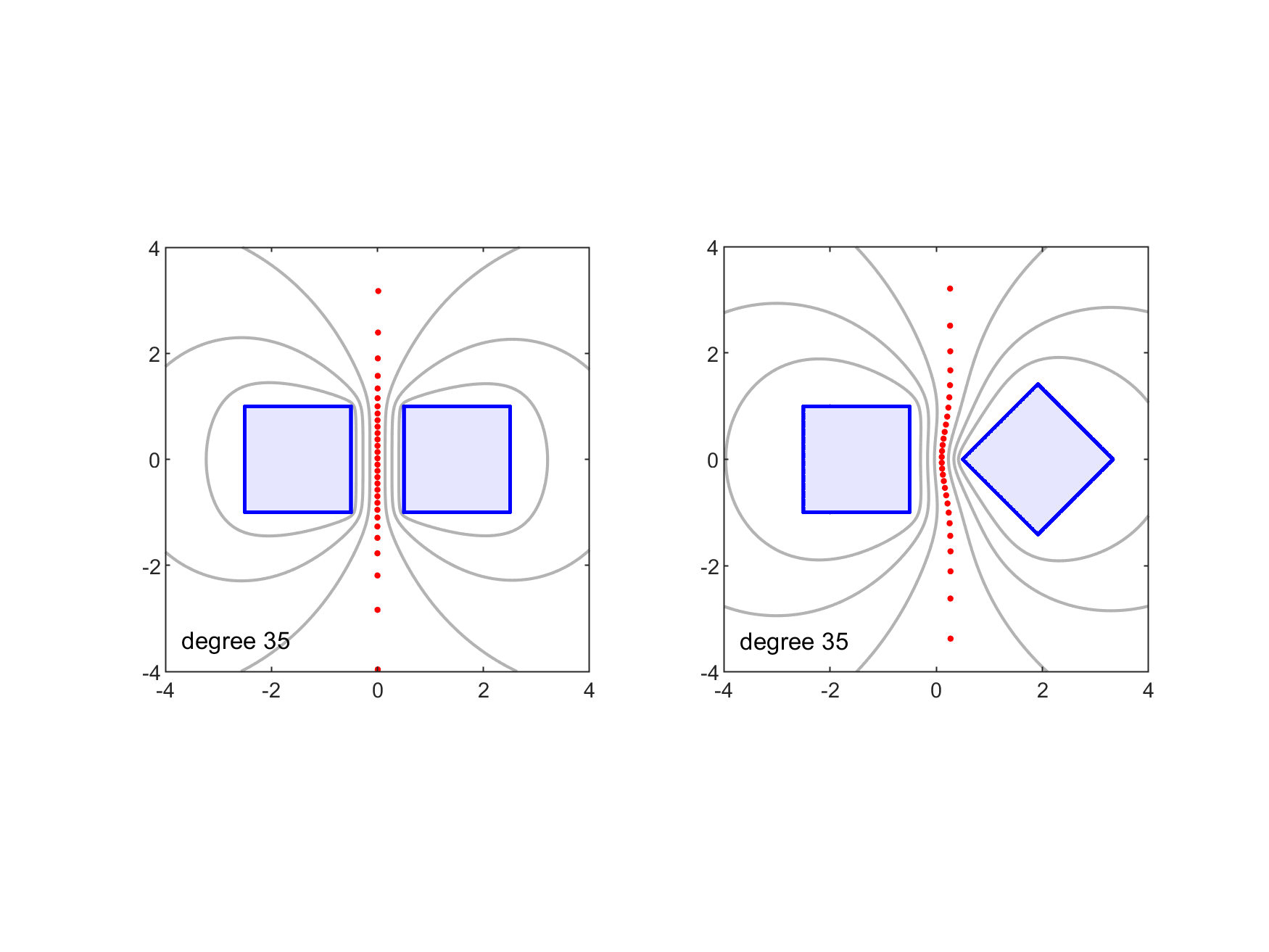}
\vspace*{-10pt}

\includegraphics[trim=30 75 30 67, clip, scale=0.85]{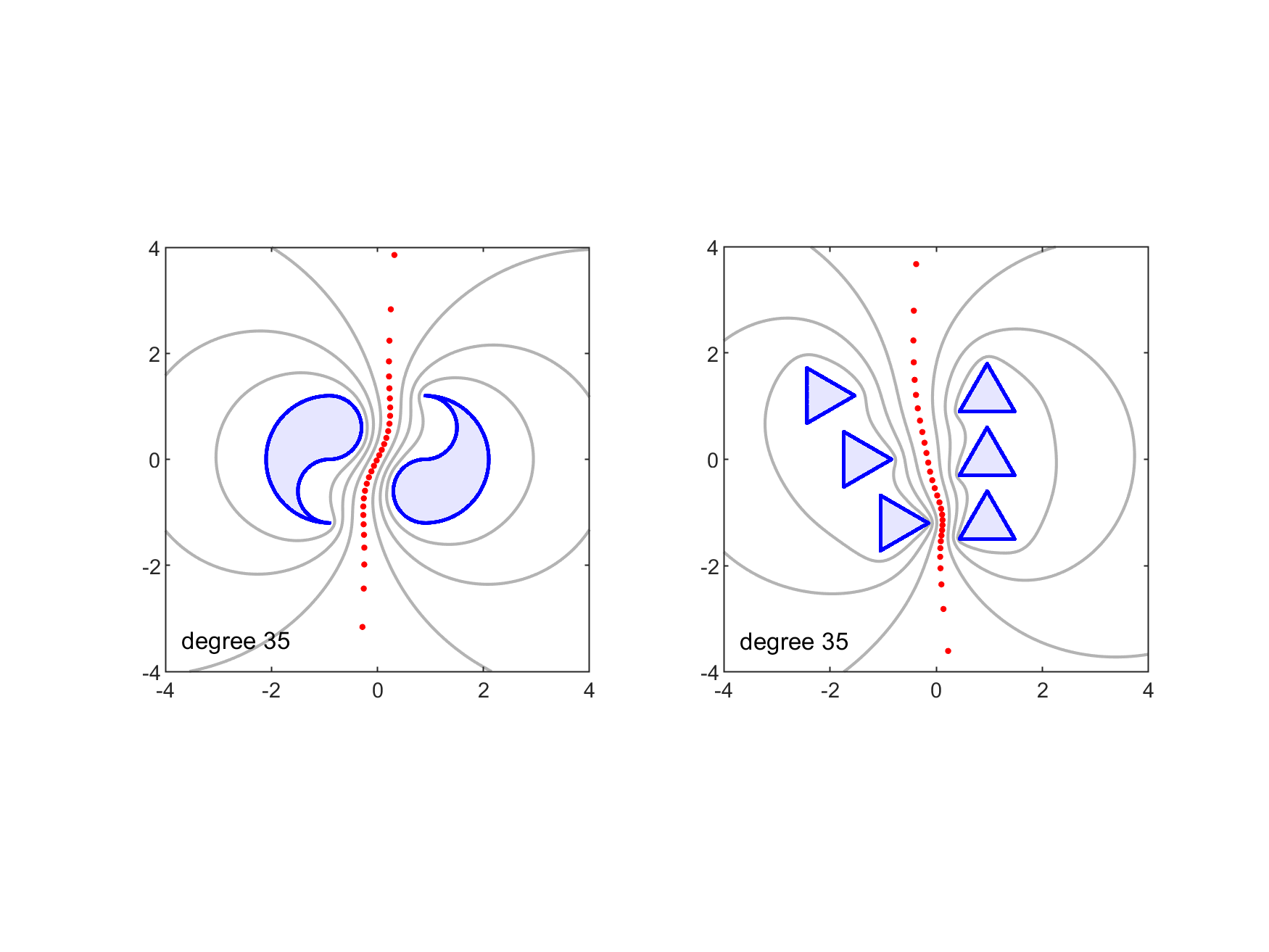}
\end{center}
\vspace*{-10pt}

\caption{\label{zolo-fig3} Four more examples of AAA-Lawson minimax solutions of
the Zolotarev sign problem, all for fixed degree $35$.  The {\tt \textquotesingle sign\textquotesingle } and
{\tt \textquotesingle damping\textquotesingle} options have been used with a damping factor $\delta = 0.95$.
The contours mark accuracies $10^{-9}$, $10^{-6}$ and $10^{-3}$.}
\end{figure}

The examples appear in Figure \ref{zolo-fig3} and illustrate, among
other things, that the sets $E$ and $F$ need not be connected.
It is interesting to take note of the grey contours in this figure
and also in Figures \ref{zolo-pairfig} and \ref{zolo-pairfiglawson}.
In all of these plots, the red dots mark the approximate branch
cut selected by AAA, and the grey contours closest to the dots show
accuracies $10^{-3}$ on two sides.  In each case we see that $r(z)$
matches $+1$ or $-1$ to three digits in large regions in the left
and right half-planes, even though {\em a priori,} it was only
tasked with approximation on $E$ and $F$.  This is the phenomenon
that Walsh called {\em overconvergence\/} \ccite{walsh}, which appears widely
in both polynomial and rational approximation and can be seen in
several of the figures of this review.

The Zolotarev sign problem arises in a number of applications where
one wants to separate one part of a system computationally from
another.  For example, classic Infinite Impulse Response digital
filters start from a rational function that is nearly constant
on one part of the real axis and nearly zero on another, and
computing coefficients for such filters is essentially a Zolotarev
sign problem~\ccite{os}.  Generalisations arise in numerical linear
algebra and computational science, for example in the computation
of functions of matrices and operators (\cite{bai,highamfA,kenneylaub}).
Divide and conquer methods for
computing eigenvalues of large matrices work by suppressing one part
of the spectrum relative to another recursively~\ccite{bdg,banks}.
This too leads to Zolotarev sign problems, typically starting from
the case where $E$ and $F$ are approximations to the left and
right complex half-planes~\ccite{nakfreund}.  In the real case,
which goes back to Zolotarev himself, they may be intervals such as
$[-a,-\varepsilon]$ and $[\varepsilon,a]$.  ``Spectral slicing''
ideas of this kind have found generalisation through algorithms
such as FEAST, which was discussed in section \ref{sec:eig} \ccite{guettelzolo}.

We just mentioned the possibility of recursion in algorithms related
to $\signEF$ approximations, and this is an important aspect of this
subject both mathematically
and philosophically.  In classical rational approximation theory, one seeks a 
rational function of a specified degree $n$ with certain properties.
(Everything we say here also applies to polynomials.)
For example, the approximations of Figure \ref{zolo-fig3} are rational
functions of degree 35.
Another possibility in many applications, however, is to compose
rational functions in structures like this:
\begin{equation}
r(z) = r_k^{}(\cdots r_3^{}(r_2^{}(r_1^{}(z)))\cdots).
\label{zolo-composite}
\end{equation}
If each element has degree $2$, for example, then this {\em
composite rational approximation\/} has degree $2^k$, which may be
far beyond the range that numerical approximation theorists are used
to imagining they can deal with~\ccite{nakfreund,gawlik}.  Though we
don't generally notice, such composite rational approximations are
exploited implicitly all the time when Newton's method is applied
\ccite{everything}.  In the case of $\signEF$ approximations
where $E$ and $F$ are disjoint, even a two-level approximation
may be powerful, for a single rational function may produce sets
$r(E)$ and $r(F)$ that are small neighborhoods of $-1$ and $+1$,
respectively, hence well separated from each other, and applying
one further polynomial or rational function can then add many more
digits of accuracy at minimal cost.

\section{\label{sec:zolo2}The Zolotarev ratio problem}

The {\em Zolotarev ratio problem\/} (also known
as {\em Zolotarev's 3\kern .3pt rd problem\/}) is another problem involving rational
functions on a pair of sets $E$ and $F$ in the complex plane.
This time the question is, how small can $r$ be on $E$ relative
to its size on $F$\kern 1pt ?  Specifically, the aim is to find a degree
$n$ rational function $\rns$ that minimises the ratio
\begin{equation}
{\max_{z\in E} |r(z)|\over \min_{z\in F} |r(z)| } ,
\label{zolo2-ratio}
\end{equation}
and we denote the minimum ratio, the {\em Zolotarev number\/} of
$E$ and $F$, by $\sigma_n$:
\begin{equation}
\sigma_n = \min_{r\in\Rn} {\max_{z\in E} |r(z)|\over \min_{z\in F} |r(z)| } .
\label{zolo2-sig}
\end{equation}
To normalise the problem in a symmetric fashion we specify
\begin{equation}
\min_{z\in F} |r(z)| =  {1\over\max_{z\in E} |r(z)|}.
\label{zolo2-norm2}
\end{equation}
(This is different from the normalisation of \ccite{istacez}
and \ccite{zolo}.)  Like the Zolotarev sign problem of the last
section, the Zolotarev ratio problem has a number of applications,
which we shall mention at the end of the section.

Like the sign problem, the ratio problem was introduced (for real
intervals $E$ and $F$) by Zolotarev \ccite{zolotarev}, and Gonchar
generalised it to the complex plane \ccite{gonchar}.  As it happens,
the two problems are equivalent, but it is not clear if Zolotarev
knew this.  The equivalence was shown by Achieser for real intervals
and by Istace and Thiran for more general sets \ccite{istacez},
and it runs as follows.  Let $E$, $F$ and $n$ be given, and let
$\rnh$ be a degree $n$ best rational approximation of the $\signEF$
function with error $\tau_n^{}$ as in (\ref{zolo1-tau}).  We state
the theorem of Istace and Thiran essentially in their words:

\begin{theorem}
\label{zolo2-thm}
Every solution $\rns, \sigma_n$ of the Zolotarev ratio
problem with
the normalisation (\ref{zolo2-norm2}) is related to the solution 
$\rnh,\tau_n^{}$ in (\ref{zolo1-tau}) of the Zolotarev sign problem by
\begin{equation}
\rnh(z) = \gamma_n \,{\rns(z) - \sigma_n
\over \rns(z) + \sigma_n}, \quad
\rns(z) = \sqrt{\sigma_n} \,\kern 1pt
{\gamma_n + \rnh(z) \over \gamma_n - \rnh(z)}
\label{zolo2-relation}
\end{equation}
with $\gamma_n = (1-\sigma_n)/(1+\sigma_n)\approx 1$.
The minimal values of the two problems satisfy
\begin{equation}
\tau_n^{} = {2\sqrt{\sigma_n}\over 1 + \sigma_n}, \quad
\sigma_n = \left({\tau_n^{}\over 1 + \sqrt{1-\tau_n^2}\kern .7pt} \right)^2,
\label{zolo2-minvals}
\end{equation}
and the set of extremal points
\begin{equation}
M = \{z\in E\cup F, ~|\rnh(z) - \signEF(z)| = \tau_n^{}\}
\label{zolo2-sets}
\end{equation}
is the union of 
$M_1 = \{z\in E, ~|\rns(z)| = \sigma_n\}$
and $M_2 = \{z\in F, ~|\rns(z)| = 1\}$.
\end{theorem}

Thanks to this theorem, once the Zolotarev sign problem has been
solved, one has in principle also solved the Zolotarev ratio problem.
Until 2024, this was true only in principle (except for a small
collection of analytic solutions), since no method was available to
compute solutions to the sign problem accurately enough.  Now that such a
method is available, the ratio problem can be handled as follows:
\medskip

(1) Solve the sign problem by AAA-Lawson as in the last section;

(2) Convert to a solution of the ratio problem by (\ref{zolo2-relation})
and (\ref{zolo2-minvals}).
\medskip

\noindent Details are presented in \ccite{zolo}.

\begin{figure}
\begin{center}
\includegraphics[trim=30 110 65 10, clip, scale=0.92]{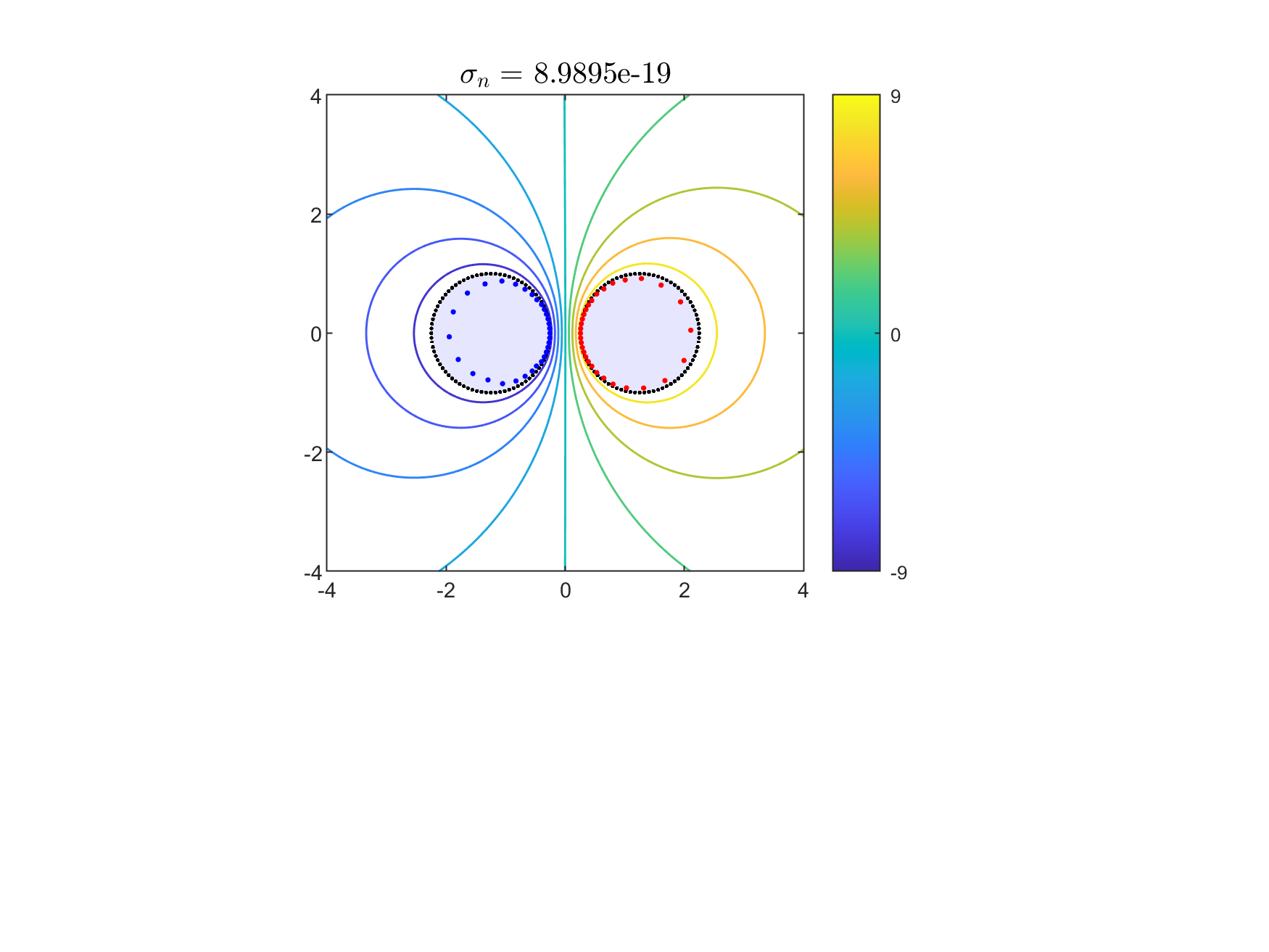}
\caption{\label{zolo2-fig1} Approximation to
the degree $n=30$ solution $\rns$ of the Zolotarev ratio
problem for the two-disks configuration of Figures \ref{zolo-pairfig}
and \ref{zolo-pairfiglawson}.  This is a transformation of
the function $\rnh$ of the right image of
Figure \ref{zolo-pairfiglawson} by the formula
(\ref{zolo2-relation}).  Blue and red dots mark zeros and poles of
the computed $\rns$ (the exact results would coalesce at two points),
and the contours mark level curves of $\log_{10}\kern -.7pt |\rns(z)|$.}
\end{center}
\end{figure}

\begin{figure}
\begin{center}
\includegraphics[trim=30 110 65 10, clip, scale=0.92]{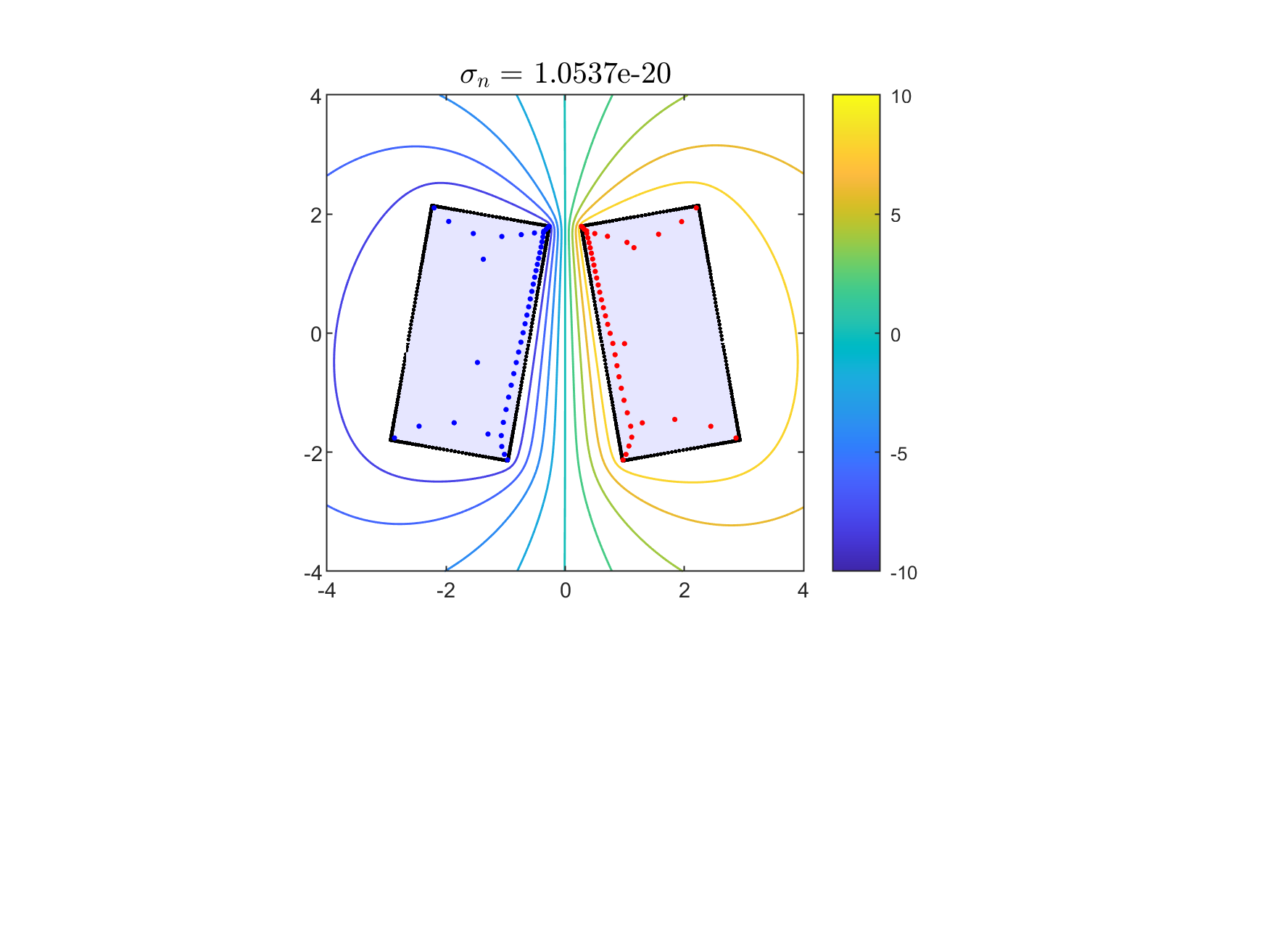}
\caption{\label{zolo2-fig2} Degree $44$ solution of the Zolotarev
ratio problem for a pair of $2\times 1$ rectangles
centred at $z = \pm 1.6$ and tilted by 10 degrees.  This
computation took 1.5 s on a laptop.}
\end{center}
\end{figure}

For example, the following code computes the solution to the
degree 30 Zolotarev ratio problem for the geometry of Figures
\ref{zolo-pairfig} and \ref{zolo-pairfiglawson}, where $E$ and $F$
are disks of radius $1$ about $z = \pm 1.25$.  This formulation
calls the Chebfun code {\tt prz} (which is also included in {\tt
aaa.m}) to find the zeros and poles of $\rns$.

{\small
\begin{verbatim}
      Z = exp(2i*pi*(1:80)'/80); Z = [1.25+Z; -1.25+Z];
      F = sign(real(Z));
      [rhat,~,~,~,zj,fj,wj] = ...
         aaa(F,Z,'degree',30,'sign',1,'lawson',200,'damping',0.95);
      tau = norm(F-rhat(Z),inf);
      sigma = (tau/(1+sqrt(1-tau^2)))^2;
      gamma = (1-sigma)/(1+sigma);
      rstar = @(z) (gamma+rhat(z))./(gamma-rhat(z));
      [~,~,zeros] = prz(zj,fj+gamma,wj);
      [~,~,poles] = prz(zj,fj-gamma,wj);
\end{verbatim}
\par}

\noindent The result is shown in Figure \ref{zolo2-fig1}.
This problem is one of those that can be solved analytically
\ccite{starke}:\footnote{The main known cases with analytic solutions involve
disks and intervals that can be reduced by
elementary maps to a symmetric pair of disks.}
\begin{equation}
\rns = \left({z+{3\over 4}\over z-{3\over 4}}\right)^{30}, \quad
\sigma_n = 2^{-60} \approx 8.67\cdot 10^{-19},
\label{zolo2-exactrns}
\end{equation}
\begin{equation}
\tau_n^{} = {2^{-29}\over 1 + 2^{-60}} \approx 1.86\cdot 10^{-9}.
\label{zolo2-exactrnh}
\end{equation}
In the numerical solution shown in the figure, $\log(\sigma_n)$ has
come out within $0.1\%$ of the analytically known optimal value.
The zeros and poles are distributed along circles close to the
boundaries of $E$ and $F$ rather than coalescing at $z = \pm 3/4$,
an instance of balayage \ccite{balayage} that evidently has had
little effect on the accuracy.

\begin{figure}
\vspace*{4pt}
\noindent\hbox{\kern96pt}\includegraphics[trim=20 80 190 0, clip, scale=0.9]{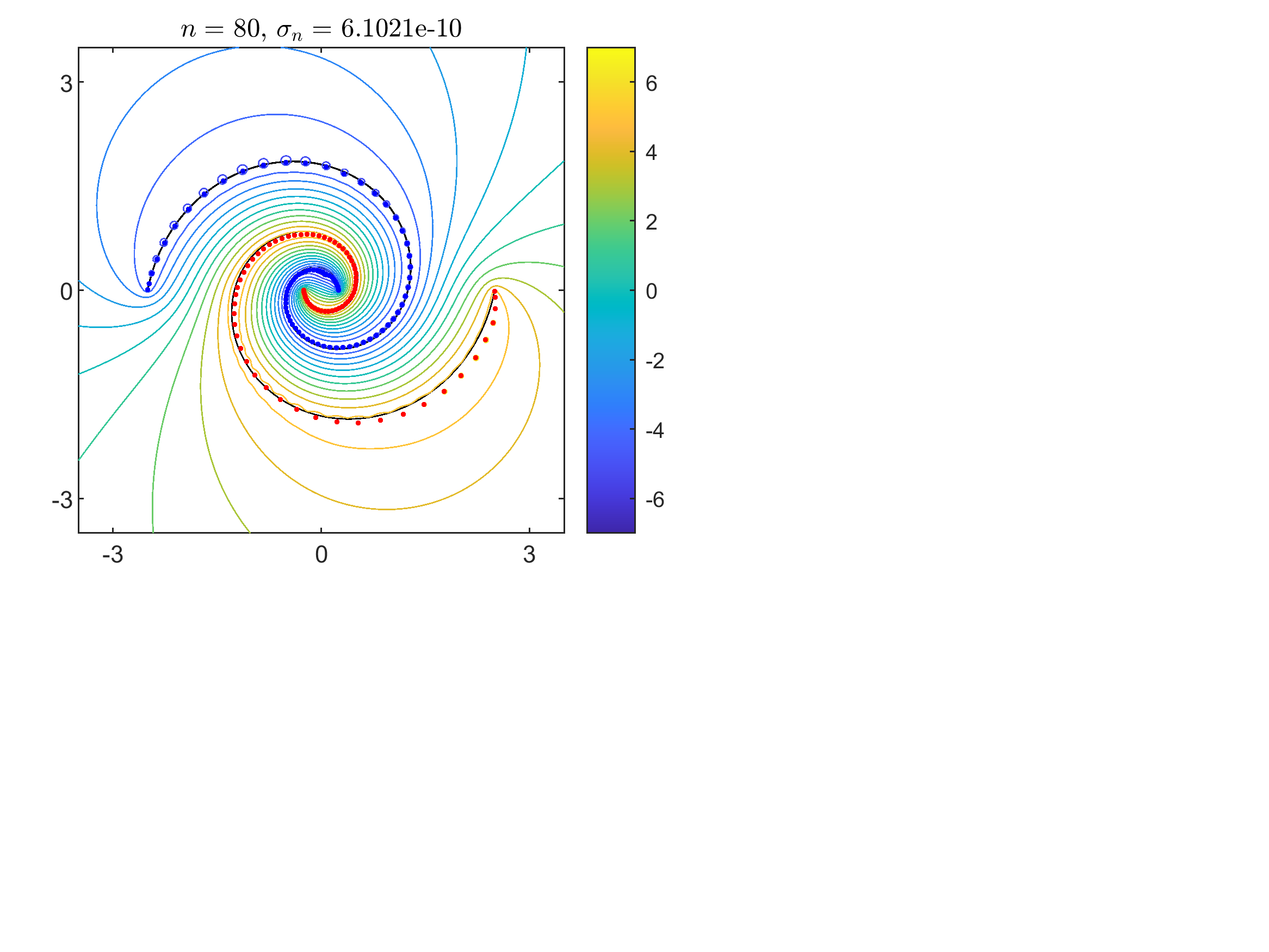}
\vspace*{-42pt}
\caption{\label{zolo2-lollifig} The Batakci lollipop.
The image shows a degree $80$ Zolotarev minimal ratio rational function for a
pair of interlocking spirals proposed by
Levent Batakci of the University of
Washington.  The sets $E$ and $F$ are marked
by black lines, and the colourbar marks levels of $\log_{10}\kern -.7pt |\rns(z)|$.
Zeros and poles are marked by blue and red dots.  The value of $\sigma_n$ is
only approximately correct, as the computation is not fully converged.}
\end{figure}

More complicated regions $E$ and $F$ can also be successfully
treated.  For example, Figure \ref{zolo2-lollifig} shows a geometry
proposed to us by Levent Batakci of the University of Washington.
Here, $E$ and $F$ are a pair of interlocking spirals, each winding $1
\half$ times around the origin.  A minimal-ratio rational function
$\rns$ of degree $n=80$ has been computed, differing by about ten
orders of magnitude on $E$ and $F$.  Our current algorithms are at
their limits for domains like this, however.  For this computation
we used \verb|'sign', 1| and \verb|'damping', 0.75| with $600$
AAA-Lawson iterations, but the convergence is only approximate.

Zolotarev, Achieser and Gonchar were motivated by approximation
theory, which in Gonchar's case was strongly allied with potential
theory, but the practical applications of the Zolotarev ratio
problem are in numerical linear algebra.  Such applications go
back to the {\em Alternating Direction Implicit (ADI)\/} method
introduced by Peaceman and Rachford in 1955 for rapid solution of
algebraic equations arising from multidimensional finite difference
discretisations of \PDEs.\ \ Later this was generalised to methods
for the solution of {\em Lyapunov equations\/} and {\em Sylvester
equations\/}, and a summary can be found in \ccite{simoncini}.
The PhD thesis of Heather Wilber was largely about numerical solution
and application of the Zolotarev problems \ccite{wilber}.

\section{\label{sec:quad}Quadrature formulas}

A quadrature formula is a discrete formula
\begin{equation}
I_n = \sum_{k=1}^n \wk f(z_k^{})
\label{quad-quadform}
\end{equation}
whose purpose is to approximate an integral
\begin{equation}
I = \int_\gamma w(z) f(z) dz,
\label{quad-int}
\end{equation}
where $\gamma$ is an integration interval or contour, $w(z)$ is an
integrable weight function, $f$ is an integrand, and $\{z_k^{}\}$
and $\{\wk\}$ are the quadrature {\em nodes\/} and {\em weights\/}.
In the case $\gamma=[-1,1]$ and $w(z)=1$ we have the context for
well-known formulas such as Gauss, Clenshaw-Curtis and Newton-Cotes
\ccite{gautschi}.  If $\gamma$ is the unit circle in $\complex$
and the nodes are roots of unity and the weights are all equal, we
have a version of the exponentially convergent periodic trapezoidal
rule \ccite{trap}.

In any of these quadrature rules, the nodes form a string of
points on a real or complex contour.   In this review, we have
seen a number of strings of poles of rational functions that delineate
approximate branch cuts: see
Figures \ref{sec:schwarz}.2--3, 
\ref{sec:potential}.4--5, 
\ref{sec:solution}.2--3 and 
\ref{sec:zolo}.1--3. 
One might speculate, is there a connection here?  The answer is
strikingly yes, as explained in \ccite{horning}.  {\em Quadrature
nodes are poles of rational approximations of a Cauchy transform,
delineating approximate branch cuts of that Cauchy transform.}
The Cauchy transform (eq.\ (\ref{cauchytrans}) below) is that
of the weight function $w(z)$ in (\ref{quad-int}): an analytic
function in $\complex\cup \{\infty\}$ except for a jump by $w(z)$
across the contour.

\begin{figure}
\vskip 5pt
\begin{center}
\includegraphics[trim=0 0 0 0, clip,scale=0.83]{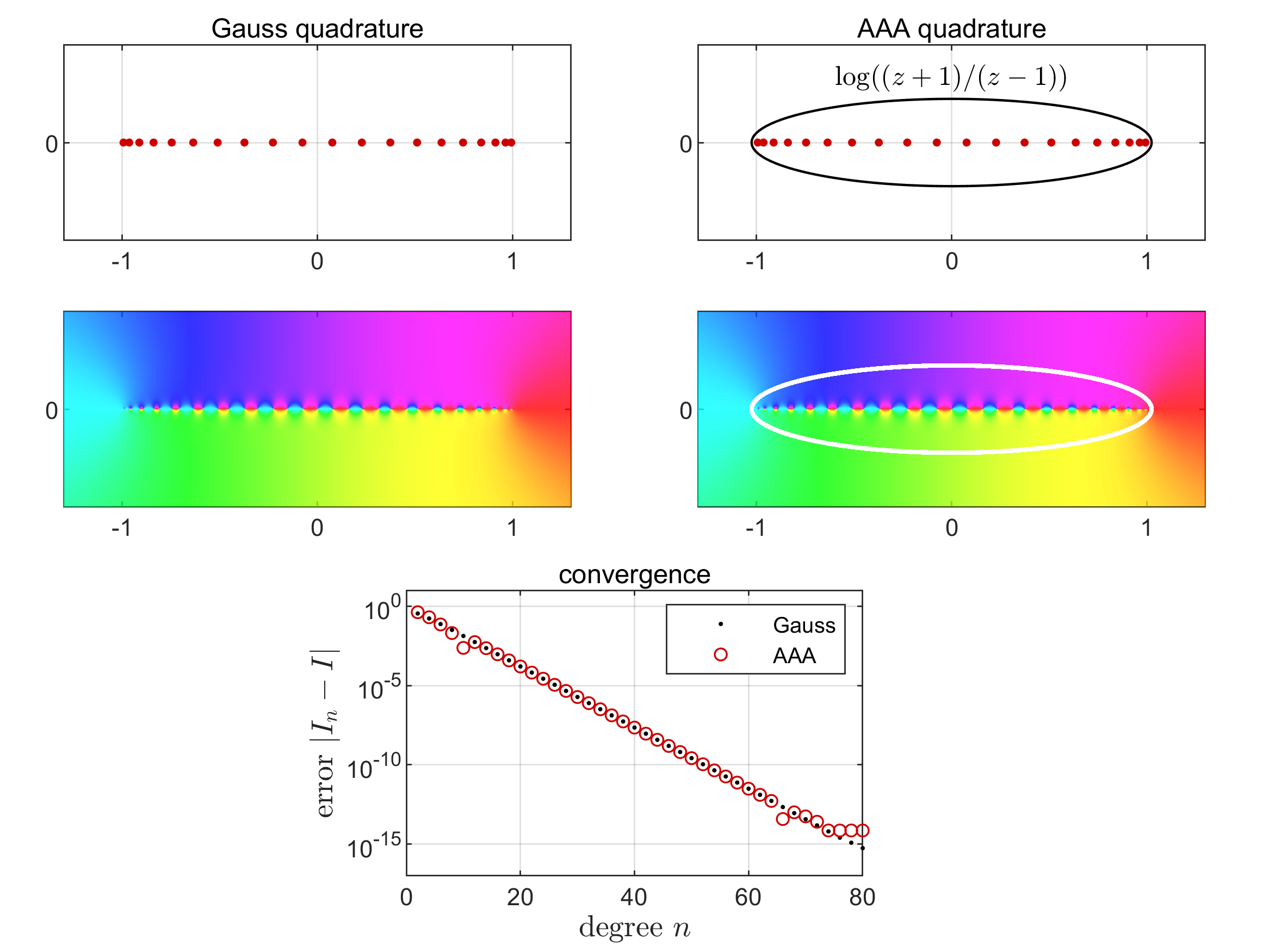}
\end{center}
\vspace*{-4pt}
\caption{\label{figgauss}On the left, nodes for the 20-point Gauss quadrature
formula and a phase portrait of the associated rational function
(\ref{quad-quadrat}).  On the right, looking virtually identical,
the same plots for a degree 20 AAA rational 
approximation of $\log((z+1)/(z-1))$ on a Bernstein ellipse.  Below,
convergence of the two quadrature formulas for the test integrand
$f(x) = 1/(1+20x^2)$.}
\end{figure}

Figure~\ref{figgauss} illustrates this effect in the familiar
case of Gauss \hbox{(-Legendre)} quadrature on $[-1,1]$.  In the
upper-left are
the nodes of the $20$-point Gauss formula in $[-1,1]$, plotted in the
complex plane.  Below this is a phase portrait of the associated rational function
\begin{equation}
r(z) = \sum_{k=1}^n {\wk\over z-z_k^{}},
\label{quad-quadrat}
\end{equation}
where $\{z_k^{}\}$ and $\{\wk\}$ are the Gauss nodes and weights.
On the right is the same pair of plots for a different rational
function $r$: the degree $20$ AAA rational approximation to
$\log((z+1)/(z-1))$ on the Bernstein ellipse plotted as a black
or white curve, with parameter chosen so that it passes through
the points $z=\pm i/\kern -1pt\sqrt{20}$.  (A Bernstein ellipse
is an ellipse in $\complex$ with foci $\pm 1$.)  Apart from the
presence of the ellipse, there is no visible difference in the plots.
This exemplifies the general principle that any quadrature formula
(\ref{quad-quadform}) can be associated with a rational function
$r$ as in (\ref{quad-quadrat}), whose approximation properties
will have much to do with the accuracy of the formula.  To see
how close Gauss and AAA are in this example, if we integrate the
test function $f(z) = 1/(1+20z^2)$, the 20-point Gauss and AAA
formulas give essentially identical errors $1.575 \cdot 10^{-4}$
and $1.560 \cdot 10^{-4}$, respectively.  The bottom plot in the
figure shows that this close agreement is maintained for all $n$
down to the usual AAA tolerance of $10^{-13}$.

The association of quadrature formulas with rational functions
goes back to Gauss in 1814, who derived his formula (as we
would now say) by Pad\'e approximation of $\log((z+1)/(z-1))$ at
$z=\infty$.\footnote{Potential theorists would explain the close
agreement in Figure \ref{figgauss} of the rational approximations
at $\infty$ and on a Bernstein ellipse
as an example of balayage \ccite{balayage}.} Many connections between
rational functions and quadrature formulas have been exploited in
the ensuing two centuries.  However, before \ccite{horning} there
was little discussion of the link with rational approximation per se.
The main example we know is \ccite{talbot}.

Following \ccite{horning}, here is a summary of
the mathematics of these approximations.
Assume the integrand $f$ of (\ref{quad-int}) is 
analytic in the closure of a Jordan region
containing $\gamma$ bounded by a Jordan curve $\Gamma$.
By the Cauchy integral formula, we have
\begin{equation}
I = \int_\gamma {1\over 2\pi i} \left[ \kern 1.5pt \int_\Gamma {f(s) \kern 1pt
ds\over s-z}\right] \kern -1pt w(z)\kern 1pt dz ,
\label{twoint1}
\end{equation}
hence by exchanging orders of integration,
\begin{equation}
I = \int_\Gamma f(s)\kern 1pt  \C(s)\kern 1pt  ds ,
\label{twoint2}
\end{equation}
where $\C(s)$ is the {\em Cauchy transform\/} of $w$,
\begin{equation}
\C(s) = {1\over 2\pi i} \int_\gamma {w(z)\kern 1pt dz\over s-z}.
\label{cauchytrans}
\end{equation}
This function is analytic throughout 
$\complex\cup \{\infty\}\backslash \gamma$ and jumps by $w(z)$
across~$\gamma$.  For the 
Gauss-Legendre case $\gamma=[-1,1]$ with $w(z) = 1$, the value is
$\C(s) = (2\pi i)^{-1} \log((s+1)/(s-1))$.
Suppose that $r$ is a rational function of degree~$n$
that approximates $\C$ on $\Gamma$.
Then (\ref{twoint2}) suggests approximating $I$ by
\begin{equation}
\In = \int_\Gamma f(s) \kern 1pt r(s) \kern 1pt ds .
\label{approx}
\end{equation}
Assume the poles of $r$ are $n$ distinct finite numbers
$z_1^{},\dots,z_n^{}$, so $r$ can be written
\begin{equation}
r(s) = C + \sum_{k=1}^n {\wk\over s-z_k^{}},
\end{equation}
where $C = r(\infty)$ and $\wk$ is the residue of $r$ at $z_k^{}$.
Then, assuming all the poles are enclosed
by $\Gamma$, residue calculus converts (\ref{approx}) to 
\begin{equation}
\In = 2\pi i \sum_{k=1}^n \wk \kern 1pt f(z_k^{}) \kern 1pt.
\label{approx2}
\end{equation}
If $\|\kern .5pt \C-r\|_\Gamma^{} \le \varepsilon$,
where $\|\cdot\|_\Gamma^{}$ is the $\infty$-norm on $\Gamma$, then
by (\ref{twoint2}) and (\ref{approx}),
\begin{equation}
| I - \In| \le \varepsilon \kern 1.5pt|\Gamma| \kern 1.5pt \|f\|_\Gamma^{},
\label{bound}
\end{equation}
where $|\Gamma|$ is the arc length of $\Gamma$.

\begin{figure}
\vskip 5pt
\begin{center}
\includegraphics[trim=0 120 0 0, clip,scale=0.83]{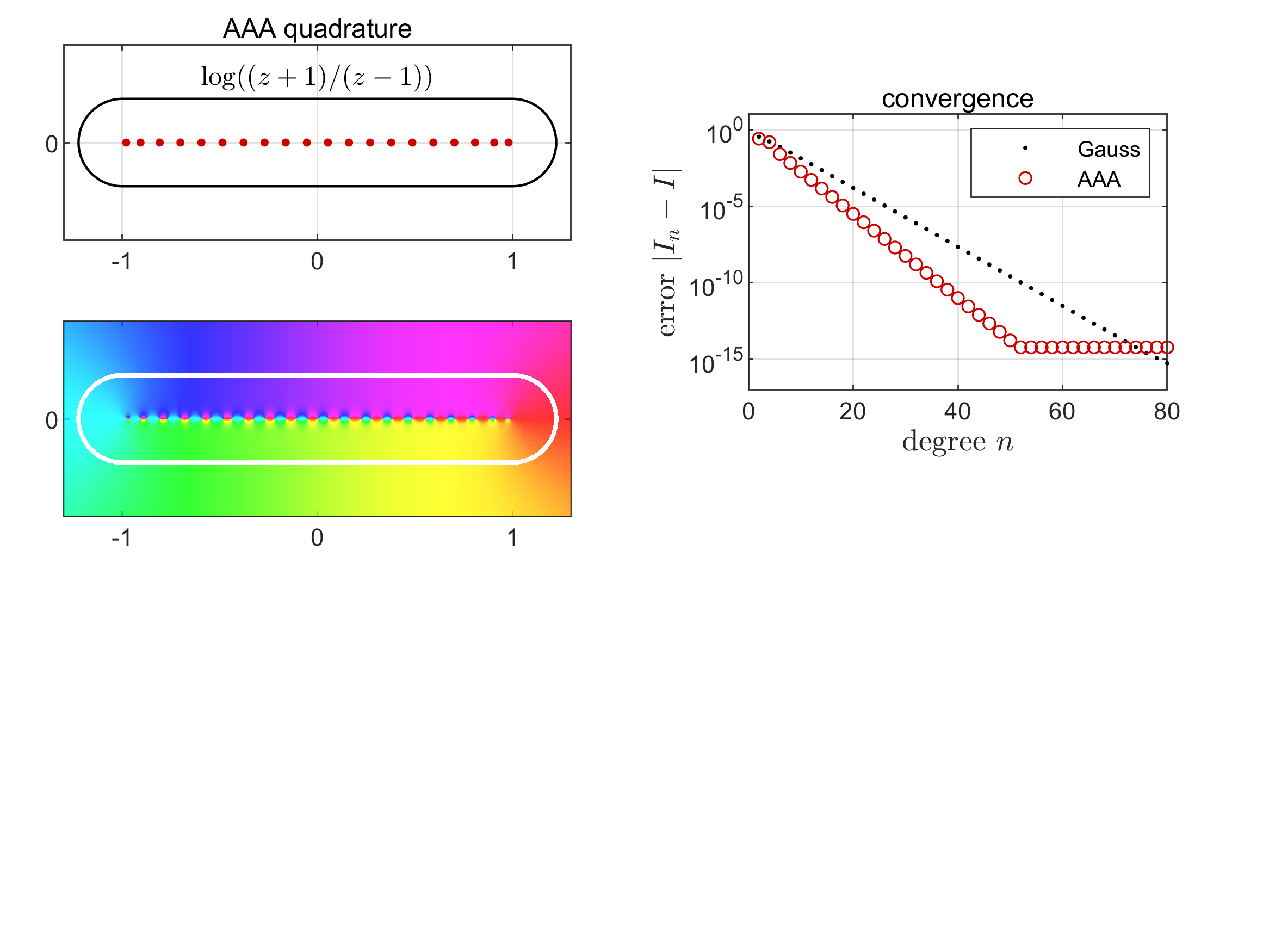}
\end{center}
\vspace*{-14pt}
\caption{\label{figgaussv}Like Figure \ref{figgauss}, but now
the quadrature formulas are based on approximation of
$\log((z+1)/(z-1))$ over the boundary of an $\varepsilon$-neighborhood
of $[-1,1]$.  For the same test integrand $f(z) = 1/(1+20z^2)$
as before, the convergence is about $\pi/2$ times faster
than with Gauss quadrature, and the quadrature nodes are much more
evenly spaced.}
\end{figure}

These relationships enable us to construct quadrature formulas
targeted to all kinds of applications by means of rational
approximation.  For example, it is known that although Gauss
quadrature is optimal from the point of view of integrating
polynomials of maximal degree, it is suboptimal for dealing,
say, with functions analytic in an $\varepsilon$-neighboorhood of
$[-1,1]$ rather than in a Bernstein ellipse.  Methods for deriving
formulas adapted to such situations have been proposed based on
conformal mapping \ccite{hale} and prolate spheroidal wave functions
\ccite{exactness}, but rational approximation gives us a simpler
new possibility, as illustrated in Figure \ref{figgaussv}.

\begin{figure}
\vskip 5pt
\begin{center}
\raisebox{17pt}{\includegraphics[trim=550 0 10 0, clip,scale=0.23]{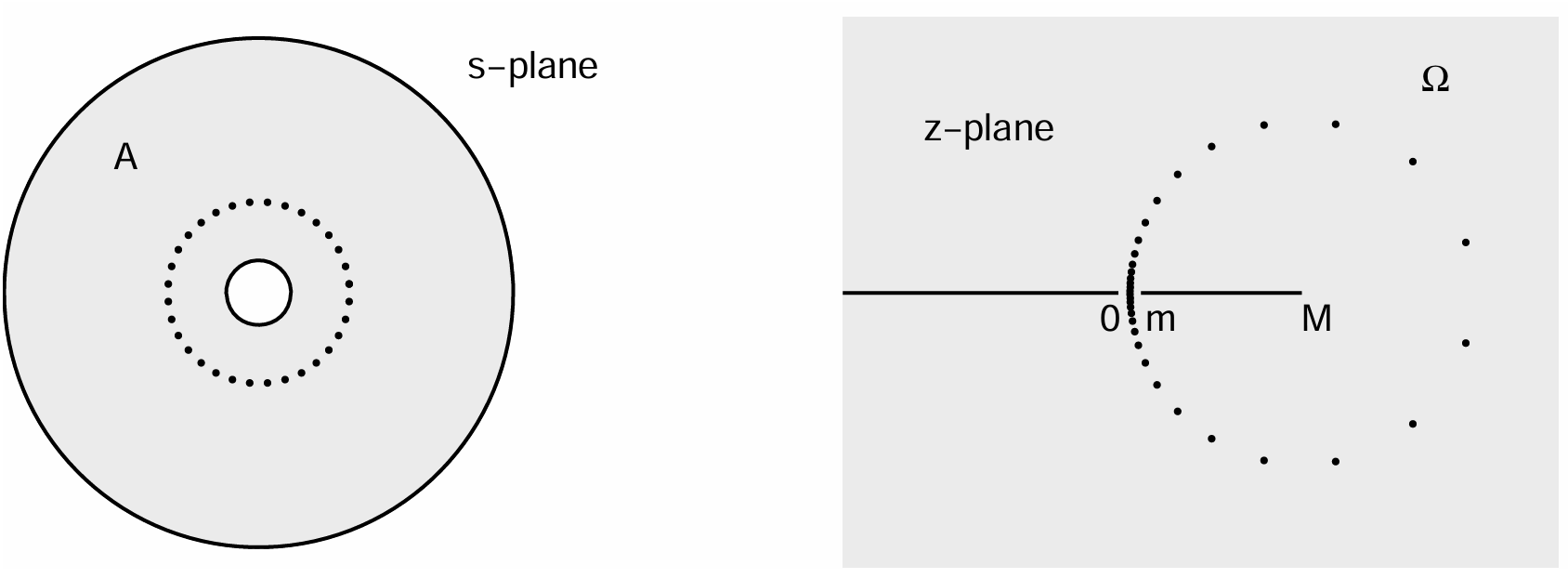}}%
\includegraphics[trim=45 190 220 0, clip,scale=1.1]{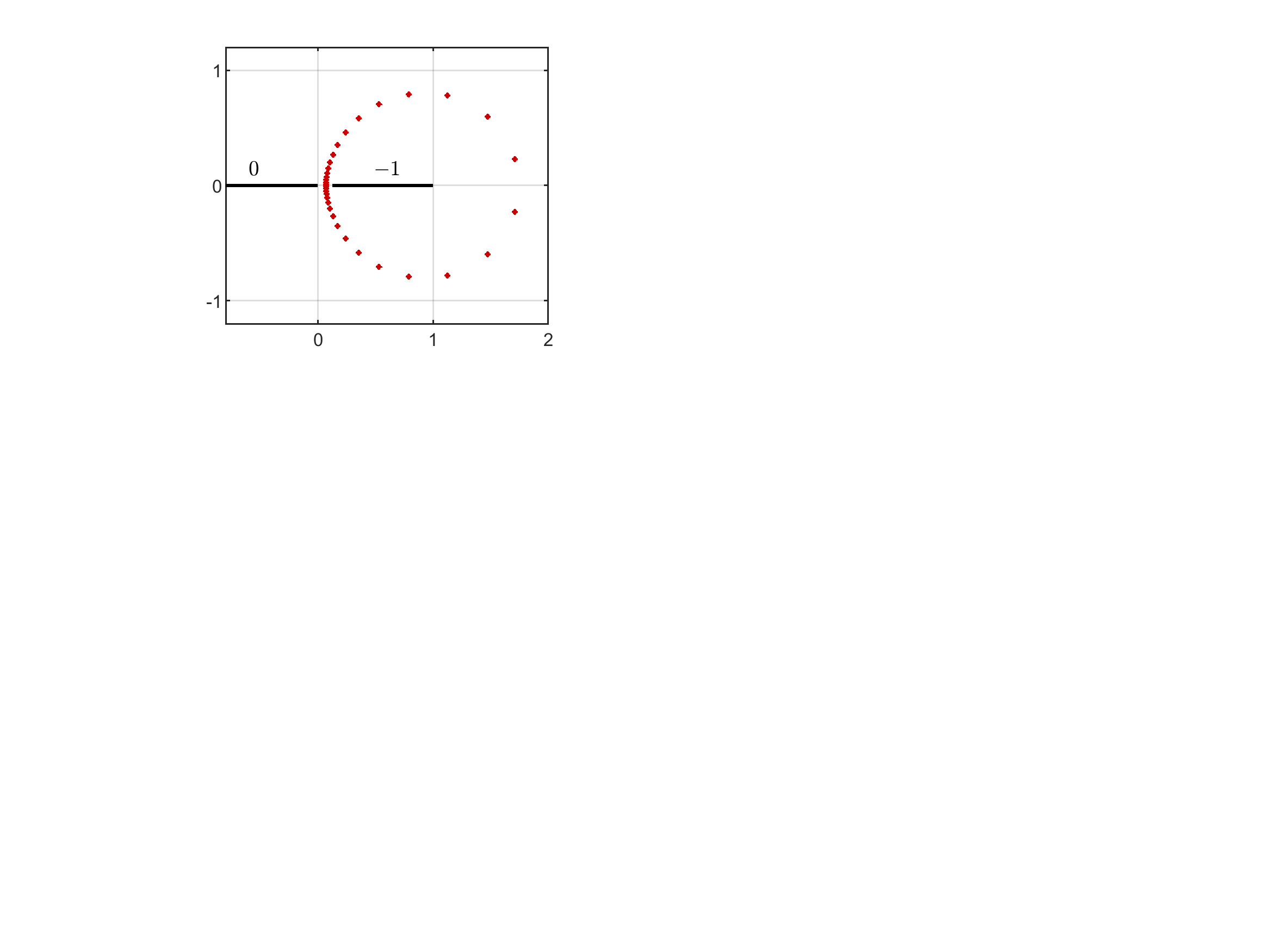}
\end{center}
\vspace*{-14pt}
\caption{\label{fignicks}On the left, a figure from \ccite{nicks} showing
quadrature nodes derived by conformal mapping for a problem of computing
functions of matrices with eigenvalues in $[m,M]$.
On the right, a figure adapted from \ccite{horning}
showing nodes derived
by a rational approximation of the Zolotarev sign type.  The accuracy of
the two formulas is comparable.}
\end{figure}

Many further examples of quadrature formulas derived from rational
approximations can be found in \ccite{horning}; we will report
just one of these, shown in Figure~\ref{fignicks}.  In the
``three Nicks'' paper published by the second author with Nick
Hale and Nick Higham \ccite{nicks}, the problem is considered of
computing matrix functions like $A^\alpha$ or $\log(A)$, where $A$
is a symmetric positive definite matrix with eigenvalues in $[m,M]$
for some $0<m<M<\infty$.  As suggested in the figure, this can be
reduced mathematically to a contour integral over a closed curve
enclosing $[m,M]$ and passing through the gap between $0$ and $M$.
From the matrix point of view, this means the matrix function is
reduced to a collection of linear systems of equations of the form
$(A-z_k^{}I\kern 1pt )x = b\kern 2pt$---the mathematics is close to that of
sections \ref{sec:eig} and \ref{sec:nonlin}.  If a simple formula
like the equispaced trapezoidal rule is used, the number of such
systems is prohibitive: $O(M/m)$ for each digit of required accuracy.
By a conformal map of $\complex \backslash ((-\infty,0\kern .3pt]
\cup [m,M])$ to a circular annulus, however, this number can be
reduced hugely to $O(\log(M/m))$.  For a problem with $M/m=8$, for
example, one gets 6-digit accuracy from as few as 10 matrix solves
(exploiting real symmetry).

As suggested in the right image of Figure \ref{fignicks}, the
same specialised quadrature form can alternatively be derived
by rational approximation without the need for conformal mapping.
The Cauchy transform for a closed curve surrounding $[m,M]$ is simply
the jump function taking the values $-1$ inside the curve and $0$
outside, so we have an approximation of exactly the Zolotarev sign
form discussed in section~\ref{sec:zolo}.  This is readily solved
numerically, and the figure shows the agreement with the result
of conformal mapping.  Convergence is exponential, as predicted
(not shown).  In this example the convergence just goes down to
about the level of $10^{-8}$ for reasons related to rounding errors
that have not yet been investigated.

We should emphasise an aspect of the mathematics of this section
that is particularly distinctive.  Normally we think of an integral
(\ref{quad-int}) as associated with a fixed contour $\gamma$,
and the same is true of a Cauchy transform (\ref{cauchytrans}).
In the analytic functions context of this section, however,
any choice of $\gamma$ that connects the same endpoints in
(\ref{quad-int}) and lies in the region of analyticity of $f$
and inside the outer contour $\Gamma$ in (\ref{twoint2}) will
give the same integral $I$ in (\ref{quad-int}) and the same Cauchy
transform $\C(s)$ along $\Gamma$ in (\ref{twoint2}).  When we find a
rational approximation $r(s)\approx \C(s)$ along $\Gamma$, its poles implicitly make a
choice of one of these equivalent contours $\gamma$.  For problems
associated with real functions on real intervals, one will hardly
notice this effect since symmetry will make $\gamma$ real too, as
in Figures \ref{figgauss} and \ref{figgaussv}.  In cases like
Figure \ref{fignicks}, however, the rational approximation is
genuinely choosing the contour as well as the nodes and weights.
More examples can be found in \ccite{horning}.

\section{\label{sec:cauchy}Cauchy, Wiener-Hopf, Riemann-Hilbert}
Many problems in complex analysis involve jumps of analytic functions
across contours.  AAA approximation usually has something to contribute,
and in this section we will illustrate the
possibilities by a Cauchy transform, a Wiener-Hopf
problem and a Riemann-Hilbert problem.

The starting point for many problems of this kind is the Cauchy integral
\begin{equation}
f(z) = {1\over 2\pi i} \int_\Gamma {h(t)\over t-z}\kern 2pt dt,
\label{ctrans}
\end{equation}
which appeared as (\ref{cauchytrans}) in the last section.
Here $\Gamma$ is a smooth bounded contour that may be open or closed,
and we assume $h$ is H\"older continuous.
The formula implies that $f$ is
analytic for $z\in\complex\cup \{\infty\} \backslash \Gamma$,
with $f(\infty)=0$.  If $\kern 1pt\Gamma$ is a closed contour, then it divides
the plane into two components, on which $f$ will normally take different analytic
branches with $\Gamma$ as a branch cut between them.  If $\kern 1pt\Gamma$ is open, then 
$\complex\cup \{\infty\} \backslash \Gamma$ is connected, and $f$ is 
a single-valued analytic function with $\Gamma$ as an internal branch cut, typically
with logarithmic singularities at the endpoints.

A key mathematical result associated with (\ref{ctrans})
goes by the name of the {\em Sokhotski-Plemelj formulas\/}.
The first of these formulas asserts that $f(z)$ jumps by
$h(t)$ as $z$ crosses $\Gamma$ at a point $t\in\Gamma$ (not
an endpoint).  This is written as
\begin{equation}
h(t) = f_+(t) - f_-(t),
\label{jump}
\end{equation}
where $f_+$ denotes $f$ as defined by (\ref{ctrans})
on the left of $\Gamma$
and $f_-$ is the same on the right.  (``Left'' and ``right'' are
defined with respect to the orientation of $\Gamma$.)
The other Sokhotski-Plemelj
formula states that for $t\in \Gamma$, we have
\begin{equation}
f(t) =  {\textstyle{1\over 2}} (\kern .7pt f_+(t)+f_-(t)),
\label{pval}
\end{equation}
where $f(t)$ is defined for $t\in\Gamma$ by the {\em principal value\/}
of (\ref{ctrans}).
A careful and detailed reference on these topics is \ccite{henrici3}.

\begin{figure}
\vskip 5pt
\begin{center}
\includegraphics[trim=0 170 0 0, clip,scale=0.84]{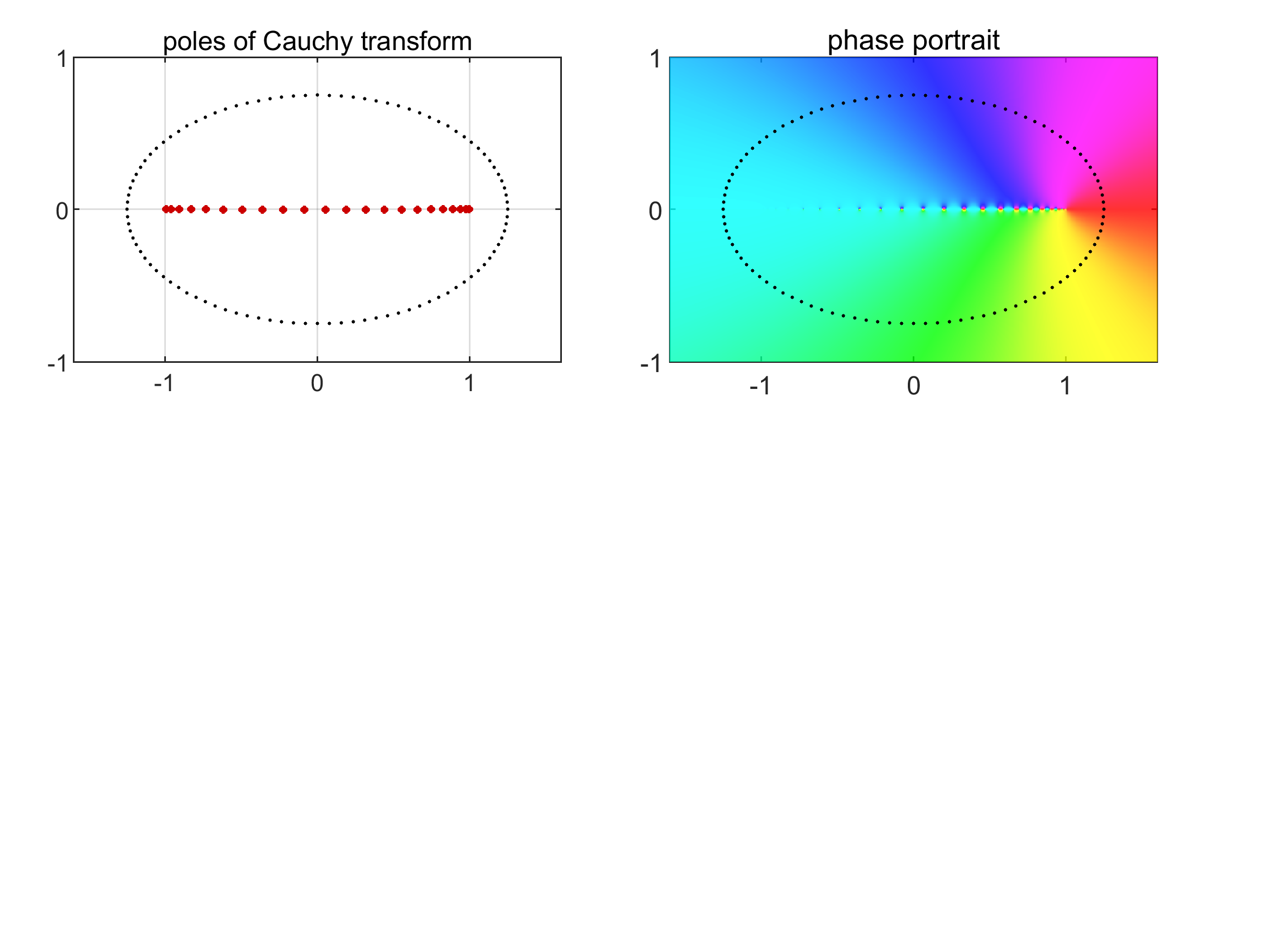}
\end{center}
\vspace*{-14pt}
\caption{\label{cauchyfig}The Cauchy transform of $h(t) = e^{3\kern .5pt t}$ on $\Gamma = [-1,1]$
is computed by evaluating the integral (\ref{ctrans}) at 100 points on an
ellipse by adaptive quadrature
and then approximating the result by a AAA rational function.
The approximation has degree 22, with little asymmetry in the pole positions
(10 poles in $[-1,0\kern .3pt]$, 12 poles in $[\kern .3pt 0,1]$) but a strong asymmetry in their
strengths as revealed by the phase portrait.}
\end{figure}

As discussed in the last section, the map from $h$ to $f$ defined by
(\ref{ctrans}) is called the {\em Cauchy transform\/}.  The application of
AAA approximation to the numerical computation of Cauchy transforms was
proposed in \ccite{horning}.  Figures \ref{figgauss}--\ref{fignicks} presented
examples from that paper involving AAA approximation of
Cauchy transforms that are known
analytically: $\log((z+1)/(z-1))$ outside an interval or $0$ and $1$ on two
sides of a closed curve (both divided by $2\pi i$).
But AAA can also treat problems that are fully numerical,
and section 5 of \ccite{horning} illustrates these too.
Suppose, for example, we want a rational approximation to the
Cauchy transform of $h(t) = e^{3\kern .5pt t}$ on $\Gamma = [-1,1]$.
The following lines evaluate (\ref{ctrans}) at 100 points on a
Bernstein ellipse surrounding $[-1,1]$ using the Matlab adaptive quadrature
code {\tt quadgk} and then compute a AAA approximation to the usual
13-digit relative accuracy.  The result is shown in Figure~\ref{cauchyfig}.

{\small
\begin{verbatim}
      c = 2*exp(2i*pi*(1:100)'/100); Z = (c+1./c)/2;
      f = @(z) arrayfun(@(zz) quadgk(@(t) exp(3*t)./(zz-t),...
                                     -1,1,'RelTol',1e-13), z);
      [r,pol] = aaa(f(Z),Z);
\end{verbatim}
\par}

We noted above that the Cauchy transform $f$ takes the value
$f(\infty) = 0$.  When $f$ is approximated by a rational function $r$
by standard AAA, the algorithm does not enforce this constraint,
though the analytic structure of the problem is likely to make it
hold approximately: in the example of Figure \ref{cauchyfig}, we have
$|r(\infty)| < 10^{-15}$.  The possibility of computing approximants
that are constrained to be zero at $\infty$ is discussed at the
end of section \ref{sec:mn}.  As a practical matter, so long as
the aim is near-best as opposed to truly best approximation, one
can often approximate a function without a condition at $\infty$
and then subtract the constant $r(\infty)$ from the result.

Many applications of (\ref{ctrans}) involve situations where
$\Gamma$ is a real interval, in which case $h$ is often assumed
to be real too.  Often $h(t)\kern .5pt dt$ is replaced by a nonnegative
measure $d\mu(t)$ that can include point masses.  
(At such points (\ref{jump}) and (\ref{pval}) will not make sense.)
With $\Gamma = (-\infty,0\kern .3pt]$ this leads to {\em Stieltjes functions\/},
and with $\Gamma=[a,b\kern .5pt]$ one gets {\em Markov functions\/}.
Finite intervals arise in quadrature problems, as
discussed in the last section.
In these cases the Cauchy transform is normalised differently and renamed
the {\em Stieltjes\/} or {\em Markov-Stieltjes transform\/}.

A context in which infinite intervals in (\ref{ctrans}) arise is
in certain problems of {\em Wiener-Hopf factorisation\/}.  In the
simplest scalar additive Wiener-Hopf problem, one has a function
$h(t)$ defined for $t\in \real$ and the aim is to write it as a
sum $h(t) = f_+(t) + f_-(t)$ of functions analytic in the upper
and lower half-planes $\complex_+$ and $\complex_-$, respectively.
We assume $h(t)\to 0$ as $|\kern .6pt t \kern .6pt|\to\infty$.
According to (\ref{jump}), this problem is solved by (\ref{ctrans})
with the adjustment that we need to flip the sign of $f_-$.  In a
multiplicative Wiener-Hopf problem, it is the same except that
one wants a representation as a product $h(t) = g_+(t) g_-(t)$.
If $h(t)>0$ for all $t\in\real$, then the multiplicative Wiener-Hopf
problem for $h$ reduces to the additive Wiener-Hopf problem for
$\log h$.

\begin{figure}
\vskip 5pt
\begin{center}
\includegraphics[trim=0 0 0 0, clip, scale=0.63]{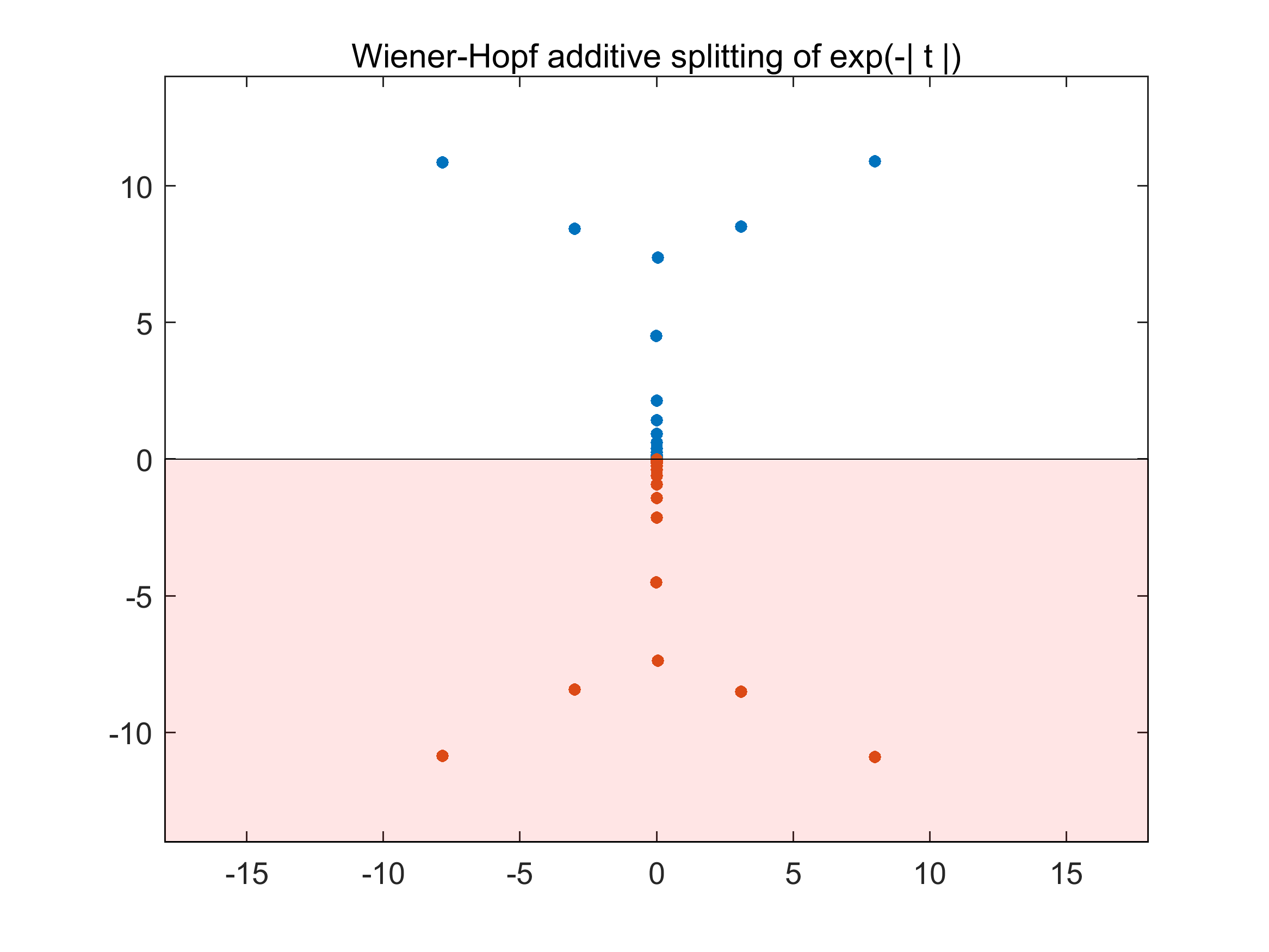}
\end{center}
\vspace*{-14pt}
\caption{\label{wienerhopffig} Poles of a AAA rational approximation to $h(t) = e^{-|\kern .6pt t \kern .6pt|}$ for
$t\in \real$.  To decompose $h$ approximately into $f_+ + f_-$,
the contributions from the poles above and below the real axis
are separated.  Here there are 41 poles in each half-plane
and the combination gives ten-digit accuracy.}
\end{figure}

Starting in 2021, we have explored AAA approximation for a
few Wiener-Hopf problems, but there is no publication as yet.
We illustrate the idea here by computing the additive Wiener-Hopf
decomposition of $h(t) = e^{-|\kern .6pt t \kern .6pt|}$ on $\real$.
First we set up a grid, exponentially clustered at $t=0$ in
view of the singularity there, and apply AAA to approximate $h$
on the grid by a rational function~$r$ to accuracy $10^{-10}$.
What we need at this stage are the poles of~$r$, which are plotted
in Figure \ref{wienerhopffig}, and the corresponding residues.
(We discard a couple of poles that AAA returns on the real line, which
are Froissart doublets with negligible residues.)

{\small
\begin{verbatim}
      h = @(t) exp(-abs(t));  
      T = logspace(-10,10,200)'; T = [T; -T];
      [~,pol,res] = aaa(h(T),T,'tol',1e-10);
      jj = find(imag(pol)==0); pol(jj) = []; res(jj) = [];
\end{verbatim}
\par}

\noindent We then need to split $r$ into its components analytic
in $\complex_+$ and $\complex_-$.  We do this by reconstructing it
from its poles and residues in each half-plane:

{\small
\begin{verbatim}
      p = pol.';
      ii = (imag(p)>0); p1 = p(ii); p2 = p(~ii);
      res1 = res(ii).'; res2 = res(~ii).';
      f1 = @(z) reshape(sum(res1./(z(:)-p1),2),size(z));
      f2 = @(z) reshape(sum(res2./(z(:)-p2),2),size(z));
\end{verbatim}
\par}

\label{whpage}
\noindent A check of accuracy shows a maximal error $|h(t_j) -
(f_+(t_j)+f_-(t_j))|$ of $2.1\cdot 10^{-10}$ over the sample grid.
The computation takes about $0.1$ s on our laptop.

It is interesting to try to explain the forked structure of the
distribution of poles in Figure \ref{wienerhopffig}.  In the terminology
of the final paragraph of section \ref{sec:potential}, we believe
this figure shows a mix of branch cut poles along
the imaginary axis near $t=0$ and approximation poles further out.

We are Wiener-Hopf amateurs, but among the experts there is
awareness of rational approximation as a useful tool \ccite{kisil13}.
A survey of constructive aspects of Wiener-Hopf problems can
be found in \ccite{kisil}, and AAA approximation is used to
treat a Wiener-Hopf problem on the unit circle in \ccite{kisil2}.
Whether rational approximations can be used to treat matrix-valued
Wiener-Hopf problems is unknown.

\begin{figure}
\vskip 5pt
\begin{center}
\includegraphics[trim=0 0 0 0, clip, scale=0.65]{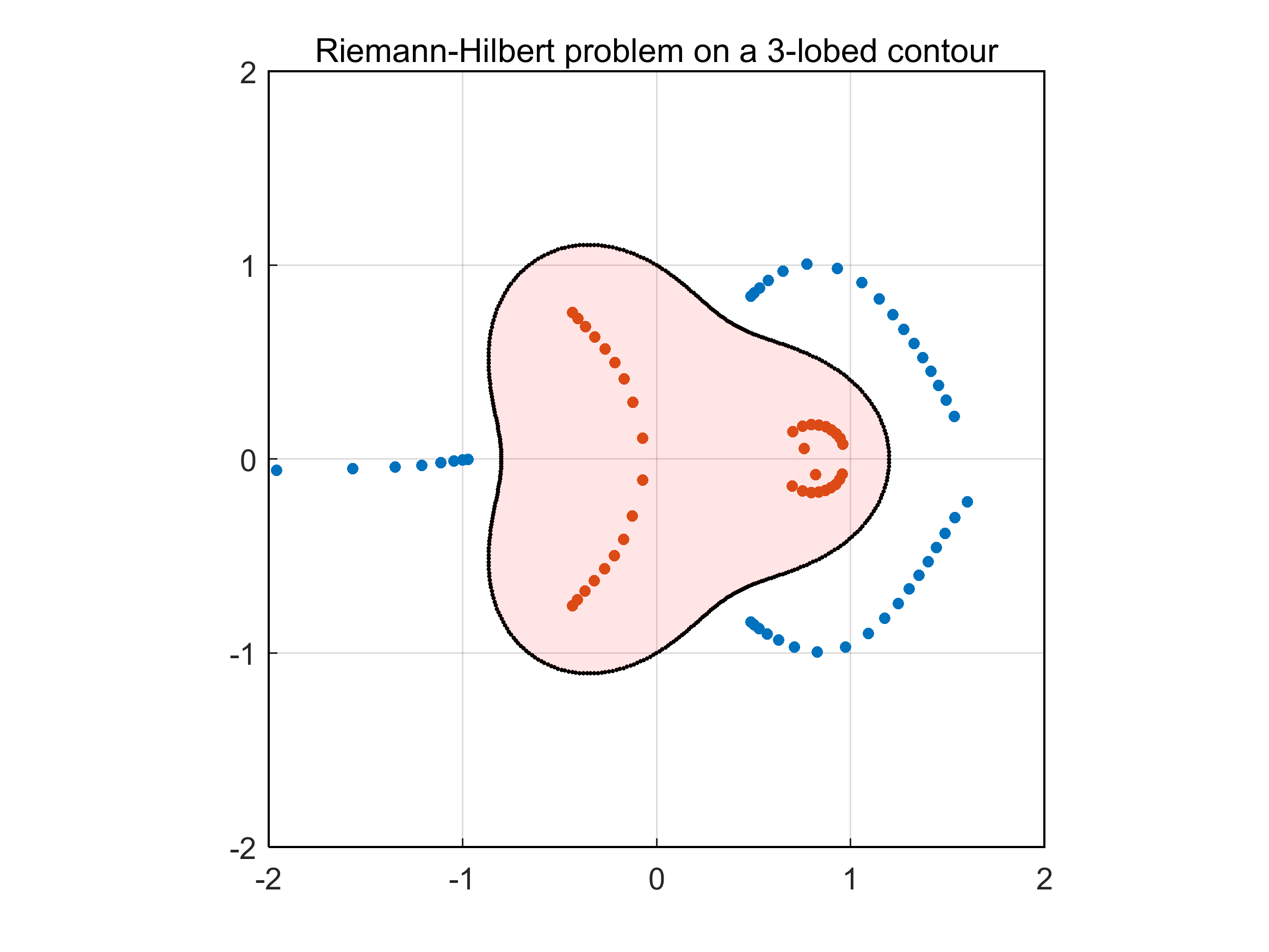}
\end{center}
\vspace*{-14pt}
\caption{\label{rhfig} Poles for a rational approximation to solve
a Riemann-Hilbert problem involving splitting of
the function (\ref{rhfun}) on a 3-lobed contour.  Piecing together
rational approximations to $f_+$ and $f_-$ from the poles outside and
inside, respectively, solves the problem to 6 digits of accuracy.}
\end{figure}

Wiener-Hopf problems are a special case of
{\em Riemann-Hilbert problems\/}, in which arbitrary
linear jump conditions are applied across a contour $\Gamma\subset\complex$
that is generally not just the real line or unit circle.
As early as 2020, Costa proposed in the final sentence
of \cite{aaalscosta} that AAA should be effective for solving
certain Riemann-Hilbert problems, but so far as we know, not much has
yet been done in this direction.  
As encouragement to further work we will show just one
example.
Suppose we want to find functions
$f_+$ and $f_-$ analytic inside and outside the 3-lobed curve $\Gamma$
of Figure~\ref{rhfig},
respectively, that satisfy the jump condition
\begin{equation}
f_+(t) - f_-(t) = g(t) = \cos(\exp(3\kern .7pt  \Re (t)))
\label{rhfun}
\end{equation}
on $\Gamma$ with the normalisation $f_-(\infty) = 0$.
We carry out AAA approximation as usual, but instead of outputting
poles and residues as in the last example we just take the poles, because
the residues will have to be computed separately as we shall explain.
The poles are plotted in Figure~\ref{rhfig}.

{\small
\begin{verbatim}
      Z = exp(2i*pi*(1:400)'/400);
      Z = Z.*(1+.2*cos(3*angle(Z)));
      g = @(z) cos(exp(3*real(z)));
      [~,pol] = aaa(g(Z),Z); pol = pol.';
\end{verbatim}
\par}

\noindent To compute residues, we solve a least-squares problem.
This is what the code {\tt aaa.m} does too, internally (see footnote 4),
but there, a constant
term is included in the least-squares fit.  Here, we must
omit the constant term to end up with a fit with $r(\infty)=0$:

{\small
\begin{verbatim}
      A = 1./(Z-pol); c = A\g(Z); res = c.';
\end{verbatim}
\par}

\noindent  Now we split $r$ into two pieces $r_+$ and $r_-$ as before. 
The accuracy is $6.8 \cdot 10^{-7}$ in this example.

{\small
\begin{verbatim}
      inpoly = @(z,w) inpolygon(real(z),imag(z),real(w),imag(w));
      ii = inpoly(pol,Z);
      fp = @(z) reshape(sum(res(~ii)./(z(:)-pol(~ii)),2),size(z));
      fm = @(z) reshape(-sum(res(ii)./(z(:)-pol(ii)),2),size(z));
      err = norm(fp(Z)-fm(Z)-g(Z),inf)
\end{verbatim}
\par}

\section{\label{sec:aaals}Laplace problems and AAA least-squares}
We are now at a transition point of this paper.  The last 21 sections
have explored problems in which data are
given that we wish to approximate by a rational function.  For the
next four sections, data are given that we wish to approximate by
the {\em real part\/} of a rational function.  This is an entirely different
problem, and the breakthrough
came with a paper published by Stefano Costa in 2020 \ccite{aaalscosta}.

It will take a couple of pages to explain.  We start with the observation
that the real part of an analytic function is harmonic, i.e., a solution
to the Laplace equation $\Delta u = 0$.  In the other direction,
given a harmonic function $u(z)$ in a planar domain $\Omega$, is
it the real part of an analytic function $f(z) = u(z) + i\kern .4pt
v(z)\kern 1pt$?  The answer is yes if $\Omega$ is simply-connected,
and the harmonic conjugate $v$ is uniquely determined up to an
arbitrary additive constant.  If $\Omega$ is multiply connected, then
one only needs to generalise the result slightly by allowing $f$ to
have a logarithmic term in each hole.  This is the {\em logarithmic
conjugation theorem,} originating perhaps in \ccite{walsh29},
which we reproduce in the words of \ccite{axler} since this basic
result is not as well known as it should be.  Its consequence is
that the methods to be discussed here and in the next three sections
generalise readily to multiply-connected geometries.

\begin{theorem}
{\bf Logarithmic Conjugation Theorem.}
Suppose $\Omega$ is a finitely connected region, with $K_1, \dots,
K_N$ denoting the bounded components of the complement of $\Omega$.
For each $j$, let $a_j$ be a point in $K_j$.  If $u$ is a real
valued harmonic function on $\Omega$, then there exist an analytic
function $f$ on $\Omega$ and real numbers $c_1^{},\dots,c_N^{}$
such that
\begin{equation} u(z) = \Re f(z) + c_1^{} \log |z-a_1^{}| +
\cdots + c_N^{} \log |z-a_N^{}|
\end{equation}
for every $z\in\Omega$.
\end{theorem}

It is also known that analytic functions can be approximated by
polynomials.  Specifically, if $f$ is analytic in the closure of
a simply-connected bounded domain, then it can be approximated by
polynomials with exponential convergence as a function of degree
$n$, and if the domain has holes, then the same is true so long
as additional polynomials are included involving reciprocal powers at
an arbitrary point in each hole.  This is {\em Runge's theorem}
\ccite{runge85}.

For simplicity we will confine the remaining
discussion to the simply-connected case of a Dirichlet problem on a
Jordan domain $\Omega$ bounded by a contour $\Gamma$.
The aim is to find a function $u$ satisfying
\begin{equation}
\Delta u(z) = 0, ~ z\in \Omega, \quad u(z) = h(z), ~ z\in \Gamma,
\label{laplaceprob}
\end{equation}
where $h$ is a given real continuous function.  A unique solution
exists, as can be proved by transplanting the problem to the unit disk by the
Riemann mapping theorem and applying the Poisson integral formula.
Variations such as Neumann or mixed boundary conditions are easily treated
by the same methods we shall describe.

In view of the results summarised above,
it is not surprising that an old idea for solving (\ref{laplaceprob})
is to represent $u$, to sufficient accuracy, as the real part of a
polynomial or a rational function that is pole-free in $\Omega$.
Thanks to the maximum principle for harmonic functions, one only
needs to approximate on $\Gamma$, and the result will be accurate
throughout $\Omega$.  This vision was pursued by Joseph Walsh and
his student John Curtiss (\cite{walsh29,curtiss62}).  Numerical
analysis and computers were not advanced enough for their attempts
to be relevant computationally, but in the polynomial case, such
methods are now eminently practical if coefficients are determined
by least-squares fitting \ccite{series}.  For complicated domains,
it is important to use a well-conditioned basis, as can
be computed by Vandermonde with Arnoldi orthogonalisation \ccite{VA}.

The approximation to $u$ is not a polynomial: it is the real part
of a polynomial.  To fit the boundary function $h$ by means of a
polynomial of degree $n$, therefore, one does not have $n+1$ real
parameters to determine, but $2\kern .3pt n+1$ real parameters corresponding to
real and imaginary parts.  If $m$ sample points $\{z_j^{}\}$ are
distributed along $\Gamma$, this will lead to an $m\times (2\kern .3pt n+1)$
matrix whose columns correspond to $1$ and to the real and imaginary
parts of the powers $z_j^{}, z_j^2, \dots,z_j^n$.  This is in the
monomial basis; on an interval one would use Chebyshev polynomials,
and on a general domain the basis would be constructed by Vandermonde
with Arnoldi orthogonalisation.

All of this is excellent except that polynomials fail in two
cases, as mentioned in the final paragraph of 
section \ref{sec:approx} and discussed at length in \ccite{inlets}.  First, if $u$
has a singularity on $\Gamma$, which will usually be the case if
$\Gamma$ has a corner or $h$ itself is singular, then the polynomial
convergence rate falls from exponential to algebraic, and obtaining
more than 2 or 3 digits of accuracy typically becomes difficult.
Of course, polynomials also have trouble
when there are singularities near $\Gamma$ but not exactly on it.
Rational approximations, by contrast, can achieve root-exponential
convergence despite the presence of branch point singularities
\ccite{newman,gt}.  
Second, if $\Omega$ is nonconvex, as illustrated in the example
of Figure \ref{potential-pot4}, then polynomials may converge
exponentially in asymptotic theory but at a rate so low that in
practice their convergence is again just algebraic.  The asymptotic
convergence rate scales as $\exp(-\pi L)$, where $L$ is the aspect
ratio of a nonconvex indentation of $\Omega$ \cite[sec.~7]{inlets}.
For example, in a Laplace problem on an inverted ellipse of parameter
$\rho=1.2$ (reciprocal of a Bernstein ellipse), the polynomial degree
will have to be increased asymptotically by about $430{,}000$ for
each additional digit of accuracy \cite[Table~1]{inlets}.  Rational
approximations, by contrast, can converge very fast because they
can place poles in the inlets, as in Figure \ref{potential-pot4} and
Table~\ref{speedtable}.

Thus there are potentially huge advantages if rational
functions rather than polynomials are used to solve Laplace
problems.\footnote{This does not seem to have been Walsh's
perspective.  He considered rational functions more because they
are a natural extension of polynomials, and because they are
needed for regions with holes, than because of their superior
approximation properties, which he may not have been very aware
of---perhaps unsurprisingly in that pre-numerical era.} Before AAA,
the main paper we are aware of that explores this possibility is
\ccite{hoch}, which presents very impressive results confirming the
power of rational functions.  The challenge is, how to compute them?
This is not rational approximation of the kind accomplished by AAA,
because only the real part of $f$ is known.  Or to put it another
way, AAA can approximate the real function $h(z)$ on $\Gamma$,
but the resulting approximation will be real on $\Gamma$, to the
approximation tolerance, and it will achieve this by having poles
inside $\Gamma$ as well as outside, making it useless for solving
the Laplace problem.

So one is left with the problem of where to put poles outside
$\Gamma$ for a good approximation, which Hochman, et al.\
interpret as an example of the Method of Fundamental Solutions
(MFS) \ccite{doicu,fk}.
(Their extensive introduction, with many references, is very much
worth reading.)  They propose an iterative process of pole-adjustment
which is successful in a number of examples, but lacks the easy
speed and reliability of AAA.

During 2016--2020, we had the AAA algorithm but could not see a
way to use it to approximate real harmonic functions as opposed
to complex analytic ones.  Instead, in this period the second
author and Gopal devised the {\em lightning Laplace method\/} for
solution of Laplace problems dominated by corner singularities,
in which poles are clustered exponentially at pre-assigned positions near each corner
in the manner established as effective by theory and experiment
\ccite{gt,clustering}.  The poles are placed in an a priori fashion,
independently of the right-hand-side $h$, making this essentially a linear method
that reduces to matrix least-squares problems.  (See Figure~\ref{wedges} for
an illustration of lightning pole clustering.)  It easily solves
Laplace problems of this kind to high accuracy, and software 
is available at \ccite{lightning}.  In \ccite{gt}, theorems
are developed establishing that exponentially clustered poles, together
with a low-degree polynomial, can
resolve general Laplace problems with any finite number of branch
point singularities with root-exponential convergence.

The challenge, however, was to find a AAA-style method to
work automatically on arbitrary regions, and the new idea
came in a paper posted on arXiv in 2020 by Stefano Costa 
\ccite{aaalscosta}.  Costa's original {\em AAA-least squares (AAALS) method\/},
described also in \ccite{aaals}, is as follows.  (In a moment we 
will modify it by changing poles to double poles and then optionally compressing the final
result by a further AAA approximation.)

\begin{enumerate}
\item Use AAA to approximate $h$ by a complex rational function $r$ on $\Gamma$\kern .7pt;
\item Compute the poles of $r$, and discard those inside $\Gamma$\kern .7pt;
\item Use linear least-squares to fit $h$ by the remaining poles outside $\Gamma$\kern .3pt.
\end{enumerate}

\noindent
If $r$ has $n$ poles $\{\pk\}$ outside $\Gamma$ and there are $m$
sample points $\{z_j^{}\}$ on $\Gamma$, the least-squares fit will
involve a matrix of dimension $m\times (2\kern .3pt n+1)$, with one column
corresponding to the constant $1$ and the others corresponding
to $\Re (1/(z_j^{}-\pk))$ and $\Im (1/(z_j^{}-\pk))$.  For better
numerical stability when poles cluster near the boundary $\Gamma$,
it is better to rescale the columns via $\Re (d_k^{}/(z_j^{}-\pk))$
and $\Im (d_k^{}/(z_j^{}-\pk))$, where $d_k{} = \min_j |z_j-\pk|$.
As mentioned above, sometimes a few columns are added to the matrix corresponding to a low-degree
polynomial, whether expressed in a monomial, Chebyshev, or
Vandermonde with Arnoldi basis.  For an analysis of the sometimes surprisingly great power
of low-degree polynomial additions to rational function approximations, see \ccite{herremans}.

Like the lightning Laplace method, the AAALS method produces a
solution to a Laplace problem in the form of the real part of a
single rational function, typically evaluable on a laptop in a matter
of microseconds per point.  The function is globally analytic
throughout $\Omega$, so there are no derivative discontinuities
of the sort one would find with finite elements or with integral
representations evaluated adaptively.  Another appealing feature
of these methods is that the harmonic conjugate, being just the
imaginary part of the same rational function, is available for free.
In applications, the harmonic conjugate is usually of physical
interest, representing flow lines of heat, charge, ideal fluid, or probability orthogonal to
the equipotential lines of the original Laplace problem.

\begin{figure}[t]
\vskip 5pt
\begin{center}
\includegraphics[trim=43 0 43 0, clip,scale=0.44]{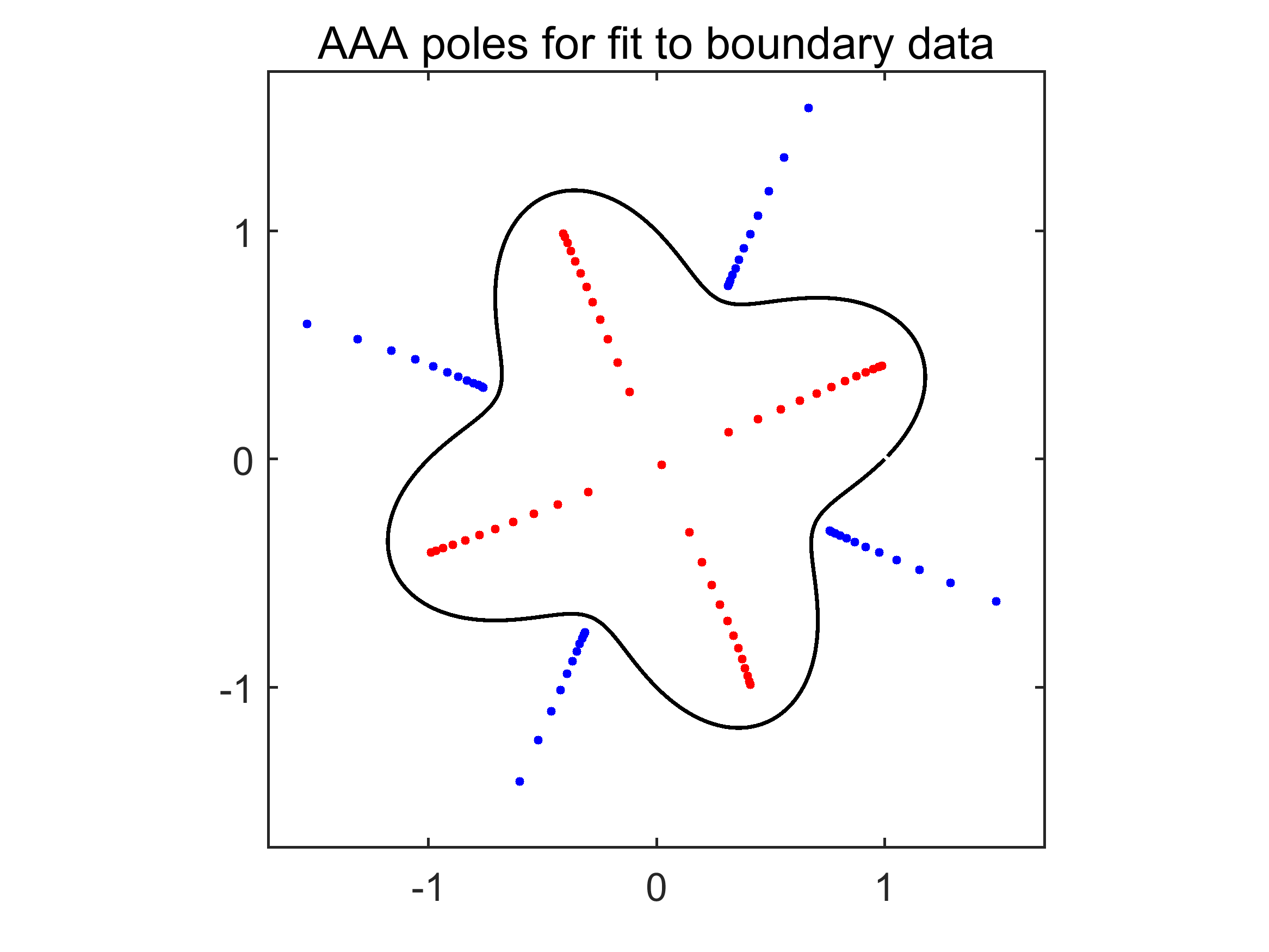}%
\includegraphics[trim=43 0 43 0, clip,scale=0.44]{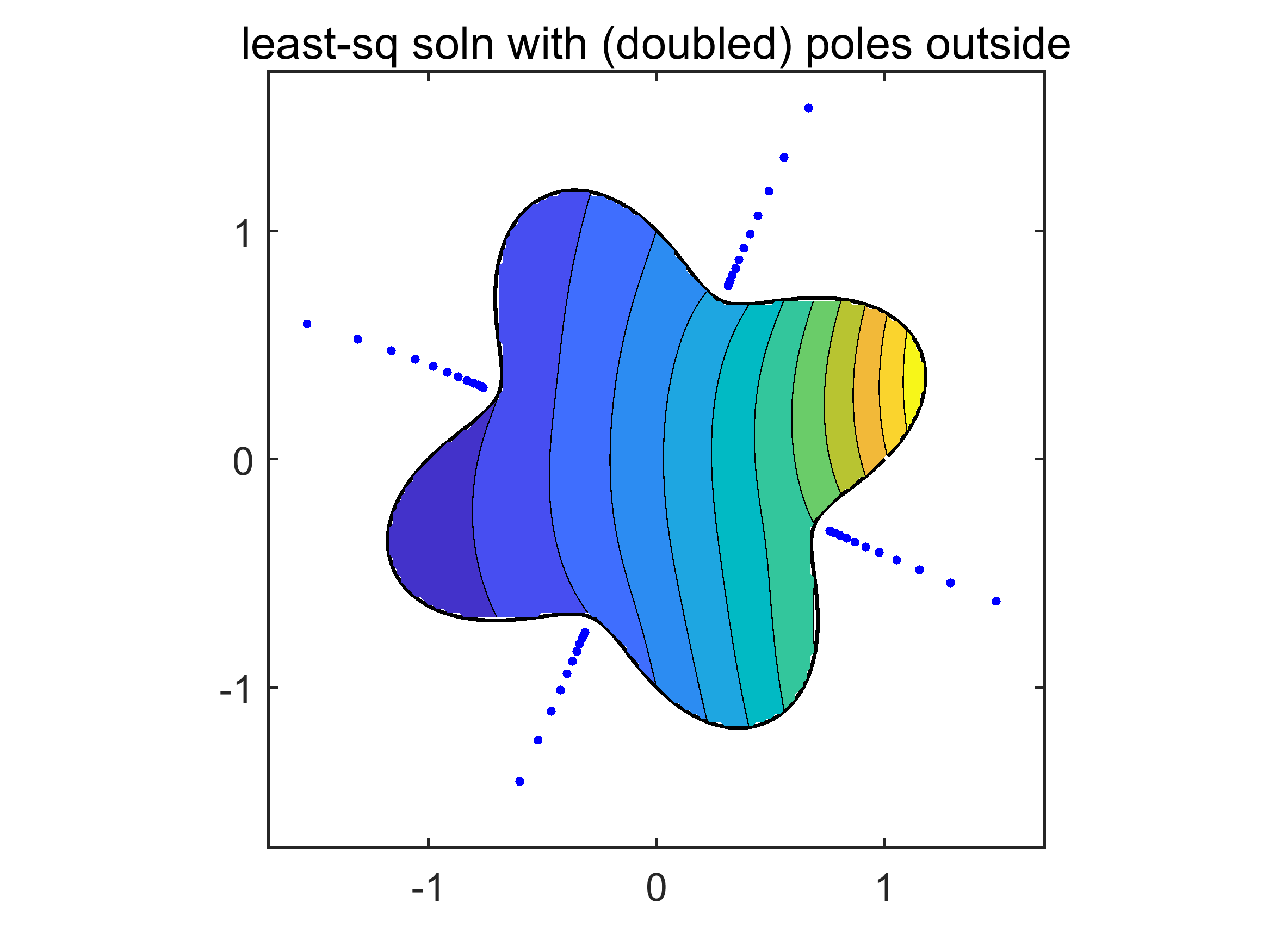}
\vskip 2pt
\includegraphics[trim=43 0 43 0, clip,scale=0.44]{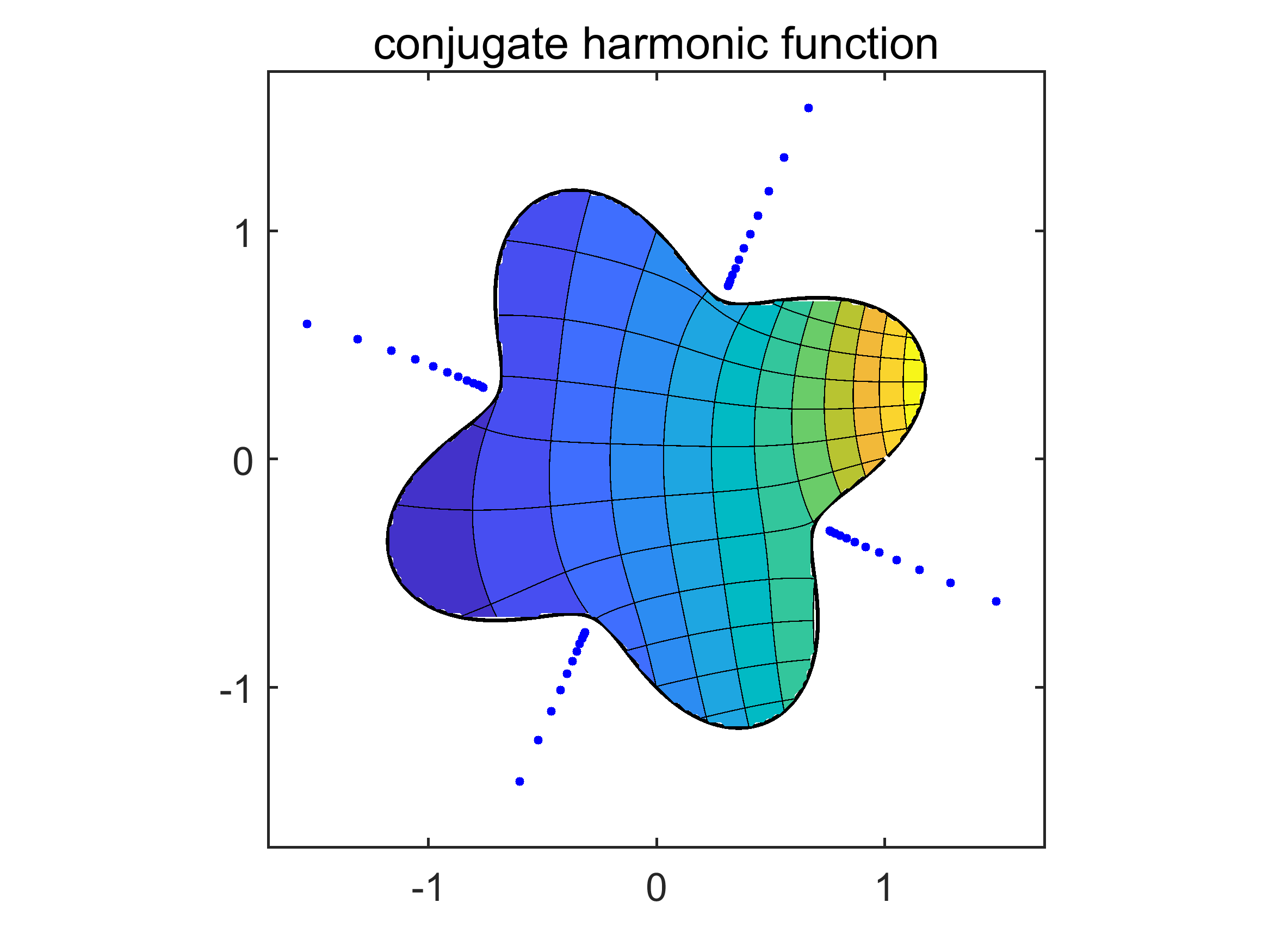}%
\includegraphics[trim=43 0 43 0, clip,scale=0.44]{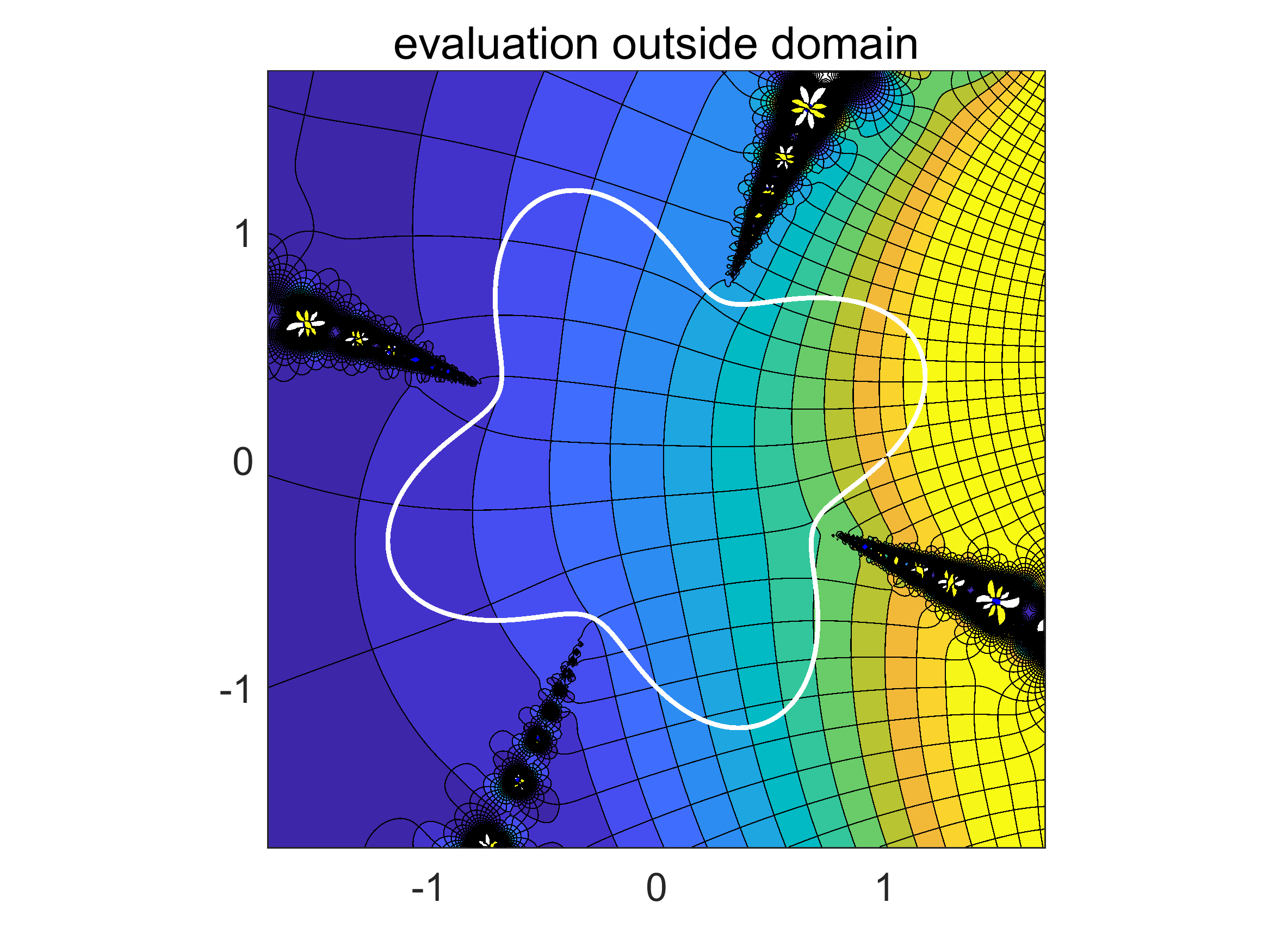}
\end{center}
\vspace*{-10pt}
\caption{\label{aaals}Illustration of AAA-least squares (AAALS) for a Laplace Dirichlet problem on
a smooth nonconvex domain $\Omega$ with boundary data $h(z) = \exp(\Re(z))$.  First
$h$ is approximated by a complex rational function $r$ with poles inside 
and outside $\Omega$.  The poles inside are discarded, and the poles outside are
taken as the locations of 2nd-order poles for a least-squares fit to $h$ by the
real part of a rational function of degree $114$.  In this example
$12$-digit global accuracy is achieved in about $1/4$ second on a laptop; we estimate
the value at the centre as $u(0)\approx 1.2157681370405$.  The same
accuracy with a polynomial fit would require a degree in the thousands, which
would increase to millions for slightly deeper indentations.}
\end{figure}

Let us illustrate how AAA works for a simple problem, shown in
Figure~\ref{aaals}.  We define the domain and the boundary
function like this:

{\small
\begin{verbatim}
      t = 2*pi*(1:500)'/500;
      r = 1 + sin(4*t)/4;
      Z = r.*exp(1i*t);
      h = @(z) exp(real(z));
\end{verbatim}
\par}

\noindent
A complex rational AAA fit is then computed to the boundary data,
which turns out to have degree $105$:

{\small
\begin{verbatim}
      [r,pol] = aaa(h(Z),Z,'mmax',300);
\end{verbatim}
\par}

\noindent
This rational function has 48 poles inside $\Omega$ and $57$ outside.
The inner poles are discarded and the outer ones are retained for
a least-squares fit.

{\small
\begin{verbatim}
      inpoly = @(z,w) inpolygon(real(z),imag(z),real(w),imag(w));
      ii = inpoly(pol,Z); p = pol(~ii).';
\end{verbatim}
\par}

\noindent
From here, if we proceed as described above, we will fit $h(z)$
by the real part of a rational function of degree 57, and the
maximum absolute error (measured a posteriori on a finer grid)
will be $6.6\cdot 10^{-8}$.  This improves to $5.0\cdot 10^{-13}$,
however, if we use a modification of the original AAALS method in
which instead of simple poles, we fit with double poles, following
a rationale to be explained in a moment along with a potential
further improvement.  The code looks as follows.
We find a vector $c$ of fitting coefficients by solving a $500\times
230$ matrix least-squares problem.  This vector is then used to
construct function handles for the approximate solution $u$, its
harmonic conjugate $v$, and the analytic combination $f = u + iv$.

{\small
\begin{verbatim}
      d = min(abs(p-Z));
      P = Z.^0; Q = [d./(Z-p) d.^2./(Z-p).^2];
      A = [real(P) real(Q) -imag(P) -imag(Q)];
      c = reshape(A\h(Z),[],2)*[1; 1i];
      f = @(z) reshape([z(:).^0 d./(z(:)-p) ...
                        d.^2./(z(:)-p).^2]*c, size(z));
      u = @(z) real(f(z)); v = @(z) imag(f(z));
\end{verbatim}
\par}

\noindent
A curious fine point is that the matrix $A$ is of dimensions
$500\times 230$, not $500\times 229$, as one of its columns (the 116th) is the imaginary
part of $1$, namely $0$.  Thus $A$ has infinite condition number,
raising questions about the associated least-squares problem.
As it happens, the Matlab backslash operation ignores such
columns \cite[sec.~5.7]{moler}, so we have not found it necessary to remove them explicitly
from our matrices.

The third panel of Figure~\ref{aaals} shows level curves of
the conjugate function $v(z)$, making an orthogonal net in the
problem domain as always with Laplace problems.  The fourth panel
extends the colours and contour lines of the rational approximation
$\Re\kern 1pt (r(z))$ outside the domain.  This is not the same
as an analytic continuation of $u(z)$ itself, which would be
multivalued, but one sees that $u(z)$ must have branch points
in the four indentations.  A plot like this highlights
the difference between solving \PDEs\ by rational approximation
and by integral equations.  With integral equations---a highly
advanced subject computationally---the name of the game is to
compute the unique single- and/or double-layer density functions on
the boundary~$\Gamma$ that generate the solution in the interior \ccite{kress}.
This necessitates accurate representation of the density functions
followed by accurate calculation of an integral whenever the solution
is to be evaluated in $\Omega$.  With rational functions, the
solution is represented by poles outside the boundary rather than
charge distributions along the boundary. The poles and residues are
numerically non-unique, and evaluation of $u(z)$ for $z\in\Omega$
is just a matter of adding them up.
Rational approximation methods for \PDEs\ are thus a variant of the
Method of Fundamental Solutions (MFS) \ccite{doicu,fk} in which the singularity locations
are chosen with exceptional efficiency, which ordinarily is
such a big challenge for these methods that it is probably the main
reason why they have not been more widely used \ccite{barnettb}.

No theory has yet been developed for AAALS beyond a couple
of theorems in \ccite{aaals} applying just on a disk or a half-plane
together with counterexamples there showing that the method does not
always work.  However, here is our understanding of why it usually
works in practice.  In typical Laplace problems involving
boundaries with corners or indentations, most of the difficulty
is geometry-driven rather than data-driven in the sense that the
singularities of the solution being computed, when analytically
continued outside $\Omega$, are essentially those of the Schwarz
function for $\Omega$, as discussed in section \ref{sec:schwarz}\kern
1pt; see \ccite{inlets}.  (Indeed, the results come out much the
same if {\tt aaa(h(Z),Z)} above is changed to {\tt aaa(conj(Z),Z)}.)
The poles chosen by the AAALS method are thus adapted
to the function that needs to be approximated.

An important challenge for the next few years is to make this
rigorous and develop a true theory of AAALS solution of Laplace
problems.

We can now explain our understanding of why it may be advantageous
to use poles of order 2 instead of 1 in the least-squares fit, while
emphasising that solid understanding of how AAALS works and how to
optimise it is not yet available.
One might imagine that this compensates for the fact that about
half the poles have been thrown away, those inside $\Gamma$, but
we believe the actual logic is different.  In fact, we think a more
accurate view may be that the
doubling of the poles compensates for the ``factor of 2'' discussed
in section~\ref{sec:potential} in connection with the transition from
(\ref{potential-potential}) to (\ref{potential-potential2}) and the
green and yellow dots in the figures of that section.  From ordinary
potential theory involving the function (\ref{potential-potential}),
one gets a certain accuracy that is on the order of half that of
the AAA approximation $r\approx h$ on $\Gamma$.  There would exist
some very slightly perturbed choice of poles, corresponding to a
best harmonic approximation of $u$, for which this accuracy could
be doubled.  However, the AAA-LS rational approximation does not provide
this appropriate perturbation
of the poles.  Instead, we explicitly double the number of poles,
yielding an approximation with more or less the accuracy associated
with an optimal approximation of degree $n$\kern 1pt ---but it is
actually a rational function of degree $2\kern .3pt n$.

\begin{figure}
\vskip 5pt
\begin{center}
\includegraphics[trim=43 0 43 0, clip,scale=0.44]{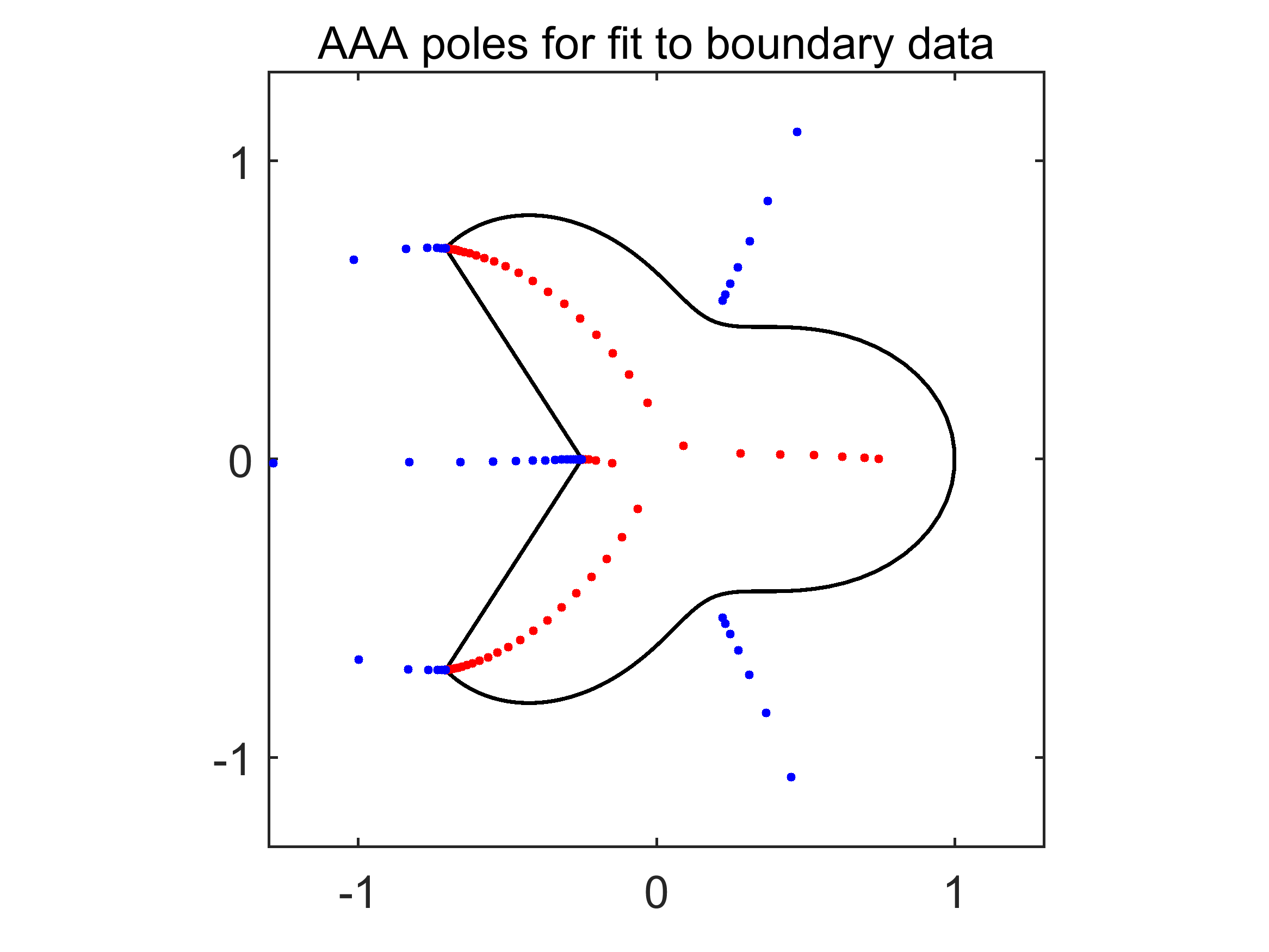}%
\includegraphics[trim=43 0 43 0, clip,scale=0.44]{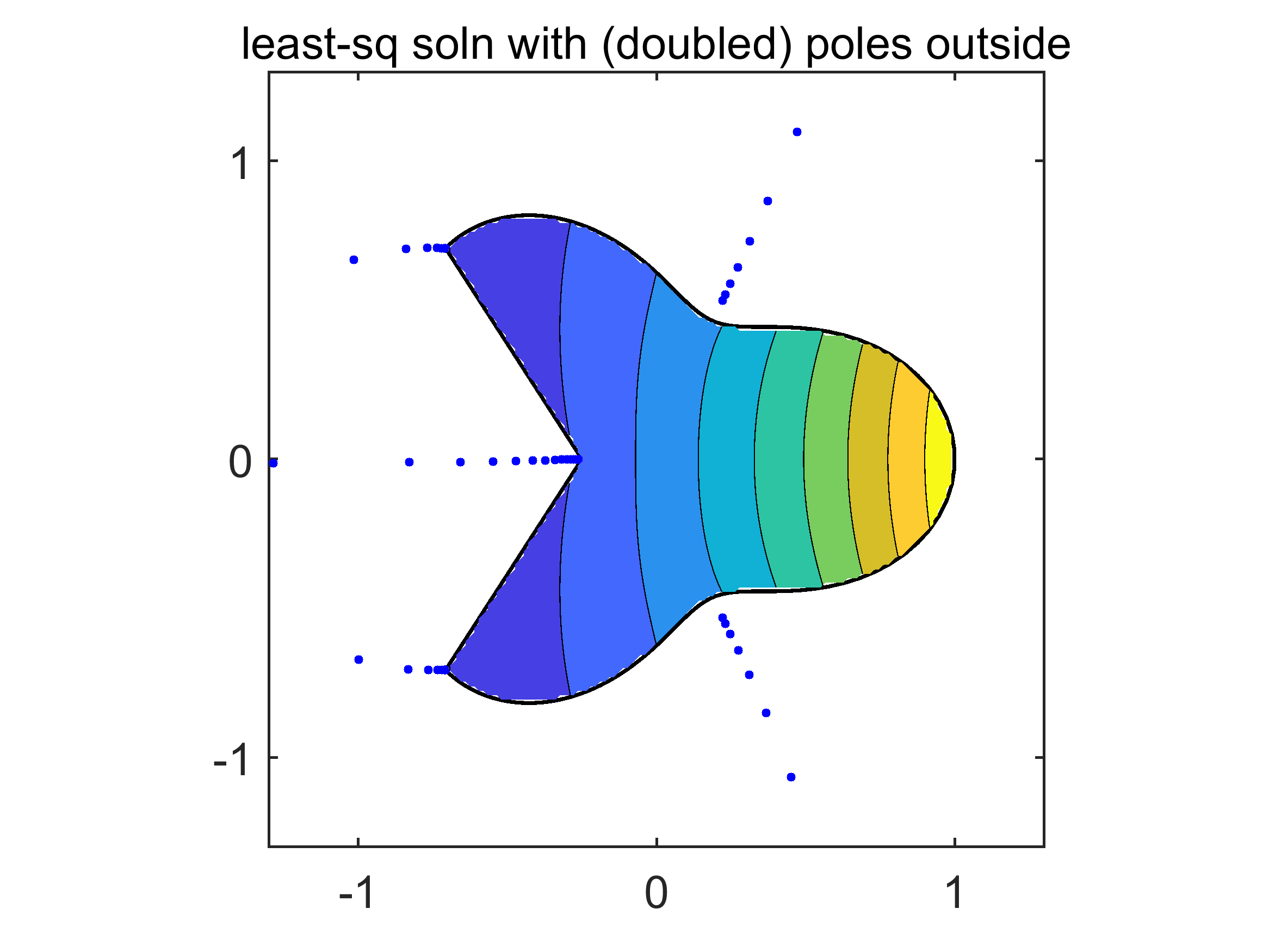}
\vskip 2pt
\includegraphics[trim=43 0 43 0, clip,scale=0.44]{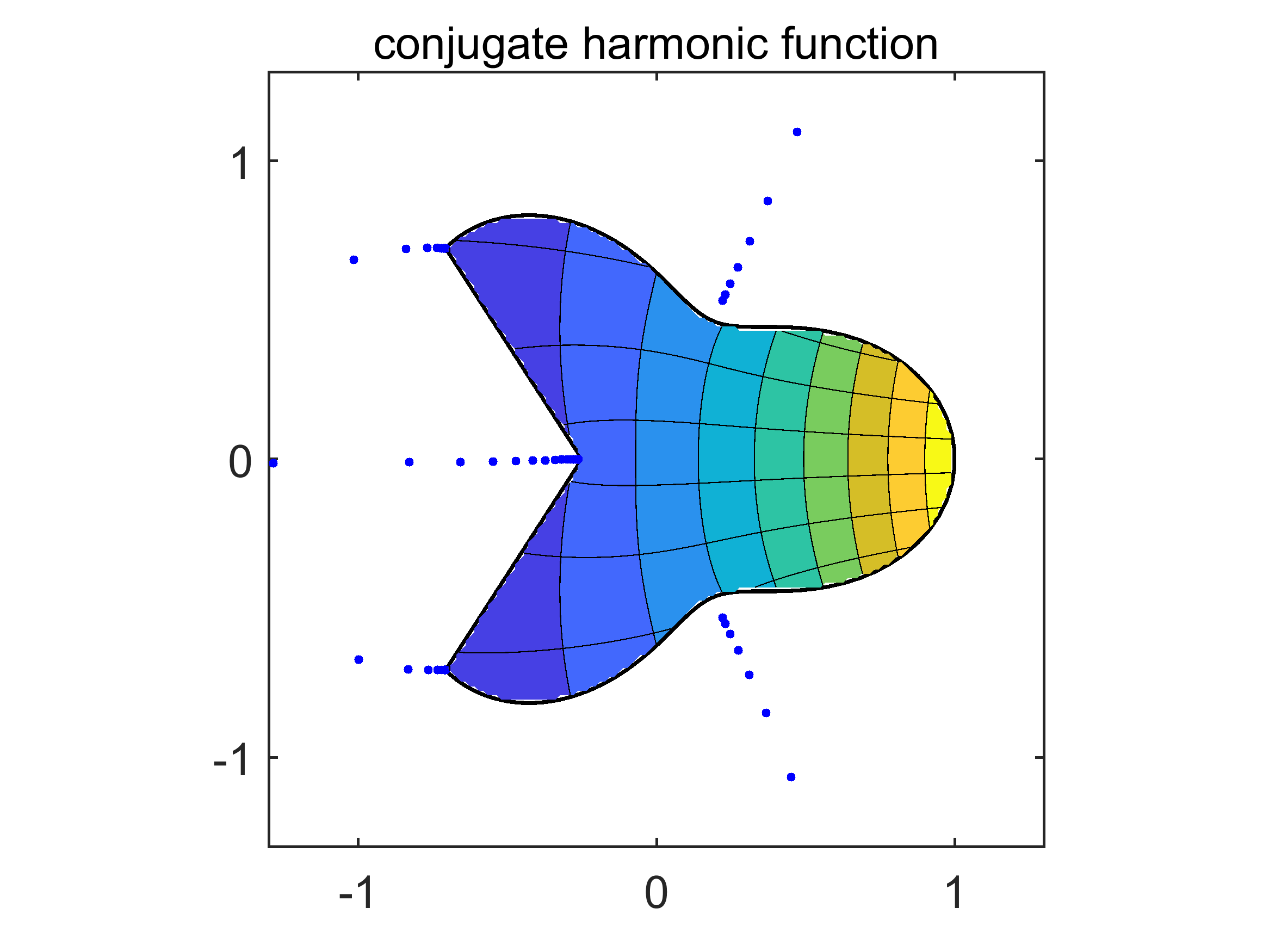}%
\includegraphics[trim=43 0 43 0, clip,scale=0.44]{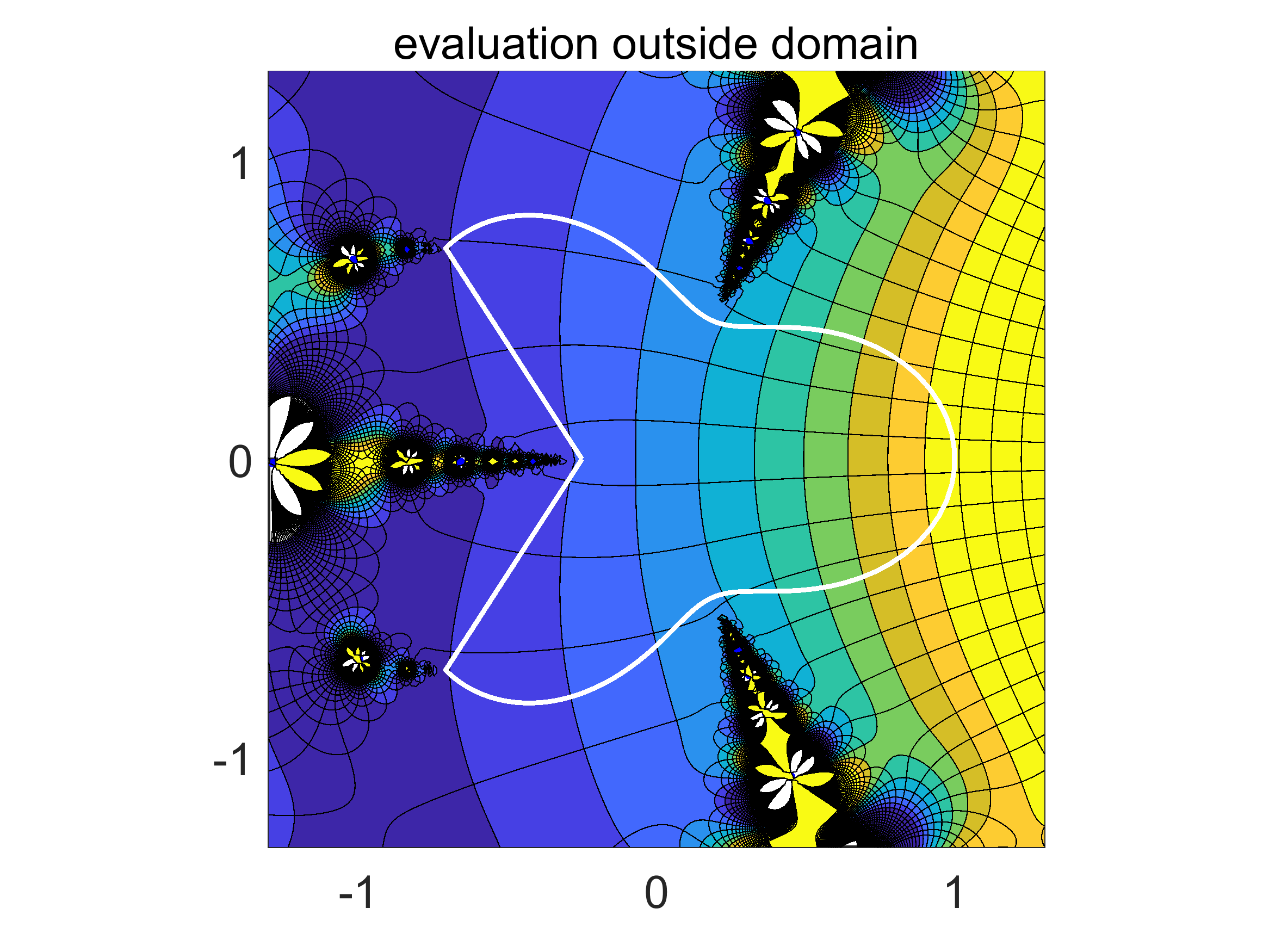}
\end{center}
\vspace*{-10pt}
\caption{\label{bigaaals}Repetition of Figure~\ref{aaals} for a harder
problem with both corner singularities and nonconvexity.  The global
accuracy is about 7 digits, and the value
at the centre is $u(0)\approx 1.081415507437$.}
\end{figure}

Thus AAA-LS with pole doubling tends to provide a solution of degree $2\kern .3pt n$
whose accuracy may approximate that of an optimal solution of degree $n$.
We can optionally exploit this connection by taking the computed solution,
together with its computed harmonic conjugate so that we have a computed
analytic function, and then re-approximating it by AAA, perhaps with the
tolerance slightly loosened.  (This idea is due to Stefano Costa.)
For the present example we can do this with the code

{\small
\begin{verbatim}
      f2 = aaa(f(Z),Z,'tol',1e-12);
      u2 = @(z) real(f(z)); v2 = @(z) imag(f(z));
\end{verbatim}
\par}

\noindent The result is an approximation with about the same accuracy
as before but degree now 54 instead of 115.

We close this section with another figure like Figure~\ref{aaals},
but for a more difficult domain mixing both corners and inlets,
as defined by this code:

{\small
\begin{verbatim}
      s = tanh(linspace(-12,12,1000)');
      r = 1 - .5*sin(pi*s).^2;
      z = -.25 + (exp(-3i*pi/4)+.25)*(s+1)/2;
      Z = [z; exp((3i/4)*pi*s).*r; flipud(conj(z))];
\end{verbatim}
\par}

\noindent The {\tt tanh(linspace(...))} construction clusters
sample points exponentially near corners.  In view of the
corner singularities, the AAA relative tolerance is loosened to
$10^{-8}$, and the absolute accuracy of the solution $u$ comes out
as about $5\cdot 10^{-8}$.  The computation takes 3 s because the
boundary grid has 3000 points.  This could probably be speeded
up by switching to the continuum AAA algorithm (section
\ref{sec:continuum}).

Since its introduction, the AAALS method has been applied by a number
of authors to quite a few problems involving interior, exterior,
and coupled domains, as well as for variations on the Laplace
equation as we will mention in the upcoming sections.  For example,
see \ccite{costamagna,harrismacfloral,harrismacvortex},
where problems with matching conditions across interfaces are treated.   
A doubly-connected geometry is treated in 
\ccite{harrismacwildfire}; another can be found in Figure~\ref{gasket}.
A particularly simple variation
concerns Poisson problems with simple right-hand sides, such as
$\Delta u = C$ for a constant $C$.
Here, one can subtract off a known function like
$u = C(x^2+y^2)/4$ with the desired Laplacian \ccite{harrismacpenguin}, and
what remains is a Laplace problem.  Computing Green's functions
can follow the same pattern with subtraction of a point singularity
$\log|z-z_0^{}|$, a method we shall exploit in the next section.
Much the same mathematics enables AAALS to compute Hele-Shaw fluid flows
as well as solutions to other moving boundary problems associated with
Laplace potentials~\ccite{hinton}.

\section{\label{sec:conformal}Conformal mapping}
Conformal mapping problems can be reduced to Laplace problems, and
the AAALS method of the last section is a good way of solving them.
Before illustrating this, however, we discuss a more distinctive
application of rational approximation to conformal mapping.
This is the matter of compactly {\em representing\/} a map once
it has been computed, as well as its inverse, an idea put forward
in \ccite{conformalrep}.

Numerical conformal mapping is an established topic in numerical
analysis; for high-level and detailed surveys, respectively,
see \ccite{ncm} (5 pages) and \ccite{wegmann} (127 pages), and
see also \ccite{henrici3}.  Many different methods appear in
this field, based often on integral equations and equally often
on other ideas.  In the case of conformal mapping of polygons,
one can take advantage of the Schwarz-Christoffel formula, whose
numerical realisation is available in Driscoll's widely-used SC
Toolbox in Matlab \ccite{sctoolbox}.

Almost every conformal mapping method maps a priori in one direction
or the other: from the problem domain to the canonical domain (say,
the unit disk), or from the canonical domain to the problem domain.
Traditionally, one first solves the mapping problem in the given
direction, typically obtaining the {\em boundary correspondence
function (BCF),} the homeomorphism that maps one boundary to
the other.  Then one applies an auxiliary method for
inverting the map, perhaps by iteration.  The new observation of
\ccite{conformalrep} was that, regardless of how the BCF may have
been computed, both the a priori map and its inverse can often
be represented most efficiently by rational functions.  If $Z$
is a vector discretising one boundary and $F$ is the corresponding
vector obtained by applying the BCF to $Z$, then the computation
of efficient representations of the conformal map and its inverse
may be as simple as
$$
\hbox{\small\verb|f = aaa(F,Z); finv = aaa(Z,F);|}
$$
at least for smooth domains.  This is an example of AAA computation
of inverse functions as discussed in section \ref{sec:inv}.
(In mapping a domain with corners, it will probably be necessary
to loosen the tolerance for reasons of ill-conditioning as well
as efficiency, and additionally, one must be careful to cluster
sample points near the corners and their preimages, or to let the
clustering be taken care of automatically by the continuum AAA
algorithm of section \ref{sec:continuum}.)

For example, the following code defines the smooth cross-shaped
domain shown in the top row of Figure~\ref{conf} and then invokes
the {\tt conformal} command of Chebfun, which computes the BCF of
a domain onto the unit disk by means of a discretisation of the
Kerzman-Stein integral equation \ccite{kerz}:

\begin{figure}[t]
\vskip 5pt
\begin{center}
\includegraphics[trim=0 90 0 10, clip,scale=0.66]{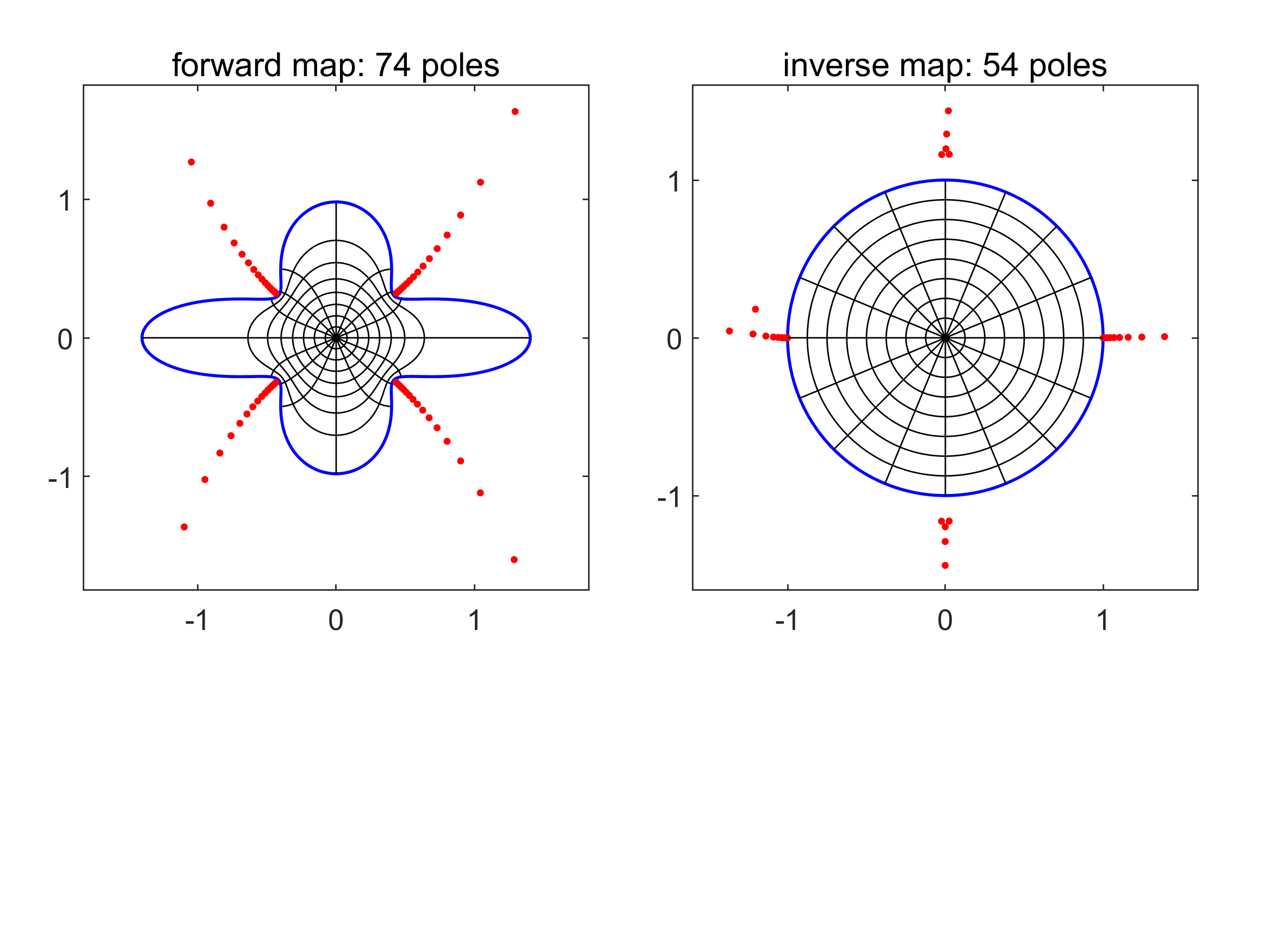}
\vskip -10pt
\hbox{\kern 40pt\includegraphics[trim=0 90 0 10, clip,scale=0.67]{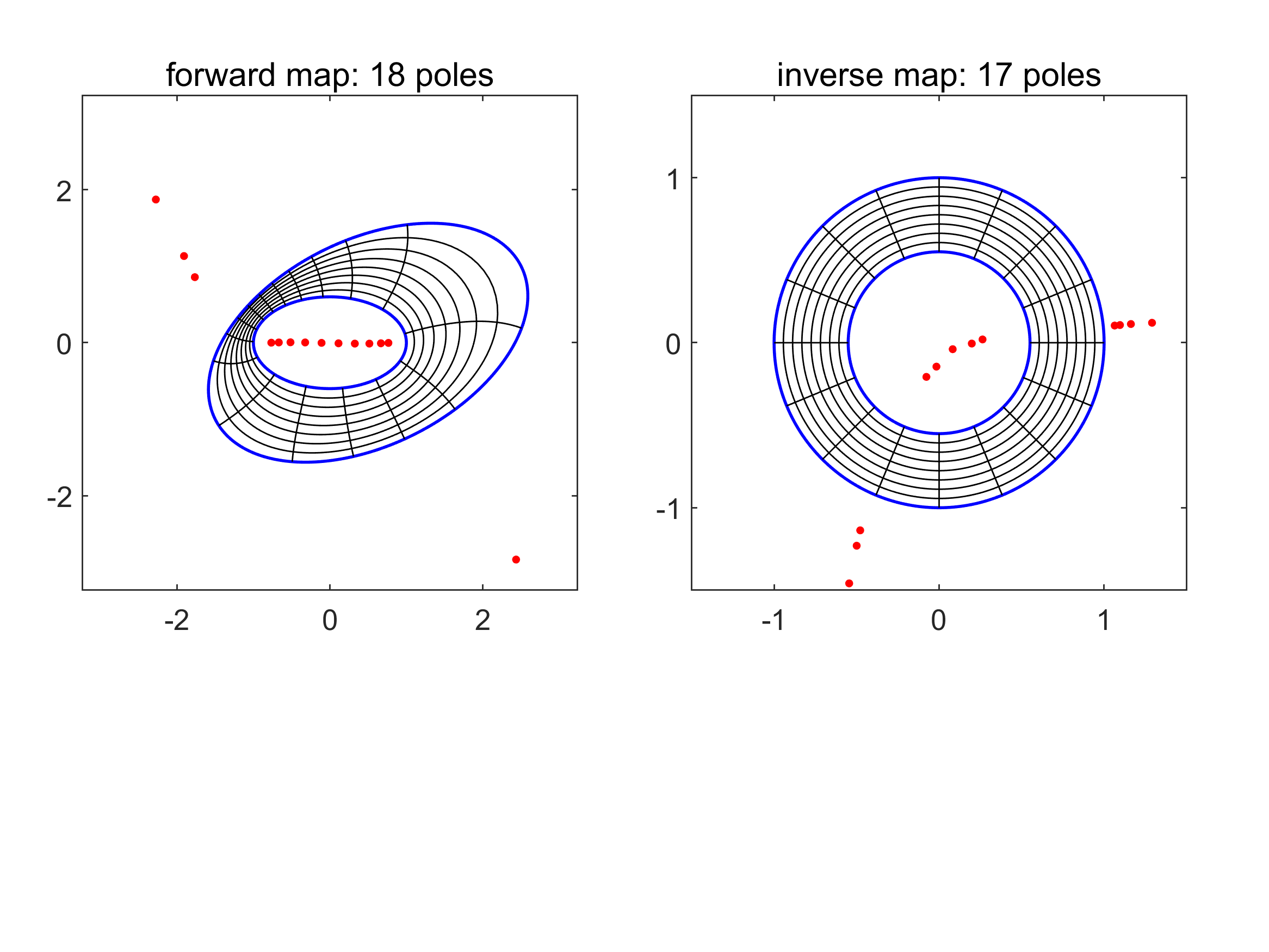}}
\end{center}
\vspace*{-16pt}
\caption{\label{conf}First row: conformal map computed by Chebfun
{\tt conformal} from a smooth cross-shaped
domain to a disk, and its inverse.  Each map
is represented to 9-digit accuracy by a rational function obtained by AAA,
and can be evaluated in about $0.5$ $\mu$\kern .3pt s per point.
Second row: doubly-connected conformal map computed by Chebfun
{\tt conformal2} and its inverse.  Each map
is represented to 7-digit accuracy by a rational function obtained by AAA,
and can be evaluated in about $0.1$ $\mu$\kern .3pt s per point.}
\end{figure}

{\small
\begin{verbatim}
      C = chebfun('exp(pi*1i*t)*(1+.4*cos(4*pi*t))','trig');
      C = real(C) + .7i*imag(C);
      tol = 1e-10;
      [f,finv] = conformal(C,'plots','tol',tol);
\end{verbatim}
\par}

\noindent 
The figure shows 
conformally equivalent nets of orthogonal contours in the two
domains.  As a sample value we have $f(1) \approx 0.984313068276404$,
accurate we believe to these 15 digits.  What is notable about
this figure for our purposes, however, are the poles plotted as
red dots outside each boundary, whose numbers are indicated in the
titles.  The {\tt conformal} code has called AAA twice and found that
the map can be represented to the prescribed tolerance by a rational
function of degree $74$ in the direction from the cross to the disk, and
by another rational function of degree $54$ in the direction from
the disk to the cross.  Each of these maps can be evaluated on our laptop
in less than 1 $\mu$\kern .4pt s per point, and since the representations
utilise global rational functions, no numerical discontinuities will
appear even if one or more derivatives are taken.  Note that although the smooth cross
contour is analytic, its nonconvexity is such that poles of $f$ come
as close as $0.024$ to the boundary; likewise the poles of $f^{-1}$
come as close as $0.0008$ to the unit circle.  It follows that
polynomial representations of these maps would require degrees in
the tens of thousands to get a few digits of accuracy \ccite{inlets}.
Conformal mappers refer to such distortions 
as the ``crowding phenomenon,'' which makes
certain numerical methods terribly inefficient but is circumvented
by the strategic placement of poles of near-best rational approximations. 
See \ccite{wegmann} and \cite[Thms.\ 2--5]{conformalrep}.

The second row of Figure \ref{conf} shows the same kind of image,
but for a doubly-connected domain mapped with the Chebfun {\tt
conformal2} command (based on series expansions as in \ccite{series})
by this code:

{\small
\begin{verbatim}
      circle = chebfun('exp(1i*pi*z)','trig');
      ellipse2 = real(circle) + .6i*imag(circle);
      ellipse1 = (2+1i)*ellipse2 + .5;
      [f,finv,rho] = conformal2(ellipse1,ellipse2,'plots'); 
\end{verbatim}
\par}

\noindent For a doubly-connected domain like this, the circular
annulus image has a uniquely determined ratio of inner to outer radii,
known as the {\em conformal modulus\/} of the domain.  Although
{\tt conformal2} with its default tolerance has pointwise accuracy
of just about $10^{-7}$ for this problem, numerical explorations
suggest that the computed modulus $\rho\approx 0.55098544815237$
is accurate to $15$ digits.

We now turn to the other connection of rational approximation with
conformal mapping, the use of AAALS to compute a map by
solving a Laplace problem, an idea that was first suggested
in \ccite{aaalscosta} and is discussed more
extensively in \ccite{conformal}.  We start with the basic case of the
conformal map $f$ of a Jordan
domain $\Omega$ containing the
origin to the unit disk, with $f(0) = 0$ and $f'(0) > 0$.
Following an old idea presented for example
in \cite[Theorem 16.5a]{henrici3}, we note that
$g(z) = \log(f(z)/z)$ is a nonzero analytic function on
$\Omega$ with
real part $-\log|z|$ on the boundary $\Gamma$.
If we write $g(z) = u(z) + i\kern .3pt v(z)$, where $u$ and
$v$ are harmonic functions, then $u$ is the solution of
the Laplace Dirichlet problem
\begin{equation}
\Delta u = 0, ~ z\in \Omega, \quad
u(z) = -\log|z|, ~ z\in \Gamma
\label{conflap}
\end{equation}
and $v$ is its harmonic conjugate in $\Omega$ with $v(0) = 0$.
This implies
\begin{equation}
f(z) = z\kern .4pt e^{u(z) + i\kern .2pt v(z)},
\label{confcomb}
\end{equation}
so a solution to (\ref{conflap}) solves the
conformal mapping problem.  Note that, following a remark at the end of the
last section, 
$u(z) + \log|z|$ is the Green's function of $\Omega$ with
respect to the point $z=0$, so $f$ is essentially the exponential of
the Green's function.

For example, the following code computes the map of
the upper row of
Figure~\ref{conf} to about 15 digits of accuracy in half a second
on our laptop, which is about ten times faster than
the Chebfun {\tt conformal} code that was used to generate that figure. 
First we apply AAA to approximate $-\log|z|$ on a 600-point
discretisation of the boundary.

{\small
\begin{verbatim}
      C = chebfun('exp(pi*1i*t)*(1+.4*cos(4*pi*t))','trig');
      C = real(C) + .7i*imag(C);
      Z = C(-1 + 2*(1:600)'/600);
      U = -log(abs(Z)); 
      [r,pol] = aaa(U,Z,'mmax',300,'tol',1e-10);
\end{verbatim}
\par}

\noindent Then we discard the poles in the domain and use
AAALS as in the last section to solve the Laplace problem.
The computed map $f$ gives the same value
$f(1) \approx 0.984313068276404$ as before.

{\small
\begin{verbatim}
      inpoly = @(z,w) inpolygon(real(z),imag(z),real(w),imag(w));
      p = pol(~inpoly(pol,Z)).'; d = min(abs(p-Z));
      P = Z.^0; Q = [d./(Z-p) d.^2./(Z-p).^2];
      A = [real(P) real(Q) -imag(P) -imag(Q)];
      c = reshape(A\U,[],2)*[1; 1i];
      g1 = @(z) reshape([z(:).^0 d./(z(:)-p) ...
                d.^2./(z(:)-p).^2]*c, size(z));
      g = @(z) g1(z) - 1i*imag(g1(0));
      f = @(z) z.*exp(g(z));
\end{verbatim}
\par}

\begin{figure}
\vskip 12pt
\begin{center}
\includegraphics[trim=0 155 10 10, clip,scale=0.85]{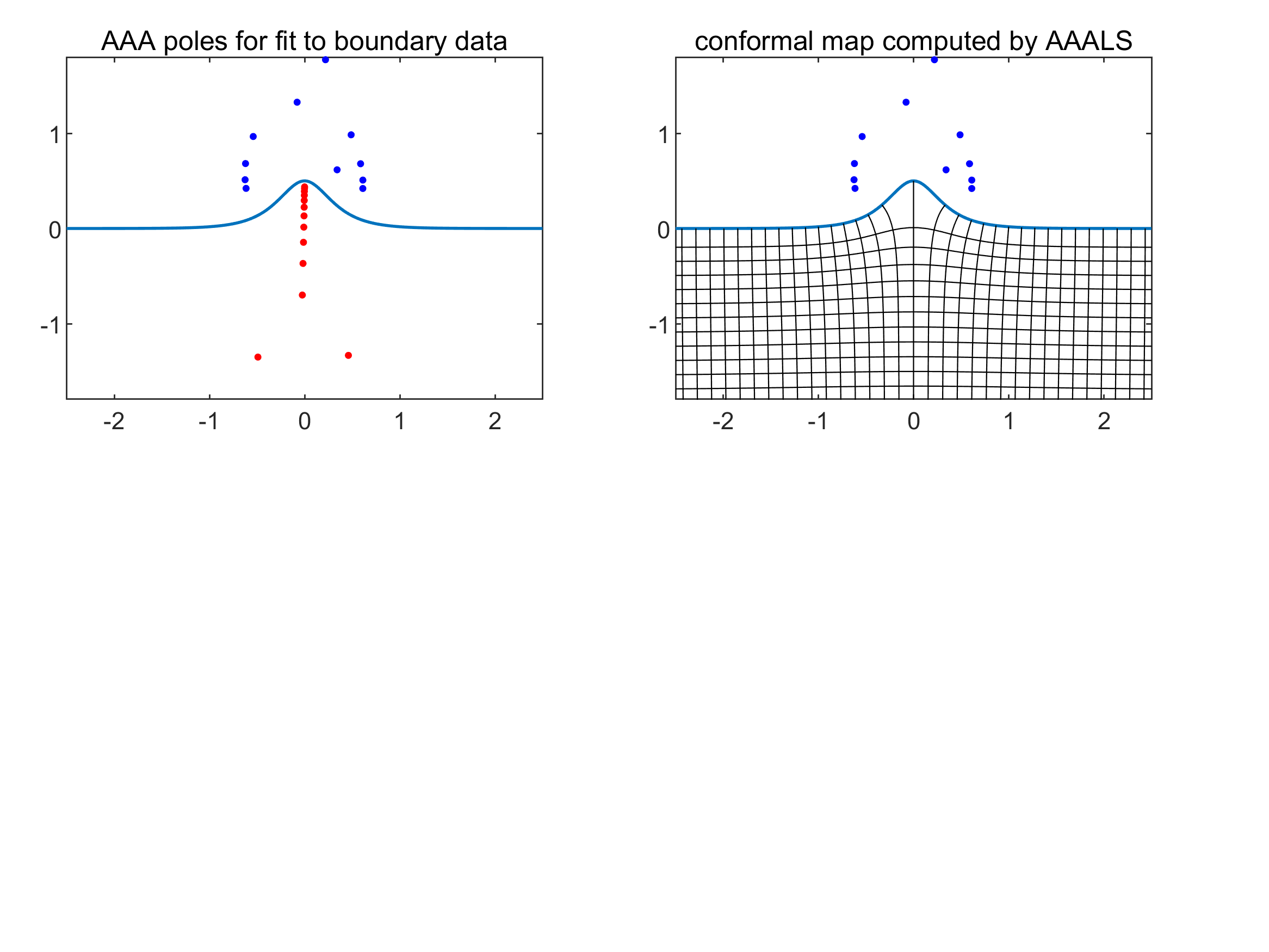}
\end{center}
\vspace*{-18pt}
\caption{\label{wavemapfig}Conformal map onto a half-plane computed
by AAALS.}
\end{figure}

\begin{figure}
\vskip 0pt
\begin{center}
\includegraphics[trim=25 67 200 57, clip,scale=0.81]{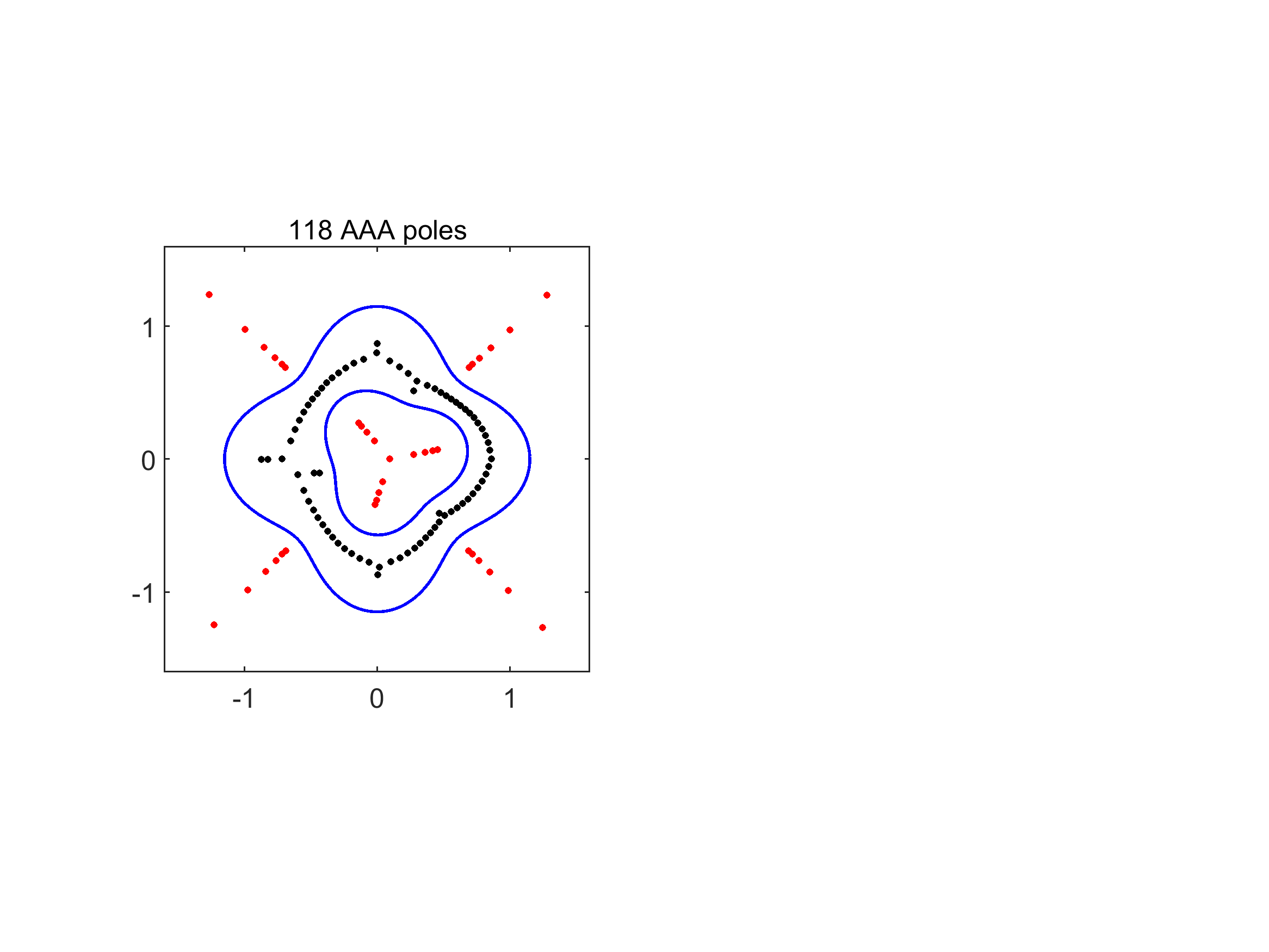}
\vskip -10pt
\hbox{\kern 15pt\includegraphics[trim=20 69 0 57, clip,scale=0.81]{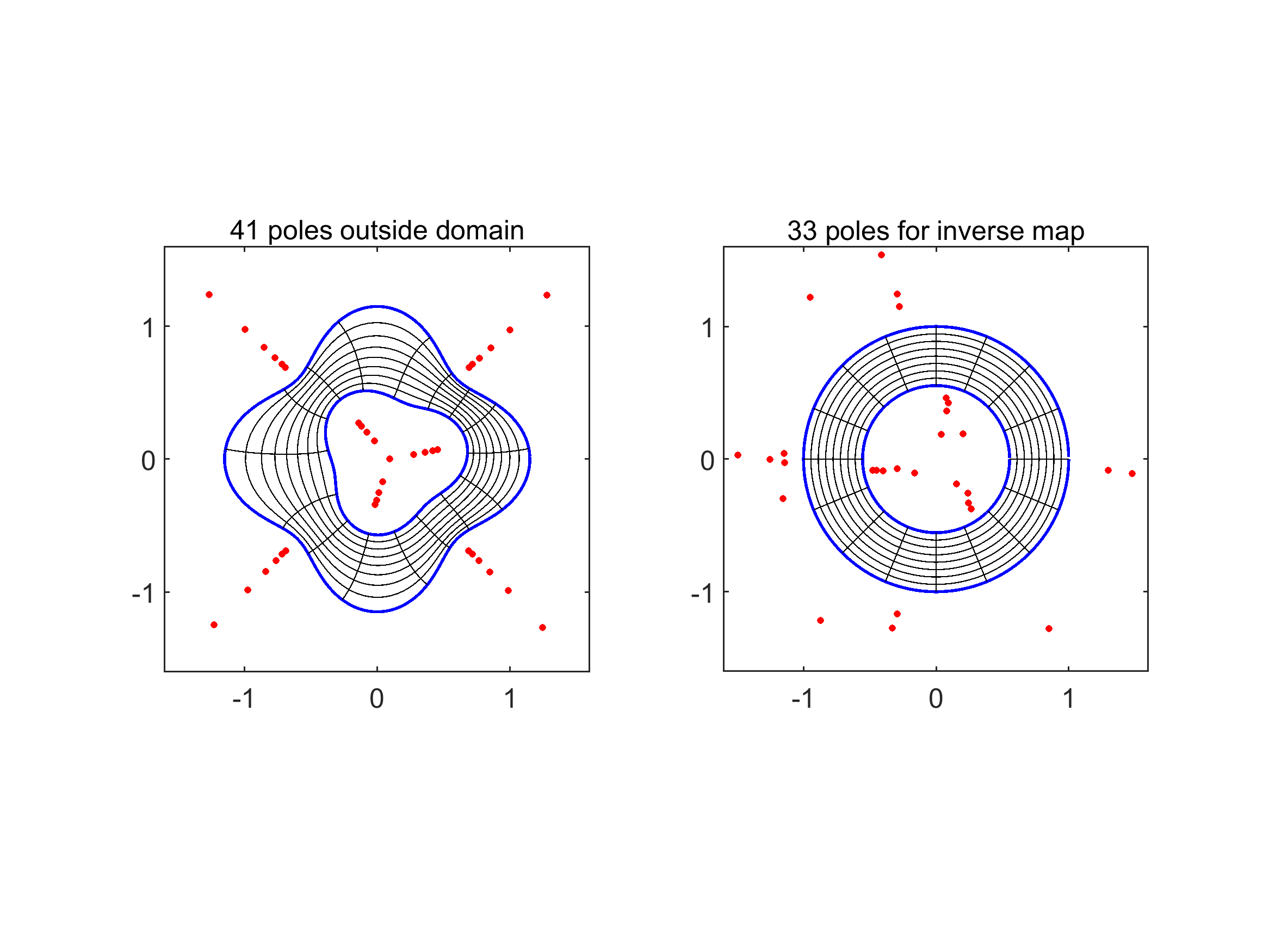}}
\end{center}
\vspace*{-11pt}
\caption{\label{gasket}Doubly-connected conformal map
computed by AAALS.}
\end{figure}

For other conformal mapping configurations, other Laplace problems become
relevant.  For example, Figure \ref{wavemapfig} shows the conformal
map to the lower half-plane of an infinite region under a curve. 
Without spelling out the mathematics, we list a code, below.
The value $f(0)$ comes out as $f(0)\approx -0.20646984\kern .3pt i$.

{\small
\begin{verbatim}
      T = tan((pi/2.000001)*(-1:.01:1)');
      Z = T + .5i*sech(4*T);
      H = -imag(Z); [~,pol] = aaa(H,Z,'tol',1e-8);
      ZZ = [Z(1)-100i; Z; Z(end)-100i];
      p = pol(~inpoly(pol,ZZ)).';
      d = min(abs(p-Z));
      P = Z.^0; Q = [d./(Z-p) d.^2./(Z-p).^2]; 
      A = [real(P) real(Q) -imag(P) -imag(Q)]; 
      c = reshape(A\H,[],2)*[1; 1i];
      w = @(z) z + 1i*reshape([z(:).^0 d./(z(:)-p) ...
               d.^2./(z(:)-p).^2]*c, size(z));
      W = w(Z); f = aaa(W,Z); finv = aaa(Z,W);
\end{verbatim}
\par}

\noindent
A doubly-connected example is shown in Figure~\ref{gasket}, corresponding
to the code below.  Again we do not explain the mathematics;
see \ccite{conformal}.  The conformal modulus comes out
as $\rho \approx 0.5544210$.
For higher-connectivity conformal mapping with AAA, see \ccite{mckee}.

{\small
\begin{verbatim}
      z = exp(2i*pi*(1:400)'/400);
      c1 = z.*abs(1+.15*z.^4); c2 = .1+(.5+.1i)*z.*abs(1+.15*z.^3);
      Z = [c1; c2];
      [r,pol] = aaa(conj(Z),Z,'tol',1e-8,'mmax',300); 
      ii = inpoly(pol,c1) & ~inpoly(pol,c2); p = pol(~ii).';
      H = -log(abs(Z));
      rvec = [zeros(np,1); ones(np,1)];
      H(np+1:2*np) = H(np+1:2*np) + 1;
      A = [Z.^0 real(1./(Z-p)) imag(1./(Z-p)) rvec];
      c = A\H;
      logr = 1 - c(end); rho = exp(logr), c(end) = [];
      u = @(z) [ ones(size(z)) real(1./(z-p))  imag(1./(z-p))]*c;
      v = @(z) [zeros(size(z)) imag(1./(z-p)) -real(1./(z-p))]*c;
      f = @(z) z.*exp(u(z)+1i*v(z));
\end{verbatim}
\par}

\section{\label{sec:helmholtz}Helmholtz equation and scattering}

As we have mentioned, the AAALS method of section \ref{sec:aaals}
can be interpreted as an application of the Method of Fundamental
Solutions (MFS) \ccite{doicu,fk} in which the locations of the
singularities are chosen by rational approximation.  This is a
crucial improvement to the MFS, since otherwise, effective placement
of singularities for problems with corners or nonconvex geometry
is difficult \ccite{barnettb,inlets,liu}.

Even before AAA was employed for the calculations,
the idea had arisen of generalising Laplace equation
ideas by using rational approximations
to place singularities for problems involving the Helmholtz equation,
\begin{equation}
\Delta u + \kappa^2 u = 0,
\label{helmeq}
\end{equation}
which is the basic equation of acoustics and scattering theory (seen already
in (\ref{res-eigprob})).  This {\em lightning-Helmholtz method\/} was proposed
in \ccite{pnas} in the context of ``lightning'' rational approximations, in which poles
are positioned
near corners according to an a priori exponential clustering formula rather
than determined implicitly by approximation \ccite{clustering}.
For further developments see (\cite{gopalthesis,ginn}).
What's different with the MFS for the Helmholtz equation
is the nature of the fundamental solutions.  For the Laplace equation, these are poles
$1/(z-\tk)$ or their monopole counterparts $\log|z-\tk|$, both of which are solutions of
the equation except at the singular point $\tk$.\footnote{AAA works with
``dipoles,'' i.e., complex poles $1/(z-\tk)$.  Although dipoles and monopoles
$\log|z-\tk|$ are mathematically almost equivalent, since a dipole can be
approximated by a pair of monopoles, so far as we know, no numerically
effective algorithm has been found for approximation by monopoles.}  For the Helmholtz equation, the
corresponding fundamental solutions are the complex functions
\begin{equation}
H_0^{}(\kappa\kern .5pt|z-\tk|)
\label{hankel0}
\end{equation}
and
\begin{equation}
H_1^{}(\kappa\kern .5pt|z-\tk|) \kern 1.5pt\Re \Bigl({z-\tk\over |z-\tk|}\Bigr), ~~
H_1^{}(\kappa\kern .5pt|z-\tk|) \kern 1.5pt\Im \Bigl({z-\tk\over |z-\tk|}\Bigr),
\label{hankel1}
\end{equation}
where $H_0^{}$ and $H_1^{}$ are Hankel functions of the first kind \ccite{doicu,rokhlin}.

For general domains as opposed to domains dominated by corner
singularities, one would like to upgrade the lightning-Helmholtz
method to a {\em AAA-Helmholtz method\/}.  Unpublished experiments by
ourselves and Stefano Costa since August 2023 have shown that such
a method can work, but there is no published literature as yet, and
certainly no theory.  To illustrate the method, let us repeat Figure \ref{aaals}
but now for (\ref{helmeq}) rather than the Laplace equation.
Instead of a Dirichlet problem in the interior of the daisy, we
consider a scattering problem in the exterior in which we want
to calculate the time-harmonic response to a plane wave with wave
number $\kappa=20$ incident from the left.  Since the computational
domain is the exterior, singularities will be placed in the interior.
Each singularity corresponds to an outgoing circular wave for the
given $\kappa$ (the Sommerfeld radiation condition, imposed by
employing $H_0^{}$ and $H_1^{}$ as above instead of their complex
conjugates), and the boundary condition is that the sum of the
contributions of these outward-radiating fundamental solutions
must cancel the incident wave on the boundary (the
so-called sound-soft boundary condition).  This is the same
principle as with the integral equation of section~\ref{sec:res}.

The opening lines of an illustrative calculation
are mostly the same as before, apart from the different boundary condition:

\begin{figure}[t]
\vskip 10pt
\begin{center}
\includegraphics[trim=49 0 43 10, clip,scale=0.48]{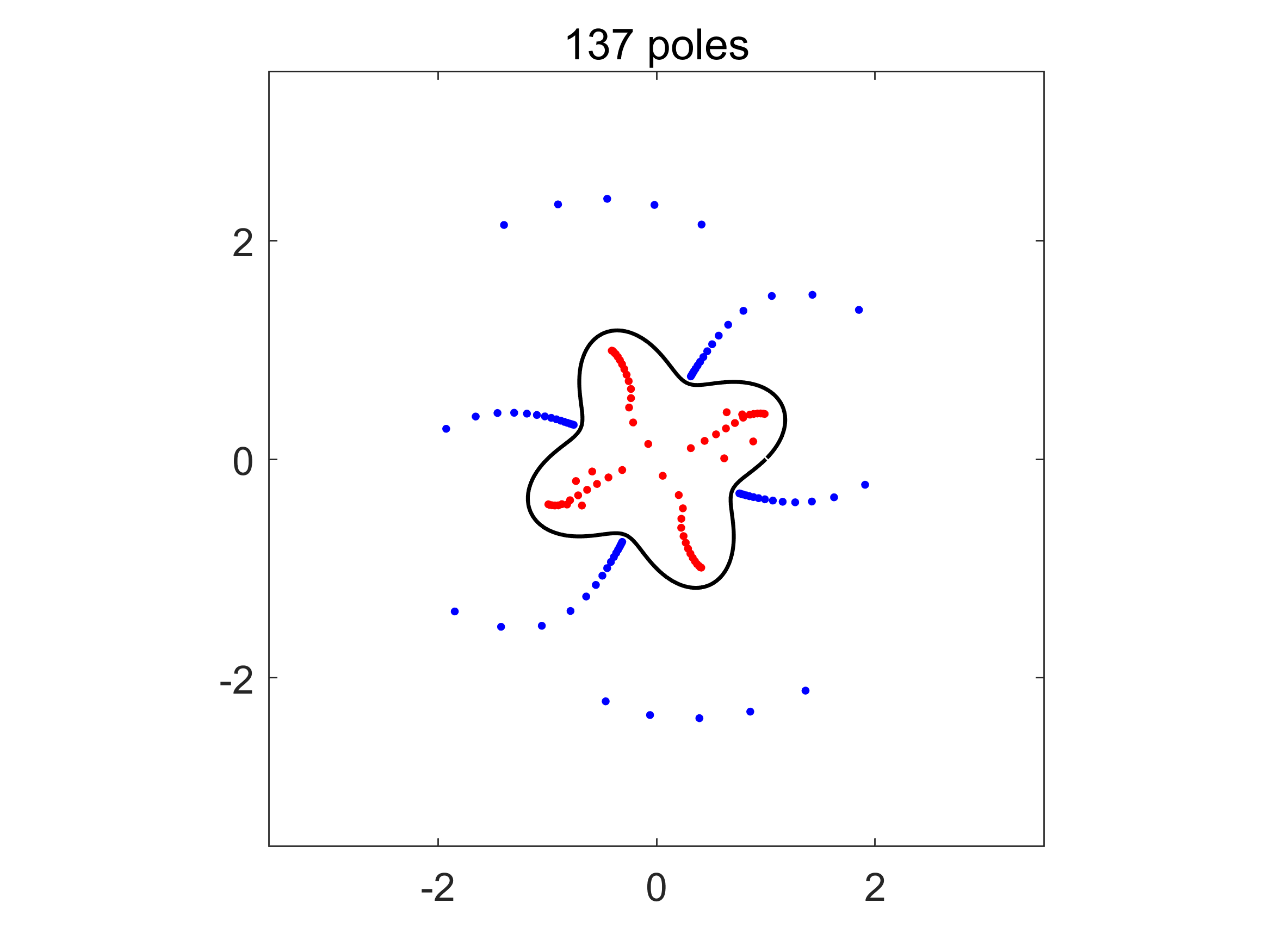}%
\includegraphics[trim=43 0 43 10, clip,scale=0.48]{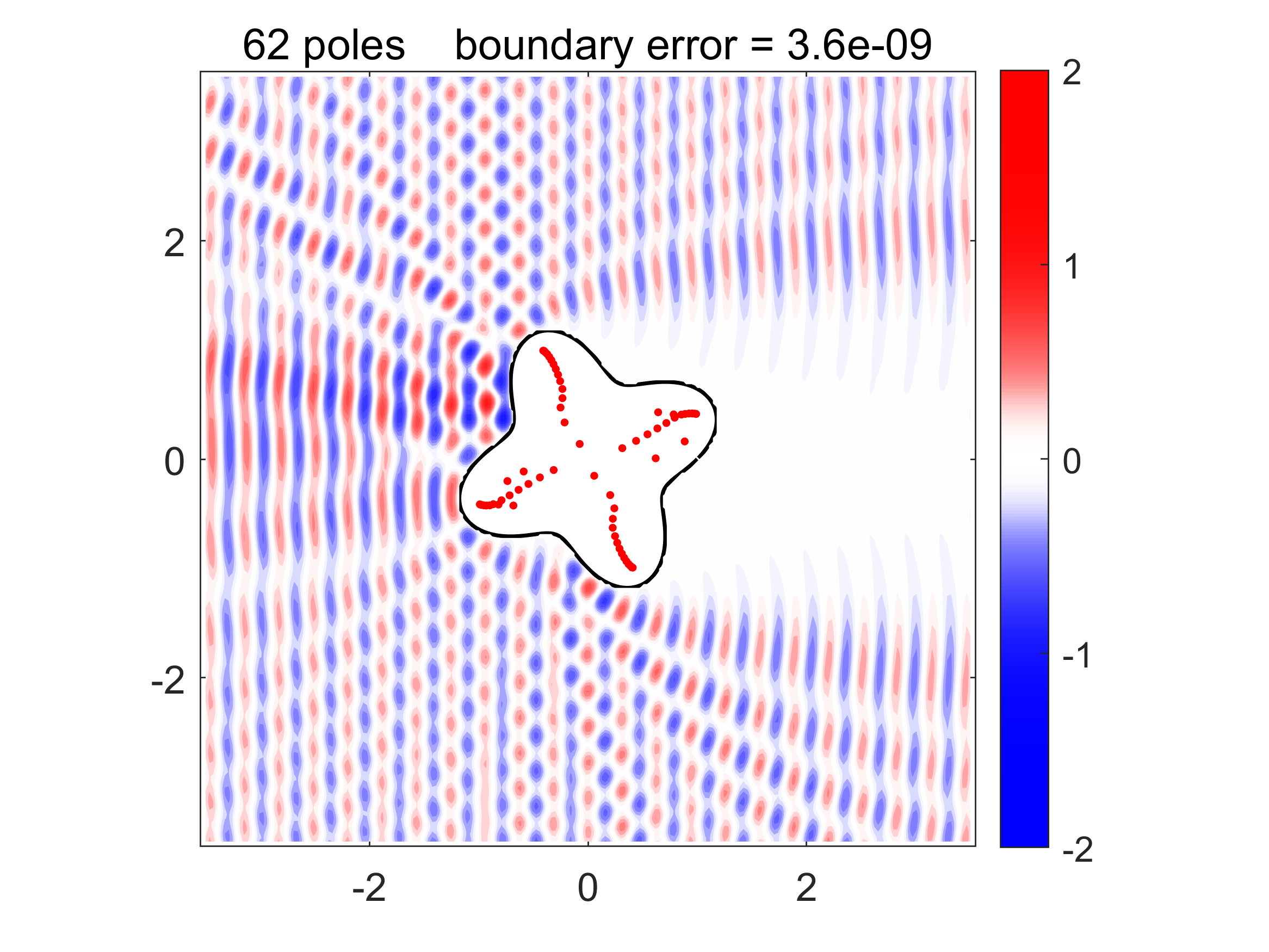}
\end{center}
\begin{center}
\includegraphics[trim=49 0 43 10, clip,scale=0.48]{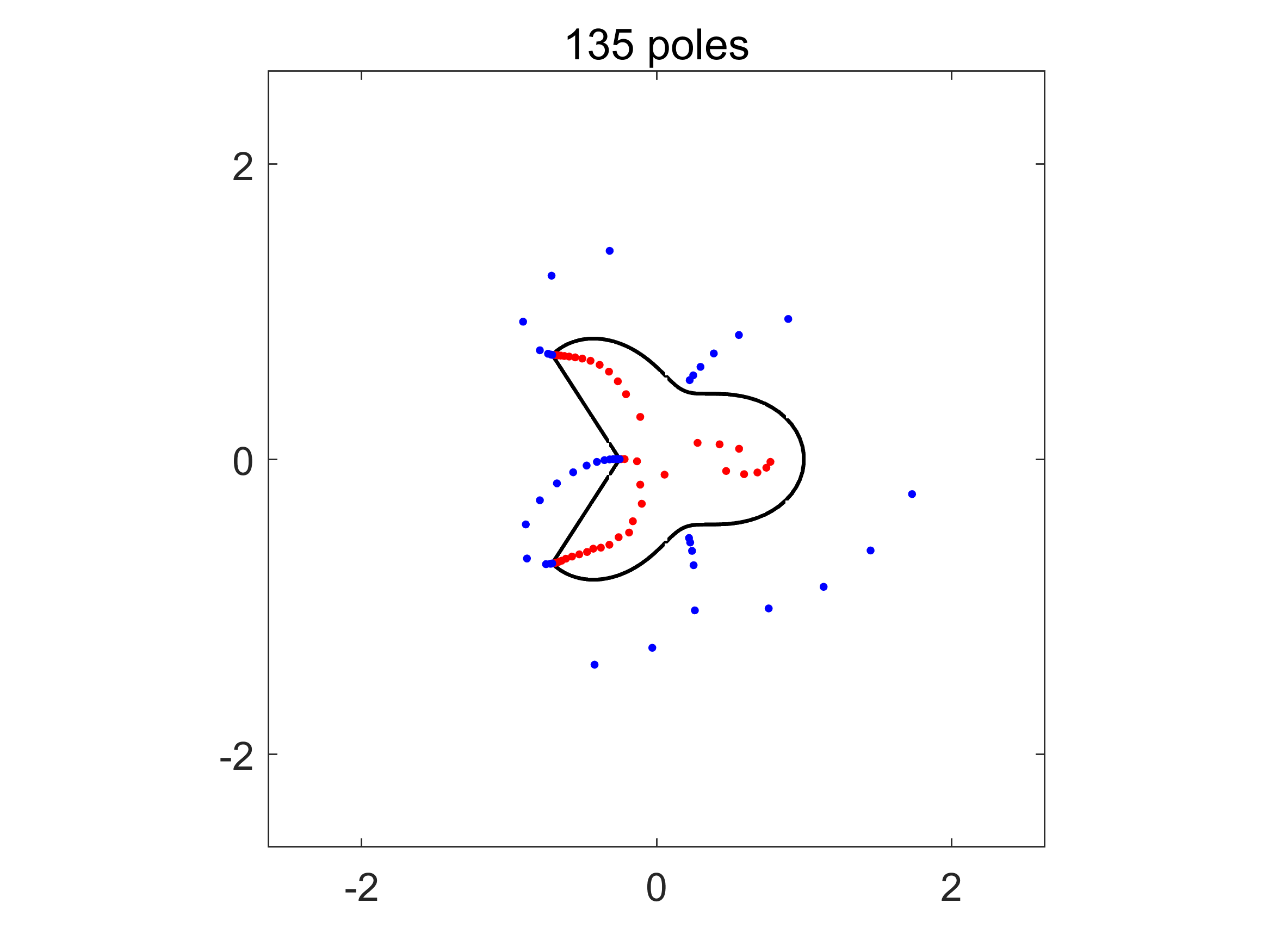}%
\includegraphics[trim=43 0 43 10, clip,scale=0.48]{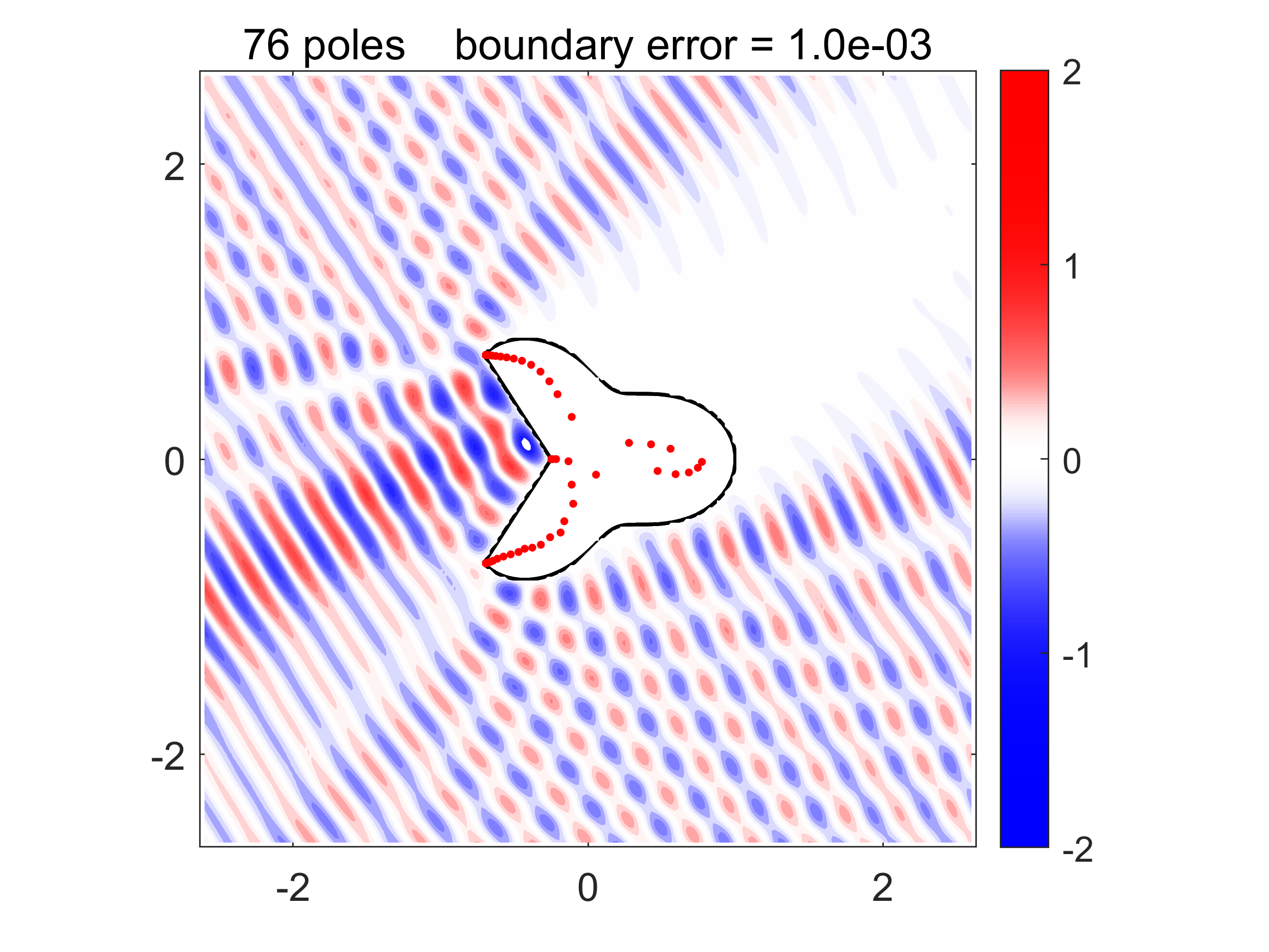}
\end{center}
\vspace*{-10pt}
\caption{\label{helmdaisy}Upper row: a AAA-Helmholtz 
scattering computation
with $k=20$ for the domain of Figure \ref{aaals}; the real part of the
complex solution $u(z)$ is plotted.
Lower row: a similar calculation for the domain
of Figure \ref{bigaaals}.  Here the incident wave has been rotated to
an angle
aligned with one of the straight boundary segments, causing strong reflection in
that direction---amplitudes $\approx \pm 2$ associated with
constructive interference.}
\end{figure}

\begin{figure}[t]
\vskip 10pt
\begin{center}
\includegraphics[trim=49 0 43 10, clip,scale=0.48]{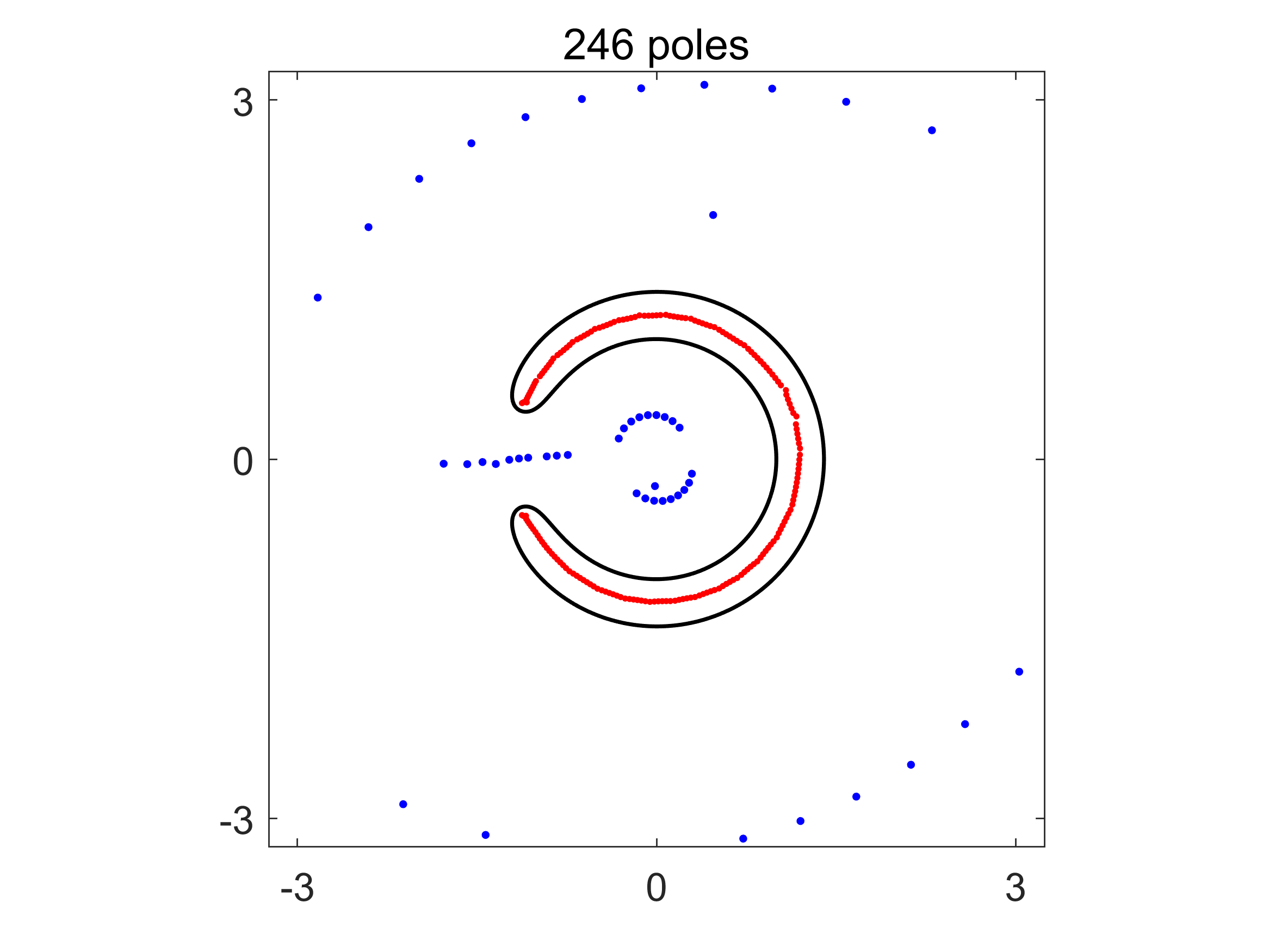}%
\includegraphics[trim=43 0 43 10, clip,scale=0.48]{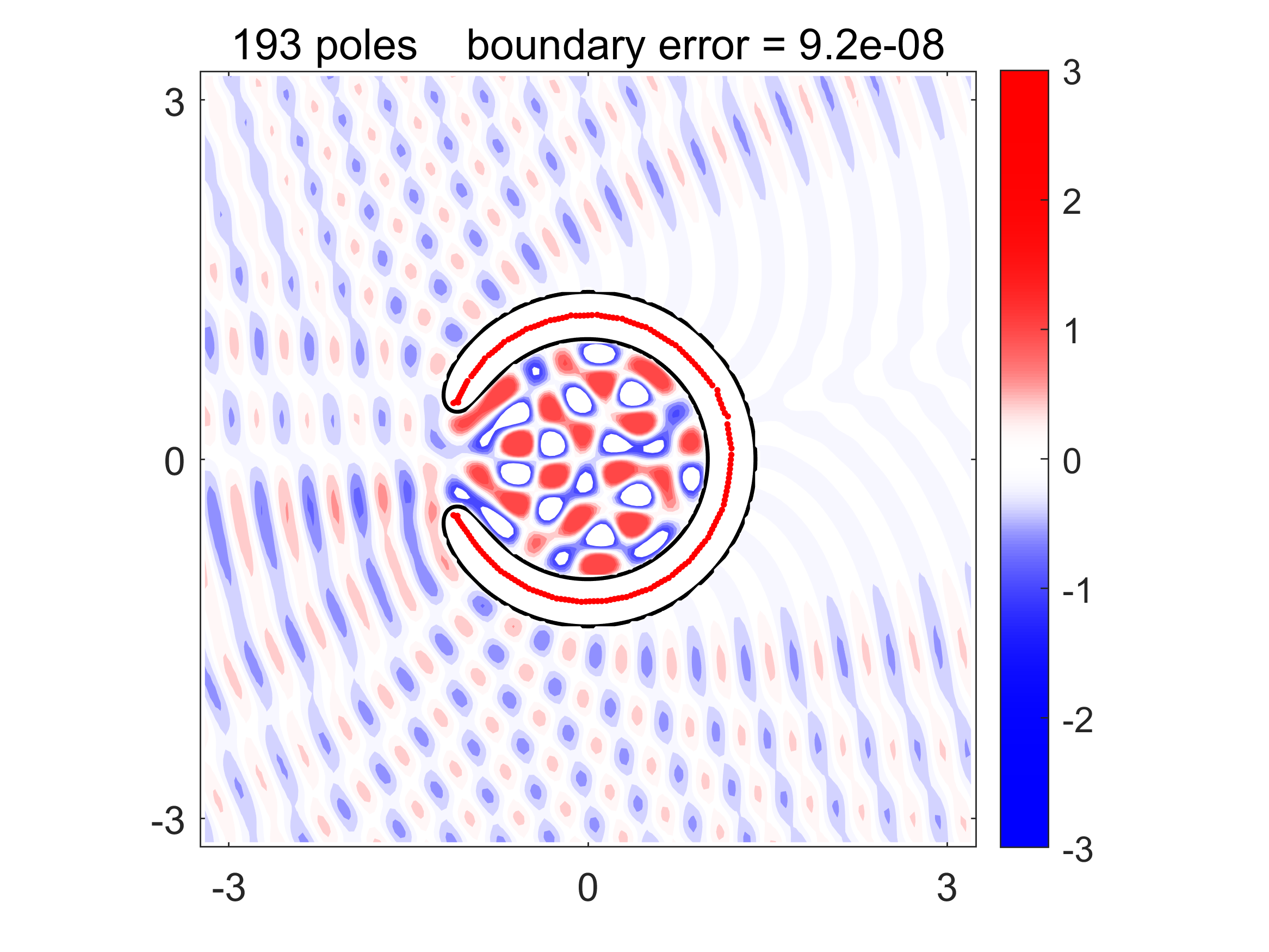}
\end{center}
\vspace*{-10pt}
\caption{\label{corral}Another AAA-Helmholtz scattering computation
involving a corral.  The wave number has been picked to be
close to a resonance for this shape,
which is why the amplitudes are bigger inside the corral
than outside.  Compare Figure \ref{res-disk}.}
\end{figure}

{\small
\begin{verbatim}
      t = 2*pi*(1:600)'/600;
      r = 1 + sin(4*t)/4;
      Z = r.*exp(1i*t);
      K = 20; h = @(z) -exp(1i*real(K*z));
      [r,pol] = aaa(real(h(Z)),Z,'mmax',300);
      ii = inpoly(pol,Z); p = pol(ii).';
\end{verbatim}
\par}

\noindent To complete the calculation, we employ $H_0^{}$
and $H_1^{}$ as indicated in (\ref{hankel0})--(\ref{hankel1}).
This example runs in about 1 s on our laptop, producing
results accurate to 9 digits as shown in the first row of Figure~\ref{helmdaisy}.
The second row shows a comparable experiment for a scatterer in
the shape of the domain of Figure~\ref{bigaaals}, but now the
accuracy is less.
Figure \ref{corral} shows a third experiment, now for a ``corral'' with
an incident wave at a near-resonant frequency in the sense of
section \ref{sec:res}.

{\small
\begin{verbatim}
      H0 = @(z) besselh(0,z); H1 = @(z) besselh(1,z);
      P = H0(K*abs(Z-p));
      Q1 = H1(K*abs(Z-p)); Q2 = (Z-p)./abs(Z-p);
      A = [P Q1.*real(Q2) Q1.*imag(Q2)];
      c = A\h(Z);
      u = @(z) reshape([H0(K*abs(z(:)-p)) ...
          H1(K*abs(z(:)-p)).*real(z(:)-p)./abs(z(:)-p) ...
          H1(K*abs(z(:)-p)).*imag(z(:)-p)./abs(z(:)-p)]*c,size(z));
\end{verbatim}
\par}

As mentioned above, there is no literature on the AAA-Helmholtz
method and no solid understanding of why it works, how well it
works, and how it can be made better.  Can rapid convergence under
certain assumptions be proved?  To get more digits of accuracy,
should one add low-degree series of Hankel functions, or use
singularities of doubled order as we did in section \ref{sec:aaals}
(involving $H_2^{}$ as well as $H_1^{}$), or use a AAA forcing function
with artificially increased frequency, or some combination of
these ideas, or something else?  How competitive can this method
ultimately be against the standard, highly-developed technologies
of integral equations?  Comprehensive research is needed
to provide answers to these and other questions.  For the moment, what we can
say is that the method certainly looks promising, and unpublished
experiments by Costa show that it can handle considerably higher frequencies than
those in the examples presented here.

\section{\label{sec:biharmonic}Biharmonic problems and Stokes flow}
The biharmonic equation is a fourth-order PDE, the square of the Laplace
equation:
\begin{equation}
\Delta^2 u(z) = 0, ~ z\in\Omega.
\label{biharm}
\end{equation}
Two boundary conditions are normally required, which, depending on
the application, would typically involve $u$ as well as its first or
second normal derivative, or perhaps both components of the gradient.
The best-known applications are in 2D
elasticity, where $u$ is the transverse deflection, and in 2D viscosity-dominated
fluid mechanics, known as {\em Stokes flow\/}, where
$u$ is the stream function \ccite{goodier}.

\begin{figure}
\begin{center}
\vskip 10pt
~~~~~~\includegraphics[trim = 0 0 0 0, clip, scale=.30]{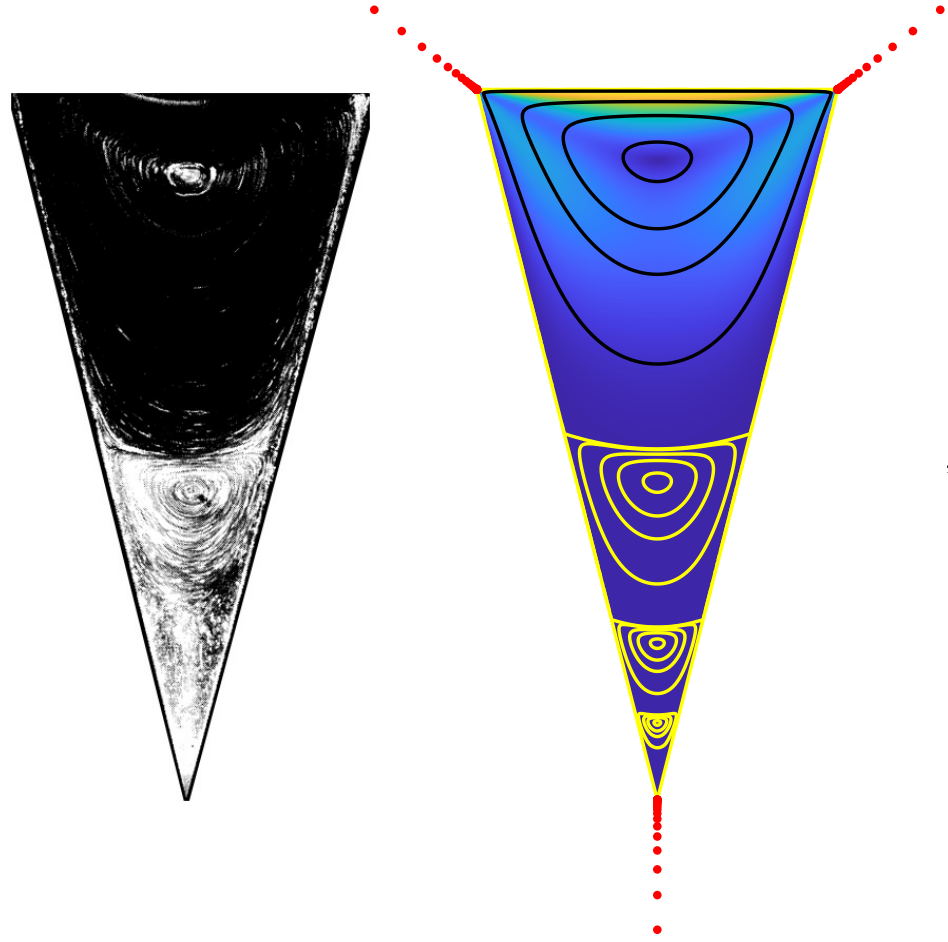}
\vspace*{-10pt}
\end{center}
\caption{\label{wedges}Figure from \ccite{stokes} showing Stokes flow in a triangular
driven cavity computed by the lightning Stokes method.
No-slip boundary conditions are applied everywhere, with the top boundary
moving at constant speed and the remainder of the boundary stationary.
The red dots are
poles placed near each corner according to the lightning a priori clustering formula.
The yellow contours mark the first three of an infinite
series of counter-rotating Moffatt eddies, each with amplitude about 830 times
less than the one above.  On the left,
a laboratory experiment in the same configuration from Taneda in 1979 reproduced
in the {\em Album of Fluid Motion\/} \ccite{vandyke}.}
\end{figure}

\begin{figure}
\begin{center}
\vskip 15pt
~~~~\includegraphics[trim = 130 35 125 25, clip, scale=.85]{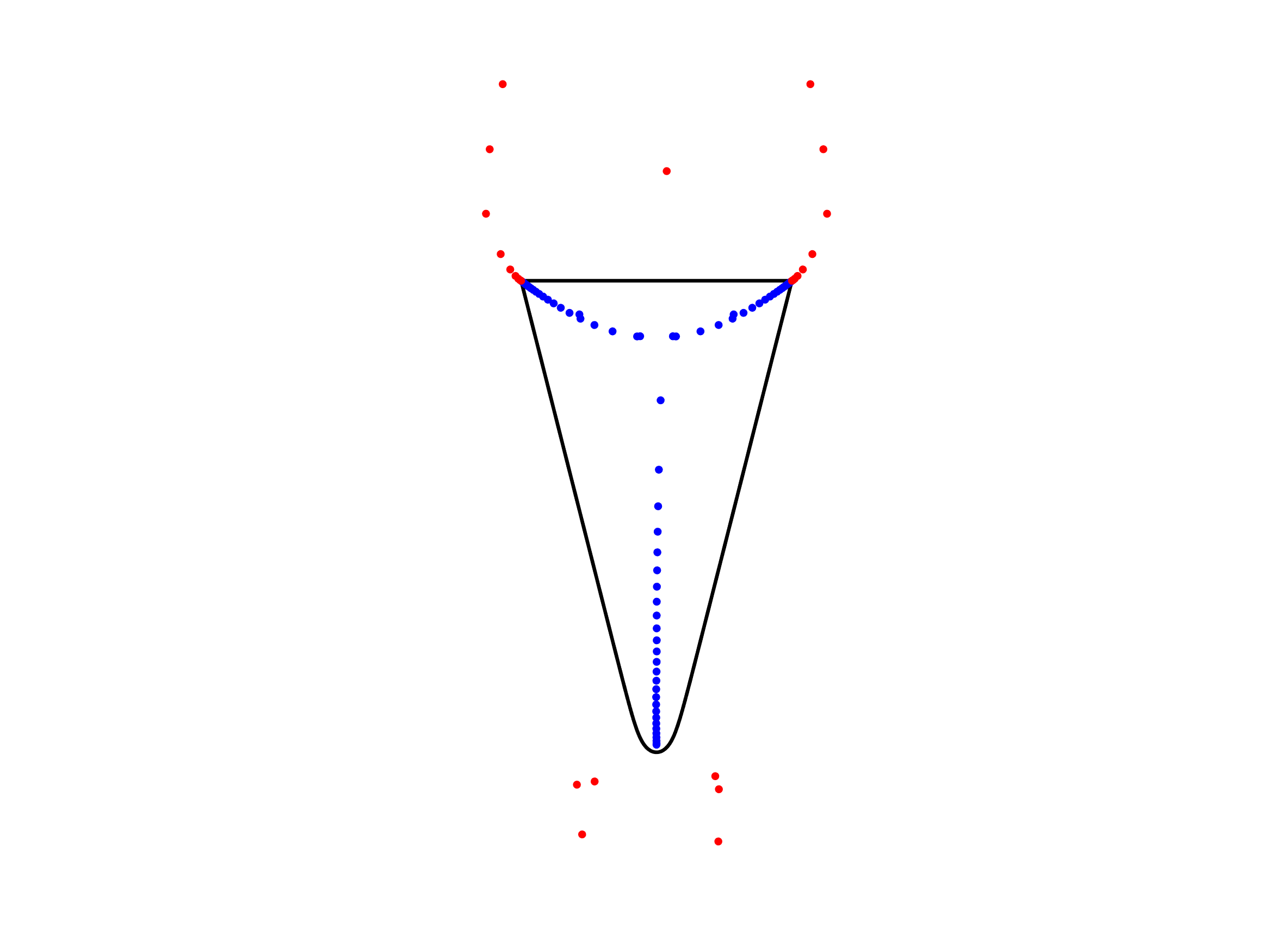}
\includegraphics[trim = 125 35 70 25, clip, scale=.85]{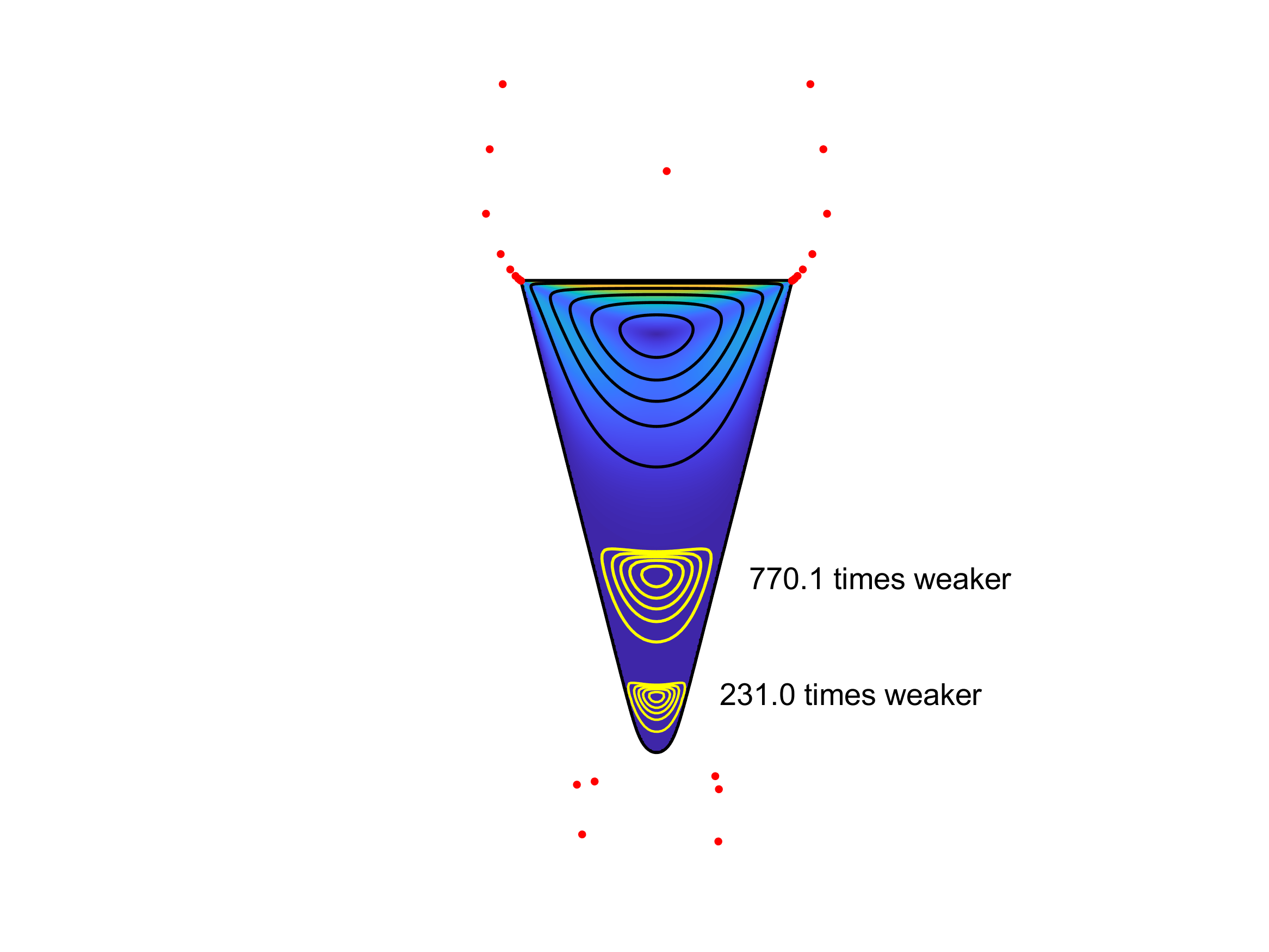}
\end{center}
\vspace*{-5pt}
\caption{\label{aaawedges}Stokes flow in another triangular
driven cavity with a curved bottom computed by the AAALS Stokes method,
following \ccite{waters}.  As in section \ref{sec:aaals}, the left image
shows all the AAA poles and the right one shows the poles outside
the problem domain together with the computed solution.
The poles placed by AAA are again clustered
near the corners, but no longer
lie on straight lines.  The poles near the bottom are
few, but still important for representing the solution to high
accuracy.  There are now two Moffatt eddies
instead of an infinite sequence.}
\end{figure}

Since 1898, it has been known that 2D biharmonic problems can be
reduced to pairs of Laplace problems.
In addition to the original method of Goursat, there are also alternative reductions due
to Alamansi and to Krakowski and Charnes;
see \ccite{stokes} for details.  In the latter paper, which was written before
AAALS was developed, Stokes problems are solved by
rational functions with the Goursat reduction
using the lightning method mentioned in section \ref{aaals}.
An example is reproduced in Figure~\ref{wedges}, showing resolution
to 10 digits of accuracy 
of a challenging driven-cavity flow by rational functions with
147 poles (17 of them off-scale in this figure).  The authors
report that this computation took half a second on a laptop for 10-digit
accuracy all the way up to the corners.

The Goursat reduction to Laplace problems can also be employed with the AAALS
method, and the leading figure here has been Yidan Xue, who, by combining
both lightning and AAA methods for placing poles, has solved a varied 
set of problems with corners and curved boundaries as well
as multiple connectivity.  Xue deals with both steady
and time-dependent flows, including in periodic geometries \ccite{waters,xue25};
his software is available at \ccite{xuesoftware}.
Figure~\ref{aaawedges} shows
an example of another driven cavity problem solved
with AAA.\ \ The bottom boundary is now curved, illustrating how
AAA places poles very differently near corners and near smooth segments.
Many more examples of Stokes flows computed with
lightning and AAA poles can be found in Xue's papers.
In addition to his published work, he has also explored further aspects of rational
function solutions of Stokes flows, including the use of continuum AAA
(see section \ref{sec:continuum}) and the treatment of unbounded domains.

There is another and possibly superior
AAA-related approach to Stokes problems that might be 
considered, which so far has apparently not yet been tried.
In the last section we solved Helmholtz problems by the Method of
Fundamental Solutions, using AAA to choose singularity locations but
then taking the actual singularities to be Hankel functions, not
poles of a rational function.  A similar approach could be
followed for the biharmonic equation, whose fundamental solutions
are known as {\em stokeslets\/} (monopoles) and {\em stresslets\/}
and {\em rotlets\/} (dipoles) \ccite{palsson,poz,poz2,zhao}.  The success of the AAA-Helmholtz method
suggests that a AAA-Stokes method based on this principle may also be successful,
and perhaps it may be simpler than the use of the Goursat representation.

So far as we know, rational approximations have not yet been used to solve
biharmonic problems arising in elasticity.

\section{\label{sec:hilbert}The Hilbert transform}
If $f_+ = u + iv$ is an analytic function in the upper half-plane
$\Cp$, the {\em Hilbert transform\/} is the map from $u$ on $\real$ to 
$v$ on $\real$.  Thus a computation of the Hilbert transform of a function
$u$ solves the Laplace Dirichlet problem in $\Cp$ with boundary
data $u$, since once you have $f_+ = u+iv$ on $\real$, 
you can extend it to $\Cp$ by a Cauchy integral or a rational approximation.
The Hilbert transform is also equivalent to 
what is known as the {\em Dirichlet-to-Neumann\/} map for the
upper half-plane.\footnote{Differentiating $v$ gives $dv/dt$ for
$t\in\real$, and by the Cauchy-Riemann conditions, this
is equal to the downward-pointing normal derivative $u_n^{}$ of
$u$ on $\real$.  Thus, once you've got $v$, you've got $u_n^{}$.}

Since the Hilbert transform solves a Laplace problem, one might expect
that to compute it with AAA, it will be necessary to use AAALS as in the last
four sections.  Indeed, this was the route followed
in \ccite{aaals}, where the effectiveness
of this method even for functions~$u$ and~$v$ with singularities was demonstrated.
However, AAALS is not actually necessary here.  Because of the real symmetry
in the Hilbert transform, it can be computed simply by additive 
Wiener-Hopf decomposition, which section \ref{sec:cauchy} showed is a matter of
AAA, not AAALS.

Before explaining and illustrating the method, we record how the Hilbert transform
can be derived
as a special case of the Cauchy transform, which was defined in 
section \ref{sec:cauchy} as the map
\begin{equation}
f(z) = {1\over 2\pi i} \int_\Gamma {h(t)\over t-z}\kern 2pt dt
\label{ctrans2}
\end{equation}
from functions $h$ on a contour $\Gamma$ to functions
$f$ in $\complex\cup \{\infty\} \backslash \Gamma$.
We defined $f_+$ and $f_-$ to be the branches of $f$ to the left and the
right of $\kern 1pt\Gamma$ near a point $t\in\Gamma$.  Then according to
the Sokhotski-Plemelj formulas (\ref{jump})--(\ref{pval}),
$h$ is equal to the jump in $f$ across the contour,
\begin{equation}
h(t) = f_+(t) - f_-(t),
\label{jump2}
\end{equation}
and if we evaluate $f$ {\em on\/} the contour we get
\begin{equation}
f(t) =  {\textstyle{1\over 2}} (\kern .7pt f_+(t)+f_-(t)),
\quad t\in\Gamma,
\label{pval2}
\end{equation}
where $f(t)$ is defined by
taking the integral (\ref{ctrans2}) in the principal value sense.
The Hilbert transform is the special case of these relationships in the
setting where $\Gamma=\real$ and $h$ is real with $\lim_{|\kern .5pt t\kern .5pt |\to\infty} h(t)=0$.
Conventionally (\ref{ctrans2}) is multiplied by $-2\kern .5pt i$,
so that (\ref{pval2}) becomes
\begin{equation}
v(s) = {1\over \pi} \kern 3pt \hbox{PV} \kern -2pt
\int_{-\infty}^\infty \kern 1pt {u(t)\over s-t} \kern 2pt dt,
\quad s\in\real 
\label{htrans}
\end{equation}
as a map from $u$ to $v$.  (The minus sign in the factor
$-2\kern .5pt i$ has been incorporated by switching
the denominator from $t-s$ to $s-t$.)  
We can explain as follows why (\ref{ctrans2}) is multiplied by $-2\kern .5pt i$.
If $f_+(z) = u(z) + iv(z)$
is analytic in the closed upper half-plane $\Cpp$ and
we define $f_-(z) = 0$ in the closed lower half-plane $\Cmp$, then by (\ref{jump2}),
$f_+$ is the Cauchy integral (\ref{ctrans2}) in $\Cp$ of $h(t) = 
f_+(t)$ on $\real$, whereas by (\ref{pval2}), the principal value on $\real$ is
$f(t) = f_+(t)/2$.  Thus the principal value Cauchy integral maps $u(t) + iv(t)$ to
${\textstyle{1\over 2}} u(t) + {\textstyle{i\over 2}}v(t)$ on $\real$.  Therefore
$-2\kern .5pt i$ times this integral maps $u(t) + iv(t)$ to
$-iu(t) + v(t)$ on $\real$.  Taking real parts, we see that
$-2\kern .5pt i$ times
the integral maps $u$ to $v$, and this is (\ref{htrans}).

Since the harmonic conjugate of $\cos(k\kern .3pt t)$ is $\sin(k\kern .3pt t)$ and the
harmonic conjugate of $\sin(k\kern .3pt t)$ is $-\cos(k\kern .3pt t)$, it follows
that the Hilbert transform acts on a Fourier transform or
series by mapping the Fourier
coefficient $\hat u(k)$ to $-\hbox{sign}(k) \kern .5pt i \hat u(k)$.
(From here, or by considering harmonic conjugates directly, one sees
that the square of the Hilbert transform is the negative of the
identity.)
Many numerical methods for the Hilbert transform have been 
based on Fourier series.
For example, this is the approach used by the {\tt hilbert} command
in the MATLAB Signal Processing Toolbox.  A different method
based on expansion in rational functions with pre-assigned poles was proposed
in \ccite{weidemanhilb}.
For some other approaches see \ccite{bilato} and \ccite{zhou}.

But let us return to AAA.\ \ Given a real function
$u(t)$ for $t\in\real$, section \ref{sec:cauchy} showed that
AAA approximation
can provide a rational Wiener-Hopf splitting
\begin{equation}
u(t) \approx r(t) = r_+(t) + r_-(t).
\label{whsplitting}
\end{equation}
(We assume that $u$ goes to $0$ at $\infty$ and likewise
$r_+$ and $r_-$.)
By the Schwarz reflection principle, $r(\bar z) = \overline{r(z)}$ for
all $z\in\complex$, implying that the poles and residues of $r_+$ are
the conjugates of those of $r_-$.
It follows that $r_+$ and $r_-$ have equal real parts for $t\in\real$ and
opposite imaginary parts.  Therefore
\begin{equation}
v(s) \approx 2 \kern 1pt \Im ( r_+(s) ).
\label{htresult}
\end{equation}
Together, (\ref{whsplitting}) and (\ref{htresult})
constitute a numerical method for computation of the Hilbert transform.

To illustrate, we return to the function
$u(t) = e^{-|\kern .5pt t\kern .5pt |}$ as in Figure~\ref{wienerhopffig}, an example that was
treated in \ccite{weidemanhilb}.
The exact solution is known analytically:
\begin{equation}
u(t) = e^{-|\kern .5pt t \kern .5pt |}, \quad
v(s) = {\hbox{sign}(s)\over\pi} \bigl[ e^{|\kern .5pt s\kern .5pt |} 
E_1(|s|) + e^{-|\kern .5pt s\kern .5pt |} Ei(|\kern .5pt s\kern .5pt |)\bigr],
\label{htexample}
\end{equation}
where $E_1$ and $Ei$ are the exponential integrals computed in MATLAB by 
{\tt expint} and {\tt ei}.  Both $u$ and $v$ have branch point singularities at the
origin, a context in which the rational representations of AAA show their power.
For a computation with tolerance $10^{-10}$ based on (\ref{htresult}), we compute {\tt f1}
and {\tt f2} exactly as on p.~\pageref{whpage} and then just execute

{\small
\begin{verbatim}
      v = @(s) 2*imag(f2(s));
\end{verbatim}
\par}

\begin{figure}
\begin{center}
\includegraphics[trim = 0 130 0 0, clip, scale=.85]{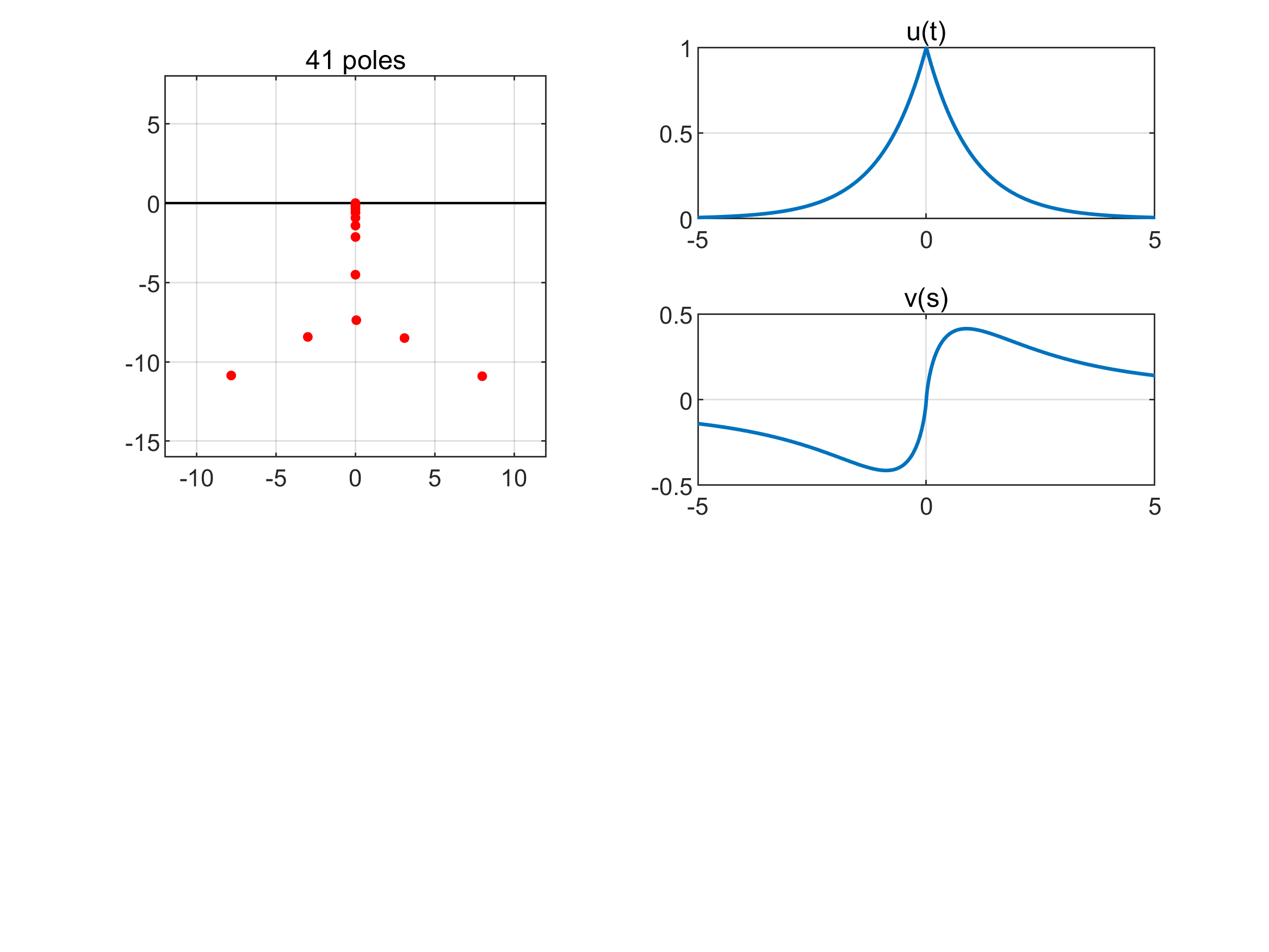}
\end{center}
\vspace*{-15pt}
\caption{\label{hilbfig}Numerical computation of the Hilbert transform
(\ref{htexample}) of $u(t) = e^{-|\kern .5pt t\kern .5pt |}$ 
by AAA (not AAALS). 
The same Wiener-Hopf splitting is used as
in Figure~\ref{wienerhopffig}, and the solution $v(s)$ is then
given by (\ref{htresult}).  A solution matching the exact result
with accuracy $10^{-9}$ is found in $0.06$ s on a laptop.}
\end{figure}

\noindent The result, shown in Figure~\ref{hilbfig}, differs from
the exact solution (\ref{htexample})
by a maximum of $9.6\cdot 10^{-10}$ for $t\in\real$.

Let us comment further on AAA vs.\ AAALS and on the
Dirichlet-to-Neumann map on $\real$
vs.\ on a general closed contour $\Gamma$.
(Here $\real$ stands for any line or circle, since these are
equivalent by a M\"obius transformation. 
The Hilbert transform on the unit circle is discussed
at length in \ccite{henrici3}.)
As we have just shown, the symmetry of $\real$ enables the
Dirichlet-to-Neumann map to be computed directly by approximation, without
the need for an additional least-squares step.  In the paper
\ccite{aaals}, this special symmetry is mirrored in the
fact that Theorems~6.1 and~6.2
only apply on lines and circles.  For a more
general domain, AAALS is needed.\ \ One might observe that a more
general domain can be reduced to a half-plane or a disk by a conformal
map, so why is the Hilbert transform setting not equivalent
to the general case?  The answer
is that while this is true theoretically, it is not relevant computationally
since the conformal map itself amounts to a Laplace problem that
must be solved somehow, as discussed in section \ref{sec:conformal}.

\section{\label{sec:lawson}AAA-Lawson for best approximation}
This completes our survey of applications of AAA approximation.
Here and in the next four sections, we will return to fundamentals
and briefly describe five extensions of the standard AAA algorithm
as presented in section~\ref{sec:AAA}.

AAA reliably produces near-best approximations of a given degree
$n$, but it cannot produce actual best approximations---``best'' in
the sense of minimising the $\infty$ norm of the error $\|r(Z)-F\|$
in the notation of section~\ref{sec:AAA}.  The reason is that AAA
represents $r$ by the barycentric formula (\ref{aaa-baryrep}),
which forces $r$ to interpolate the data vector $F$ at the support
points $\{\tk\}$ (assuming the weights $\bk$ are nonzero).  This
constraint precludes the representation of a best approximation,
since in general a best approximation on a discrete set need not
interpolate the data anywhere.\footnote{In theory, best approximations on
continua exist but need not be unique, and on discrete sets or other
sets with isolated points, they need not exist~\ccite{istace}.
(The case of real approximation on a real interval is special: here we have
both existence and uniqueness, and the Remez algorithm for computations.)
However, nonuniqueness
and nonexistence seem rarely to be issues in practice, at least for
an algorithm like AAA-Lawson that does not attempt to enforce a characterisation
of optimality.}

From the beginning of work on AAA, therefore, the idea has been
pursued of generalising (\ref{aaa-baryrep}) to 
an ``alpha-beta barycentric representation,''
\begin{equation}
r(z) = {N(z)\over D(z)} = \sum_{k=0}^n {\ak\over z-\tk} \left/
\sum_{k=0}^n {\bk\over z-\tk} \right..
\label{alphabeta}
\end{equation}
This differs from (\ref{aaa-baryrep}) in having arbitrary constants
$\{\ak\}$ rather than interpolation data $\{\fk\kern .5pt \bk\}$
in the numerator.  For any choice of support points $\{\tk\}$,
any rational function of degree $n$ can be represented in the form
(\ref{alphabeta}), and thus best approximations should be in reach
if an algorithm can be found for computing $\{\ak\}$ and $\{\bk\}$.
If $r^*$ is a best approximation to $F$ on $Z$, and the points
$\{\tk\}$ are fixed arbitrarily, then there will be a choice
of $\{\ak\}$ and $\{\bk\}$ for
which $r$ in (\ref{alphabeta}) is equal to $r^*$; the choice is
unique up to a normalisation constant if\/ $r^*$ is of degree $n$ but not $n-1$
\cite[Thm.~3.1]{aaaL}.  Arbitrary choices of $\{\tk\}$ will be numerically awful,
but the choices made by AAA iteration are generally excellent.

The algorithm that has been explored from the beginning,
{\em AAA-Lawson,} first runs AAA as usual and then switches to
the formulation (\ref{alphabeta}) with the support points $\{\tk\}$
chosen by AAA, aiming to find $\{\ak\}$ and $\{\bk\}$ by a process
known as {\em iteratively reweighted least-squares\/}.
This algorithm was introduced in \ccite{aaa} and \ccite{minimax} and became
fully practical with various improvements introduced in \ccite{aaaL}.
As mentioned in (\ref{zolo-lawson}),
the idea is to solve the weighted least-squares problem
\begin{equation}
\min_{\gamma,\, \|\gamma\|_2^{}=1} \; 
\sum_{j=1}^m {}' \; w_j \left(\kern .7pt f_j^{} \sum_{k=0}^n {\beta_k\over z_j^{}-\tk} -
\sum_{k=0}^n {\alpha_k^{}\over z_j^{}-\tk}\right)^2 ,
\label{law-ls}
\end{equation}
where $\gamma$ denotes the concatenation of the weight vectors $\alpha$ and
$\beta$,
and then adjust the weights iteratively by the formula
\begin{equation}
w_j^{\rm (new)} = w_j |e_j|.
\label{law-lawson}
\end{equation}
The prime on the summation sign in (\ref{law-ls}) reflects special treatment when
$z_j^{}$ is one of the support points $\tk$; see \ccite{aaaL} for details.
For linear approximation problems, the Lawson iteration was
introduced by Lawson in 1961 \ccite{lawson}, and convergence
was proved by Rice (real) and Ellacott and Williams (complex); see
\ccite{aaaL} for references.  The AAA-Lawson iteration is nonlinear, however,
and there is no theory proving that it will converge.  Indeed,
it does not always converge.

\begin{figure}
\vskip 7pt
\begin{center}
\includegraphics[trim=0 150 90 10, clip, scale=1]{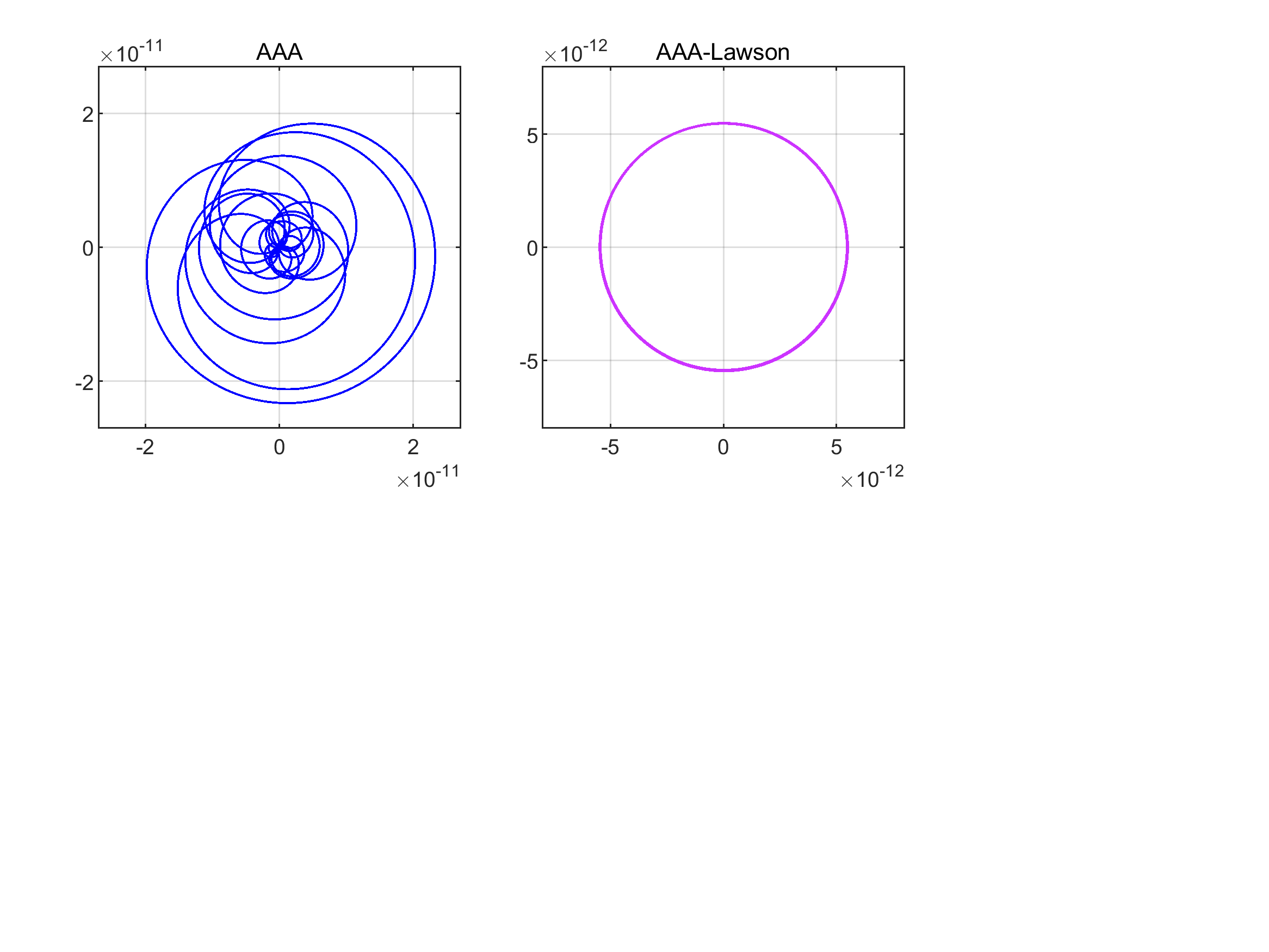}
\end{center}
\vskip -10pt
\caption{\label{gammalawfig}Repetition of Figure \ref{fig3.3}, but
with the degree $n$ of the
AAA approximation of $\kern 1pt \Gamma(z)$ on the circle $|z|=1.5$ reduced from $12$
to $11$ to stay well above the level of rounding errors.  
On the left, the error curve for AAA near-best approximation passes through
the origin $n+1 = 12$ times, at each of the support points.  On the right,
the error curve for AAA-Lawson best approximation is a
near-circle of winding number $2\kern .3pt n-3 = 19$.  Note the different scales on the two sides,
reflecting the pattern that best approximations typically have errors smaller than
those of AAA approximations by a factor less than $10$.}
\end{figure}

To illustrate AAA-Lawson approximation, we begin with the example of
Figure~\ref{fig3.3}, approximation of the gamma function $\Gamma(z)$
on the circle $|z|=1.5$.  Both sides of Figure \ref{gammalawfig}
plot error curves defined as the image of the circle under $e(z) =
\Gamma(z) - r(z)$, with AAA near-best approximation on the left and
AAA-Lawson best approximation on the right.  In Figure \ref{fig3.3}
the error curve was erratic because it was close to the level of
rounding error, but here, the degree has been reduced from $12$
to $11$ and we get a clean result.  The computation of the curve
on the left is effected by the commands

{\small
\begin{verbatim}
      Z = 1.5*exp(2i*pi*(1:100)'/100);
      r = aaa(g(Z),Z,'degree',11,'lawson',0);
      E = g(Z) - r(Z);
\end{verbatim}
\par}

\noindent with {\tt g} being the complex gamma function routine from
\ccite{godfrey}.  (Actually a finer grid than just 100 points
is used for plotting.)  For the curve on the right, the \verb|'lawson',0|
flag has been removed so that the iteration (\ref{law-ls})--(\ref{law-lawson}) is
applied for the default number of 20 steps.  The computations on
a laptop take about 4 and 10 ms, respectively.  The error curve is
a near-circle of winding number $2\kern .3pt n-3 = 19$, as we can verify with
the commands

{\small
\begin{verbatim}
      theta = unwrap(angle(E));
      windingnumber = (theta(end)-theta(1))/(2*pi)
\end{verbatim}
\par}

\noindent
If $\kern 1pt \Gamma(z)$ were analytic inside the circle, the winding
number would be $2\kern .3pt n+1$, but the two poles inside reduce it by $4$,
as explained in the footnote on p.~\pageref{footpage}.

\begin{figure}
\vskip 7pt
\begin{center}
\includegraphics[trim=0 150 90 10, clip, scale=1]{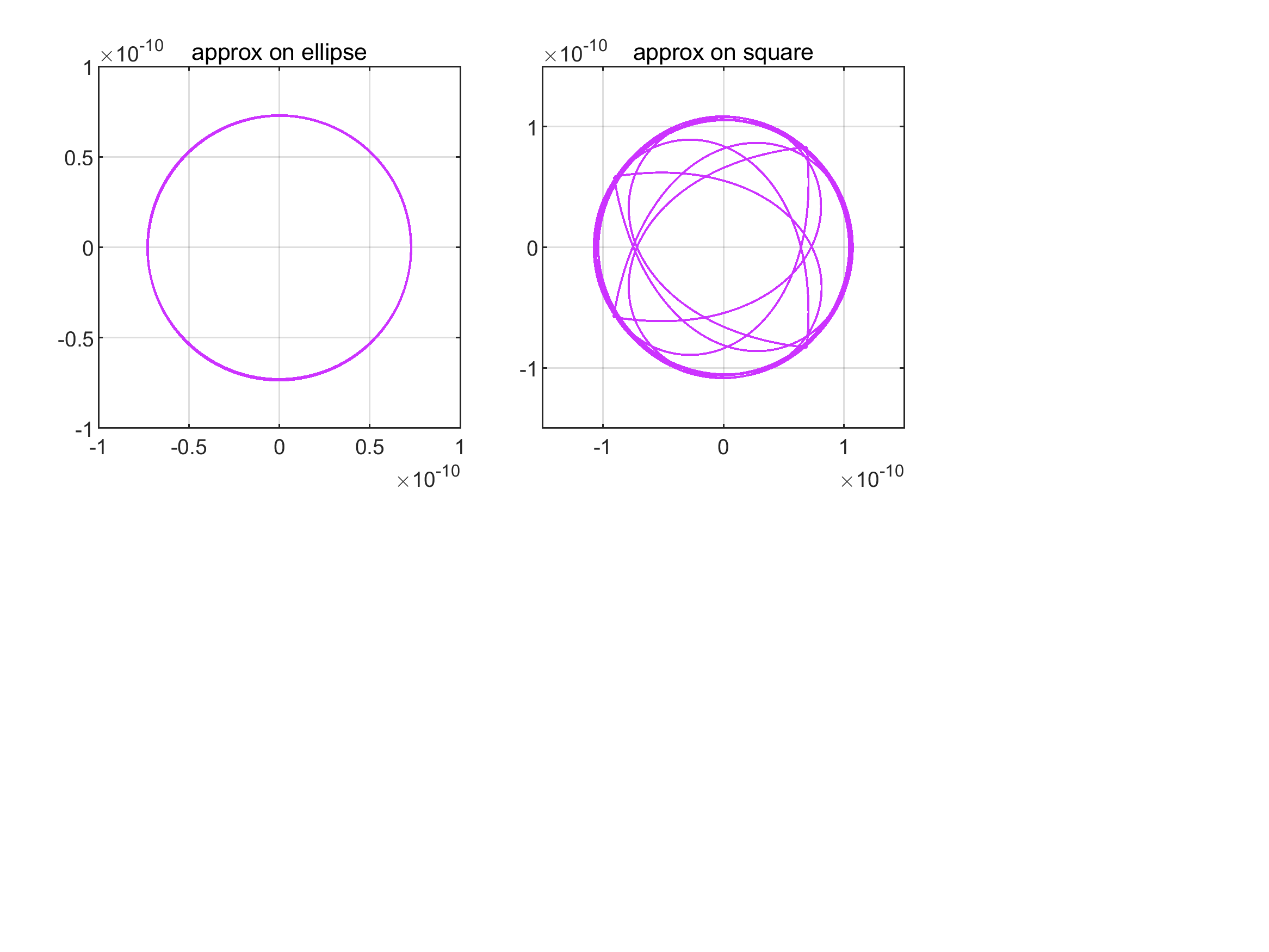}
\end{center}
\vskip -10pt
\caption{\label{gamma2fig}Error curves for
two more degree $11$ AAA-Lawson approximations of $\kern 1pt \Gamma(z)$.
On the left, the domain
is the ellipse with extreme points $\Re z = \pm 1.5$ and
$\Im z = \pm 2$ and the error curve is very nearly circular.
On the right, the domain is the square $-1.5\le \Re z, \Im z \le 1.5$, and
the four corners are visible in the error curve.  The
windings numbers are 18 on the left and 19 on the right.}
\end{figure}

The near-circularity phenomenon exemplified in the right panel
of Figure \ref{gammalawfig} was explained in \ccite{nearcirc},
and it applies not just for approximation on a disk but also on
other smooth approximation domains.  On domains with corners,
the error function $e(z)$ will have to preserve the corner angles
since it is locally a conformal map, but the rest of the error curve
will tend to come close to a circle.  To illustrate these effects,
Figure~\ref{gamma2fig} repeats the degree $11$ best approximation
of Figure~\ref{gammalawfig} but for approximation on a $3\times 4$
ellipse and on a square of side length $3$.  The winding number is
19 again on the square, but just 18 on the ellipse for reasons we
have not investigated.

Before AAA-Lawson, no method was available for computing complex
rational best approximations that was effective enough to play a role
in scientific computing.
With AAA-Lawson, most problems are now easy, and many more examples
can be found in \ccite{aaaL}, including cases with disconnected
approximation domains.  AAA-Lawson is less important for
real approximations on intervals, on the other hand,
since the Remez algorithm has been around for many years,
with robust barycentric implementations in \ccite{ionita} and
Chebfun {\tt minimax} \ccite{minimax}.
Still, it may sometimes be
advantageous for real problems, and an example of a real
AAA-Lawson computation
appears in Figure \ref{figcontex}.

As described in \ccite{aaaL},
AAA-Lawson converges often but not always, and with the
\verb|'damping'| variation described in section \ref{sec:zolo}, its
success rate improves but is still not $100\%$.  As a particularly
challenging problem of a fundamentally discrete rather than
continuous nature, what if we try to fit random complex data in
30 roots of unity by an optimal rational function of degree $10$?
Six trials generated by the code
\begin{figure}
\begin{center}
\includegraphics[trim=0 65 0 0, clip, scale=.85]{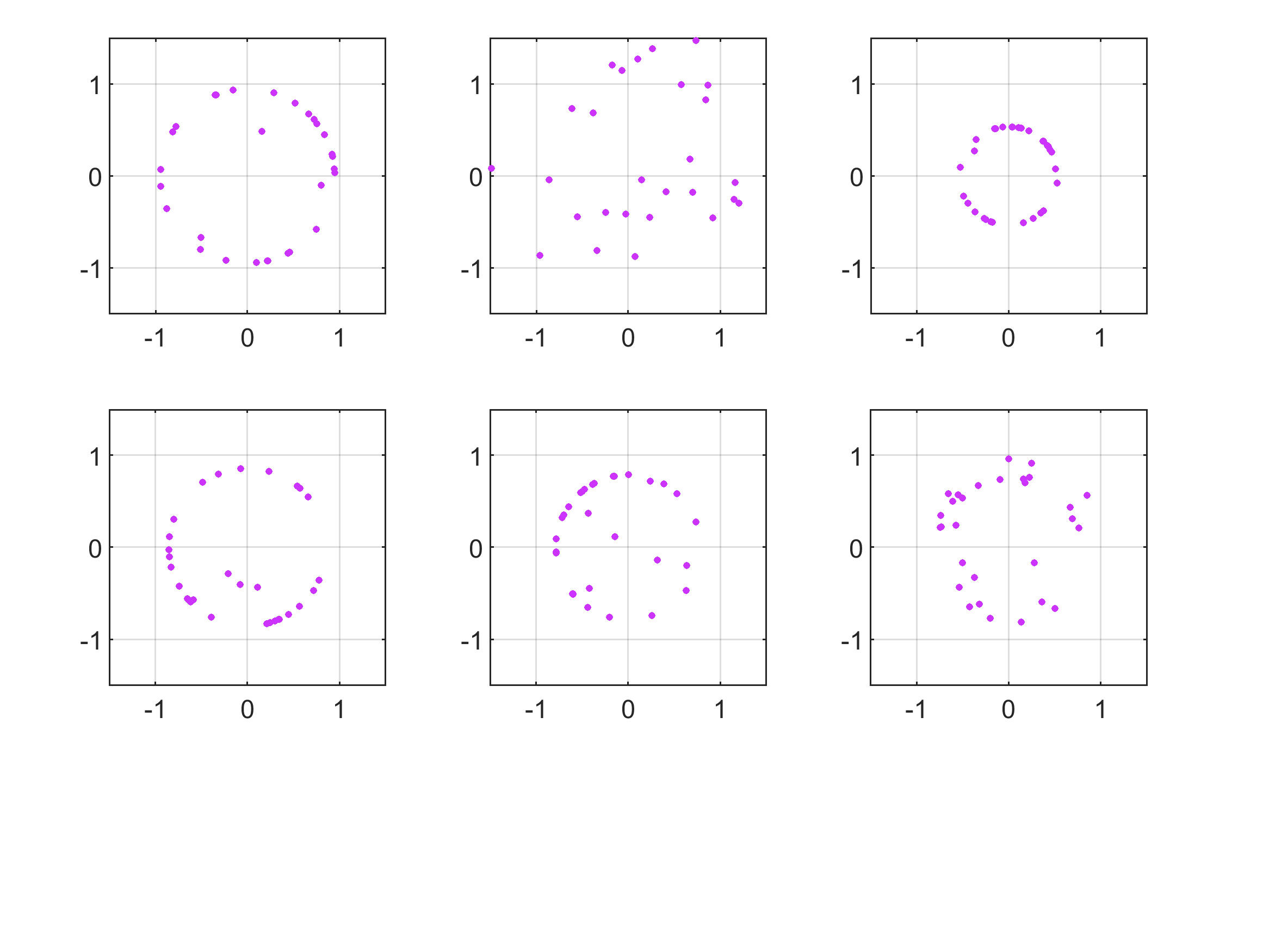}
\end{center}
\vskip -13pt
\caption{\label{discretefig} The error ``curve'' (discrete dots)
for six trials of degree 10 AAA-Lawson best
approximation of complex random data in 30 roots of unity.  The first, third,
fourth and fifth images appear to have achieved optimality, whereas the second one
is not converging.  The sixth image is close but is also not converging.
More research is needed to find a AAA-Lawson algorithm that always succeeds.}
\end{figure}

{\small
\begin{verbatim}
      rng(0)
      Z = exp(2i*pi*(1:30)'/30);
      for k = 1:6
         F = randn(30,1) + 1i*randn(30,1);
         r = aaa(F,Z,'degree',10,'lawson',500,'damping',.5);
         subplot(2,3,k), plot(r(Z)-F,'.')
      end
\end{verbatim}
\par}

\noindent lead to the images shown in Figure \ref{discretefig}.
Four of the six trials appear successful, but the failure rate rapidly
grows if the problem size is increased.

If a variant of AAA-Lawson with guaranteed convergence could be
found, would we recommend that best approximation should replace
near-best as the default?  Of course it would depend on the details of the
new algorithm, but our guess is that we would not, for two reasons.
One is that even an algorithm with provable convergence in exact
arithmetic may prove fragile in floating point for
approximations close to machine precision.  The other is that the
cost of AAA-Lawson is considerable:
an increase in computing time over standard AAA typically by at
least a factor of 2 (see Figure 1.1 of \ccite{aaaL})
and sometimes by much more, when many Lawson steps are needed.
Note that inevitably, a best-approximation algorithm like
AAA-Lawson must sample the function at more points than a near-best
algorithm like AAA, since it depends on accurate information near extrema.

\section{\label{sec:continuum}Continuum AAA}

Some approximation problems are intrinsically discrete: one may
genuinely have $m$ points $\{z_k^{}\}$ with $m$ associated data
values $\{\fk\}$ to approximate.  Most problems, however, involve
continuous data on continuous domains.  Throughout this survey,
and throughout most of the work by us and our collaborators,
we have addressed such situations by discretising the domains a
priori without worrying too much about the discretisation details.
Since the computational effort only grows linearly with the
number of points, not much is lost if, say, 1000 points are used
in a setting that might have been treated successfully with 300.
Admittedly, such constant factors can get aggravating in problems
with expensive function evaluations.

There is a context where discretisation challenges become more
extreme, however, and that is when there are singularities
on or near the boundary of the approximation region.
If poles are going to cluster near a point
$z_0$, then it is crucial that sample points cluster there too,
or one will end up with a discrete approximation that fails to
fit the continuous problem.  We are experienced at choosing grids
with exponential clustering near known singularity locations, as in
the computation for Figure \ref{bigaaals} for example, but this is
hardly an ideal way to proceed in cases with several singularities or
singularities at unknown locations.
Here the
cost of discretising ``by hand'' may be much worse than a constant
factor, and with this in mind,
a continuum AAA algorithm was introduced in \ccite{continuum}
and has been subsequently explored and extended by Driscoll
\ccite{driscolljuliacon}.
We have already illustrated some continuum AAA results in
Figure~\ref{sixfigs}.

The idea of continuum AAA is to work with
a discrete sample grid that is extended adaptively as the AAA
iteration proceeds.  Throughout the iteration, there are at least
three equispaced sample points between each pair of support points
(somewhat more in the early steps).  This means that the Loewner
matrices have only about three times as many rows as columns, and
the total number of function evaluations (including at the support
points) is about four times the degree.  At each greedy step, a new
support point is chosen where the error is largest, and then new
sample points are added on either side of the new sample point to
maintain the 3-to-1 ratio.  This leads to an algorithm that with good
though not perfect reliability distributes sample and support points
where they are needed, including clustering them near singularities.

A helpful corollary of approximation on a continuum is that one must
specify that continuum.  This means that when a AAA approximation is
computed, one can compute its poles and find out if any of them lie
where they shouldn't.  Such ``bad poles'' render an approximation
useless,\footnote{Useless in principle, at any rate, making the
approximation error $\infty$ in the $\infty$-norm.  In practice, bad
poles often have residues near machine epsilon---numerical Froissart
doublets---in which case their
effects are so localised that they may be harmless after all.} and
the continuum AAA algorithm simply proceeds to the next greedy step,
where with luck there will be no bad poles.  Usually---though not
always---one eventually reaches an approximation to the required
tolerance that is known to be free of bad poles.

\begin{figure}
\vskip 2pt
\begin{center}
\includegraphics[trim=0 0 0 0, clip,scale=0.8]{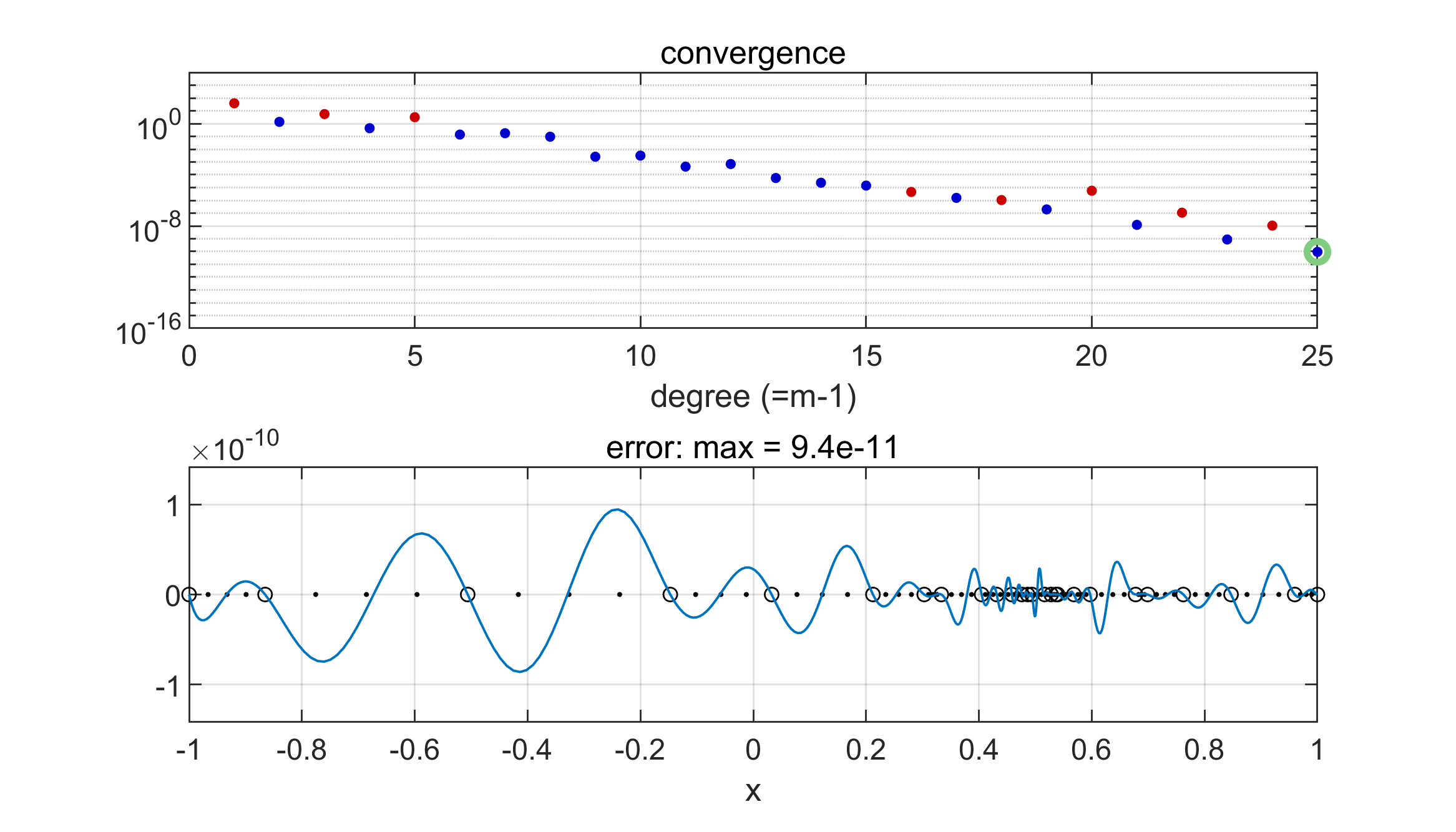}
\includegraphics[trim=0 0 0 123, clip,scale=0.8]{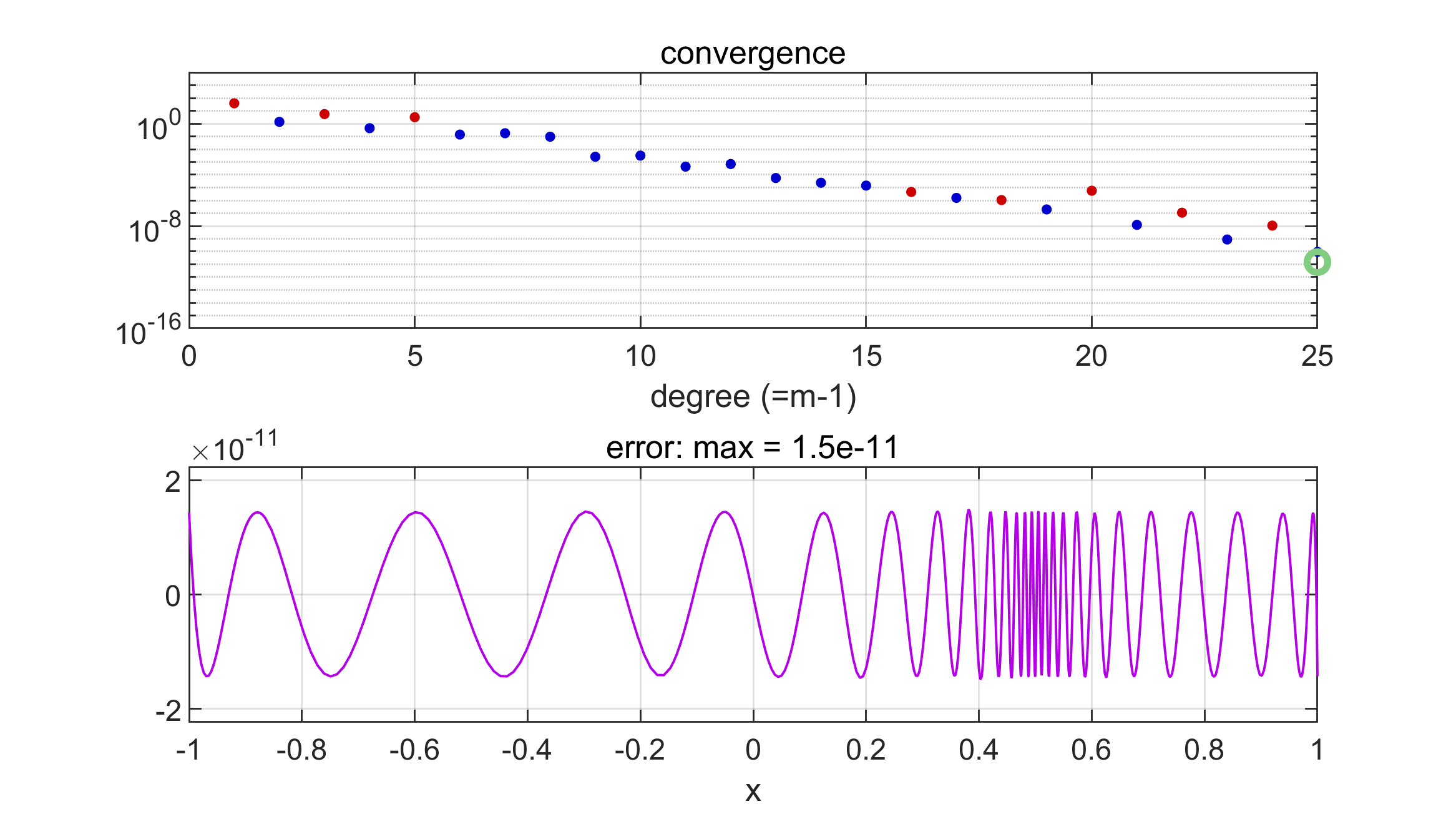}
\end{center}
\vspace*{-13pt}
\caption{\label{figcontex}Continuum AAA near-best rational
approximation of $f(x)=\log(1+1000(x-\half)^2)$ on $[-1,1]$.  The convergence
curve in the first image shows red dots at degrees with ``bad poles'' in $[-1,1]$ and blue
dots at degrees without such poles.  The second plot shows the final degree $25$ AAA error curve,
with adaptively determined support points marked by circles and other sample points by dots.  The final plot
shows the error curve for a best approximation of the same degree computed
by continuum AAA-Lawson.}
\end{figure}

The codes introduced in \ccite{continuum} are {\tt aaax} for
approximation on $[-1,1]$, {\tt aaaz} for approximation on the unit
circle or the unit disk and {\tt aaai} for approximation on the
imaginary axis or the right half-plane.  (The difference between
the circle and the disk lies in which poles are rejected as bad
ones; likewise the difference between the imaginary axis and the
half-plane.)  Many examples can be found in \ccite{continuum} and
\ccite{driscollzhou,driscolljuliacon}, and here we will just show one to give the idea.
To approximate $\log(1+1000 (x-\half)^2)$ by a degree $25$ rational
function on $[-1,1]$, one can execute

{\small
\begin{verbatim}
      aaax(@(x) log(1+1000*(x-.5).^2),25);
\end{verbatim}
\par}

\noindent This produces the upper two plots of
Figure~\ref{figcontex}, showing a near-best approximation with
error about $10^{-10}$ computed with 26 support points and just
75 additional sample points.  For a true best approximation of the
same degree, one can specify that 60 AAA-Lawson iterations should
be taken with the command

{\small
\begin{verbatim}
      aaax(@(x) log((1+1000*(x-.5).^2)),25,60); 
\end{verbatim}
\par}

\noindent This gives the third plot, with error six times smaller,
showing the characteristic equioscillation between $2\kern .3pt n+2 = 52$
extrema.  The error curve shows that much of the action is near
$x=0.5$.\footnote{A reminder of the point of rational approximations:
polynomials can't do this.} To get good equioscillation over the
whole continuum, {\tt aaax} by default works with a grid five
times finer.

Although this article has almost exclusively dealt with discrete AAA
approximation, this is partly due to our limited experience.
It is likely that many of the applications presented here could
have been carried out successfully in continuum mode, often with
a computational speedup.
Good candidates would be Figures
\ref{schwarz-squiggles3}, \ref{fig:mor}, \ref{bigaaals} and \ref{helmdaisy}.
We expect that continuum versions of AAA
will be increasingly prominent in the years ahead.

\section{\label{sec:trig}Periodic AAA}

Often a geometry is periodic, and one would like to apply
an appropriate analogue of rational approximation.   With
$2\pi$-periodicity, for example, one would like to replace
algebraic rational functions $p(z)/q(z)$ by trigonometric rational
functions $p(\theta)/q(\theta)$, where $p$ and $q$ are trigonometric
polynomials.  The approximation is then not rational but meromorphic.

There are two approaches to such problems.  One can
derive trigonometric variants of the usual algebraic formulas and
algorithms and work explicitly in the variable $\theta$.
Or one can make the change of variables $z = e^{i\theta}$ and recover
the usual algebraic case on the unit circle.  The choice between
these two approaches
long predates AAA.\ \ For example, in section~\ref{sec:extending} we
mentioned the idea of Fourier-Pad\'e approximation \ccite{bgm} and
commented that some authors call this simply Pad\'e approximation,
probably because what they do in practice is change variables and
then use standard Pad\'e.

To date, there have been three main
publications concerned with periodic AAA approximation.  One of these,
a paper by Xue mentioned in section~\ref{sec:biharmonic},
follows the transplantation approach and calls {\tt aaa.m} \ccite{xue25}.
The other two, working in the $\theta$ variable,
are by Wilber et al.\ \ccite{wdt}, who wrote a code called PronyAAA,
and Peter Baddoo\footnote{Peter Baddoo grew up in Britain, with
roots also in Ghana, and studied and worked
at Oxford, Cambridge, Imperial College and MIT.\ \  He was one
of the brightest and most appealing people we have ever encountered,
and it was devastating news when we learned in 2023 that he had died
as the result of a heart abnormality after a basketball game at
MIT.\ \ He was 29 years old.} \ccite{baddoo}, who provided the code {\tt aaatrig} in
Chebfun \ccite{baddoo}.
These projects required a generalisation of the
algebraic barycentric formula (\ref{aaa-baryrep}) to a trigonometric
barycentric formula, and on this topic there is a literature going
back to Salzer, Henrici and Berrut (\cite{salzer,henrici,berrut}).
We will not give
details, which are complicated by a distinction that must be observed
between the cases of an odd or an even number of support points.
This distinction also arises in the more elementary context of
computation with trigonometric
polynomials, as in Fourier spectral methods for \ODEs\ and \PDEs\
and the ``trig'' mode in Chebfun \ccite{wjmt,austintrig}.

\begin{figure}
\vskip 2pt
\begin{center}
\includegraphics[trim=0 0 0 0, clip,scale=0.75]{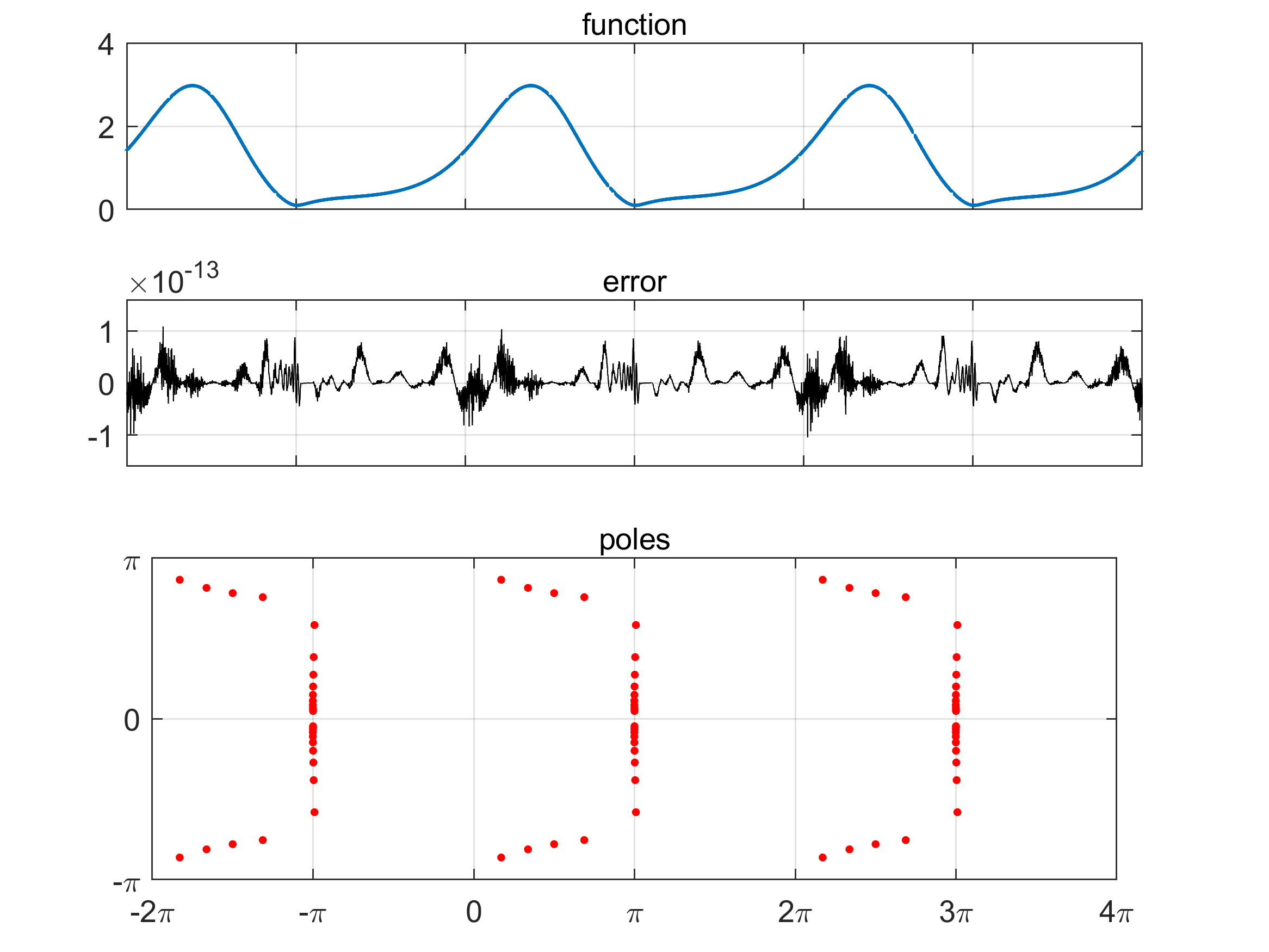}
\end{center}
\vspace*{-8pt}
\caption{\label{figtrig}Periodic AAA approximation of $(\ref{trigexample})$ computed
by {\tt aaatrig}. Three periods are shown.}
\end{figure}

To illustrate periodic AAA approximation, consider the $2\pi$-periodic function
\begin{equation}
f(\theta) = e^{\sin(\theta)} \sqrt{1.01 + \cos(\theta)} ,
\label{trigexample}
\end{equation}
which has branch points at about $\pi\pm 0.14i$ and its translates by
multiples of $2\pi$.  The code

{\small
\begin{verbatim}
      f = @(t) exp(sin(t)).*sqrt(1.01+cos(t));
      T = 2*pi*(1:300)/300;
      [r,pol] = aaatrig(f(T),T);
\end{verbatim}
\par}

\noindent calls {\tt aaatrig} to compute a periodic rational
approximant, and the results are shown in Figure \ref{figtrig}.
The function handle {\tt r} maintains
real symmetry, giving exactly real values for real arguments,
but the vector {\tt pol} of poles deviates from conjugate pairs at the level of
about $10^{-8}$ as a transformation to an algebraic
eigenvalue problem has been used in the computation. 

Which is the better approach to periodic AAA approximation,
explicitly periodic in $\theta$ or transformed to the unit circle
with $z$?  We do not have a definitive answer to this question,
but it is interesting to note that in mathematics more generally,
nobody would recommend eliminating Fourier series
by transplanting every problem to Laurent series on the unit circle.
If we ask ourselves why not, perhaps the first answer is that such a
transplantation converts real problems to complex ones, which adds
computational and conceptual complexity.  On a computer,
transplantation also breaks the symmetry across the real axis,
for whereas the set of complex floating numbers is invariant under
complex conjugation $\theta \mapsto \bar \theta$, it is not invariant under
reflection in the unit circle $z \mapsto \bar z^{\kern 1pt -1}$.  A practical
consequence in the AAA context is that a periodic approximation
to a periodic real function may not come out exactly real, if it is
computed via transplantation to the unit circle,  nor will its
poles necessarily lie in exactly conjugate pairs.

That said, there is no denying the simplicity of computing periodic AAA approximations
by transplantation.
Essentially the same results as in Figure \ref{figtrig} are produced
in this fashion by the following
code calling {\tt aaa}, not {\tt aaatrig}.  The indented lines
effect the transformations between $\theta$ and $z$.

{\small
\begin{verbatim}
      f = @(t) exp(sin(t)).*sqrt(1.01+cos(t));
          g = @(z) f(log(z)/1i);
      T = 2*pi*(1:300)/300;
          Z = exp(1i*T);
      [r0,pol] = aaa(g(Z),Z);
          r = @(t) r0(exp(1i*t)); pol = log(pol)/1i;
\end{verbatim}
\par}

\section{\label{sec:mn}Type \boldmath $(m,n)$ AAA}

A rational function is said to be of {\em type $(m,n)$} if it
can be written as a quotient $p/q$ where $p$ is a polynomial of degree
at most $m$ and $q$ is a polynomial of degree at most $n$.  It is of
{\em exact type $(\kern .6pt \mu,\nu)$} if it is of type $(\kern .6pt \mu,\nu)$ and these
numbers are as small as possible.  Thus
$z/(1+z)^2$ is of exact type $(1,2)$, for example, and it is of type $(m,n)$ for
any $m\ge 1$ and $n\ge 2$.  A constant is of exact type $(0,0)$, except
that the special case of the zero
constant is said to be of exact type $(-\infty,0)$.
A rational function of exact type $(\kern .6pt\mu,\nu)$ has a zero of
order $\nu-\mu$ at $\infty$ if $\mu \le \nu$ and a pole of order
$\mu-\nu$ if $\mu \ge \nu$. 
In the extended complex plane $\complex \cup \{\infty\}$ (the Riemann
sphere), a rational
function of exact type $(\kern .6pt \mu,\nu)$ has exactly $\max\{\mu,\nu\}$ poles and 
also $\max\{\mu,\nu\}$ zeros in total, counted with multiplicity.

Note that according to these definitions,
type $(n,n)$ is the same as degree $n$, and we denote the space
of rational functions of degree~$n$ by $\Rn$.

In approximation theory, it is common for theorists
to consider all values
of~$m$ and~$n$; we denote the corresponding
spaces of rational functions by $\Rmn$.
The {\em Pad\'e table} of a function $f$, for example,
is the array of all of its Pad\'e approximants of types $(m,n)$ for
$m\ge 0$ and $n\ge 0$, each defined by matching the Taylor series of
$f$ as far as possible.\footnote{This is the ``mathematician's
definition'' of Pad\'e approximation \ccite{atap}.
The ``physicist's definition''
\ccite{bgm} is slightly different.}  A notable property of the Pad\'e table is
that degeneracies may occur in which several pairs $(m,n)$ have the
same approximant---for example, this always
happens if $f$ is odd or even---and such
degeneracies always occupy exactly square blocks of the table 
(except in the special case of an identically zero Pad\'e approximant).
Similarly, the {\em \kern -2pt Walsh table} of a function $f$ is the array of
its best $\infty$-norm approximants on a given real or complex
domain $E$ for $m,n \ge 0$.  If $f$ and its rational approximants
are real and $E$ is a real interval,
then the equioscillation theorem implies
that degeneracies in the Walsh table also fill square blocks, with the
same exception for the identically zero approximation.
These approximations can be computed even in quite challenging cases
by the Chebfun {\tt minimax} command, which is based on the Remez
algorithm implemented with a barycentric representation \ccite{minimax}.
For robust computation of Pad\'e approximations, see \ccite{pade}.
All these matters are described in \ccite{atap}.

In approximation practice, on the other hand, usually
the interest is in type $(n,n)$ or sometimes
\hbox{$(n-1,n)$} approximations.
(These classes are coupled in some procedures for
generating continued fractions, which follow a staircase sequence
$(0,0), (0,1), (1,1), (1,2),\dots$ or
$(0,0), (1,0), (1,1), (2,1),\dots$ \ccite{henrici2,salazar}.)
This is one reason why the standard AAA algorithm,
and its Chebfun implementation {\tt aaa.m},
treat just the case of type $(n,n)$.  Nevertheless, a method
for extending AAA to type $(m,n)$ was proposed in section~9 of the
original AAA paper \ccite{aaa}.

Consider the barycentric representation (\ref{aaa-baryrep}) that
AAA is based upon.  As is shown in Theorem 2.1 of \ccite{aaa}, this
representation covers not just rational functions of exact type
$(n,n)$ but {\em all\/} rational functions of type $(n,n)$.  This statement
applies for any fixed set of support points $\{\tk\}$\kern 1pt---having a 
good choice matters numerically but not algebraically---and it 
applies even if attention is restricted to weights $\{\bk\}$ that
are all nonzero.  Thus, algebraically at least, it should be possible
to retain the form (\ref{aaa-baryrep}) for computing
type $(m,n)$ approximations with $m<n$; for $m>n$ we would increase
the limit on the sum to $m$.  The method for doing this proposed in
\ccite{aaa} is
to restrict the
weight vector $\{\bk\}$ to lie in a subspace where the numerator
degree is constrained to be less than the denominator degree if
$m<n$, and the reverse if $m>n$.  This is a condition of linear algebra
that can be imposed by
modifications of the computation that we shall not go into.
These computations implicitly work with polynomials,
and for numerical stability, it is important that a good basis for
the polynomials is utilised.  This can be achieved by employing
Vandermonde with Arnoldi orthogonalisation as described
in \ccite{VA}.  All these elements are combined successfully in
Chebfun {\tt minimax} \ccite{minimax}, but
to this date they have not yet been implemented
in a general-purpose AAA code.
Therefore, unusually, we shall conclude this discussion without
a numerical illustration.
$$
*~~*~~*
$$

However, a related matter must be brought up; we focus on the case
$(n-1,n)$ because of its practical interest.
{\em To say that\/ $r$ is of type $(n-1,n)$ rather than $(n,n)$ is
not the same as saying that it takes the value $0$ at $\infty$.}
For example, $r(z) = (z+1)/(z-1)$ is of type $(1,2)$, but
$r(\infty)=1$.  
In an application, one must be clear which of the two is wanted, and if
it is $r(\infty) = 0$, then a new algorithm will be needed from the one
mentioned above.

The condition $r(\infty)=0$ is important in
model order reduction (section \ref{sec:mor}), where the interest is
usually in approximations that tend to zero on the imaginary axis.
In the MathWorks
RF Toolbox \ccite{rf}, there is a flag {\tt TendsToZero} that the user
can set to impose this condition.  In the Control System and System Identification
toolboxes, similarly, one can specify a value at $\infty$ that may be
zero or nonzero.  Mathematically, interpolation conditions like these are
quite special, as $z=\infty$
is just one point among many of $\complex\cup \{\infty\}$.  Whereas
one could in principle impose a condition at $z=\infty$ in any rational approximation problem,
it is only likely to serve much purpose in a problem involving
approximation on a domain
that extends to $\infty$, like the imaginary axis.

Best $\infty$-norm approximation with $r(\infty)=0$ may be
problematic, for it is ill-posed.  (More generally, best approximation of
exact type $(\kern .6pt \mu,\nu)$ as opposed to type $(m,n)$ is
ill-posed.)  We can show this
with the example of approximation of $f(x) = 2 \kern .5pt x e^{1-x}$
for $x\in [\kern .5pt 0,\infty)$.  This function takes the extreme values
$0$ at $x=0$ and $\infty$ and $2$ at $x=1$, so by the equioscillation theorem,
$r(x) = 1$ is its unique best approximant of type $(1,1)$, with error
$\|\kern .7pt f-r\| = 1$.  Now suppose we ask for a best approximant of
type $(1,1)$ with
$r(\infty)=0$.  By uniqueness of best approximations, 
such a function will have to
have $\|\kern .7pt f-r\|>1$. 
It can come arbitrarily close to
$\|\kern .7pt f-r\|=1$, as we see by considering
$r(x) = 1/(1+\varepsilon x)$
as $\varepsilon \to 0$, but the limit cannot be attained.  

Fortunately, the ill-posed cases are non-generic, and approximation
with $r(\infty)=0$ often makes sense.  A method for computing
such approximants is
to replace the usual barycentric formula 
(\ref{aaa-baryrep}) by
\begin{equation}
r(z) =
\sum_{k=1}^n \left. {\fk\kern .5pt \bk\over z- \tk} \kern 1pt \right/
\kern -2pt \left( 1 + \smash{\sum_{k=1}^n}\kern 2pt  {\bk \over z- \tk} \right)
\label{baryagain}
\end{equation}
and adjust the AAA iteration accordingly.  The big change is
that it is no longer appropriate to impose the normalisation
(\ref{aaa-normalize}) on the weights $\{\bk\}$, so instead of an SVD,
we get a least-squares problem at each step.  The fact that $\{\bk\}$
may be arbitrarily large, diverging to $\infty$ in special cases,
reflects the risk of ill-conditioning and in the limit ill-posedness.
Nevertheless (\ref{baryagain}) is often an effective basis for AAA
approximation with $r(\infty)=0$, and this variant of the barycentric
formula has been used by a number of authors going back at least
to \cite[sec.~1.5.2]{ionita}.

For approximation with more than one zero at $\infty$, Dave Darrow
(private communication) has proposed a similar algorithm
with the constant~$1$ of (\ref{baryagain}) replaced by a monic
polynomial of degree $\nu-\mu-1$.  Presumably in general it would
be desirable to use Vandermonde with Arnoldi orthogonalisation to
represent that polynomial, but we have no experience with this.

\section{\label{sec:sets}Vector-valued AAA}
Instead of a single function $f(z)$, in some applications one
has several functions $\{f_j^{}(z)\}$ that one would like to
fit by rational functions $\{r_j^{}\}$ all with the same poles.
This could happen in an engineering context, for example, when
various signals are available that ultimately derive from the same
dynamical process, as in a MOR or ROM problem of the kind discussed
in section~\ref{sec:mor}.  We also encountered such a situation in
sections~\ref{sec:res}--\ref{sec:nonlin}, where a linear or nonlinear
resolvent $A(z)^{-1}$ is scalarised by conversion to $u^* A(z)^{-1}
v$, with $u$ and $v$ taken as arbitrary vectors, and the results
may be more accurate if these are replaced by matrices with several
independent columns.  If $u$ and $v$ have $p$ and $m$ columns,
respectively, as in (\ref{eq:Happrox}), then we have a vector-valued
approximation of dimension $p\kern .6pt m$ to construct.

Rational fits of this kind, with shared poles, have not only
arisen in the context of AAA.\ \ An interesting example is
\ccite{wrightfornberg}, which presents an algorithm for shared-pole
rational approximation for applications to radial basis functions
(\RBFs).  In their application, an RBF $f_\varepsilon(\vec x\kern
1.2pt)$ depends analytically on the shape parameter $\varepsilon$,
and values $\varepsilon \approx 0$ are impossible to determine by
direct evaluation but are accessible via rational approximation based
on data at larger values of $\varepsilon$.  By using several values
of $\vec x\kern 1.2pt$ to construct a fit, the authors obtain a
more reliable and accurate answer than if they used just one value
of $\vec x\kern 1.2pt$.

For AAA, generalisations of the basic algorithm to handle vector-valued data include
``FastAAA'' \ccite{hochman}, the MathWorks {\tt rational} code \ccite{rf},
and ``setAAA'' \ccite{setAAA} (first posted on arXiv in 2018; software
is available at \ccite{meerbergen}.
The idea is simple, as expressed in the words of Hochman:
\vskip 6pt

{\leftskip 1.2cm
\small\noindent
We replace the L\"owner matrix$\dots$ by a block-L\"owner analogue.
\hfill\break
Its size is now $p(N-n)\times n$.
\par}
\vskip 6pt

\noindent
(In the notation of the present paper around (\ref{aaa-Aentries}),
this becomes $p(m-n-1)\times (n+1)$.)  In other words, one fits
all the functions simultaneously, and the AAA least-squares
problem at each step minimises a residual based on all of them.
The result is an approximation with a single set of support points
and barycentric weights and hence a single set of poles, but with $p$
different function values interpolated at each each support point.
Since the number of rows of the Loewner matrix is multiplied by $p$
but the number of columns is not, the work scales just linearly
with $p$.  In applications involving scalarisation vectors $u$ and
$v$ as described above, for example, one can increase accuracy and
robustness by replacing these vectors by tall-skinny matrices with
$3$ or $4$ columns, say, and thereby increase the computational
cost by just a factor on the order of $9$ or $16$.

For an idea of what dimensions may be relevant in
applications, it is interesting to look at the SLICOT collection
of model problems for MOR/ROM mentioned in section \ref{sec:mor}
\ccite{slicot}.  Of the 56 problems, 19 involve scalars and 37
involve vectors.  Most of these 37 have dimensionality just 2
or 3, with only one problem having dimensionality bigger than 7
(namely 20).  Thus all of these problems would be readily accessible
to vector-valued AAA.\ \ On the other hand, the \hbox{SLICOT}
collection dates to 2002, and much larger problems are also of
interest nowadays.

In this section we have discussed vector-value approximations,
making no distinction between vectors and matrices since the
approximations in question treat components individually though
with shared poles.  FastAAA, setAAA and the MathWorks RF Toolbox
all operate in this mode.  In \ccite{goseaalgs}, a different
approach called blockAAA is introduced, which deals intrinscially
with matrix-valued approximations.  This can be advantageous in
MIMO applications, and the method is implemented in the MathWorks
Control System and System Identification toolboxes.
\
\section{\label{sec:discussion}Discussion}

Rational approximations are changing numerical analysis, and more
broadly, applied mathematics.  The change was sparked by AAA,
which has been the tool used in this review, but other methods
may also prove competitive.  As mentioned in the introduction,
the greedy Thiele Continued Fractions (TCF) method is impressive
\ccite{salazar,driscollzhou,driscolljuliacon}.  It is faster than
AAA, and while it does not take advantage of least-squares fitting,
it may perhaps prove to be equally reliable.  One limitation of TCF
is that, so far at least, it does not have a variant analogous to
AAA-Lawson for upgrading near-best to best approximations.

We have discussed many applications, but there are potentially
many more, for rational functions have a role to play in almost
all areas of computational continuous mathematics.  One topic
we are particularly curious about is the theory and computation
of {\em composite\/} rational approximations, as discussed in
the final paragraph of section \ref{sec:zolo}.  As argued in
\ccite{everything}, almost everything one can do numerically on a
computer comes down to piecewise composite rational approximations,
and the machine learning revolution has been driven by composite
approximations---not necessarily rational, although as shown in
\ccite{boulle}, they may be.

For AAA, TCF and potentially other methods, an important question is
to what extent continuum variants as in section \ref{sec:continuum}
should replace the discrete approximations we have emphasised here.
Most applica\-tions---in fact, {\em every\/} example considered
in this paper apart from that of Figure~\ref{discretefig}---ultimately
concern continuous domains, and
Driscoll has shown that continuum AAA is extremely
effective.  With more experience, it is possible that it will
emerge that continuum rational approximation algorithms should be
the default.

Many research questions have arisen in this survey, including the following.
\begin{enumerate}[label=\arabic*.,itemsep=1pt]

\item Section \ref{sec:AAA}:
How might AAA linear algebra be modified to speed up $O(n^4)$
to $O(n^3)$ in a numerically stable fashion?

\item Section \ref{sec:AAA}:
Can theorems be developed to explain the rapid convergence of AAA to
near-best approximations?  (Question~13
below may play into this investigation.)

\item Section \ref{sec:AAA}: What are the prospects for developing
an AAA-style algorithm that avoids producing poles in unwanted 
regions of $\real$ or $\complex$?  A method based on semidefinite
programming has been put forward
in \ccite{bradde}.  See also question 19.

\item Sections \ref{sec:AAA} and \ref{sec:equi}:
Can methods be developed to help AAA stop at the right moment if
noise is encountered, or to otherwise reduce the impact of noise?
For polynomial approximations, a method is proposed 
in \ccite{matsuda}.

\item Sections \ref{sec:AAA} and \ref{sec:poles}: Should AAA be modified to measure errors
in the Riemann sphere metric $|w-w'|/\sqrt{(1+|w|^2)(1+|w'|^2)}$
rather than the complex plane metric $|w-w'|$?  This might
improve performance in problems with poles near the approximation domain,
and it might make polefinding and zerofinding more nearly equivalent.

\item Section \ref{sec:branch}:
How can rational approximations be used to identify branch point structure
for functions with several branch points, as in
Figures \ref{schwarz-squiggles}, \ref{schwarz-squiggles3},
and \ref{DEs-lorenzfig}?

\item Section \ref{sec:branch}:
Can a AAA-style algorithm be developed for reciprocal-log
approximations \ccite{baddooT,reciplog}?  This would have the potential
of speeding up root-exponential to exponential convergence, and also
of offering new methods for dealing with functions near branch points and on
Riemann surfaces, and for approximation on domains with slits.

\item Section \ref{sec:equi}:
What theorem might be proved to explain the success of rational
approximations in interpolating equispaced data?

\item Section \ref{sec:analcont}:
What's going on with the ``one-wavelength principle'' in analytic
continuation?

\item Section \ref{sec:schwarz}: How can the branch structure of Schwarz functions
be computed numerically?  See question 6.

\item Section \ref{sec:potential}: Can a Hermite integral theory for
rational approximations be developed based on zeros of $r(z)$
instead of poles, or based on both zeros and poles?  See question~5 above.

\item Section \ref{sec:solution}:
How can adaptive rational approximations be used on the fly
to solve \ODEs\ and \PDEs\ whose solutions have
singularities or near-singularities?  The vision here is to
create a kind of ``Ratfun'' in analogy to Chebfun.

\item Section \ref{sec:zolo}:
How can the \verb|'sign'| variant of AAA be modified into a
bulletproof procedure that might replace the current default
selection of the smallest singular vector of the Loewner matrix?

\item Section \ref{sec:cauchy}:
Can AAA be useful for solving more advanced Riemann-Hilbert problems?
Very little has been attempted so far.

\item Sections \ref{sec:aaals}--\ref{sec:biharmonic}:
What theoretical support can be developed for the AAALS method,
and for its extensions from the Laplace equation to other \PDEs\
such as the Helmholtz and biharmonic equations?

\item Section \ref{sec:biharmonic}:
Can AAA placement of stokeslets be the basis of a new MFS
method for solving biharmonic problems?

\item Sections \ref{sec:aaals}--\ref{sec:biharmonic}:
Can AAALS-style methods be developed
based on monopoles $\log|z-\tk|$ rather than dipoles $1/(z-\tk)\kern 1pt$?

\item Section \ref{sec:lawson}:
Can the AAA-Lawson algorithm be improved to give guaranteed
convergence, at least in exact arithmetic?
Presumably this will require a reconsideration of the
\verb|'damping'| option discussed in section \ref{sec:zolo}.

\item Section \ref{sec:continuum}:
Can the continuum AAA algorithm be modified to avoid ``bad poles''
more unfailingly, such as arise in many problems especially with
real functions on real intervals?

\item Section \ref{sec:mn}:
Can type $(m,n)$ AAA approximation
be realised in a fashion as robust numerically as what we have for
type $(n,n)$?
\end{enumerate}

These questions are all important, but let us distinguish two further
questions which perhaps are even more so.  One is, what can be done
in multiple dimensions?  At present, neither multivariate rational
approximation per se nor applications to \PDEs\ in more than two
dimensions have tools comparable to univariate AAA (see section
\ref{sec:bivariate}).  Is this restriction to a single variable---or to
the plane---unbreakable?

Secondly, now that greedy algorithms are changing rational
approximation so profoundly, what about analogous algorithms
for other problems that might seem to have a similar flavor?
In particular, how about fitting a function $f(x)$ near-optimally
by a linear combination of Gaussians, or radial basis functions,
of unknown widths and centres?  If there were a AAA-style algorithm
for such problems, applications would be found everywhere.  So far,
we have been stymied by the lack of a barycentric representation
for Gaussians or \RBFs, but is there an idea here we haven't spotted?

\section*{Acknowledgments}
\addcontentsline{toc}{section}{Acknowledgments}
Every section of this paper has been read and reread by a number of
our colleagues during the writing process, and the combined help
from all of these people has added up to a major contribution.
We would like to thank Anthony Austin, Alex Barnett, Jean-Paul
Berrut, Felix Binkowski,
Jake Bowhay, Oscar Bruno, Sven Burger, Keaton Burns, Anil Damle,
Dave Darrow, Tom DeLillo, Mark Embree, Marco Fasondini,
Eamonn Gaffney, Abi Gopal,
Victor Gosea, Serkan Gugercin, Stefan G\"uttel, Sam Harris, Astrid Herremans,
Sam Hocking, Andrew Horning, Daan Huybrechs, Cosmin Ionita, Anastasia
Kisil, Andrei Martinez-Finkelshtein, Robb
McDonald, Kyle McKee, Karl Meerbergen, Andrea Moiola, Hadrien
Montanelli, Mark Reichelt, Oliver Salazar Celis, Manuel Santana,
Olivier S\`ete, Nian Shao, Shuai Shao, Anna-Karin Tornberg,
Alex Townsend, Tom Trogdon, Michael Tsuk, Kieran Vlahakis, Andr\'e Weideman,
Nicholas West, Heather Wilber, Grady Wright, Tony Xu, Fei Xue,
Yidan Xue and Maxim Yattselev.  A few errors in a late draft of the paper
were also caught by an AI system: by the Refine software created by
Yann Calv\'o L\'opez and Ben Golub.

In addition to these wonderfully helpful friends and colleagues, we wish to single
out two others whose contributions to the AAA enterprise have been particularly
extensive and are visible in many sections of this article.  One is Stefano Costa, the originator
of the AAALS method and the first to propose that it could be applied not just
to Laplace problems on simple domains but also to \PDEs\ with interfaces,
conformal mapping, Riemann-Hilbert problems and other applications.  As we have
worked on AAALS in the past five years, many of the best ideas for improvements
and understanding
of the method have come from Costa.  The other is Toby Driscoll, who was already a close
colleague of ours through Schwarz-Christoffel mapping,
the Chebfun project and the book {\em Exploring \ODEs}.  
For AAA, Driscoll has been a leader in Julia implementations, high-precision analytic
continuation,
exploration of the alternative method of Thiele continued fractions,
and the development and implementation of the continuum AAA algorithm.  

\clearpage 
\addcontentsline{toc}{section}{References}
\bibliography{refs}
\clearpage 
\label{lastpage}
\end{document}